\documentclass[11pt,leqno]{amsart}
\usepackage{a4wide}
\usepackage{amssymb}
\usepackage{amsmath}
\usepackage{amsthm}
\usepackage{mathrsfs}
\usepackage{stmaryrd}
\usepackage{enumerate} 

\usepackage{euscript}
\usepackage[all]{xypic}
\usepackage[usenames]{color}
\usepackage{hyperref}
\usepackage{mathabx} 

\usepackage[T1]{fontenc}
\usepackage[utf8]{inputenc}
\usepackage{comment}
\usepackage{tikz-cd} 

\newcommand{\Footnote}[1]{\footnote{{\color{Sepia} #1}}}  
\renewcommand{\Footnote}[2][]{\relax}  

\theoremstyle{plain}

\newcommand{\id}{\operatorname{id}}
\newcommand{\im}{\operatorname{im}}
\newcommand{\ev}{\operatorname{ev}}
\newcommand{\gr}{\operatorname{gr}}

\newcommand{\Sp}{\operatorname{Sp}}

\newcommand{\Hom}{\operatorname{Hom}}
\newcommand{\End}{\operatorname{End}}

\newcommand{\Fil}{\operatorname{Fil}}

\newcommand{\Lie}{\operatorname{Lie}}
\newcommand{\Ind}{\operatorname{Ind}}

\newcommand{\Rep}{\operatorname{Rep}}

\newcommand{\Spec}{\operatorname{Spec}}

\newcommand{\GL}{\operatorname{GL}}

\newcommand{\Gal}{\operatorname{Gal}}
\newcommand{\rk}{\operatorname{rk}}

 \newcommand{\CC}{\mathbb{C}}

 \newcommand{\NN}{\mathbb{N}}
 \newcommand{\QQ}{\mathbb{Q}}
 
 \newcommand{\ZZ}{\mathbb{Z}}

 \newcommand{\bD}{\mathbf{D}}

 \renewcommand{\cD}{\mathcal{D}}
 
 \renewcommand{\cL}{\mathcal{L}}

 \newcommand{\cC}{\mathcal{C}}

\newcommand{\cT}{\mathcal{T}}

\renewcommand{\cR}{\mathcal{R}}

\newcommand{\cP}{\mathcal{P}}

 \newcommand{\Qp}{\QQ_p}

 \newcommand{\be}{\begin{equation}}
\newcommand{\ee}{\end{equation}}

\newcommand{\fpbar}{\ifmmode {\overline{\mathbb{F}_p}}\else$\mathbb{F}_p$\ \fi}

\newcommand{\fp}{\ifmmode {\mathbb{F}_p}\else$\mathbb{F}_p$\ \fi}
\newcommand{\zp}{\ifmmode \mathbb{Z}_p\else$\mathbb{Z}_p$\ \fi}
\newcommand{\zpur}{\ifmmode \widehat{\zp^{ur}}\else $\widehat{\zp^{ur}}$\ \fi}

\renewcommand{\d}{\mathbf{d}}
\renewcommand{\u}{\mathbf{1}}

 \newcommand{\TLT}{T_{\pi}}

 \newcommand{\cris}{\mathrm{cris}}
 \newcommand{\dR}{\mathrm{dR}}

 \newcommand{\perf}{\mathrm{perf}}
 \newcommand{\pst}{\mathrm{pst}}
\renewcommand{\H}{\mathrm{H}}

\newcommand{\Kl}{\{\!\!\{}
\newcommand{\Kr}{\}\!\!\}}

\newcommand{\La}{\ifmmode\Lambda\else$\Lambda$\fi}

\newcommand{\q}{\ifmmode {\Bbb Q}\else${\Bbb Q}$\ \fi}
\newcommand{\qp}{\ifmmode {\Bbb Q}_p\else${\Bbb Q}_p$\ \fi}
\newcommand{\z}{\mathbb{Z}}

\newcommand{\Q}{\ifmmode {\Bbb Q}\else${\Bbb Q}$\ \fi}
\newcommand{\ql}{\ifmmode {{\Bbb Q}_l}\else${\Bbb Q}_l$\ \fi}
\renewcommand{\d}{\mathbf{d}}
\renewcommand{\u}{\mathbf{1}}

\newcommand{\betr}[1]{\lvert #1\rvert}

\newcommand{\nach}{\longrightarrow}
\newcommand{\auf}{\longmapsto}

\newcommand{\nachinj}{\lhook\joinrel\relbar\joinrel\rightarrow}
\newcommand{\iso}{\overset{\sim}{\nach}}

\DeclareMathOperator{\Ddif}{\mathbf D_{\mathrm{dif}}}
\DeclareMathOperator{\Ddifp}{\mathbf D_{\mathrm{dif}}^+}

\DeclareMathOperator{\Frob}{Frob}
\DeclareMathOperator{\dif}{dif}

\newcommand{\PDif}{\mathfrak D}
\DeclareMathOperator{\cone}{cone}

\DeclareMathOperator{\Tr}{Tr}

\newcommand{\phigam}{\mathfrak}
\DeclareMathOperator{\Max}{Max}

\newcommand{\mf}{\mathfrak}
\newcommand{\tn}{\textnormal}

\DeclareMathOperator{\expnM}{\exp^{(n)}_M}

\DeclareMathOperator{\expnfM}{\exp^{(n)}_{f,M}}

\newtheorem{theorem}{Theorem}[section]
\newtheorem{corollary}[theorem]{Corollary}
\newtheorem{lemma}[theorem]{Lemma}
\newtheorem{remark}[theorem]{Remark}
\newtheorem{proposition}[theorem]{Proposition}
\newtheorem{conjecture}[theorem]{Conjecture}
\newtheorem{definition}[theorem]{Definition}

\newtheorem{example}[theorem]{Example}

\newtheorem*{bconj}{Conjecture}

\theoremstyle{remark}

\author{Milan Malcic, Rustam Steingart,  Otmar Venjakob and Max Witzelsperger}
\begin{document}

\title{\textbf{$\epsilon$-isomorphisms for rank one $(\varphi,\Gamma)$-modules over Lubin-Tate Robba rings}}
\maketitle
\begin{abstract}
Inspired by Nakamura's work \cite{NaANT} on $\epsilon$-isomorphisms for $(\varphi,\Gamma)$-modules over (relative) Robba rings with respect to the cyclotomic theory, we formulate an analogous conjecture for $L$-analytic Lubin-Tate $(\varphi_L,\Gamma_L)$-modules over (relative) Robba rings for any finite extension $L$ of $\Qp.$ In contrast to Kato's and Nakamura's setting, our conjecture involves $L$-analytic cohomology instead of continuous cohomology within the generalized Herr complex. Similarly, we restrict to the identity components of $D_{cris}$ and $D_{dR},$ respectively. For rank one modules of the above type or slightly more generally for trianguline ones, we construct $\epsilon$-isomorphisms for their Lubin-Tate deformations satisfying the desired interpolation property.
\end{abstract}
\tableofcontents

\section{Introduction}

In \cite{NaANT} Nakamura generalized Kato's $p$-adic local $\epsilon$-conjecture \cite{Kat, FK} to the framework of $(\varphi,\Gamma)$-modules over the Robba ring (over $\Qp$-affinoid algebras) and proved the essential parts of it for rigid analytic families of trianguline  $(\varphi,\Gamma)$-modules. The technical foundations for this had been laid by the work of Kedlaya, Pottharst and Xiao \cite{KPX} who had established the fundamental theorems concerning their cohomology (finiteness, base change property, Tate duality, Euler-Poincar\'{e} formula) and Nakamura's work \cite{Na}, in which he generalized the theory of Bloch-Kato exponential maps and Perrin-Riou's exponential maps in that framework.

Recently there has been much progress concerning $(\varphi_L,\Gamma_L)$-modules over Lubin-Tate extensions \cite{Fou,KR,BF,FX,SV15, SV20}. In particular, the results by Steingart \cite{Ste1,Ste2} regarding such $(\varphi_L,\Gamma_L)$-modules over families (finiteness, base change property, Euler-Poincar\'{e} formula, perfectness of Iwasawa cohomology) make it possible to study a version of Nakamura's approach  for $L$-analytic trianguline modules.

 Let $L \subseteq \CC_p$ be a finite extension of $\QQ_p$ and $L_\infty$ a Lubin-Tate extension of $L$ with Galois group $\Gamma_L = \operatorname{Gal}(L_\infty/L)$ corresponding to a uniformiser $\pi_L$ of the ring of integers $o_L$ of $L.$ A continuous representation of $G_L$ on a finite dimensional $L$-vector space $V$ is called $L$-analytic, if the semi-linear representation $\CC_p\otimes_{\QQ_p}V \cong \prod_{\sigma\colon L \to \CC_p}\CC_p\otimes_{L,\sigma}V$ is trivial at the components where $\sigma \neq \id.$ By a theorem of Berger the category of $L$-analytic representations is equivalent to the category of \'{e}tale $L$-analytic $(\varphi_L,\Gamma_L)$-modules over the Robba ring $\cR_L$ (cf. \cite{Be16}). Analyticity means here, that the action of the Lie group $\Gamma_L$ is differentiable and the action of $\operatorname{Lie}(\Gamma_L)$ is (not only $\Qp$-, but even) $L$-bilinear.
For analytic $(\varphi_L,\Gamma_L)$-modules one can define analytic cohomology (see Section \ref{sec_analcoho} for a precise definition).  	Finiteness of analytic cohomology allows us to attach to a family $M$ of analytic $(\varphi_L,\Gamma_L)$-modules over $A$ a graded invertible line bundle $\Delta_A(M)$ over $A$ which is essentially the determinant of the analytic cohomology of $M$. Note that, for an  $L$-analytic \'{e}tale $(\varphi_L,\Gamma_L)$-module attached to some $L$-analytic Galois representation $V$ of $G_L$ with coefficients in $L$, these analytic cohomology groups in general do {\it not} coincide with the Galois cohomology groups $H^i(L,V)$ of $V$ for $i>0.$ Nonetheless they behave similarly to Galois cohomology and allow us to study certain invariants of $V$ ``at the identity component''. If $M$ is the $(\varphi_L,\Gamma_L)$-module attached to an $L$-analytic de Rham representation $V,$ then one can also attach an $\varepsilon$-constant to the ``identity component'' of $D_{pst}(V),$ i.e., the $G_L$-smooth vectors in  $\mathbf{B}_{st}\otimes_{L_0}V$ (which injects into the full $\mathbf{B}_{st}\otimes_{{\QQ_p}}V$). This can be generalised to the non-\'{e}tale case as well (see Section \ref{sec:deRham} for details).
	 The content of the analytic variant of the $\varepsilon$-conjecture is a trivialisation of $\Delta_A(M)$ which interpolates these $\varepsilon$-constants at the de Rham points, i.e., the points $x \in \operatorname{Sp}(A)$ where the specialisation $M_x$ is de Rham.

We formulate the following conjecture in a more general setting  (and indicate in Remark \ref{rem:compNakamura} (ii) how to formulate a version of this conjecture for $L$-analytic $(\varphi_L,\Gamma_L)$-modules over the character variety $\mathfrak{X}_{o_L}$ in the sense of Schneider-Teitelbaum).
\begin{bconj}[See Conjecture \ref{conj}] \label{conj:intro}
	Choose a compatible system $u=(u_n)$ of $[\pi_L^n]$-torsion points of the Lubin-Tate group and a generator $t'_0$ of the Tate module of its Cartier dual.
	Let $K$ be a complete field extension of $L$ containing $L^{ab}$, and
	$A$ an affinoid algebra over $K$.
	 For each $L$-analytic $(\varphi_L,\Gamma_L)$-module $M$ over $\cR_A$  satisfying condition \eqref{cond:charactertype}
	there exists a unique trivialisation $$\varepsilon_{A,u}(M): \u_A \xrightarrow{\cong}\Delta_A(M)$$ satisfying the following axioms:
	\begin{enumerate}
		\item For any affinoid algebra $B$ over $A$ we have $$\varepsilon_{A,u}(M)\otimes_A \id_B = \varepsilon_{B,u}(M\hat{\otimes}_AB)$$ under the canonical isomorphism $\Delta_A(M)\otimes_AB \cong \Delta_B(M\hat{\otimes}_AB).$
		\item $\varepsilon_{A,u}$ is multiplicative in short exact sequences.
		\item For any $a \in  o_L^\times$ we have $$\varepsilon_{A,a \cdot u}(M)=\delta_{\det M}(a)\varepsilon_{A,u}.$$
		\item $\varepsilon_{A,u}(M)$ is compatible with duality in the sense that  for the dual module $\tilde{M}$ (see section \ref{sec:duality}) we have
		$$\varepsilon_{A,u}(\tilde{M})^* \otimes h(\chi^{r_M}) = (-1)^{\dim_KH^0(M)}\Omega_{t'_0}^{-r_{M}} \varepsilon_{A,-u}(M)$$ under the natural isomorphisms $\mathbf{1}_A \cong \mathbf{1}_A \otimes \mathbf{1}_A$ and $\Delta(M)\cong \Delta(\tilde{M})^* \otimes (A(r_M),0),$ where $h(\chi^{r_M})  \colon A(r_M) \to A$ maps $e_{\chi^{r_M}}$ to $1$ and     $r_M$ denotes the rank of $M$ over $\cR_K$.
		\item For $L=\Qp$, $\pi_L=p$ and $u=(\zeta_{p^n}-1)_n$ the trivialisation coincides with that of Nakamura, in the sense of Proposition \ref{prop:precisecomp}.
		\item Let $F/L$ be a finite subextension of $K,$ $M_0$ be a de Rham $(\varphi_L,\Gamma_L)$-module over $\cR_F$ and $M = K\hat{\otimes}_FM_0.$ Then $$\varepsilon_{K,u}(M)=\varepsilon^{dR}_{F,u}(M_0),$$
 where the isomorphism $\varepsilon^{dR}_{F,u}(M_0): \mathbf{1}_K\xrightarrow{\cong} \Delta_K(M)$ is called the de Rham $\varepsilon$-isomorphism which is defined in \eqref{f:epsdR} unconditionally using a generalized  Bloch–Kato
exponential and   dual exponential map as well as the  $\varepsilon$-constant associated to $M_0$ in section \ref{sec:deRham}.
	\end{enumerate}
\end{bconj}
While in the cyclotomic setting the  $\epsilon$-constants depend on the choice of a norm compatible system of $p$-power roots of unity, in the Lubin-Tate setting this is replaced by a compatible system of $\pi_L$-power torsion points of the Lubin-Tate formal group, see Remarks \ref{rem:psiu}, \ref{rem:epsind} for a comparison of both. We also fix a generator $t'_0$ of the Tate module of the Cartier dual of the Lubin Tate group which determines a certain period $\Omega_{t'_0} \in \CC_p$ (cf. \cite{ST}).
We prove parts of this conjecture for $L$-analytic trianguline $(\varphi_L,\Gamma_L)$-modules. More precisely, we construct the $\varepsilon$-isomorphism for the Lubin-Tate deformation of a rank one, de Rham $L$-analytic $(\varphi_L,\Gamma_L)$-module $M$ over some finite extension $F$ of $L$
\[\varepsilon_{D(\Gamma_L),u}(\mathbf{Dfm}(K\hat{\otimes}_FM)):\mathbf{1}_{D(\Gamma_L)}\xrightarrow{\cong} \Delta_{\mathfrak{X}_{\Gamma_L}}(\mathbf{Dfm}(K\hat{\otimes}_FM)),\]
see Theorem \ref{thm:LTDeform}. This lives over the rigid analytic character variety $\mathfrak{X}_{\Gamma_L}$  over $L$.  The $\mathbb{C}_p$-points of this variety correspond to locally $L$-analytic characters $\Gamma_L\to \mathbb{C}_p^\times.$  We refer to subsection \ref{sec:Iwasawa} for the precise definition of the Lubin-Tate deformation $\mathbf{Dfm}(N)$ of a $(\varphi_L,\Gamma_L)$-module $N$ over $\cR_K.$ Heuristically one can think of it as the base changed $(\varphi_L,\Gamma_L)$-module $D(\Gamma_L,K)\hat{\otimes}_K N$ over the relative Robba ring $D(\Gamma_L,K)\hat{\otimes}_K\cR_K.$ But due to the complicated behaviour of completed tensor products over LF-spaces which are not Fr\'{e}chet, it requires a more technical treatment. The correct point of view, which is used for the cyclotomic setting in earlier articles of Pottharst (but apparently neither consequently pursued nor carefully explained  in \cite[Def.\ 4.4.7, Thm.\ 4.4.8]{KPX} unfortunately), consists of viewing this deformation as a sheaf of $(\varphi_L,\Gamma_L)$-modules over $\mathfrak{X}_{\Gamma_L},$  which is not affinoid and hence does not strictly speaking fit into the above Conjecture.
Instead, the isomorphism $\varepsilon_{D(\Gamma_L),u}$ is a trivialisation of a line bundle over $\mathfrak{X}_{\Gamma_L}$ which restricts to an isomorphism of the conjectured type on each affinoid subdomain.

{\it Philosophically,} the $L$-analytic theory over Lubin-Tate extensions is {\it one-dimensional} and thus very similar to the cyclotomic case in the sense that $\Gamma_L$ is - although $[L:\Qp]$-dimensional over $\Qp$   - one-dimensional as $p$-adic Lie group over $L$. Nevertheless, {\it technically} we have had to overcome serious difficulties. We are going to describe these differences compared to Nakamura's work in the following.

In the cyclotomic setting, {\it Herr-complexes} are formed with respect to the two operators $\varphi$ and $\gamma-1$ for a topological generator $\gamma$ of the torsion-free part of $\Gamma;$ moreover, one can directly go over to the complex consisting of the fixed part under the torsion subgroup $\Delta$ of $\Gamma$. In the Lubin-Tate setting (with $L\neq \Q$) there is no intrinsic counterpart of $\gamma$ as one needs at least $[L:\Qp]$ elements to generate the (torsion-free part of) $\Gamma_L$ topologically. So instead we make use of Fourier theory and the Lubin-Tate isomorphism \`{a} la Schneider and Teitelbaum \cite{ST}
\[D(\Gamma_n,K)\cong \mathcal{O}(\mathfrak{X}_{\Gamma_n})\cong  \mathcal{O}({\mathbf{B}})\]
over a huge field extension $K$ of $L$, over which the character variety $\mathfrak{X}_{\Gamma_n}$ for the subgroup of $n$-th higher units $\Gamma_n\cong o_L$ of $\Gamma_L$ can be identified with the open unit disk $\mathbf{B}$ for $n$ sufficiently big. Via this  isomorphism
we can now choose ${\mf Z_n}\in D(\Gamma_n,K) $ corresponding to the choice of a coordinate of $\mathbf{B}$. The generalized Herr-complex in the Lubin-Tate setting can thus be formed using the two operators $\varphi_L$ and $\mf Z_n.$ Unfortunately, in contrast to $\Delta\subseteq \Gamma_{\Qp}$, the remaining quotient $\Gamma_L/\Gamma_n$ in general cannot be identified with a subgroup of $\Gamma_L$, whence we cannot take $\Gamma_L/\Gamma_n$-invariants as before, but have to circumvent this problem.

 An important step for our approach consists of establishing the analogue of local Tate duality for analytic cohomology, see subsection \ref{sec:duality}. In contrast to \cite{NaANT} we find an intrinsic way to normalize our trace map without any comparison to Galois cohomology (which is not available anyway as we indicated); nevertheless for $L=\Qp$ our choice coincides with that of Nakamura (for an appropriate choice of period $ \Omega$).

Another price we have to pay is the fact that even the minimal choice for $K$ is no longer spherically complete, which means that the functional analysis requires some additional care. For the explicit descent calculation Lemma \ref{lem:gamma} we make use of the explicit reciprocity law from \cite{SV15}.

Contrary to the cyclotomic case, it seems difficult to establish integral results in the analytic case. On the one hand the ``dualizing character'' $\chi$ used to establish Tate duality has Frobenius action given by $\frac{\pi_L}{q}$ and hence does not make sense integrally (unless $L=\QQ_p$), on the other hand the period $\Omega$ is not a unit  (unless $L=\QQ_p$). The $L$-analytic distribution algebra $D(\Gamma_L,L)$ contains the ring $\Lambda_{\mathfrak{X}_{\Gamma_L}}$ of power-bounded functions on the character variety. It is not known whether $\Lambda_{\mathfrak{X}_{\Gamma_L}} = o_L\llbracket \Gamma_L\rrbracket.$ Paradoxically, the Iwasawa algebra $o_L\llbracket \Gamma_L \rrbracket[1/p]$ is dense inside both, the $d$-dimensional $\QQ_p$-analytic distribution algebra and the $1$-dimensional $L$-analytic distribution algebra making it unclear how to descend to integral results even under the assumption  $\Lambda_{\mathfrak{X}_{\Gamma_L}} = o_L\llbracket \Gamma_L\rrbracket.$

The structure of the paper is as follows: In section \ref{sec:phigamma} we introduce (analytic) $(\varphi_L,\Gamma_L)$-modules.
In section \ref{sec_analcoho} we introduce and study analytic cohomology of analytic $(\varphi_L,\Gamma_L)$-modules and recall the main results of \cite{Ste1} while providing some generalisations suited to our needs. Furthermore we develop an analogue of Tate duality for analytic cohomology. In section \ref{sec:BKexponential} we develop an analogue of the Bloch-Kato (dual) exponential map for analytic cohomology. We recall classical $\varepsilon$-constants in section  \ref{sec:varepsilon-constants} and state the conjecture in Section
\ref{sec:EpsilonStatement}. Section \ref{sec:RankoneCase} is dedicated to proving the main result.
In the Appendix we adapt Nakamura's density argument to the Lubin-Tate setting.

\textbf{Acknowledgements:} We are grateful to L\'{e}o Poyeton  for discussions about analytic B-pairs and to Kentaro Nakamura for answering generously questions concerning his work.   The project has been funded by the Deutsche Forschungsgemeinschaft (DFG, German Research Foundation) under   TRR 326, {\it Geometry and Arithmetic of Uniformized Structures}, project-ID 444845124.

\section{Notation}

We denote by $\mathbb{N}$ the natural numbers including $0.$

Let $\mathbb{Q}_p \subseteq L \subset \mathbb{C}_p$ be a field of finite degree $d$ over $\mathbb{Q}_p$, $o_L$ the ring of integers of $L$, $\pi_L \in o_L$ a fixed prime element, $k_L = o_L/\pi_L o_L$ the residue field,  $q := |k_L|$ and $e$ the absolute ramification index of $L$. We always use the absolute value $|\ |$ on $\mathbb{C}_p$ which is normalized by $|\pi_L| = q^{-1}$. This choice of normalisation is consistent with \cite{Co2} and \cite{SV20}.  We normalize the reciprocity map of local class field theory such that $\pi_L$ is sent to  the geometric Frobenius.

We fix a Lubin-Tate formal $o_L$-module $LT = LT_{\pi_L}$ over $o_L$ corresponding to the prime element $\pi_L$. We always identify $LT$ with the open unit disk around zero, which gives us a global coordinate $Z$ on $LT$. The $o_L$-action then is given by formal power series $[a](Z) \in o_L[[Z]]$. For simplicity the formal group law will be denoted by $+_{LT}$.

The power series $\frac{\partial (X +_{LT} Y)}{\partial Y}_{|(X,Y) = (Z,0)}$ is a unit in $o_L[[Z]]$ and we let $g_{LT}(Z)$ denote its inverse. Then $g_{LT}(Z) dZ$ is, up to scalars, the unique invariant differential form on $LT$ (\cite{Haz} \S5.8). We also let \begin{equation}\label{f:tLT}
  \log_{LT}(Z) = Z + \ldots
\end{equation}
denote the unique formal power series in $L[[Z]]$ whose formal derivative is $g_{LT}$. This $\log_{LT}$ is the logarithm of $LT$ in the sense of \cite[\S 8.6]{Lan} and converges on the maximal ideal in ${o_{\mathbb C_p}}$ (by \S 8.6, Lemma 3 (ii) ibid.).
By $\exp_{LT}:=\log_{LT}^{-1}$ in $L[[Z]]$ we denote the inverse power series of $\log_{LT} $, i.e., satisfying $\log_{LT}\circ \exp_{LT}(Z)=\exp_{LT}\circ\log_{LT}(Z)=Z.$\footnote{$\exp_{LT}$ converges on $D:=\{z\in\mathbb{C}_p| v_{\pi_L}(z)>\frac{1}{q-1}\}$ and induces on $D$ the inverse of $\log_{LT}$ respecting the valuation, see \cite[\S 8.6, Lem.\ 4]{Lan}}

In particular, $g_{LT}dZ = d\log_{LT}$. The invariant derivation $\partial_\mathrm{inv}$ corresponding to the form $d\log_{LT}$ is determined by
\begin{equation*}
  f' dZ = df = \partial_\mathrm{inv} (f) d\log_{LT} = \partial_\mathrm{inv} (f) g_{LT} dZ
\end{equation*}
and hence is given by
\begin{equation}\label{f:inv}
  \partial_\mathrm{inv}(f) = g_{LT}^{-1} f' \ .
\end{equation}
For any $a \in o_L$ we have
\begin{equation}\label{f:dlog}
  \log_{LT} ([a](Z)) = a \cdot \log_{LT} \qquad\text{and hence}\qquad ag_{LT}(Z) = g_{LT}([a](Z))\cdot [a]'(Z)
\end{equation}
(\cite{Lan} 8.6 Lemma 2).

Let $\TLT$ be the Tate module of $LT$. Then $\TLT$ is
a free $o_L$-module of rank one and we choose a generator $u=(u_n)_{n\in\mathbb{N}}$ where $u_0=0$, $u_1\neq 0$ and, for all $n,$ we have $u_n\in \mathfrak{m}_{\mathbb{C}_p}$ as well as $[\pi_L](u_{n+1})=u_n.$ Then the action of
$G_L := \Gal(\overline{L}/L)$ on $\TLT$ is given by a continuous character $\chi_{LT} :
 G_L \longrightarrow o_L^\times$.   Let $\TLT'$ denote the Tate module of the $p$-divisible group Cartier dual to $LT$ with period $\Omega_{t_0'}\in\widehat{L^{ab}}$, which again is a free $o_L$-module of rank one and where $t_0'$ is a generator. The Galois action on $\TLT'\cong\TLT^*(1)$ is given by the continuous character $\tau := \chi_{cyc}\cdot\chi_{LT}^{-1}$, where
$\chi_{cyc}$ is the cyclotomic character.
As mentioned in \cite[\S 1]{boxall} and \cite[\S3]{ST2} it follows from the work of Tate on $p$-divisible groups that we have natural $o_L$-linear isomorphisms
\begin{align}\label{f:Tate}
\TLT'\cong \Hom_{o_{\mathbb{C}_p}}(LT,\hat{\mathbb {G}}_m)\cong \Hom_{\mathbb{Z}_p}(T_\pi,\mathbb{Z}_p(1))\cong \Hom_{\mathbb{Z}_p, cts}(T_\pi\otimes_{o_L}L/o_L,\mu(p)),
\end{align}
where the last isomorphism   is induced by Pontrjagin duality and the adjunction between $Hom$ and $\otimes.$ According to the proof of \cite[Lem.\ 13]{boxall} the above composite sends $at_0'$ to  the map sending $u\otimes \frac{1}{\pi_L^n}$ to $\eta_{t_0'}(a,u_n)$, where, for $x\in o_L$, we define $\eta_{t_0'}(x,Z):=\exp\left(\Omega_{t_0'} x\log_{LT}(Z)\right)\in 1+Zo_{\widehat{L_\infty}}[[Z]]$; when the choice of $t_0'$ is clear from the context, we often omit this index from $\Omega_{t_0'}$ or $\eta_{t_0'}(x,Z)$.

Our constructions will depend crucially on the choices of $u$ \emph{and} $t_0'$, which determines the period $\Omega= \Omega_{t_0'}.$  By \eqref{f:Tate} these two choices  automatically determine a system \[\eta(1,T)_{\mid T=u_n} = \exp(\Omega \log_{LT}(T))_{\mid T=u_n}\]  of compatible $p$-power roots of unity.\footnote{E.g.\ if $L = \QQ_p$ and $LT$ is the special group corresponding to $pX+X^p$ and $\Omega=1,$ then $\eta(1,T)$ is the Artin-Hasse exponential $\exp(X + X^p/p+\dots)$.} In the cyclotomic case where $LT = \mathbb{G}_m$ it suffices to fix a choice of compatible $p$-power roots of unity because one can then take the identity as a canonical generator $t_0'$ of $\TLT'=\Hom(\mathbb{G}_m,\mathbb{G}_m)$.

For $n \geq 0$ we let $L_n/L$ denote the extension (in $\mathbb{C}_p$) generated by the $\pi_L^n$-torsion points of $LT$, and we put $L_\infty := \bigcup_n L_n$. The extension $L_\infty/L$ is Galois. We let $\Gamma_L := \Gal(L_\infty/L)$ and $H_L := \Gal(\overline{L}/L_\infty)$. The Lubin-Tate character $\chi_{LT}$ induces an isomorphism $\Gamma_L \xrightarrow{\cong} o_L^\times$. Note that by \cite[Rem.\ 1.17]{BSX} we have $N_{L/\qp}\circ \chi_{LT}=\chi_{cyc}$ if and only if $N_{L/\qp}(\pi_L)\in p^\z.$

Note that we have homomorphisms $o_L\to  1+Zo_{\widehat{L_\infty}}[[Z]],\; x\mapsto \eta(x,Z),$ and $LT\to \hat{\mathbb{G}}_m,\; Z \mapsto \eta(x,Z)$, respectively. For a  $\pi_L^n$-torsion point $a$ (whence $p^m$-torsion with $m=\lceil\frac{n}{e}\rceil$ being the smallest integer greater or equal to $\frac{n}{e}$) we thus obtain a character $o_L\to \zp[\zeta_{p^m}]^\times, x\mapsto \eta(x,a),$ of finite order. In particular $\eta(x,u_n)$ belongs to $\mu_{p^m}$ for any $x\in o_L.$ If $\gamma\in\Gamma_L$, we have $\gamma \eta(x,Z)=\eta(\chi_{LT}(\gamma)x,Z), $ while $\varphi(\eta(x,Z))= \eta(\pi_Lx,Z).$
\begin{remark}\label{rem:eta}
Since for $\sigma$   in $G_L, $  one has $\sigma(\Omega ) =\Omega\tau(\sigma)$ by \cite[Lem.\ 4.1.24]{SV20}, it follows that  ${^\sigma \eta(x,Z)}=\eta(x\tau(\sigma),Z)=\eta(x,[\tau(\sigma)](Z))$, if we let act $G_L$ on the coefficients only, and $\sigma ( \eta(x,Z))=\eta(x\tau(\sigma),[\chi_{LT}(\sigma)](Z))=\eta(x,[\chi_{cyc}(\sigma)](Z))$, if we let act $G_L$ on the coefficients and on the variable. In particular, $\sigma(\eta(x,u_n))=\eta(x\tau(\sigma),[\chi_{LT}(\sigma)](u_n))=\eta(x,[\chi_{cyc}(\sigma)](u_n))=\eta(x\chi_{cyc}(\sigma),u_n)=\eta(x,u_n)^{\chi_{cyc}(\sigma)}.$ Moreover, for a fixed choice $\zeta_{p^n}$ of  a primitive $p^n$th root of unity, there is a unique homomorphism $\beta_{u_n}:o_L\to \mathbb{Z}/p^n\z$ such that the following diagram is commutative
\[\xymatrix{
                & \mathbb{Z}/p^n\z \ar[dr]^{ \zeta_{p^n}^-}             \\
o_L \ar[ur]^{\beta_{u_n}} \ar[rr]^{\eta(-,u_n)} & &     \mu_{p^n},        }\] i.e., $\eta(x,u_n)=\zeta_{p^n}^{\beta_{u_n}(x)}.$ One easily checks that $\beta_{u_n}(\chi_{cyc}(\sigma)x)=\chi_{cyc}(\sigma)\cdot\beta_{u_n}(x).$
\end{remark}

Henceforth   we use the same notation as in \cite{SV15}.
In particular, the ring endomorphisms induced by sending $Z$ to $[\pi_L](Z)$ are called $\varphi_L$ where applicable; e.g.\  for the ring $\mathscr{A}_L$ defined to be the $\pi_L$-adic completion of $o_L[[Z]][Z^{-1}]$, or $\mathscr{B}_L := \mathscr{A}_L[\pi_L^{-1}]$ which denotes the field of fractions of $\mathscr{A}_L$.  Recall that we also have introduced the unique additive endomorphism $\psi_L$ of $\mathscr{B}_L$ (and then  $\mathscr{A}_L$) which satisfies
\begin{equation*}
  \varphi_L \circ \psi_L = \pi_L^{-1} \cdot \mathrm{Tr}_{\mathscr{B}_L/\varphi_L(\mathscr{B}_L)} \ .
\end{equation*}
Moreover, the projection formula
\begin{equation*}
  \psi_L(\varphi_L(f_1)f_2) = f_1 \psi_L(f_2) \qquad\text{for any $f_i \in \mathscr{B}_L$}
\end{equation*}
as well as the formula
\begin{equation*}
  \psi_L \circ \varphi_L = \frac{q}{\pi_L} \cdot \id \
\end{equation*}
hold. An  \'{e}tale $(\varphi_L,\Gamma_L)$-module $M$ comes with a Frobenius operator $\varphi_M$ and an induced operator  denoted by $\psi_M$. \\
For a perfectoid field extension $F$ of $L$  in the sense of \cite[Section 1.4]{GAL} let $o_{F^\flat} := \varprojlim o_{F}/p o_{F}$ with the transition maps being given by the Frobenius $\varphi(a) = a^p$. We may also identify $o_{F^\flat}$ with
$\varprojlim o_{F}/\pi_L o_{F}$ with
the transition maps being given by the $q$-Frobenius
$\varphi_q (a) = a^q$. We recall that $\widehat{L_\infty}$ and $\CC_p$ are perfectoid and that $o_{\CC_p^\flat}$ is a complete valuation ring with residue field $\overline{\mathbb{F}_p}$ and its field of fractions $\CC_p^\flat = \varprojlim \mathbb{C}_p$ is algebraically closed of characteristic $p$ (cf. \cite[Lemma 1.4.6,Proposition 1.4.7 and Lemma 1.4.10, Proposition 1.4.12]{GAL}). Let $\mathfrak{m}_{\CC_p^\flat}$ denote the maximal ideal in $o_{\CC_p^\flat}$.
The $q$-Frobenius $\varphi_q$ first extends by functoriality to the rings of the Witt vectors $W(o_{\CC_p^\flat}) \subseteq W(\CC_p^\flat)$ and then $o_L$-linearly to $W(o_{\CC_p^\flat})_L := W(o_{\CC_p^\flat}) \otimes_{o_{L_0}} o_L \subseteq W(\CC_p^\flat)_L := W(\CC_p^\flat) \otimes_{o_{L_0}} o_L$, where $L_0$ is the maximal unramified subextension of $L$. The Galois group $G_L$ obviously acts on $\CC_p^\flat$ and $W(\CC_p^\flat)_L$ by automorphisms commuting with $\varphi_q$.  This $G_L$-action is continuous for the weak topology on $W(\CC_p^\flat)_L$ (cf.\ \cite{GAL} Lemma 1.5.3).

Sometimes we omit the index $q, L,$ or $M$ from the Frobenius operator, but we always write $\varphi_p$ when dealing with the $p$-Frobenius.

Evaluation of the global coordinate $Z$ of $LT$ at $\pi_L$-power torsion points induces a map (not a homomorphism of abelian groups) $\iota: T_\pi \longrightarrow o_{\CC_p^\flat}$. Namely, if $t=  (z_n)_{n\geq 1} \in T_\pi$ with $[\pi_L](z_{n+1}) = z_n$ and $[\pi_L](z_1) = 0$, then $z_{n+1}^q \equiv z_n \bmod \pi_L$ and hence $\iota(t) := (z_n \bmod \pi_L)_n \in o_{\CC_p^\flat}$. As before we fix an $o_L$-generator $u$ of $T_\pi$ and put $\omega_u:=\iota(u).$ Then there exists a (unique) lift $Z_u \in W(o_{\CC_p^\flat})_L$  of $\omega_u$
satisfying (cf.\  \cite[Lem.\ 4.1]{SV15})
\begin{enumerate}
  \item if $u'=au$ with $a \in o_L^\times$ denotes another generator of $T_\pi$, then $ Z_{u'}=[a](Z_u ) $ is the corresponding lift;
  \item $\phi_q(Z_u) =  [\pi_L](Z_u)$;
  \item $\sigma(Z_u)=[\chi_{LT}(\sigma)](Z_u)$ for any $\sigma \in G_L$.
\end{enumerate}
By sending $Z$ to $Z_u\in W(o_{\CC_p^\flat})_L$  we obtain an $G_L$-equivariant, Frobenius compatible  embedding of rings
\begin{equation}\label{f:Zu}
     o_L\llbracket Z\rrbracket \longrightarrow W(o_{\CC_p^\flat})_L \ .
\end{equation}

Let $K\subseteq\mathbb{C}_p$ be a complete subfield containing $L_\infty$ and $\Omega,$ i.e., the minimal choice is the completion of the extension $L_\infty(\Omega)$ of $L_\infty$; by an observation of Colmez the completion $\widehat{L^{ab}}$ would be a possible choice, where we write $L^{ab}=L^{nr}L_\infty$ and $L^{nr}$ for the maximal abelian and for the maximal unramified extension of $L$, respectively. If $L\neq \qp$, such $K$ cannot be discretely valued even if we replace $L_\infty$ by $L,$ see \cite[Lem.\  3.9]{ST2}.
Following Colmez we define $K_n:=L_n\otimes_LK=\prod_{(o_L/\pi_L^n)^\times}K,$ where the latter identification is given by mapping $l\otimes_L k$ to $(\sigma_a(l)\cdot k)_{a\in  (o_L/\pi_L^n)^\times }$, and have the maps
\[\mathrm{Tr}_{K_n/K}:\prod_{(o_L/\pi_L^n)^\times}K\to K,\; (l_a)_{a\in  (o_L/\pi_L^n)^\times }\mapsto \sum_{a\in  (o_L/\pi_L^n)^\times } l_a.\]
Note that  we have $v_p(\Omega)=\frac{1}{p-1}-\frac{1}{e(q-1)}$ and, for $n\geq 1,$ $r_n:=v_p(u_n)=\frac{1}{e(q-1)q^{n-1}}$.\\

For any ring $R$, let $\mathbf{D}^{[a,b]}_{\mathrm{perf}}(R)$ (respectively $\mathbf{D}^b_{\mathrm{perf}}(R)$, $\mathbf{D}^-_{\mathrm{perf}}(R)$) denote the   triangulated   subcategory of the derived category $\mathbf{D}(R)$ of (cochain) complexes of $R$-modules consisting of the complexes of $R$-modules which are quasi-isomorphic to complexes of finitely generated projective $R$-modules concentrated in degrees $[a,b]$ (respectively bounded degrees, degrees bounded above).

Furthermore, if $R$ is a commutative ring, $X$ an $R$-module and $t\in R$
 a non-zerodivisor,
 we write $X_{t}
 :=X[\frac{1}{t}]$ for the localisation at the multiplicatively  closed set $\{1,t,t^2,\ldots\}$.

For a locally $L$-analytic group $G$ and a complete field $F\subseteq \mathbb{C}_p$ containing $L$ we write  $D(G,F)$ for the locally $L$-analytic distribution algebra with coefficients in $F$; if the coefficients are clear from the context we often abbreviate this as $D(G)$.  Dirac distributions associated with group elements $g\in G$ are denoted by $\delta_g$ or $[g]$.

\section{\texorpdfstring{$(\varphi_L,\Gamma_L)$}{(phi,Gamma)}-Modules over the Robba ring}
\label{sec:phigamma}

 For the entire section, fix a complete intermediate field $F$ of the extension
$\mathbb{C}_p/L.$

\subsection{Definition of the Robba ring $\cR$}\label{subsec:Robba}
For any
interval $I\subseteq (0,\infty)$ that is either compact or of the form
$(0,r]$, $r>0$, we define
\begin{align*}
 \cR^{I}_F:=\left\{\sum_{k\in \ZZ}a_k\cdot Z^k\mid a_k\in F,\
 \lim_{\betr k\to\infty}v_{F}(a_k)+kt=\infty\text{ for all }t\in I \right\}.
\end{align*}
We always assume that the boundary points of $I$ are in the value group of
$v_F$, so that $\cR_F^I$ is the ring of rigid analytic functions on the
annulus
\begin{align*}
 \{x\in F\mid v_F(x)\in I\}.
\end{align*} 
Furthermore, for $r>0$, let
\begin{align*}
 \cR_F^r:=\cR_F^{(0,r]}.
\end{align*}
i.e. the ring of rigid analytic functions on the annulus with outer radius $1$
and inner radius depending on $r$.
For any $s\in (0,r]$, one has $\cR_F^{[s,r]}\subseteq\cR_F^r$, and
$\cR_F^{[s,r]}$ is a Banach algebra over $F$ with the norm
\begin{align*}
 V_{[s,r]}\left(f \right)=\min_{t\in [s,r]}\left(\inf_{k\in\ZZ}
 (v_F(a_k)+kt)\right),\ \ \
 \text{ where }f=\sum_{k\in\ZZ}a_k Z^k\in \cR_F^{[s,r]}.
\end{align*}
Thus $\cR_F^r=\bigcap_{0<s\leq r}\cR_F^{[s,r]}$ is a Fr\'{e}chet space.
There are natural inclusions $\cR_F^r\subseteq \cR_F^s$ for $s\leq r$.
Now the \textbf{Robba ring} over
$F$ in the variable $Z$ is defined by
\begin{align*}
 \cR_F:=\bigcup_{r>0}\cR_F^r.
\end{align*}
We endow $\cR_F$ with the locally convex direct limit topology of the
$\cR_F^r$, making it an LF-space.\\
Moreover, let
\begin{align*}
 \cR_F^+:=\cR_F\cap F[[ Z]].
\end{align*}
This is the ring of power series with coefficients in $F$ that are
convergent on the open unit disk. In particular,
we have $\cR_F^+\subseteq\cR_F^r$ for all $r>0$.
For a complete field extension $F\subset F' \subset \CC_p$ we have $$F' \hat{\otimes}_{F,i}\cR_F \cong \cR_{F'}$$ (see. \cite[Corollary 2.1.8]{BSX}). Their proof also shows $F' \hat{\otimes}_{F,\pi}\cR_F^{r}\cong \cR_{F'}^r.$
 \\
 Inside $\cR_F$, we have the subring $\cR_F^b$
 of \textbf{bounded elements}, i.e., those Laurent series
 $f=\sum_{k\in\ZZ}a_k Z^k$ where the coefficients
 $a_k$ are bounded in $F$.
 It is well-known that  $\cR_F^\times=(\cR_F^b)^\times$.
 Furthermore, the map $f\mapsto \lVert f\rVert_1:=\sup_k\lvert a_k\rvert$ defines a multiplicative norm on
 $\cR_F^b$, see \cite[\S 1.3]{BSX}.

\subsection{Frobenius and $\Gamma_L$-action on $\cR$}
On $\cR_F$, we define a Frobenius $\varphi_L$ and a commuting
$\Gamma_L$-action by
\begin{align*}
 \varphi_L(Z):=[\pi_L](Z)\ \ \ \text{ and }\ \ \ \gamma(Z):=[\chi_{LT}(
 \gamma)](Z)\ \text{for }\gamma\in \Gamma_L
\end{align*}
on the variable
and trivial actions on the coefficients. For $r>0$, the Frobenius
$\varphi_L$ and each $\gamma\in\Gamma_L$ restrict to maps
\begin{align*}
 \varphi\colon\cR_F^r\nach \cR_F^{r/q}\ \ \ \text{ and }\ \ \ \gamma\colon\cR_F^r
 \iso\cR_F^r.
\end{align*}
For $r$ small enough, there is a left inverse
$$\psi_L\colon\cR_F^{r/q}\nach\cR_F^r$$ of $\varphi_L$, given by
$\Psi=\frac{\varphi_L^{-1}}{q}\circ \Tr_{\cR_F^{r/q}/\varphi_L(\cR_F^r)}$,
see \cite[§2]{FX}. We have $\Psi= \frac{\pi_L}{q}\psi_L.$

\subsection{\texorpdfstring{$(\varphi_L,\Gamma_L)$}{(phi,Gamma)}-Modules}
\begin{definition}
 A $\varphi_L$-\textbf{module} over $\cR_F$ is a finitely generated
 free $\cR_F$-module $M$, equipped with a continuous\footnote{where
 $M$ is, of course, endowed with the product topology from $\cR_F$.},
 $\varphi_L$-semilinear endomorphism $\varphi_M$,
 such that the induced $\cR_F$-linear map
 \begin{align*}
  \cR_F\otimes_{\cR_F,\varphi_L}M\nach M,\ f\otimes x\auf f\cdot
  \varphi_M(x)
 \end{align*}
 is an isomorphism. Note that in
 the above tensor product, $\cR_F$ is viewed as
 a left-module over itself in the usual way and as a right module via $
 \varphi_L$.\\
 We will often simply write $\varphi$ instead of $\varphi_M$.
\end{definition}

\begin{proposition}\label{phiGammaModel}
 Let $M$ be a $\varphi_L$-module over $\cR_F$. Then there exists an
 $r(M)>0$ such that, for each $0<r\leq r(M)$, there exists a unique
 finitely generated free
 $\cR_F^r$-submodule $M^r\subseteq M$ satisfying the following
 properties:
 \begin{enumerate}
  \item $M=\cR_F\otimes_{\cR_F^r}M^r$.
  \item $\varphi_M$ induces an isomorphism
   $\cR_F^{r/q}\otimes_{\cR_F^r,\varphi_L}M^r\iso
   \cR_F^{r/q}\otimes_{\cR_F^r}M^r$.
  \end{enumerate}
  In particular, for $0<s\leq r\leq r(M)$, one has
  \begin{align*}
   M^s=\cR_F^s\otimes_{\cR_F^r}M^r.
  \end{align*}
\end{proposition}
\begin{proof}
 See Thm. I.3.3 in \cite{Be04}.
\end{proof}

\begin{remark}
 Let $M$ be a $\varphi_L$-module over $\cR_F$. Then for $0<s\leq r
 \leq r(M)$ and $?\in \{s,[s,r],\varnothing\}$ we write
 \begin{align*}
  M^?:=\cR_F^?\otimes_{\cR_F^{r(M)}}M^{r(M)}.
 \end{align*}
 Composing the canonical map
 $M^r\nach \cR_F^{r/q}\otimes_{\cR_F^r,\varphi_L}M^r,\
   m\auf 1\otimes m$
 with the isomorphism $\cR_F^{r/q}\otimes_{\cR_F^r,\varphi_L}M^r\cong
  M^{r/q}$ from Prop. \ref{phiGammaModel}(ii) above, we obtain
  $\varphi_L$-semilinear maps
 $$\varphi\colon M^r\nach M^{r/q}.$$
 There is also an operator
 \begin{align*}
  \Psi_M\colon M^{r/q}\cong \cR_F^{r/q}\otimes_{\cR_F^r,\varphi_L}M^r
  \nach M^r
 \end{align*}
 given by $f\otimes m\auf \Psi(f)\cdot m$.
\end{remark}

\begin{definition}
 A $(\varphi_L,\Gamma_L)$-\textbf{module} over $\cR_F$ is a
 $\varphi_L$-module $M$ over $\cR_F$
 which carries a continuous, semilinear action
 of $\Gamma_L$ that commutes with $\varphi_M$. We shall write $\phigam{M}(\cR_F)$ for the category of $(\varphi_L,\Gamma_L)$-modules over $\cR_F.$
\end{definition}

\begin{remark}
 If $M$ is a $(\varphi_L,\Gamma_L)$-module over $\cR_F$ and $0<r\leq
 r(M)$, then from the uniqueness in Prop. \ref{phiGammaModel} it follows
 that $\gamma(M^r)=M^r$ for all $\gamma\in\Gamma_L$.
\end{remark}

\begin{definition}
\begin{enumerate}
 \item For $n\geq 1$ we put $r_n:=v_p(u_n)=\frac{1}{e(q-1)q^{n-1}}$.
 \item Let $M$ be a $(\varphi_L,\Gamma_L)$-module over $\cR_F$.
 For any $n$ such that $r_n\leq r(M)$, define
 \begin{align*}
  M^{(n)}:=M^{r_n}.
 \end{align*}
 Observe that for the Frobenius we then have $\varphi\colon M^{(n)}
 \nach M^{(n+1)}$ for $n\gg 0$.
\end{enumerate}
\end{definition}

Let $M$ be a $(\varphi_L,\Gamma_L)$-module over
  $\cR_F$. After fixing a basis of $M$, consider
  the matrix
  $P\in\GL_n(\cR_F)$ representing
   $\varphi_M$. Then
   we have $\det(P)\in \cR_F^\times=(\cR_F^b)^\times
   $ and may thus take the norm $\lVert \det(P)
   \rVert_1$ introduced at the end of Subsection
   \ref{subsec:Robba}.
   Define the \textbf{degree} $\deg(M)$ of $M$ as
   the number satisfying $\lVert\det(P)\rVert_1=q^{-\deg(M)}$; one checks that this is
   independent of the initial choice of basis of
   $M$, see \cite[\S 3.3]{BSX} for details.
   Furthermore, the \textbf{slope} of $M$ is
   defined as $\mu(M):=\deg(M)/\rk(M)$.

 \begin{definition}
  A $(\varphi_L,\Gamma_L)$-module $M$ over $\mathcal{R}_F$ is called \textbf{\'{e}tale},
if it has degree $0$ and every $(\varphi_L,\Gamma_L)$-submodule has slope $\geq0$.
 \end{definition}

	\begin{definition}
		For an affinoid algebra $A$ over $F$ we define $\cR_A^I:= A \hat{\otimes}_F\cR_F^I$ (with the projective tensor product topology) and similarly $\cR_A^r$ and $\cR_A.$ We can extend $A$-linearly the actions of $\varphi_L$ and $\Gamma_L.$ By a $(\varphi_L,\Gamma_L)$-\textbf{module} over $\cR_A$ we mean a $\cR_A$-module $M$ which arises as a base change of a projective $\cR_A^r$-module $M^r$ for some $r \gg0,$ together with a continuous $\cR_A^r$-semilinear action of $\Gamma_L$ on $M^r$ and a $\varphi_L$-semilinear map $\varphi_M \colon M^r \to M^{r/q},$ which commutes with $\Gamma_L.$ We can analogously extend the definition of $\Psi.$
\end{definition}
If $F$ is not spherically complete, we do not know if there exist non-free, projective $(\varphi_L,\Gamma_L)$-modules over $\cR_F.$ In all cases considered by us, we will only need free modules.
According to \cite{BSX} Prop.\ 2.25 the $\Gamma_L$-action on a $(\varphi_L,\Gamma_L)$-module $M$ is differentiable so that the derived action of the Lie algebra $\Lie(o_L^\times)$ on $M$ is available.

\begin{definition}\label{def:L-analytic}
  A $(\varphi_L,\Gamma_L)$-module $M$ over $\cR \in \{\cR_F,\cR_A\}$ is called $L$-\textbf{analytic}, if the derived action $\Lie(\Gamma_L) \times M \rightarrow M$ is $L$-bilinear, i.e., if the induced action $\Lie(\Gamma_L)\to \End(M)$ of the Lie algebra $\Lie(\Gamma_L)$ of $\Gamma_L$ is $L$-linear (and not just $\Qp$-linear). We shall write $\mathfrak{M}^{an}(\cR)$ for the category of $L$-analytic $(\varphi_L,\Gamma_L)$-modules over $\cR.$
  \\
  In the case $\cR=\cR_F$,
   we write $\mathfrak{M}^{an,\acute{e}t}(\cR_F)$
for the category of \'{e}tale, $L$-analytic $(\varphi_L,\Gamma_L)$-modules over
$\mathcal{R}_F.$
\end{definition}

For the relation with $L$-analytic continuous Galois representations $Rep_L^{an}(G_L)$ of $G_L$ on finite dimensional vector spaces $V$, which are {\em analytic}, i.e., satisfying that, if $D_{dR}^{\Qp}(V) := (V \otimes_{\mathbb{Q}_p} B_{dR})^{G_L}$, the filtration on $D_{dR}^{\Qp}(V)_\mathfrak{m}$ is trivial for each maximal ideal $\mathfrak{m}$ of $L \otimes_{\mathbb{Q}_p} L$ which does not correspond to the identity $\id : L \to L$,    Berger's theorem is crucial.

\begin{theorem}\label{thm:BergerEquiv}{
  There is an equivalences of categories
  \begin{align*}
    Rep_L^{an}(G_L) & \longleftrightarrow \mathfrak{M}^{an,\acute{e}t}(\cR_L) \\
    V & \mapsto D_{\mathrm{rig}}^\dagger(V).
  \end{align*}}
\end{theorem}

\begin{proof}{
Thm.\  D in \cite{Be16}  }
\end{proof}

The embedding  $o_L\llbracket Z\rrbracket  \to W(\CC_p^\flat)_L$ in \eqref{f:Zu} depends by construction on the choice of $u$. Any other choice does not change the image of the embedding  $o_L\llbracket Z\rrbracket  \to W(\CC_p^\flat)_L$ because $Z_{au} = [a](Z_u)$ for $a \in o_L^\times$ by property (i) above \eqref{f:Zu}.  As   explained in \cite[\S 8]{SVcomp}   the image $Z_u$ of the variable $Z$ already lies in $W({\hat{L}_\infty}^\flat)_L$, so that we actually have an embedding $o_L[[Z]]\to  W({\hat{L}_\infty}^\flat)_L$.
Similarly as in \cite[Def.\ 4.3.1]{KLI} for the cyclotomic situation one shows that the latter  embedding   extends to a $\Gamma_L$- and $\varphi_L$-equivariant topological monomorphism  $\mathcal{R}_L\to \tilde{\mathcal{R}}_L$ into the perfect Robba ring, see \cite[\S 5]{SVcomp} for a definition and \cite[Konstruktion 1.3.27]{W} for a proof in the Lubin-Tate setting.
\begin{remark}\label{rem:hiddendependence}In order to trace the choice of $u$ in our constructions, we should view $\cR_L$ as a subring of $\widetilde{\cR}_L$  via the embedding induced by $Z \mapsto Z_u$ and define $(\varphi_L,\Gamma_L)$-modules over this (isomorphic) subring. We will ignore this dependence for the most part by working with a fixed $Z = Z_u.$ This ``hidden'' dependence on $u$ is only relevant if an element of a $(\varphi_L,\Gamma_L)$-module is explicitly defined in terms of power series in the variable $Z$, see e.g.\  \eqref{f:Mellinsigmaminus1}, \eqref{eq:dependenceonu}, \eqref{f:phif}.
\end{remark}

\subsection{Rank one modules and characters}
\label{sec:rankonecharacters} Let $A$ be an affinoid algebra over $F.$
To each continuous character $\delta \colon L^\times \to A^\times$ we can attach a $(\varphi_L,\Gamma_L)$-module of rank one $\cR_A(\delta) := \cR_A\mathbf{e}_\delta$ by setting $\varphi_L(\mathbf{e}_\delta) = \delta(\pi_L)\mathbf{e}_\delta$ and $\gamma(\mathbf{e}_\delta) = \delta(\chi_{LT}(\gamma))\mathbf{e}_\delta$ for $\gamma \in \Gamma_L.$
We say a module is \textbf{of character type} if it arises in this way. A $(\varphi_L,\Gamma_L)$-module of character type is $L$-analytic (in the sense of Definition \ref{def:L-analytic}) if and only if $\delta$ is locally $L$-analytic (or equivalently $\delta_{|o_L^\times}$ is locally $L$-analytic). Over $\cR_L$ any rank one module is of character type (cf. \cite[Proposition 1.9]{FX}). We write $\Sigma=\Sigma(A)$ for the set of continuous characters $\delta:L^\times\to A^\times$. We denote by $\Sigma_{an}:=\Sigma_{an}(A)$  the set of locally $L$-analytic characters $\delta \colon L ^\times \to A^\times$.
Consider the following characters $\delta_{LT},\chi=x|x|, \delta^{un}_c:L^\times \to L^\times$  for $c\in L^\times$ given by
\begin{align*}
	\delta_{LT}(\pi_L)=1, &\; \delta_{LT|o_L^\times}=\id _{o_L^\times}, \\
	\chi(\pi_L)=\frac{\pi_L}{q} ,&\; \chi_{|o_L^\times}=\id _{o_L^\times}, \\
	\delta^{un}_c(\pi_L)=c, & \; (\delta^{un}_c)_{|o_L^\times}\equiv1.
\end{align*}
In particular, $\chi=\delta^{un}_{\frac{\pi_L}{q}}\delta_{LT}.$ Then $\delta_{LT}$ corresponds via class field theory to the character $\chi_{LT}\colon G_L\to o_L^\times.$
 Let $\delta\colon L^\times\to L^\times$
be any continuous character;
setting $\delta_0:=\delta^{un}_{\delta^{-1}(\pi_L)}\delta$ we may always decompose $\delta=  \delta^{un}_{\delta(\pi_L)}\delta_0$  satisfying $\delta_{|o_L^\times} = (\delta_0)_{|o_L^\times}$ and $\delta_0(\pi_L)=1.$
 If   $|\delta(\pi_L)|=1$, the character $\delta$ corresponds  to a Galois character $\chi_\delta$ via local class field theory. Then $D^\dagger_{rig}(L(\chi_\delta))=\cR_L(\delta)$ and we call $\delta$  \'{e}tale.

Later, for descent calculations we will have to select out the sets of \textit{special} characters $\Sigma_1:=\{x^{-i}| i\in\mathbb{N}\}$ and $\Sigma_2:=\{x^{i}\chi| i\in\mathbb{N}\}$ from the \textit{generic} ones $\Sigma_{gen}:=\Sigma_{an}\setminus (\Sigma_1\cup\Sigma_2).$

Note that we have two $\psi$-operators. While $\psi$ satisfies the identity $\psi\circ\varphi=\frac{q}{\pi_L}\id$ and makes sense even integrally, $\Psi$ denotes the left inverse of $\varphi,$ i.e., satisfying $\Psi\circ\varphi= \id$. In particular, we have $\psi=\frac{q}{\pi_L}\Psi.$ Note that $\psi(\mathbf{e}_\delta)=\frac{q}{\pi_L}\delta^{-1}(\pi_L)\mathbf{e}_\delta.$

If $\delta \in \Sigma_{an}(K)$ and $a \in o_L^\times$ such that $\log(a) \neq 0,$ then one defines the \textbf{weight} of $\delta$ as $\omega_\delta:= \log(\delta(a))/\log(a)$ (which is independent of $a$).
   We shall say that $\delta$ is \textbf{de Rham}, if the attached $(\varphi_L,\Gamma_L)$-module $R_K(\delta)$ is de Rham in the sense that will be introduced in subsection \ref{sec:dualexp} below. As shown in the Appendix \ref{App:density}, Remark \ref{rem:deRhamchar}, $\delta$ is de Rham if and only if there exist  some locally constant character  $\delta_{\mathrm{lc}}$ and $k(=\omega_\delta)\in \mathbb{Z,}$  such that
\[\delta=\delta_{\mathrm{lc}}x^k  \mbox{ (or equivalently }\delta=\delta_{\mathrm{lc}}\delta_{LT}^k \mbox{ for some other $\delta_{\mathrm{lc}}$)},   \]
  see also \cite[Rem.\ 3.2.3/4]{SV20} for the \'{e}tale case.

We fix some notation for the remainder of the article. Consider the differential operator $\partial:= \partial_{\mathrm{inv}}= \frac{1}{\log_{LT}'(Z)}\frac{d}{dZ}$ acting on $\cR_K.$ (This differs from \cite{Co2} by a constant.) Let  $\nabla\in Lie(\Gamma_L)\cong L$ be the element corresponding to  $1 \in L.$
\begin{remark}\label{rem:FColmezFormeln}
	We obtain the following properties (cf. \cite[1.2.4]{Co2}):
	\begin{enumerate}
		\item $\partial \circ \varphi = \pi_L \varphi \circ \partial.$
		\item $\partial \circ \gamma  = \chi_{LT}(\gamma) \gamma \circ \partial.$
		\item $\nabla f = t_{LT}\partial f$ for $f\in \cR_K.$
		\item $\nabla(f\mathbf{e}_\delta) = (\nabla f+\omega_{\delta}f)\mathbf{e}_\delta$ for $\delta  \in \Sigma_{an}.$
		\item $\partial \eta(x,T) = \Omega x \eta(x,T)$
	\end{enumerate}
\end{remark}

\subsection{The modules $\mathbf D_{\mathrm{dif}}^{(+)}(M)$}

We set $t_{LT}=\log_{LT}(Z) \in L[[Z]]$, so that
\begin{align*}
 \varphi(t_{LT})=\pi_L\cdot t_{LT}\ \ \ \text{ and }\ \ \ \gamma(t_{LT})=\chi_{LT}(\gamma)\cdot t_{LT}
 \text{ for all }\gamma\in\Gamma_L
\end{align*}
by (\ref{f:dlog}). For $n\geq 1$, we set
\[[\pi_L^{-n}](Z):=u_n+_{LT}\exp_{LT}(\frac{t_{LT}}{\pi_L^n}) \in L_n[[Z]]. \] Then $[\pi_L^n](u_n+_{LT}\exp_{LT}(\frac{t_{LT}}{\pi_L^n}))=Z,$ which is how Colmez justifies this notation in \cite[1.4.2]{Co2}. Note that the constant term of $[\pi_L^{-n}](Z)$ is equal to $u_n$ and hence is non-zero, so $[\pi_L^{-n}](Z)$ is a unit in $L_n[[Z]]$.

Furthermore, let $\theta:K_n[[t_{LT}]]\to K_n$ denote the $K_n$-linear map sending $t_{LT}$ to $0,$ i.e., the reduction modulo $t_{LT}.$ This is the completed base change to $K$ of the restriction of $\theta:\mathbf{B}_{dR}^+\to \mathbb{C}_p$ to $L_n[[t_{LT}]]\to L_n.$

In the following, let
$F$ be a complete nonarchimedean field containing $L$.

\begin{definition} \label{def:iota} The group $\Gamma_L$ acts diagonally on
$$F_n:= L_n\otimes_LF$$ (trivially on the right factor and naturally on the left), and we extend this to an action on $F_n[[Z]]$ via its usual action on $Z$.\footnote{Note that $F_n[[t_{LT}]]=F_n[[Z]]$ because the map $F_n[[Z]]/Z^k \nach F_n[[Z]]/Z^k, Z \auf t_{LT}$ is an isomorphism for all $k$, a consequence of $t_{LT}$ being an element of $Z+Z^2L[[Z]]$.} Now define
\begin{align*}
  \iota_n=\iota_n^{(F)}:\cR_F^{r_n} & \to F_n[[t_{LT}]]=F_n[[Z]], \\
  \sum_{k\in\z} a_kZ^k & \mapsto \sum_{k\in\z} a_k([\pi_L^{-n}](Z))^k,
\end{align*}
where on the right-hand side $a_k$ denotes (by abuse of notation) the image under the canonical embedding $F \nachinj L_n\otimes_LF$ and $[\pi_L^{-n}](Z)$ is viewed as a power series over $L_n\otimes_LF$ via the embedding $L_n \nachinj L_n\otimes_LF$.
\end{definition}


 \begin{remark}\label{iota t}\phantom{section}
\begin{enumerate}
\item The map $\iota_n$ is well-defined.
\item For the power series $t_{LT}=\log_{LT}(Z
  )\in \cR_F^+$, we have
  \begin{align*}
   \iota_n(t_{LT})=\frac{t_{LT}}{
    \pi_L^n}.
  \end{align*}
  	\item $\iota_n$ is injective for every $n.$ 	\footnote{The injectivity of $\iota_n$ in the cyclotomic case is \cite[Proposition 2.11, Proposition 2.25]{Be02}, but the map $\iota_n$ is defined in terms of Witt vectors. The argument given by us is in a similar spirit as \it{(loc.cit.)}.}
\end{enumerate}
 \end{remark}
\begin{proof}
	To see that $\iota_n$ is well-defined for $F=L$ first recall that by \cite[Prop. 8.10]{Co4}, the ring $\mathbf{B}_{\mathrm{dR}}^+$ contains a period $t_L$ for the Lubin-Tate character, i.e. we have $g(t_L)= \chi_{LT}(g)t_L$ for all $g \in G_L$ and $t_L$ differs from the usual $t$ by a unit. Thus $L_n[[t_{LT}]]$ embeds into $\mathbf{B}_{\mathrm{dR}}^+$ via $t_{LT} \auf t_L$ and we endow it with the subspace topology, making it a closed subspace of $\mathbf{B}_{\mathrm{dR}}^+$.
	It hence suffices to show that $\iota_n(f)$ converges in $\mathbf{B}_{\mathrm{dR}}^+.$
	 A series of the form $x=\sum_{k \gg -\infty} p^k [x_k] \in W(o_{\mathbb C_p^\flat})[1/p]$ converges in $\mathbf{B}_{\mathrm{dR}}^+$ if and only if $\theta(x)$ converges in $\mathbb C_p$ (which is the case precisely when $k+v(x_k) \to \infty$ for $k \to \infty$).
	 As in \cite[Prop. III.2.1 (i)]{CC99}, the condition $x=\sum_{k \in \mathbb Z}a_k Z_u^k$ with $a_k \in o_L$ and  $v_p(a_k) + k \cdot r_n \to \infty$ for  $k \to -\infty$  implies that $\iota_n(x)$ converges in $\mathbf{B}_{\mathrm{dR}}^+$.
	Even though the coefficients of an element $x \in \cR_L^{r_n} $ are not bounded, they do satisfy the same growth condition, which  ensures the convergence of  $\iota_n(x).$
	 The case of general $F$ is obtained via completed base change $\cR_F^{r_n} = F \hat{\otimes}_{L,\pi}\cR_L^{r_n} \to F_n[[t_{LT}]] = F \hat{\otimes}_{L,\pi}L_n[[t_{LT}]].$ \\
	For the second point we compute  \begin{align*}
		\iota_n(t_{LT})=\log_{LT}([\pi_L^{-n}](Z))
		=\underbrace{\log_{LT}(u_n)}_{=0}+\log_{LT}\exp_{LT}\left(
		\frac{t_{LT}}{\pi_L^n}\right)=\frac{t_{LT}}{
			\pi_L^n}.
	\end{align*}
For the injectivity of $\iota_n$ we can assume $F=L$ because completed base change $F \hat{\otimes}_{L,\pi}-$ preserves injectivity by \cite[1.1.26]{Eme}.
Consider $\theta\colon L_n \llbracket t_{LT} \rrbracket \to L_n$ the reduction modulo $t_{LT}$ such that $\theta\circ\iota_n(f)=f(u_n).$
If $\iota_n(f)=0$ then $f(u_n)=0$ and hence $f$ is divisible by $Q_n(Z):=\frac{[\pi_L^n](Z)}{[\pi_L^{n-1}](Z)}.$ Because $\iota_n(t_{LT})\neq 0$ and $Q_n \mid t_{LT} = Z\prod_{\mu\geq 1}\frac{Q_\mu}
{\pi_L}$ we conclude $\iota_n(Q_n)\neq0$ and hence $\iota_n(f/Q_n)=0.$ Inductively $Q_n^k \mid f$ for every $k\geq 0.$ The choice of $n$ ensures that $Q_n \in \cR_L^{r_n}$ is a non-unit. By considering the image of $f$ in the noetherian domain $\cR_L^{[r_n,r_n]}$ under the inclusion $\cR_L^{r_n}\subset \cR_L^{[r_n,r_n]}$ we conclude $f=0$ by Krull's Intersection Theorem.
\end{proof}

The map $\iota_n$ commutes with the action of $\Gamma_L$.
 Writing $\mathrm{Tr}=\id_F\otimes\frac{1}{q}\mathrm{Tr}_{L_{n+1}/L_n}$ we obtain the commutative diagrams
\[
\begin{tikzcd}
\cR_F^{r_n} \arrow[d, "\varphi"] \arrow[r, "\iota_n"] & F_n[[t_{LT}]] \arrow[d, hook] \\
\cR_F^{r_{n+1}} \arrow[r, "\iota_{n+1}"]                  & F_{n+1}[[t_{LT}]]
\end{tikzcd}
\]
and
\[
\begin{tikzcd}
\cR_F^{r_n}  \arrow[r, "\iota_n"] & F_n[[t_{LT}]]  \\
\cR_F^{r_{n+1}} \arrow[r, "\iota_{n+1}"]  \arrow[u, "\psi"]                & F_{n+1}[[t_{LT}]]  \arrow[u, "\mathrm{Tr}"].
\end{tikzcd}
\]
\begin{definition}\label{Ddif.Def}
Let $M$ be a $(\varphi_L,\Gamma_L)$-module over $\cR_F$. Viewing $F_{n}[[t_{LT}]]$ as an $\cR_F^{r_n}$-module via $\iota_n$, we define the $F_{n}[[t_{LT}]]$- and $F_{n}((t_{LT}))$-modules
\begin{align*}
\bD_{\textnormal{dif},n}^+(M):= F_{n}[[t_{LT}]] \otimes_{\cR_F^{r_n}} M^{(n)} \quad \textnormal{  and  } \quad \bD_{\textnormal{dif},n}(M):=\bD_{\textnormal{dif},n}^+(M)_{t_{LT}},
\end{align*}
respectively, where $(-)_{t_{LT}}$ means localising
an $F_n[[t_{LT}]]$-module at the multiplicative
subset generated by $t_{LT}$. Furthermore,
$\bD_{\textnormal{dif},n}^+(M)$ carries the diagonal action of $\Gamma_L$, which also extends to $\bD_{\textnormal{dif},n}(M)$. Under the isomorphism $\varphi^{*}(M^{(n)}) = \cR_F^{r_{n+1}} \otimes_{\varphi,\cR^{r_n}_F} M^{(n)} \cong M^{(n+1)}$, the map $\varphi \colon M^{(n)} \to M^{(n+1)}$ corresponds to the canonical map $\operatorname{can_{n,n+1}}\colon M^{(n)}\to  \varphi^{*}(M^{(n)}); x\mapsto 1 \otimes x.$ The above diagrams then induce the diagrams (see \cite[\S 2.B]{NaANT} for details)
\[
\begin{tikzcd}
M^{(n)} \arrow[d, "\varphi"] \arrow[r, "\iota_n"] & \bD_{\textnormal{dif},n}^{(+)}(M) \arrow[d, hook, "\mathrm{can}_{n,n+1}"] \\
M^{(n+1)} \arrow[r, "\iota_{n+1}"]                  & \bD_{\textnormal{dif},n+1}^{(+)}(M)
\end{tikzcd}
\]
where the map $\mathrm{can}_{n,n+1}$ is given by $f(t)\otimes x\mapsto f(t)\otimes \varphi(x)  \corresponds f(t)\otimes \operatorname{can}_{n,n+1}(x) = f(t)\otimes 1 \otimes x$, and $\iota_n$ by $m\mapsto 1\otimes m$,
as well as, for $n\geq 1,$
\[
\begin{tikzcd}
M^{(n)} \arrow[r, "\iota_n"] & \bD_{\textnormal{dif},n}^{(+)}(M)  \\
M^{(n+1)} \arrow[r, "\iota_{n+1}"]  \arrow[u, "\psi"]                & \bD_{\textnormal{dif},n+1}^{(+)}(M)  \arrow[u, "\mathrm{Tr}"]
\end{tikzcd}
\]
 with transitions maps $f(t)\otimes x\mapsto \mathrm{Tr}(f(t))\otimes \psi(x)$ on the right hand side.
Finally, we define
\begin{align*}
\bD_{\textnormal{dif}}^{(+)}(M):= \varinjlim_{n \gg 0}\bD_{\textnormal{dif},n}^{(+)}(M)
\end{align*}
 with $\mathrm{can}_{n,n+1}$ as transition maps.
\end{definition}
As in \cite{NaANT}, we have $\bD_{\textnormal{dif},n}^{+}(M) \otimes_{F_n[[t_{LT}]]} F_{n+1}[[t_{LT}]] \iso \bD_{\textnormal{dif},n+1}^{+}(M)$ and hence
\begin{align*}
\bD_{\textnormal{dif}}^{(+)}(M) = \bD_{\textnormal{dif},n}^{(+)}(M) \otimes_{F_n[[t_{LT}]]} (\bigcup_{m \geq n} F_m[[t_{LT}]])
\end{align*}
for $n \gg 0$.

\begin{remark}
	\label{rem:DdifFrechet}
 Since $M^{(n)}$ is a free module over $\cR_F^{(n)}$, say of rank $d$, we have $\bD_{\mathrm{dif},n}^+(M)\cong F_n[[t_{LT}]]^d$.
 The Fr\'{e}chet-space-structure on $F_n[[t_{LT}]]=\varprojlim F_n[[t_{LT}]]/(t_{LT}^k)$ (with
 base field $F,$ where each factor is a finite-dimensional $F$-vector space endowed with it's canonical topology) thus induces one on $\bD_{\mathrm{dif},n}^+(M)$, which is
 of course independent of the choice of the isomorphism above.
 Furthermore, $\bD_{\mathrm{dif},n}(M)=\varinjlim_{k}\bD_{\mathrm{dif},n}^+(M)\cdot t_{LT}^{-k}$ becomes an LF-space over $F$ in this way.
 Finally, the modules $\bD_{\mathrm{dif}}^+(M)$ and $\bD_{\mathrm{dif}}(M)$ are also LF-spaces\footnote{{ Note that this topology is not the norm topology on $L_\infty  $ because a strict LF-space is complete.}} over $F$.
\end{remark}

Later on it will be crucial to form the cohomology groups $H^i_{\varphi, \mathfrak{Z}}(\bD_{\mathrm{dif}}^+(M))$ from Section \ref{sec_analcoho}. For this we need a $D(\Gamma_L,F)$-module-structure on $\bD_{\mathrm{dif}}^{(+)}(M)$, which we get from Proposition \ref{ana:to:distr} below after showing that the action is pro-$L$-analytic.
Let us first recall this notion.

\begin{definition}\label{defi:proAn}
 Let $G$ be an $L$-analytic group (the main example to have in mind is $G=\Gamma_L$).
\begin{enumerate}[(a)]
 \item Let $V$ be a Banach space over $F$ equipped with
 a continuous linear $G$-action. We say that a
 vector $v\in V$ is \textbf{locally $L$-analytic}
 if there exists an open subgroup
 $\Gamma_n\subseteq G$ together with a chart
 $\ell\colon\Gamma_n\overset\sim\to\pi_L^n o_L$ (for some $n\geq 0$)
 such that the orbit map of $v$ restricted to
 $\Gamma_n$ is given by a power series
 \begin{align*}
  \gamma(v)=\sum_{k\geq 0}\ell(\gamma)^k v_k
 \end{align*}
 where $v_k\in V$ is a sequence of vectors
 satisfying $\pi_L^{nk}v_k\to 0$ for $k\to\infty$.
 We say that the action of $G$ is \textbf{locally $L$-analytic} if all $v\in V$ are locally $L$-analytic.
 \item Let $W=\varinjlim_m W_m$ be an LF-space
 over $F$,
 with Fr\'{e}chet spaces $W_m=\varprojlim_n W_{m,n}$
 and Banach spaces $W_{m,n}$, such that $G$ acts
 linearly and continuously on $W$. We say that
 a vector $w\in W$ is \textbf{pro-$L$-analytic}
 if its orbit map $G\to W$ factors over some $W_m$
 and the induced maps $G\to W_{m,n}$ are locally
 $L$-analytic for all $n$. We denote by $W^{L-pa}$ the subset of pro-$L$-analytic vectors of $W.$ We say the action is \textbf{pro-$L$-analytic} if $W^{L-pa}=W.$
\end{enumerate}
\end{definition}

\begin{proposition}\label{ana:to:distr}
 Let $W$ be an LF-space over $F$ carrying a pro-$L$-analytic
 action of $\Gamma_L$.
 Then this action extends uniquely to a separately continuous action of $D(\Gamma_L,F)$ on $W$.
\end{proposition}
\begin{proof}
 This follows from the proof of \cite[Proposition 4.3.10]{SV20}.
\end{proof}

\begin{lemma}\label{frechmodproan}
 Let $B$ be a Fr\'{e}chet $\Gamma_L$-ring over $F$ and $W$ a finitely generated free $B$-module with a compatible $\Gamma_L$-action.
 Assume there is a basis $\mathfrak A:=(e_1,\ldots,e_d)$ for $W$ such that the map
 \begin{align*}
   \Gamma_L\nach \mathrm{GL}_d(B),\ \gamma\auf \mathrm{Mat}_{\mathfrak A}(\gamma)
 \end{align*}
 is pro-$L$-analytic.  $W^{L-\mathrm{pa}}=\bigoplus_{j=1}^dB^{L-\mathrm{pa}}\cdot e_j$.
\end{lemma}
\begin{proof}
 This is proven for $F=L=\QQ_p$ in \cite[Prop. 2.4]{Be16} and the identical proof applies for general $F$ and $L$.
\end{proof}

\begin{proposition}
 For an $L$-analytic $(\varphi_L,\Gamma_L)$-module $M$ over $\cR_F$, the $\Gamma_L$-action on the LF-spaces $\Ddifp(M)$ and $\Ddif(M)$
 is pro-$L$-analytic.
\end{proposition}
\begin{proof}
 We start with $\Ddifp(M)=\varinjlim_{n\gg 0}\bD_{\mathrm{dif},n}^+(M)$. By definition, it suffices to check that the
 $\Gamma_L$-action on the Fr\'{e}chet space $\bD_{\mathrm{dif},n}^+(M)$ is pro-$L$-analytic for $n\gg 0$.\\
 We wish to apply Lemma \ref{frechmodproan} with $B:=F_n[[t_{LT}]]$
 and $W:=\bD_{\mathrm{dif},n}^+(M)=B\otimes_{\cR_F^{(n)}}M^{(n)}$: Choose any $\cR_F^{(n)}$-module basis $x_1,\dots,x_d$ of $M^{(n)}$. Then
 $\mathfrak A:=(1\otimes x_1,\ldots,1\otimes x_d)$ is a basis of the free $B$-module $W$, and the map
 $\gamma\auf \mathrm{Mat}_{\mathfrak A}(\gamma)$ is given by the composite
 \begin{align*}
  \Gamma_L\nach \GL_d(\cR_F^{(n)})\overset{\iota_n}{\nach} \GL_d(B)
 \end{align*}
 where the first map is pro-$L$-analytic because $M^{(n)}$ is pro-$L$-analytic by assumption. Moreover, since $\iota_n$
 is a continuous homomorphism of $F$-algebras, we conclude that \ref{frechmodproan} is applicable. Thus we obtain
 \begin{align*}
  W^{L-\mathrm{pa}}=\bigoplus_{j=1}^d B^{L-\mathrm{pa}}\cdot (1\otimes x_j)=
  F_n[[t_{LT}]]^{L-\mathrm{pa}}\otimes_{\cR_F^{(n)}}M^{(n)}.
 \end{align*}
 Finally, from \cite[Prop. 2.6 2.]{Por} it follows that $F_n[[t_{LT}]]^{L-\mathrm{pa}}=F_n[[t_{LT}]]$, which completes
 the proof for $\Ddifp(M)$.\\
 Moving on to $\Ddif(M)$, we write $\Ddif(M)=\varinjlim_{n,k}\bD_{\mathrm{dif},n}^+(M)\cdot t_{LT}^{-k}$ as a direct limit of
 Fr\'{e}chet spaces. By what we have just shown,
 one can express $\bD_{\mathrm{dif},n}^+(M)$ for $n\gg 0$
 as inverse limit
 $\bD_{\mathrm{dif},n}^+(M)=\varprojlim_r B_{n,r}$
 for certain $F$-Banach spaces $B_{n,r}$ on which
 $\Gamma_L$ acts $L$-analytically. So for any $k$
 one has $$\bD_{\mathrm{dif},n}^+(M)\cdot t_{LT}^{-k}
 =\varprojlim_r B_{n,r}\cdot t_{LT}^{-k},$$ where by
 $B_{n,r}\cdot t_{LT}^{-k}$ we denote the $\Gamma_L
 $-module $B_{n,r}$ whose $\Gamma_L$-action is
 twisted by $\chi_{LT}^{-k}$. Since the inversion in $\Gamma_L$ is an
 $L$-analytic map, we see that the twisted action
 $b\auf  \chi_{LT}(\gamma^{-k})\cdot \gamma(b)$
 on $B_{n,r}$ is again $L$-analytic. Thus $\Gamma_L
 $ acts pro-$L$-analytically on $\bD_{\mathrm{dif},n}^+(M)\cdot t_{LT}^{-k}$ for
 $n\gg 0$ and $k\geq 1$, so the claim follows.
\end{proof}

Note that $\bD_{\mathrm{dif}}(M)$ depends on the coefficient field of $\cR_F.$ For a complete field extension $F'/F$ and an $L$-analytic $(\varphi_L,\Gamma_L)$-module $M$ over $\cR_F$ one checks that $F' \hat{\otimes}_{F,i}M$ is an $L$-analytic $(\varphi_L,\Gamma_L)$-module over $\cR_{F'}.$ Here $\hat{\otimes}_{F,i}$ denotes the inductive tensor product topology.
\begin{remark}
	\label{rem:DdifBasechange}
	Let $F'/F$ be a complete field extension and let $M$ be an $L$-analytic $(\varphi_L,\Gamma_L)$-module over $\cR_F.$ The natural maps $$F'\hat{\otimes}_{F,i}\bD_{\mathrm{dif}}(M) \to \bD_{\mathrm{dif}}(M\hat{\otimes}_{F,i}F')$$
	and
	$$F' \hat{\otimes}_{F,i}\bD_{\mathrm{dif},n}(M) \to \bD_{\mathrm{dif},n}(M\hat{\otimes}_{F,i}F')$$
	are $\Gamma_L$-equivariant isomorphisms.
	
\end{remark}
\begin{proof}
	The completed inductive tensor product commutes with strict locally convex inductive limits by \cite[Theorem 1.1.30]{Eme} together with the argument in the proof of \cite[2.1.7(i)]{BSX}. Hence the first statement follows from the second.
	For Fr\'{e}chet spaces inductive and projective tensor products agree and commute with projective limits (of Hausdorff spaces) with dense transition maps (cf. \cite[2.1.4]{BSX} and \cite[17.6]{NFA}). This allows us to first reduce to the corresponding statement for $\bD_{\dif,n}^+(M)t_{LT}^{-k_0}$ since $\bD_{\dif,n}(M) = \varinjlim_k \bD_{\dif,n}^+(M)t_{LT}^{-k}$ and
by \ref{rem:DdifFrechet} we have $\bD^+_{\mathrm{dif},n}(M) \cong \varprojlim_k \bD^+_{\mathrm{dif},n}(M)/(t_{LT}^k)$, hence we even have surjective transition maps which allow us to reduce to the corresponding statement for $\bD^+_{\mathrm{dif},n}(M)/(t_{LT}^k)$ (assuming for simplicity $k_0=0$, the general case being treated analogously). Since each $\bD^+_{\mathrm{dif},n}(M)/(t_{LT}^k)$ is finite dimensional over $F$, we may omit the completion and see that $$F' \otimes_{F}\bD^+_{\mathrm{dif},n}(M)/(t_{LT}^k) \to \bD^+_{\mathrm{dif},n}(M\hat{\otimes}_{F,i}F')/(t_{LT}^k)$$ is an isomorphism of finite dimensional $F'$-vector spaces, which follows from the fact that any basis of $M^{(n)}$ gives rise on the one hand to a $F_n[[t_{LT}]]$ basis of $\bD^+_{\mathrm{dif},n}(M)$ and on the other hand to a basis of $F'\hat{\otimes}M^{(n)}$ and thus to a $F'_n[[t_{LT}]]$ basis of $\bD^+_{\mathrm{dif},n}(F'\hat{\otimes}M).$ Note that $F_n[[t_{LT}]]/(t_{LT})^k \otimes_F F' \cong F'_n[[t_{LT}]]/(t_{LT})^k$ by a dimension argument.
\end{proof}
\begin{lemma}
	\label{lem:InvariantsTensor}
Let $V$ be a $F$-Banach space and let $G$ be a group acting on $V$ via continuous $F$-linear maps. Let $W$ be an $F$-Banach space of countable type endowed with the trivial $G$-action. Then $$(V \hat{\otimes}W)^G = V^G \hat{\otimes}W$$
\end{lemma}
\begin{proof}
	Assume without loss of generality, that $W$ is infinite dimensional (the finite dimensional case being simpler).
	By \cite[Corollary 2.3.9]{PGS} $W$ is isomorphic to $c_0(F),$ the space of zero sequences in $F$ indexed by $\NN$. We obtain a $G$-equivariant isomorphism $V \hat{\otimes}W \cong c_0(V)$ by first identifying $c_0(F)$ (resp. ($c_0(V)$) ) with the completion of $\bigoplus_{n \in \NN}F$ (resp. $\bigoplus_{n \in \NN}V$) and using the $G$-equivariant isomorphism $(\bigoplus_{n \in \NN} F)\otimes_F V \cong \bigoplus_{n \in \NN} V$ and passing to completions. Note that $g \in G$ acts via continuous automorphisms with respect to the product topology and hence extends to an automorphism of the completions with $g$ acting on a sequence $(v_1,v_2,\dots)$ via $g(v_1,v_2,\dots)= (gv_1,gv_2,\dots).$ It is clear that any such sequence is $G$-invariant if and only if each component is $G$-invariant.
\end{proof}
\begin{corollary}
	Let $F'/F$ be a complete field extension contained in $\CC_p$ and let $M$ be an $L$-analytic $(\varphi_L,\Gamma_L)$-module over $\cR_F.$ We have
	$$\bD^{(+)}_{\mathrm{dif}}(M)^{\Gamma_L} \hat{\otimes}_{F,i}F'= \bD^{(+)}_{\mathrm{dif}}(M\hat{\otimes}_{F,i}F')^{\Gamma_L}.$$
\end{corollary}
\begin{proof}
	Like in the proof of \ref{rem:DdifBasechange} we reduce to the corresponding statement for the Banach spaces $\bD_{\mathrm{dif},n}(M)^+t_{LT}^{-k_0}/\bD_{\mathrm{dif},n}(M)^+t_{LT}^{k-k_0}$. The field $F'$ is of countable type over $F$ since $F'\cap \overline{\QQ_p}$ is dense in $F'$ by \cite[Theorem 1]{IovitaZaharescu} (and of at most countable dimension over $\QQ_p$) and hence also $F(F'\cap \overline{\QQ_p})$ is a dense $F$-subspace of at most countable $F$-dimension. Because the action on $F'$ is trivial, we can deduce the result from \ref{lem:InvariantsTensor}.
\end{proof}

\section{(Analytic) Cohomology groups}\label{sec_analcoho}

For the moment let $F$ be any field extension of $L$ and $G$   be any $L$-analytic group (of dimension one); we shall reserve the letter $U$ for  a (sub)group isomorphic to $o_L$. If $K$ is big enough such that $D(U):=D(U,K)\cong\cR_K^+=: \cR^+$ then we denote by $\mf Z\in D(U)$
the element corresponding to the variable $Z\in \cR^+$.  Let $V$ be any (abstract) $D(G,F)$-module.  We define cohomology groups $H^\bullet_{\clubsuit,\spadesuit}(V)$ for $\clubsuit\in\{\varphi,\psi\}$ and $\spadesuit\in \{D(G,F),\mf Z, Lie(G),\nabla\}$  as follows: By $\mathrm{RHom}_{D(G)}(F,V)$ we denote any (bounded) complex of $F$-vector spaces whose cohomology gives $\mathrm{Ext}_{D(G)}^\bullet(F,V)$ (extensions  as abstract $D(G)$-modules). Let $f$ be any   endomorphism of $V$ which commutes with the $D(G)$-action   inducing an operator on $\mathrm{RHom}_{D(G)}(F,V)$
and  we denote by \[K_{f,D(G)}(V):=\mathrm{cone}\left( \mathrm{RHom}_{D(G)}(F,V)\xrightarrow{f-\id}\mathrm{RHom}_{D(G)}(F,V) \right)[-1]\] the induced mapping fibre.
For $U\cong o_L$ and  $K$ being big enough
\[\xymatrix@C=0.5cm{
  0 \ar[r] & D(U) \ar[rr]^{\mf Z} && D(U) \ar[rr]^{} && K \ar[r] & 0 }\]
  is a projective resolution of the trivial representation $K$ and we can choose $V\xrightarrow{\mf Z}V$ (functorially) for $\mathrm{RHom}_{D(U)}(K,V)$. In this context we shall also use  the notation $K_{f,\mf Z}(V)$ for  $K_{f,D(U)}(V).$
  Note that
\[K_{f,\mf Z}(V)\cong\mathrm{cone}\left(V\otimes^\mathbb{L}_{D(G)}K \xrightarrow{f-\id}V\otimes^\mathbb{L}_{D(G)}K \right)[-2]\]
as $ \mathrm{RHom}_{D(G)}(K,V)\cong V\otimes^\mathbb{L}_{D(G)}K[-1] .$ Analogous isomorphisms exist for  $K_{f,D(G)}(V)$ for any $G$ of dimension one, since in our context taking $G/U$-invariants and -coinvariants coincide and form exact functors by Maschke's theorem.

 Following \cite{Ko} we write $D^\infty(G)$ for the algebra of locally constant distributions, i.e., the quotient of $D(G)$ by the  ideal generated by $Lie(G)\subseteq D(G).$ We then obtain  isomorphisms by \cite[p. 306]{ST3}
\begin{equation}
\label{f:Ext}\mathrm{Ext}_{D(G)}^\bullet(D^\infty(G),V)\cong H^\bullet (Lie(G),V),
\end{equation}
where the latter denotes Lie algebra cohomology. Since the reference does not cover coefficient fields such as our $K$, which is not spherically complete, we would like to briefly justify this isomorphism:  For $Lie(G)=L\nabla$ we have a strict exact sequence of Hausdorff locally convex vector spaces over $L$
\begin{equation}\label{f:resolution}
   \xymatrix@C=0.5cm{
     0 \ar[r] & D(G,L) \ar[rr]^{\nabla} && D(G,L) \ar[rr]^{pr} && D^\infty(G,L) \ar[r] & 0 }
\end{equation}
by \cite[\S 3]{ST3}, i.e., a resolution of $D^\infty(G,L)$ by free $ D(G,L)$-modules. Moreover, it arises by base change $D(G,L)\otimes_{U_L(Lie(G))}-$ from the following resolution of $L$ by free $U_L(Lie(G))$-modules, where the latter denotes the enveloping algebra of $Lie(G):$
\begin{equation}\label{f:resolutionL}
\xymatrix@C=0.5cm{
     0 \ar[r] & U_L(Lie(G)) \ar[rr]^{\nabla} && U_L(Lie(G)) \ar[rr]^{pr} && L \ar[r] & 0 }
\end{equation}
see \cite[Rem.\ 1.1]{ST3}. Base change $K\hat{\otimes}_L-$ of \eqref{f:resolution}  leads to the strict exact sequence of Hausdorff locally convex $K$-vector spaces
\begin{equation}\label{f:resolutionBaseChange}
   \xymatrix@C=0.5cm{
     0 \ar[r] & K\hat{\otimes}_L D(G,L) \ar[rr]^{\nabla} && K\hat{\otimes}_L D(G,L) \ar[rr]^{pr} && K\hat{\otimes}_L D^\infty(G,L) \ar[r] & 0 }
\end{equation}
by \cite[Lem.\  4.3.6]{SV20}. Since $K\hat{\otimes}_L D(G,L)\cong D(G,K)$ by the proof of \cite[Lem.\ 4.1.2]{SV20}, we also obtain $K\hat{\otimes}_L D^\infty(G,L)\cong D^\infty(G,K)$, i.e., this sequence is the analogue of \eqref{f:resolution} for $K$ replacing $L$ and visibly it arises again by base change $D(G,K)\otimes_{U_K(Lie(G))}-$ from the analogue of \eqref{f:resolutionL}
\begin{equation}\label{f:resolutionK}
\xymatrix@C=0.5cm{
     0 \ar[r] & U_K(Lie(G)) \ar[rr]^{\nabla} && U_K(Lie(G)) \ar[rr]^{pr} && K \ar[r] & 0 }.
\end{equation}
Since $\Hom_{D(G,K)}(D(G,K),V)\cong \Hom_{U_K(Lie(G))}(U_K(Lie(G),V)$ the isomorphism \eqref{f:Ext} follows.

If $\mathrm{RHom}_{D(G)}(D^\infty(G),V)$ denotes any (bounded) complex of $K$-vector spaces having the  groups \eqref{f:Ext} as cohomology, we again write \[K_{f,Lie(G)}(V):=\mathrm{cone}\left( \mathrm{RHom}_{D(G)}(D^\infty(G),V)\xrightarrow{f-\id}\mathrm{RHom}_{D(G)}(D^\infty(G),V) \right)[-1]\] for the induced mapping fibre.

Assume   $\nabla\in Lie(U)=K$ corresponds to $1.$ Then  $V\xrightarrow{\nabla}V$ is a valid (functorial) choice for $\mathrm{RHom}_{D(U)}(D^\infty(U),V)$ and we shall also use the notation $K_{f,\nabla}(V)$ instead.

Finally, we set \[H^\bullet_{\clubsuit,\spadesuit}(V):=h^\bullet(K_{\clubsuit,\spadesuit} (V) ).\]

Note that we have isomorphisms (see proof of \cite[Thm.\ 4.8]{Ko} or \cite[\S 10.8.2]{Wei})
\[\mathrm{RHom}_{D^\infty(G)}(K,\mathrm{RHom}_{D(G)}(D^\infty(G),V))\cong \mathrm{RHom}_{D(G)}(K,V)\]
and, for $G_0\subseteq G$ any $L$-analytic normal subgroup, (see \cite[Exc.\ 10.8.5]{Wei})
\[\mathrm{RHom}_{D (G/G_0)}(K,\mathrm{RHom}_{D(G_0)}(K,V))\cong \mathrm{RHom}_{D(G)}(K,V)\]
in the derived category, therefore inducing the spectral sequences
\[\mathrm{Ext}_{D^\infty(G)}^i(K,\mathrm{Ext}_{D(G)}^j(D^\infty(G),V))\Rightarrow \mathrm{Ext}_{D(G)}^{i+j}(K, V)\] and
\[H^i(G/G_0, \mathrm{Ext}_{D(G_0)}^{j}(K, V))\Rightarrow \mathrm{Ext}_{D(G)}^{i+j}(K, V).\] They both degenerate by the projectivity of $K$ as $D^\infty(G)$- and $D(G/G_0)=K[G/G_0]$-module (cf.\ the proof of \cite[Thm.\ 4.10]{Ko} for the first claim over $L$, from which the general case again follows by complete base change to $K$, and using Maschke's theorem for the second claim).
Moreover, note that $\mathrm{Hom}_{D^\infty(G)}(K,W)\cong W^G$, for any $D^\infty(G,K)$-module $W$, because the Dirac measures $\delta_\gamma\in D^\infty(G,K)$ induce the elements $\delta_\gamma-1$ in the augmentation ideal, which is the kernel of $D^\infty(G,K)\twoheadrightarrow K$ and which is a finitely generated ideal by Cor.\ 4.6 of (loc.\ cit.) plus complete (exact) base change; using this,
the above spectral sequences induce the isomorphisms
\begin{equation}\label{f:LieExt}H^i(Lie(G),V)^G=\mathrm{Ext}_{D(G)}^{i}(K, V)
\end{equation}
and
\begin{equation}
\label{f:groupExt}H^0(G/G_0, \mathrm{Ext}_{D(G_0)}^{j}(K, V))= \mathrm{Ext}_{D(G)}^{j}(K, V).
\end{equation}

\begin{remark}\label{rem:analCoh}
In \cite{Co2} the pro-$L$-locally analytic cohomology groups $H^i_{an}(A^+,M)$ for the $L$-analytic semi-group $A^+\cong\Gamma_L\times \{\varphi^\mathbb{Z}\}$ with $M$ being specified below are defined. By \cite[3.7.6]{Th} they are isomorphic to the cohomology groups $H^i_{\varphi_L,\Gamma_L, an}(M)$  which arise as follows:  Following \cite[\S 5]{Co2} we write $  \cC^\bullet_{an}(G,M)$ for the locally $L$-analytic cochain complex of an $L$-analytic group $G$ with coefficients in $M$ and $H^i_{an}(G,M) := h^i(\cC^\bullet_{an}(G,M))$ for  locally $L$-analytic group cohomology.  More precisely, let $M=\varinjlim_s\varprojlim_r M^{[r,s]}$ with Banach spaces $M^{[r,s]} $  be an LF space  with a pro-$L$-analytic action of $G$ (cf. Definition \ref{defi:proAn}).
If $\mathrm{Maps}_{loc L-an}(G,M^{[r,s]})$ denotes the space of locally $L$-analytic maps from $G$ to  $M^{[r,s]}$, then \[C_{an}^n(G,M)=\varinjlim_s\varprojlim_r\mathrm{Maps}_{loc L-an}(G^n,M^{[r,s]})\] is the space of locally $L$-analytic functions (locally with values in  $\varprojlim_r M^{[r,s]} $ for some $s$ and such that the composite with the projection onto $M^{[r,s]}$ is locally $L$-analytic for all $r$). Then $H^i_{\varphi_L,G, an}(M):=h^i(K_{\varphi_L,G, an}(M))$ is the cohomology of the mapping fibre $K_{\varphi_L,G,an}(M)$ of  $\cC^\bullet_{an}(G,\varphi_L)$ and analogously for $\psi$ instead of $\varphi_L.$
By \cite[Corollary 4.2.6]{Ste1} we have natural isomorphisms
\begin{equation}
 H^i_{ an}(G,M)\cong  \mathrm{Ext}^i_{D(G)}(K,M)
\end{equation}
and hence, for $ \clubsuit\in\{\varphi,\psi\},$
\begin{equation}
 H^i_{\clubsuit,G, an}(M)\cong  H^i_{\clubsuit,D(G)}(M).
\end{equation}
\end{remark}
	For $n \gg 0$ we have that $1+\pi_L^n o_L$ is isomorphic to $\pi_L^n o_L$ via $\log_p$. In particular, we have the chain of isomorphisms
	\begin{align}\label{f:elln}
		\ell_n:\Gamma_n \xrightarrow{\chi_{LT}} 1+\pi_L^n o_L \xrightarrow{\log_p} \pi_L^n o_L \xrightarrow{\cdot \pi_L^{-n}} o_L
	\end{align} which yields
	\begin{align}\label{f:isoGamma_no_L}
		D(\Gamma_n, K) \cong D(o_L, K) \cong \mathcal{R}_K^+,
	\end{align}
	the last isomorphism being the Fourier isomorphism.
	Since $\Gamma_n$ is clopen in $\Gamma$, every locally analytic function on $\Gamma_n$ is the restriction of a locally analytic function on $\Gamma_L$. Hence, by considering the restriction of functions from $\Gamma_L$ to $\Gamma_n$ and taking its dual, we obtain an injective map $D(\Gamma_n, K) \nachinj D(\Gamma_L,K)$.
\begin{definition}
	For $n\gg0$ such that $\Gamma_n = \operatorname{Gal}(L_\infty/L_n)\cong o_L$ we denote by $\mathfrak{Z}_n \in D(\Gamma_n,K)\cong \cR_K^+$ the element corresponding to $Z \in \cR_K^+.$ If the precise choice of $n$ is not relevant we frequently write $(U,\mathfrak{Z})$ instead of $(\Gamma_n,\mathfrak{Z}_n).$
\end{definition}

\begin{remark} \label{rem:Zproperties}
Let $n \gg 0$ such that $\Gamma_n \cong o_L.$ Under the  natural inlcusions $D(\Gamma_{n+1},K) \subseteq D(\Gamma_n,K)$ and $\Lie(\Gamma_L)=\Lie(\Gamma_n) \subseteq D(\Gamma_n,K)$ we have:
\begin{enumerate}
	\item $\mathfrak{Z}_{n+1} = [\pi_L](\mathfrak{Z}_n)$ in the ring $D(\Gamma_n,K).$
	\item $\nabla = \frac{\Omega}{\pi_L^n} \log_{LT}(\mathfrak{Z}_n)$ in $D(\Gamma_n,K).$ In particular, $\nabla$ is divisible by $\mathfrak{Z}_n.$
\end{enumerate}
\end{remark}
\begin{proof}
	For \it{(i)} see \cite[Definition 1.2.10]{Ste1}.
	For \it{(ii)} see (the proof of) \cite[Remark 4.4.8]{SV20}.
\end{proof}
\begin{lemma} \label{lem:Zinvertible}
	Let $A$ be $K$-affinoid and let $M \in \mathfrak{M}^{an}(\mathcal{R}_A).$  Then:
	\begin{enumerate}
		\item For $r \in (0,1)$ large enough the action of $\mathfrak{Z}$ on $(M^{r})^{\psi=0}$ is invertible.
		\item The action of $\mathfrak{Z}$ on $M^{\psi=0}$ is invertible.
		\item Analogous results hold for $M$ (resp. $M^r$) replaced by $M_{t_{LT}}$ (resp. $M^r_{t_{LT}}).$
	\end{enumerate}
\end{lemma}
\begin{proof}
The case $A=K$ is originally treated in \cite[Theorem 4.3.21]{SV20}. For \textit{(i)} in the general case see  \cite[Theorem 2.4.5]{Ste1}. The second point follows by passing to the colimit. The third point also follows by passing to the colimit, where for $mt_{LT}^{-k} \in Mt_{LT}^{-k}$ one extends $\psi$ by setting $\psi(mt_{LT}^{-k}):=\pi_L^{k}\psi(m)t_{LT}^{-k}$ (cf. \cite[Lemma 4.5.23]{SV20} for details in the case $A=K$).
\end{proof}

\begin{lemma}\label{lem:compcoh}\phantom{section}
\begin{enumerate}
\item $H^i_{\clubsuit,\spadesuit}(V)=0$ for $i\neq 0,1,2.$
\item $H^\bullet_{\varphi,D(G)}(M) \cong H^\bullet_{\varphi,Lie(G)}(M)^{ G }$ for $M$ in $\phigam{M}^{an}(\cR)$.\footnote{In \cite[Thm.\ 5.6]{Co2} the analogous statement for $H^i_{an}(A^+,M)$ and $H^i_{Lie}(A^+,M)$, as defined in (loc.\ cit.), is claimed  referring to \cite[Thm.\ 4.2]{FX}, but this only  covers $i=0,1$. }
\item $H^\bullet_{\varphi,D(G)}(M)\cong H^\bullet_{\Psi,D(G)}(M)$ for $M$ in $\phigam{M}^{an}(\cR)$. \footnote{ Cf. \cite[Thm.\ 5.5]{Co2} and \cite[Cor.\ 2.2.3]{BF} for related statements in cohomological degrees $0,1.$}
\end{enumerate}
\end{lemma}

\begin{proof}
Part (i) holds due to the length the of total complex. (ii) follows immediately from \eqref{f:LieExt} upon considering one of the spectral sequences attached to the double complexes arising from the defining mapping fibres.  By \eqref{f:groupExt}, (iii) is reduced to the case $H^\bullet_{\varphi,\mf Z}(M)\cong H^\bullet_{\Psi,\mf Z}(M)$, which is a consequence of Lemma \ref{lem:Zinvertible}.

\end{proof}

\subsection{Finiteness of analytic Cohomology}
\begin{theorem}
	\label{thm:CohomologyFinite}
	Let $A,B$ be $K$-affinoid and let $M$ be an $L$-analytic $(\varphi_L,\Gamma_L)$-module over $\cR_A.$ Let $f\colon A  \to B$ be a morphism of $K$-affinoid algebras. Then:
	\begin{enumerate}[(1)]
		\item $K_{\varphi_L,\mathfrak{Z}}(M) \in \mathbf{D}^{[0,2]}_{\mathrm{perf}}(A).$
		\item The natural morphism $K_{\varphi_L,\mathfrak{Z}}(M) \otimes_{A}^\mathbb{L} B \to K_{\varphi_L,\mathfrak{Z}}(M \hat{\otimes}_A B)$ is a quasi-isomorphism.
	\end{enumerate}
\end{theorem}
\begin{proof}
	See \cite[Theorem 3.3.12]{Ste1}.
\end{proof}
For a commutative ring $R$ and an object $C \in \mathbf D^b(R)$ whose cohomology groups are of finite rank over $R,$ we denote by $\chi_R(C) = \sum_i (-1)^i \operatorname{rank}H^i(C)$ the Euler-Poincar\'{e}-characteristic of $C.$ Recall that a $(\varphi_L,\Gamma_L)$-module is called \textbf{trianguline} if it can be written as a successive extension of rank one modules of character type in the sense of section \ref{sec:rankonecharacters}.
\begin{remark}\label{rem:EP-formula} Let $A/K$ be affinoid and
	let $M $ be a trianguline $L$-analytic $(\varphi_L,\Gamma_L)$-module over $\cR_A.$ Then the Euler-Poincar\'{e} Formula holds, i.e.,$$\chi(M):= \chi(K_{\varphi_L,\mathfrak{Z_n}}(M))=\sum (-1)^i\mathrm{rk}_{\cR_A}(H^i_{\varphi_L,\mf Z_n}(M)) = [\Gamma_L:\Gamma_n]\mathrm{rk}_{\cR_A}(M).$$
\end{remark}
\begin{proof} Without loss of generality we may assume that $M = \cR_A(\delta)$ is an $L$-analytic module of character type (attached to an $A$-valued locally $L$-analytic character $\delta\colon L^\times \to A^\times).$ Then the case $A=K$ is treated in \cite[Remark 6.3]{Ste2}. The validity of the formula can be checked at each maximal ideal of $A.$ Note that $\cR_A(\delta)/\mathfrak{m}$ is a $(\varphi_L,\Gamma_L)$-module of character type over $\cR_{K'}$ for some finite extension $K'/K$ for each $\mathfrak{m} \in \Max(A)$ by the Nullstellensatz and the claim hence follows from the previous case.
\end{proof}
We will require a slight generalization of \ref{thm:CohomologyFinite}. Recall that $K_{\varphi_L,\mathfrak{Z}}(M)$ is (up to shift) quasi isomorphic to the cone of $1-\varphi$ on $\operatorname{RHom}_{D(U,K)}(K,M).$ As a consequence of \cite[Lemma 2.5]{Ste2} $K$ admits a finite projective resolution consisting of finitely generated projective $D(\Gamma_L,K)$-modules. In particular the complex computing $\operatorname{RHom}_{D(U,K)}(K,M)$ (and hence also $K_{\varphi_L,\mathfrak{Z}}(M)$) can be represented by a complex of $A[\Gamma_L/U]$-modules the terms of which are all of the form $\Hom_{D(U,K)}(P,M),$ where $P$ is the restriction of scalars of a projective $D(\Gamma_L,K)$-module with $\Gamma_L$ acting via $(\gamma f)(x) = \gamma (f(\gamma^{-1}x))$ and $A$ acting by multiplication on $M.$
\begin{remark}
	\label{rem:CohomologyRefinement}
	In the situation of \ref{thm:CohomologyFinite}, if we view $K_{\varphi_L,\mathfrak Z}(M)$ as an object in $\mathbf{D}(A[\Gamma_L/U])$  we have
$$K_{\varphi_L,\mathfrak{Z}}(M) \in \mathbf{D}^{[0,2]}_{\mathrm{perf}}(A[\Gamma_L/U]).$$
\end{remark}
\begin{proof}
	The finiteness of the cohomology groups over $A$ already implies that $K_{\varphi_L,\mathfrak{Z}}(M)$ belongs to $ \mathbf{D}^{-}_{\perf}(A[\Gamma_L/U]).$
Choosing a complex of bounded above projective $A[\Gamma_L/U]$-modules representing $K_{\varphi_L,\mathfrak{Z}}(M)$, truncating and using \cite[Lemma 4.1.3]{KPX}, we can conclude that the complex in question is quasi isomorphic to a bounded complex of finitely generated projectives outside of perhaps degree $0,$ where the module is finitely generated over $A[\Gamma_L/U]$ and its underlying $A$-module is flat.
But then it is projective as an $A$-module and by \cite[Lemma 2.5]{Ste2} also projective as an $A[\Gamma_L/U]$-module, hence the claim.
\end{proof}

\subsection{Perfectness of analytic Iwasawa cohomology and the Lubin-Tate deformation}\label{sec:Iwasawa}

For $M$ any $(\varphi_L,\Gamma_L)$-module over any basis consider the complex
\[\cT_{\Psi}(M):=[M\xrightarrow{\Psi-1}M]\]
concentrated in degrees $1$ and $2$, whose cohomology we call {\it (analytic) Iwasawa cohomology} due to Fontaine's classical result, which relates these groups in the \'{e}tale case to usual Iwasawa cohomology defined in terms of Galois cohomology. We set $D:=D(\Gamma_L,K).$ The following result \cite[Thm.\ 4.8]{Ste2}   will be central for the whole article:

\begin{theorem}\label{thm:perf}
For $M\in\phigam{M}^{an}(K)$ trianguline, $\cT_{\Psi}(M) $ is a perfect complex of $D$-modules, i.e., belongs to $\mathbf{D}^b_{\mathrm{perf}}(D)$.
\end{theorem}

For the rest of this subsection we assume that $M\in\mathfrak{M}^{an}(K)$ is   trianguline.

Later for our approach it will be important to interpret Iwasawa cohomology as  analytic cohomology of a deformation $\mathbf{Dfm}(M)$ of $M$ via generalized Herr complexes. This deformation lives over the character variety $\mathfrak{X}_{\Gamma_L}$ (base changed to $K$)  of the locally $L$-analytic group $\Gamma_L$ (\cite{ST2}) and will allow to use density arguments to deduce many   properties of the Epsilon-isomorphism for rank one modules just from properties over its de Rham points.

We pick an affinoid  cover $X_n\cong\Sp(D_n)$ of $\mathfrak{X}_{\Gamma_L}$ with $D_n:=K[\Gamma_L]\otimes_{K[U]}D_{r_n}(U,K)$ for a decreasing sequence $r_n$ such that each $D_{r_n}(U,K) $ corresponds to the ring of rigid analytic functions on the annulus $[r_n,\infty]$ via the Fourier isomorphism for $D(U,K)$.  Over the space $\mathfrak{X}_{\Gamma_L}$ we have the sheaf of Robba rings $\cR_{\mathfrak{X}_{\Gamma_L}}$ given by mapping $X_n$ to $\cR_{\mathcal{O}_{\mathfrak{X}_{\Gamma_L}}(X_n)}$ and $\mathbf{Dfm}(M)$ should be thought of as a $(\varphi_L, \Gamma_L)$-module (sheaf) over $\cR_{\mathfrak{X}_{\Gamma_L}}$ (but unfortunately, Schneider's and Teitelbaum's formalism of coadmissible modules does not apply here as $\cR_{\mathfrak{X}_{\Gamma_L}}(\mathfrak{X}_{\Gamma_L})$ does not form a Fr\'{e}chet-Stein algebra in any obvious sense):

For an $L$-analytic $(\varphi,\Gamma_L)$-module $M$ over $\cR_L$ we define
	$$\mathbf{Dfm}(M)(X_n):=\mathbf{Dfm}_n(M):=\mathcal{O}_{\mathfrak{X}_{\Gamma_L}}(X_n) \hat{\otimes}_L M,$$ where $\Gamma_L$ acts diagonally, on the left factor via the inversion and on $M$ via its  given action. For each $n$  this is a $(\varphi,\Gamma_L)$-module $M$ over $\cR_{\mathcal{O}_{\mathfrak{X}_{\Gamma_L}}(X_n)}$ by \cite[Prop.\ 3.2]{Ste2}.


As definition for the generalized Herr complex for the sheaf $\mathbf{Dfm}(M)$, philosophically, we would  like to take the complex in $\mathbf{D}(D(\Gamma_L,K))$ \footnote{instead of e.g.\ forming the generalized Herr complex attached to  the  global sections $\mathbf{Dfm}(M)(\mathfrak{X}_{\Gamma_L})!$ } defined as total derived sheaf cohomology of the complex of sheaves
\[K_{\Psi,D(\Gamma_L,K)}(\mathbf{Dfm}(M))  =\cT_{\Psi}(\mathbf{Dfm}(M))\otimes^\mathbb{L}_{D(\Gamma_L,K),\mathrm{diag}}K
\ ("\cong K_{\Psi,\mf Z}(\mathbf{Dfm}(M))\otimes_{K[\Gamma_L/U],\mathrm{diag}}^\mathbb{L}K"),\]
where for the last (quasi-)isomorphism in quotation marks we used implicitly the free resolution
\begin{equation}\label{f:resZ}
  \xymatrix@C=0.5cm{
  0 \ar[r] & D(\Gamma_L,K) \ar[rr]^{\mathfrak Z} &&  D(\Gamma_L,K) \ar[rr]^{ } && K[\Gamma_L/U] \ar[r] & 0 }
\end{equation}
which induces an isomorphism $\cT_{\Psi}(\mathbf{Dfm}(M))\otimes^\mathbb{L}_{D(\Gamma_L,K),\mathrm{diag}}K[\Gamma_L/U]\cong K_{\Psi,\mf Z}(\mathbf{Dfm}(M)).$\\
   But strictly speaking one needs a resolution of $D(\Gamma_L,K)\otimes_K K[\Gamma_L/U]$-modules in order to define the $(D_n,K[\Gamma_L/U])$-bimodule structure on $K_{\Psi,\mf Z}(\mathbf{Dfm}_n(M))\cong  \cT_{\Psi}(\mathbf{Dfm}_n(M))\otimes^\mathbb{L}_{D(\Gamma_L,K)}K[\Gamma_L/U]$. To this end we can formally work with the resolution\be \label{eq:bimodules}\xymatrix@C=0.5cm{
    0 \ar[r] & C :=\ker\ar[rr]^{ } && D(\Gamma_L,K)\otimes_K K[\Gamma_L/U] \ar[rr]^-{} && K[\Gamma_L/U] \ar[r] & 0 }\ee for an explicit construction of $K_{\Psi,\mf Z}(\mathbf{Dfm}(M))$ in $\mathbf{D}((D_n,K[\Gamma_L/U])\mathrm{-bimod}),$ with the last non-trivial map given by $a\otimes b\mapsto H(a)\cdot b,$ where  $H\colon D(\Gamma_L,K)\to K[\Gamma_L/U]$ denotes the augmentation map sending the Dirac distributions of $u\in U$ to $1.$ Indeed, the kernel $C$ is projective (hence flat) as a $D$-module by the same reasoning as for \cite[Lem.\ 2.5]{Ste2}.  This sequence is related to (the direct sum of) the sequences
    \[\xymatrix@C=0.5cm{
      0 \ar[r] & I_{\chi_i}:=\ker \ar[rr]^{ } && D(\Gamma_L,K) \ar[rr]^-{\gamma \mapsto \chi_i^{-1}(\gamma) } && K(\chi_i) \ar[r] & 0 }\] for the characters $\chi_i$ of $\Gamma_L$ which factor through $\Gamma_L/U.$

Instead of verifying that we really have a complex of (coherent) sheaves we  just use the facts as a motivation  that on a Stein space $\Gamma (\mathfrak{X}_{\Gamma_L},-)=\varprojlim_n \Gamma(X_n,-)$ and     that higher sheaf cohomology of coherent sheaves vanishes on affinoids. Thus we rather take the total derived inverse limit as a formal definition, i.e.,
\[C^\bullet:=R\Gamma_{\Psi,D(\Gamma_L,K)}(\mathfrak{X}_{\Gamma_L}, \mathbf{Dfm}(M)):=\mathbf{Rlim}\bigg( K_{\Psi,D(\Gamma_L,K)}(\mathbf{Dfm}_n(M)) \bigg).\]
The following results are variants of those in \cite[\S 3.3]{Ste2}; among others they are based on the observation that for the sheaf of cohomology groups sending $X_n$ to $H^i_{\Psi, D(\Gamma_L,K)}(\mathbf{Dfm}_n(M))$
 the formalism of coadmissible modules over $D(\Gamma_L,K)$ does apply.
\begin{theorem}\label{thm:def}\phantom{section}
\begin{enumerate}
\item For all $i$, the cohomology groups $H^i_{\Psi,D(\Gamma_L,K)}(\mathfrak{X}_{\Gamma_L}, \mathbf{Dfm}(M))$ of the complex\linebreak $R\Gamma_{\Psi,D(\Gamma_L,K)}(\mathfrak{X}_{\Gamma_L}, \mathbf{Dfm}(M)) $ coincide with the global sections \[ \varprojlim_n H^i_{\Psi,D(\Gamma_L,K)}(\mathbf{Dfm}_n(M))\]  of the sheaf of cohomology groups sending $X_n$ to $H^i_{\Psi, D(\Gamma_L,K)}(\mathbf{Dfm}_n(M)). $
\item There is an isomorphism in $\mathbf{D}^b_{\mathrm{perf}}(D(\Gamma_L,K))$ \[R\Gamma_{\Psi,D(\Gamma_L,K)}(\mathfrak{X}_{\Gamma_L}, \mathbf{Dfm}(M))\cong \cT_{\Psi}(M).\]
\end{enumerate}
\end{theorem}

\begin{remark}\label{rem:normalization}
 In accordance with \eqref{f:normalization} the isomorphism in (ii) only becomes independent of the choice of $\mf Z$ if we insert the scalar factor $C_{Tr}({\mf Z_n})$ (see \eqref{f:CTr} below) in the identification  $ \mathbf{Dfm}_n(M)/\mf Z\mathbf{Dfm}_n(M)\cong D_n\tilde{\otimes}_{D(U,K)}M$ in the proof of Lemma \ref{lem:dfmTPsi}, compare with  \cite[(32), p. 369]{NaANT}.
\end{remark}

For the proof of Theorem \ref{thm:def} we need the following lemma for which we recall some  notation from \cite[Def.\ 3.20]{Ste2}: We define $D_n\hat{\otimes}_D M^r$ as the completion of $D_n\otimes_D M^r$ with respect to the quotient topology of the projective tensor product $D_n\otimes_{K,\pi}M^r$. Then we set $D_n\tilde{\otimes}_D M:=\varinjlim_r D_n\hat{\otimes}_D M^r.$
\begin{lemma}\label{lem:dfmTPsi}
\begin{enumerate}
  \item The natural map $D_n\otimes_D M\to D_n\tilde{\otimes}_D M$ induces a quasi-isomorphism
  \[D_n\otimes_D \cT_\Psi(M)=\cT_\Psi(D_n\otimes_D M)\to \cT_\Psi(D_n\tilde{\otimes}_D M).\]
  \item  Viewing $\mathbf{Dfm}_n(M)$ as $D_n$-module via the left tensor factor, there is a natural isomorphism in $\mathbf{D}(D_n)$
  \[ \mathbf{Dfm}_n(M) \otimes^\mathbb{L}_{D(\Gamma_L,K),\mathrm{diag}}K\cong D_n\tilde{\otimes}_{D}M[0],\] where the latter module is considered as complex concentrated in degree $0.$
\end{enumerate}
\end{lemma}
\begin{proof}
For (i) the same proof as for \cite[Lem.\ 3.23]{Ste2} works and the assumptions are satisfied by Theorem \ref{thm:perf}, but note that there $D_n$, $D$ have a slightly different meaning. The augmentation ideal $I_{\Gamma_L}$ is a finitely generated submodule of $D(\Gamma_L,K)$ and thus projective as a $D(U,K)$-module since the latter is a Pr\"{u}fer Domain. Using \cite[Lemma 2.5]{Ste2} one can conclude projectivity as a $D(\Gamma_L,K)$-module. Using that $D\mf Z$ is contained in $I_{\Gamma_L}$ gives rise to the projective resolution of $D$-modules
\[\xymatrix@C=0.5cm{
  0 \ar[r] & I_{\Gamma_L} \ar[rr]^{ } && D \ar[rr]^{ } && K \ar[r] & 0. }\] We can represent the complex in question in (ii) by the complex
  \[\mathbf{Dfm}_n(M)\otimes_{D,\mathrm{diag}} I_{\Gamma_L}\to \mathbf{Dfm}_n(M)\otimes_{D,\mathrm{diag}} D\]
with cokernel
\begin{align*}
\mathbf{Dfm}_n(M)/I_{\Gamma_L}\mathbf{Dfm}_n(M)&\cong (\mathbf{Dfm}_n(M)/\mf Z\mathbf{Dfm}_n(M))/(I_{\Gamma_L}/D\mf Z)\\
&=(D_n\tilde{\otimes}_{D(U,K)}M)_{\Gamma_L/U}\\
&=D_n \tilde{\otimes}_DM,
\end{align*}
 where  $ \gamma \in \Gamma_L\subseteq D(\Gamma_L,K)$ acts diagonally (via $\gamma(a \otimes b) = \delta_{\gamma^{-1}} a \otimes \gamma b)$) on $D_n\tilde{\otimes}_{D(U,K)}M$ and this action factors over $\Gamma_L/U$. For the second equality we use an obvious variant of \cite[(31)]{Ste2}, while by the exactness of colimits the last one is easily reduced to the claim that on the level of models $M^r$ we have
\[(D_m\hat{\otimes}_{D(U,K)} M^r)_{\Gamma/U}\cong D_m\hat{\otimes}_{D(\Gamma_L,K)} M^r.\] Since $\Gamma/U$ is finite and taking $\Gamma_L/U$-invariants in this situation is exact by Maschke's theorem, this follows in the context of Fr\'{e}chet spaces by completion from the well-known fact that
\[(D_m{\otimes}_{D(U,K)} M^r)_{\Gamma/U}\cong D_m{\otimes}_{D(\Gamma_L,K)} M^r.\]
 The injectivity of the non-trivial differential in the above complex can be checked by  calculating instead the cohomology in degree $-1$ of $\mathbf{Dfm}_n(M) \otimes^\mathbb{L}_{D(\Gamma_L,K),\mathrm{diag}}K[\Gamma_L/U]$, because taking $\Gamma_L/U$-(co)invariants is exact and leads to the original complex
\[\bigg(\mathbf{Dfm}_n(M) \otimes^\mathbb{L}_{D(\Gamma_L,K),\mathrm{diag}}K[\Gamma_L/U]\bigg)\otimes_{K[\Gamma_L/U]}K\cong\mathbf{Dfm}_n(M) \otimes^\mathbb{L}_{D(\Gamma_L,K),\mathrm{diag}}K.\]
For this composition of functors it is crucial that $\mathbf{Dfm}_n(M) \otimes^\mathbb{L}_{D(\Gamma_L,K),\mathrm{diag}}K[\Gamma_L/U]$ belongs to $\mathbf{D}((D_n,K[\Gamma_L/U])\mathrm{-bimod})$  as in \eqref{eq:bimodules} in order to allow an action by $\Gamma_L/U$. But then the vanishing in degree $-1$ can be checked just as complex of $K$-vector spaces and therefore it suffices to calculate the derived functor by a projective resolution of $D$-modules (instead of bi-modules).  To this end we use the resolution \eqref{f:resZ}, which leads to the complex
\[\mathbf{Dfm}_n(M)\xrightarrow{\mf Z} \mathbf{Dfm}_n(M),\]
which is left exact by an obvious variant of \cite[(31)]{Ste2}, again.
\end{proof}

\begin{proof}[Proof of  Theorem \ref{thm:def}] Using Lemma \ref{lem:dfmTPsi} we obtain  isomorphisms in $\mathbf{D}(D_n)$
\begin{align*}
  K_{\Psi,D(\Gamma_L,K)}&(\mathbf{Dfm}_n(M))\\
  & \cong \mathrm{cone} \bigg( \mathbf{Dfm}_n(M)\otimes^\mathbb{L}_{D(\Gamma_L,K),\mathrm{diag}}K\xrightarrow{\Psi-1} \mathbf{Dfm}_n(M)\otimes^\mathbb{L}_{D(\Gamma_L,K),\mathrm{diag}}K \bigg)[-2] \\
    & \cong \cT_\Psi(D_n\tilde{\otimes}_{D}M ) \\
   & \cong D_n{\otimes}_{D}\cT_\Psi(M )
\end{align*}
compatible for the variation in $n$  by an obvious variant of Theorem \ref{thm:CohomologyFinite} (2). Thus, combining \cite[Prop.\ 3.15]{Ste2} with Theorem \ref{thm:perf} we obtain in $\mathbf{D}^b_{\perf}(D)$ an isomorphism
\[\cT_{\Psi}(M)\cong\mathbf{Rlim} \big(D_n\otimes_D\cT_{\Psi}(M)\big)\cong \mathbf{Rlim} \big(K_{\Psi,D(\Gamma_L,K)}(\mathbf{Dfm}_n(M))\big).\]
This proves (ii) while (i) follows by the same arguments as in \cite[Rem.\ 3.16]{Ste2} using that
the projective system $( K_{\Psi,D(\Gamma_L,K)}(\mathbf{Dfm}_n(M)) )_m $ defines a consistent object in $\mathbf{D}(\mathrm{mod}(\mathbb{N},D))$ (using the notation of (loc.\ cit.))  together with the fact that $D$ is a Fr\'{e}chet-Stein algebra.
\end{proof}

\subsection{Replacing Local Tate duality}\label{sec:duality}

In this subsection we develop local duality analogous to local Tate duality for Galois cohomology, see \cite{Her,liu2007cohomology} for an approach purely in terms of $(\varphi,\Gamma)$-modules.
We focus  technically on the complexes $K_{f,\mf Z}$ and shall then apply \eqref{f:groupExt} to deal with $K_{f,U}$. Assume henceforth that $M$ is an analytic $(\varphi_L,\Gamma_L)$-module over $\cR=\cR_K$. For an analytic character $\delta: L^\times \to K^\times$   we define the  twisted module $M(\delta ) \in \mathfrak{M}^{an}(\cR),$ where $M(\delta) := M\otimes_{\cR}\cR(\delta)$ endowed with the diagonal  $\varphi_L$- and $\Gamma_L$-action.
Recall the residue map (at $Z$)
\[Res: \Omega^1_{\cR} := \cR dZ\otimes_{\cR}\cR(\delta^{un}_{\chi(\pi)}) \to K,\;\;\; \sum_i a_iZ^idZ\otimes \mathbf{e}_{\delta^{un}_{\chi(\pi)}} \mapsto a_{-1},\]
and that  the $(\varphi_L,\Gamma_L)$-action on $ \cR dZ $  with respect to the basis $d\log_{LT} = g_{LT}dZ$ is given by the character $\chi_{LT}.$  \footnote{The action on $\Omega^1_{\cR}$ differs by $\delta^{un}_{\chi(\pi)}$ from the action considered in \cite[Section 4]{SV20} and agrees with the action from \cite[1.3.5]{Co2}.} As a formal consequence, we have the following:
\begin{lemma}\label{Omega-as-twist}
The map
\begin{align*}
  \cR(\chi) & \xrightarrow{\;\cong\;} \Omega^1_\cR, \\
                        f {\mathbf{e}_{\chi}} & \longmapsto fd\log_{LT} \otimes \mathbf{e}_{\delta^{un}_{\chi(\pi)}} ,
\end{align*}
is an isomorphism of $(\varphi_L,\Gamma_L)$-modules.
\end{lemma}
Setting $\tilde{M}:= \Hom_{\cR}(M,{\cR})(\chi) \cong \Hom_{\cR}(M,\Omega^1_{\cR}) $, for any finitely generated projective $\cR$-module $M$, we obtain more generally the pairing
\begin{equation}\label{f:res-pair-general-Psi}
   \langle\;,\;\rangle:=\langle\;,\;\rangle_M: \tilde{M} \times M \to K,\;\;\; (g,f)\mapsto Res(g(f)),\footnote{Note that Colmez defines $\Omega Res(\sigma_{-1}(g)(f ) )$ instead.   }
\end{equation}
 where by abuse of notation we also write $Res:\cR(\chi)\to K$ for the map sending $\sum_i a_iZ^i \otimes \mathbf{e}_{\chi}$ to $a_{-1}.$
This map identifies $M$ and $\tilde{M}$ with the   (strong)   topological duals of $\tilde{M}$ and $M$, respectively. Moreover, the isomorphism $\tilde{M} \cong \Hom_{K,cts}(M,K)$ (induced by $ \langle\;,\;\rangle$) is $D(\Gamma_L,K)$-linear by \cite[Corollary 4.5.4]{SV20}.

\begin{lemma}
	\label{lem:SerreDual}
The residuum map induces an isomorphism $Res: H_{\varphi_L,\mf Z_n}^2(\Omega^1_{\cR} )^{\Gamma_L}\cong K.$
\end{lemma}
\begin{proof}
  We know from Lemma \ref{lem:dim} that $\dim_K H_{\varphi_L,\mf Z_n}^2(\Omega^1_{\cR} )^{\Gamma_L}=1$ while $Res$ is a surjection   $\Omega^1_{\cR} \twoheadrightarrow K$ which factorizes over $(\varphi_L-\id)\Omega^1_{\cR}$ and $\mf Z\cdot \Omega^1_{\cR}$ by \cite[Lemma 4.5.1]{SV20} or \cite[Prop.\ 1.5]{Co2}. The claim follows as $H_{\varphi_L,\mf Z_n}^2(\Omega^1_{\cR} )^{\Gamma_L}=H_{\varphi_L,\mf Z_n}^2(\Omega^1_{\cR} )_{\Gamma_L}.$
\end{proof}
 For compatibility questions we renormalise the residuum map to obtain the trace map $Tr=C_{Tr}({\mf Z_n})Res:H_{\varphi_L,D(\Gamma_L)}^2(\Omega^1_{\cR} )=H_{\varphi_L,\mf Z_n}^2(\Omega^1_{\cR} )^{\Gamma_L}\cong K$ by setting
\begin{equation}\label{f:CTr}
C_{Tr}({\mf Z_n}):=\frac{q}{q-1}\frac{\Omega}{\pi_L^n}.
\end{equation}
 Note that for $L=\qp$ and $\pi_L=p$ this trace map is compatible with Tate's trace map in Galois cohomology by \cite[Prop.\ 5.2]{NaANT}. Independence of $n$ follows by the same argument as for Definition \ref{def:dualexp}  below. The principle is explained as follows:

 The map of complexes, for $m\geq n$,
 \[\xymatrix{
   M \ar@{=}[d]_{ } \ar[r]^{{\mf Z_n}} & M \ar[d]^{\mathfrak{Q}_{m-n}({\mf Z_n})} \\
   M \ar[r]^{{\mf Z_m}} & M  }\]
induces the restriction maps $res^n_m:H^i_{{\mf Z_n}}(M)\to H^i_{{\mf Z_m}}(M),$ where   $\mathfrak{Q}_{m-n}({\mf Z_n}):=\frac{{\mf Z_m}}{{\mf Z_n}}=\frac{\varphi_L^{m-n}({\mf Z_n})}{{\mf Z_n}}$ with $\mathfrak{Q}_{m-n}(0)=\pi_L^{m-n}.$ Since $C_{Tr}({\mf Z_n})= \pi_L^{m-n}C_{Tr}({\mf Z_m})$ by \eqref{f:derivative} the isomorphism
\begin{align}\label{f:normalization}
 \vartheta_n: H^1_{{\mf Z_n}}(M)^{\Gamma_L}& \xrightarrow{\cong}M_{\Gamma_L}, [x]\mapsto [C_{Tr}({\mf Z_n})x]
\end{align}
into the $\Gamma_L$-coinvariants is compatible with $res^n_m,$ i.e., the diagram
\[\xymatrix{
  H^1_{{\mf Z_n}}(M)^{\Gamma_L} \ar[d]_{res^n_m} \ar[r]^{\vartheta_n} & M_{\Gamma_L}   \\
  H^1_{{\mf Z_m}}(M)^{\Gamma_L} \ar[ru]^{\vartheta_m}     }\]
commutes.

For a complex $(X^\bullet, d_X )$ of topological $K$-vector spaces we define its $K$-dual $((X^*)^\bullet, d_{X^*})$ to be the complex with \[(X^*)^i:=\mathrm{Hom}_{K,cts}(X^{-i},K)\] and \[d_{X^*}(f):=(-1)^{\mathrm{deg}(f)-1} f\circ d_X.\]

The following lemma is taken from \cite[Lemma 5.2.4 and Remark 5.2.6]{SV20}.

\begin{lemma}\label{lem:strict}
Let $(\cC^\bullet,d^\bullet)$ be a complex in the category of  locally convex topological $F$-vector spaces.
\begin{enumerate}
  \item   If $\cC$ consists of Fr\'{e}chet spaces and  $h^i(\cC^\bullet)$ is finite-dimensional over $F$, then  $d^{i-1}$ is strict and has closed image.
  \item  If  $d^{i}$ is strict and either $F$ is spherically complete or the spaces involved in degree $i$ are of countable type\footnote{From \cite{PGS} we recall that a locally convex vector space $V$ is said to be of {\it countable type}, if for every continuous
seminorm $p$ on $V$ its completion $V_p$ at $p$ has a dense subspace of countable algebraic dimension. They are  stable under forming subspaces, linear images,
projective limits, and countable inductive limits, cf.\    theorem 4.2.13 in (loc.\ cit.). By  corollary 4.2.6 in (loc.\ cit.) for such vector spaces the Hahn-Banach theorem holds, too. By \cite[Prop.\ 5.4.3]{Th} the Robba ring over any complete intermediate field $\mathbb{Q}_p\subseteq K\subseteq \mathbb{C}_p$ (and hence also finitely generated modules over it) is of countable type as $K$-vector space.}, then $h^{-i}(\cC^*)\cong h^{i}(\cC)^*.$
  \item If $\cC^\bullet$ consists of $LF$-spaces, $\cC^{i+2}=0$ and $h^i(\cC^\bullet)$ is finite dimensional, then $d^i$ is strict.
      \item If $V\xrightarrow{\alpha} W$ is a continuous linear map of Hausdorff $LF$-spaces over $F$ with finite dimensional cokernel, then $\alpha$ is strict and has closed image.
\end{enumerate}
\end{lemma}

The translation  $X[n]$ of a complex $X$ is given by $X[n]^i:=X^{i+n}$ and $d^i_{X[n]}:=(-1)^nd_X^{i+n}.$

Let $\iota$ denote the involution on $D(o_L,K)$ induced by the inversion on the group $o_L$. We observe that
\begin{equation}\label{f:lambda}
 \mf Z^\iota=\lambda \mf Z
\end{equation}
for  a unit   $\lambda\in D(o_L,K)$ as they both generate the augmentation ideal: more explicitly,  $\mf Z^\iota=[-1]({\mf Z}),$ $\lambda^{-1}=\lambda^\iota$.

\begin{theorem}\phantomsection\label{lem:explicitPairing}
\begin{enumerate}
\item There is a canonical quasi-isomorphism
\begin{equation}
\xymatrix{{K_{\varphi,\mf Z}(M ):} &  0 \ar[r]  & M \ar@{=}[d] \ar[rr]^{\begin{pmatrix}
                                        \varphi-1\\    {\mf Z}\end{pmatrix}
    }  & & {M}\oplus {M}\ar[d]^(0.35){-\Psi\oplus \lambda}  \ar[rr]^(0.55){\begin{pmatrix}
                          {\mf Z}& 1-\varphi  \\
                               \end{pmatrix}
      }    &&   {M} \ar[d]^{-\lambda\Psi}  \ar[r]^{ } & 0 \\
 {K_{\Psi,{\mf Z}^\iota}(M ):} &  0 \ar[r]  & M    \ar[rr]^{\begin{pmatrix}
                                        \Psi-1\\    {\mf Z}^\iota\end{pmatrix}
    } & & M\oplus M   \ar[rr]^(0.55){\begin{pmatrix}
                          {\mf Z}^\iota & 1-\Psi  \\
                               \end{pmatrix}
      }  &&   M   \ar[r]^{ } & 0. }
 \end{equation}
\item
Via the pairing \eqref{f:res-pair-general-Psi} there are canonical isomorphisms of complexes in the derived category $\mathbf{D}(K)$
\[K_{\varphi_L,\mf Z}(M)^*\cong K_{\Psi_L,{\mf Z}^\iota}({M})^*{ \cong} K_{\varphi_L,\mf Z}({M}^*)[2]{ \cong} K_{\varphi_L,\mf Z}(\tilde{M})[2].\]
given by the following diagram of quasi-isomorphisms
\begin{equation}\label{f:dualitycomplex}
{\footnotesize \xymatrix{{K_{\varphi,\mf Z}(M )^*[-2]:} &  0 \ar[r]  & M^*    \ar[rr]^{- \begin{pmatrix}
                          {\mf Z}& 1-\varphi  \\
                               \end{pmatrix}^*
 } & & ({M}\oplus {M})^*    \ar[rr]^{\begin{pmatrix}
                                       \varphi-1\\    {\mf Z}\end{pmatrix}^*
  } &&   {M}^* \ar@{=}[d]  \ar[r]^{ } & 0 \\
 {K_{\Psi,{\mf Z}^\iota}(M )^*[-2]:} &  0 \ar[r]  & M^* \ar[u]_{(-\lambda\Psi)^*}\ar@{=}[d]  \ar[rr]^{- \begin{pmatrix}
                          {\mf Z}^\iota& 1-\Psi  \\
                               \end{pmatrix}^*
 } & & ({M}\oplus {M})^*    \ar[u]^{\Upsilon}_\cong \ar[rr]^{\begin{pmatrix}
                                       \Psi-1\\    {\mf Z}^\iota\end{pmatrix}^*
  } &&   M^*    \ar[r]^{ }\ar@{=}[d]  & 0 \\
{K_{\varphi,{\mf Z}}({M}^*)}: &  0 \ar[r]^{ }   & { {M}^*  } \ar@{=}[d] \ar[rr]^{\begin{pmatrix}
                                       \varphi-1\\    {\mf Z}\end{pmatrix}
    } && {{M}^*\oplus  {M}^* }\ar[u]^{\Xi}\ar@{=}[d] \ar[rr]^(0.55){\begin{pmatrix}
                          {\mf Z}& 1-\varphi  \\
                               \end{pmatrix}
      } & & { {M}^* } \ar@{=}[d] \ar[r]^{ } & 0   & \\
      {K_{\varphi,{\mf Z}}(\tilde{M})}: &  0 \ar[r]^{ }   & { \tilde{M}  }  \ar[rr]^{\begin{pmatrix}
                                        \varphi-1\\    {\mf Z}\end{pmatrix}
    } && {\tilde{M}\oplus  \tilde{M} } \ar[rr]^(0.55){\begin{pmatrix}
                          {\mf Z}& 1-\varphi  \\
                               \end{pmatrix}
      } & & { \tilde{M} }  \ar[r]^{ } & 0
   }}\end{equation}
with $\Upsilon=(-\Psi\oplus \lambda)^*: (M\bigoplus M)^*\to (M\bigoplus M)^*$ and $\Xi(x,y)= y\oplus - x$. In particular, we obtain  isomorphisms
\begin{align}\label{f:dualityisos}
  H^i_{\varphi_L,\mf Z}(M)^* &\cong   H^{2-i}_{\Psi_L,\mf Z}(\tilde{M}){\cong} H^{2-i}_{\varphi_L,\mf Z}(\tilde{M}).
\end{align}
induced by the perfect pairings, denoted by  $\langle-,-\rangle:=\langle-,-\rangle_M$,
\[H^1_{\varphi, \mf Z}(M) \times H^1_{\varphi, \mf Z}(\tilde{M}) \to K, (\overline{(m,n)},\overline{(f,g)})\mapsto -Res\Big(\varphi(g)(m)+({\lambda}^\iota f)( n)\Big),\]
\[H^2_{\varphi, \mf Z}(M) \times H^0_{\varphi, \mf Z}(\tilde{M}) \to K, (\overline{m},\tilde{n})\mapsto -Res \Big(\tilde{n}(\lambda^\iota\varphi(m))\Big) ,\]
\[H^0_{\varphi, \mf Z}(M) \times H^2_{\varphi, \mf Z}(\tilde{M}) \to K, ({m},\overline{\tilde{n}})\mapsto Res \Big(\tilde{n}(m)\Big).\]
\end{enumerate}
\end{theorem}

\begin{remark}\label{rem:dualePaarung} Identify $M$ with $\tilde{\tilde{M}}$ via $m \mapsto m^{**}$ and consider the pairing in degree $(1,1)$ from Theorem \ref{lem:explicitPairing} obtained by exchanging the roles of $M$ and $\tilde{M},$ i.e.,
	
\begin{align*} 	\langle -,-	\rangle_{\tilde{M}} \colon   H^1_{\varphi,\mathfrak{Z}}(\tilde{M}) \times H^1_{\varphi,\mathfrak{Z}}(\tilde{\tilde{M}})&\to K,\\
((f,g),(m^{**},n^{**})) &\mapsto  -\operatorname{Res}(\varphi(n^{**})(f)+(\lambda^{\iota}m^{**})(g)).
\end{align*}
We have
	\begin{equation}
	 \langle \overline{(m,n)},\overline{(f,g)}\rangle_{M} = -\langle \overline{(f,g)},\overline{(m^{**},n^{**})} \rangle_{\tilde{M}}.
	\end{equation}
	In the other degrees consider
	\begin{align*}
		\langle -,-	\rangle_{\tilde{M}}\colon H^2_{\varphi,\mathfrak{Z}}(\tilde{M}) \times H^0_{\varphi,\mathfrak{Z}}(\tilde{\tilde{M}})&\to K,\\
	{	(\overline{f}, m^{**})} &\mapsto   \operatorname{Res}(m^{**}(-\lambda^\iota(\varphi(f)))),
	\end{align*}
	satisfying {$\langle \overline{f}, {m^{**}}\rangle_{\tilde{M}} = \langle  {m},\overline{f} \rangle_{M}$}
	and
	\begin{align*}
	H^0_{\varphi,\mathfrak{Z}}(\tilde{M}) \times H^2_{\varphi,\mathfrak{Z}}(\tilde{\tilde{M}})&\to K,\\
(g,\overline{n^{**}}) &\mapsto   \operatorname{Res}({ n^{**}(g)}),
	\end{align*}
satisfying
$\langle g,\overline{n^{**}}\rangle_{\tilde{M}} = \langle \overline{n},g \rangle_{M}.$ \footnote{ In the cyclotomic case $L= \QQ_p$ and $\mathfrak{Z} = \gamma-1$ one has $\lambda^\iota = -\gamma$ because $\mathfrak{Z}^\iota = \gamma^{-1}-1 = (-\gamma^{-1})(\gamma-1).$ We see that the pairing from \cite[Definition 2.13]{NaANT} agrees with our $\langle -,- \rangle_{\tilde{M}}.$}
\end{remark}
\begin{proof}
	 By viewing $K_{\varphi,\mathfrak{Z}}(\tilde{M})$ as a Koszul complex attached to the automorphisms $\varphi-1,\mathfrak{Z}$ of $M$ one can see that $\mathfrak{Z}$ and $\varphi_L-1$ act as $0$ on the cohomology groups. Since $$\lambda = -1 + \text{ terms divisible by } \mathfrak{Z}$$ we see that the class $\overline{(f,g)} \in H^1_{\varphi,\mathfrak{Z}}(\tilde{M})$ is equal to the class of $\overline{(-\lambda f,- \lambda g)}.$ Now let $ \overline{(m,n)} \in H^1_{\varphi,\mathfrak{Z}}(M).$ Using $\mathfrak{Z}f = (\varphi-1)g$ and $\mathfrak{Z}m = (\varphi-1)n$ we compute
	 \begin{align*}
	 	\langle  \overline{(m,n)},\overline{(f,g)}\rangle_{M}	 &=  Res\Big(-\varphi(g)(m)-(\lambda^\iota f)(n)\Big)\\
	 		 &=	Res\Big(-[g+\mathfrak{Z}f](m)-(\lambda^\iota f)(\varphi(n)-\mathfrak{Z}m)\Big)\\
	 		 &= Res\Big(-g(m)-(\lambda^\iota f)\varphi(n))- \underbrace{(\mathfrak{Z}f)(m) + (\lambda^\iota \mathfrak{Z}^\iota f)(m)}_{=0}\Big)\\
	 		 & = Res\Big((\lambda g)(m)+(\lambda \lambda^\iota f)(\varphi(n))\Big)\\
			 & =  Res\Big(m^{**}(\lambda g)+\varphi(n^{**})( f)\Big)\\
	 		 & = - \langle \overline{(f,g)},\overline{(m^{**},n^{**})} \rangle_{\tilde{M}},
	 \end{align*}
	 where in the fifth equation we replace $\overline{(f,g)}$ with $(-\lambda f,-\lambda g).$ Now consider the degree $(0,2)$ case  with regard to $\langle-,-\rangle_M$. Since $\overline{\varphi(f)} = \overline{f}$ we get $\operatorname{Res}(f(m)) = \operatorname{Res}(\varphi(f)(m)) =  - \operatorname{Res}(\varphi(f)(\lambda m))$ using that $\mathfrak{Z}{m} = 0$ in $H^0$ and hence $\lambda {m} =-{m}.$
The computation in degree $(2,0)$ is similar.
\end{proof}


Later in explicit calculations we will need to work partly with $\Psi$-versions, which we therefore establish in the next remark.
\begin{remark}\label{rem:asymmetricII}
As a obvious variant of (i) in Theorem \ref{lem:explicitPairing} there is also a
canonical quasi-isomorphism
\begin{equation}
\xymatrix{{K_{\varphi,\mf Z}(M ):} &  0 \ar[r]  & M \ar@{=}[d] \ar[rr]^{\begin{pmatrix}
                                        \varphi-1\\    {\mf Z}\end{pmatrix}
    }  & & {M}\oplus {M}\ar[d]^(0.35){-\Psi\oplus \id }  \ar[rr]^(0.55){\begin{pmatrix}
                          {\mf Z}& 1-\varphi  \\
                               \end{pmatrix}
      }    &&   {M} \ar[d]^{-\Psi}  \ar[r]^{ } & 0 \\
 {K_{\Psi,{\mf Z} }(M ):} &  0 \ar[r]  & M    \ar[rr]^{\begin{pmatrix}
                                        \Psi-1\\    {\mf Z} \end{pmatrix}
    } & & M\oplus M   \ar[rr]^(0.55){\begin{pmatrix}
                          {\mf Z}  & 1-\Psi  \\
                               \end{pmatrix}
      }  &&   M   \ar[r]^{ } & 0. }
 \end{equation}
 In particular, we obtain an isomorphism $\Upsilon'_M: H^1_{\varphi_L,\mf Z}(M)\cong H^1_{\Psi_L,\mf Z}(M)$ sending a class $[(x,y)]$ to the class $[(-\Psi(x),y)].$

Using this one derives from $\langle -, - \rangle_M$ in Theorem \ref{lem:explicitPairing} the (asymmetric) perfect pairings, denoted by $\Kl -, - \Kr_M,$
\[H^1_{\Psi, \mf Z}(M) \times H^1_{\varphi, {\mf Z} }(\tilde{M}) \to K, (\overline{(m,n)},\overline{(f,g)})\mapsto Res\Big( g( m)- (\lambda^\iota f)( n)\Big),\]
\[H^2_{\Psi, \mf Z}(M) \times H^0_{\varphi, {\mf Z} }(\tilde{M}) \to K, (\overline{m},\tilde{n})\mapsto  Res \Big( \tilde{n}(\lambda m)\Big) ,\]
\[H^0_{\Psi, \mf Z}(M) \times H^2_{\varphi, {\mf Z} }(\tilde{M}) \to K, ({m},\overline{\tilde{n}})\mapsto Res \Big( \tilde{n}(m)\Big),\]
for which by construction we have
\[\langle x,y \rangle_M=\Kl \Upsilon'_M(x), y \Kr_{M}.\]
Moreover, we obtain, for $x\in H^i_{\varphi, {\mf Z} }({M}),\ y\in H^{2-i}_{\varphi, {\mf Z} }(\tilde{M}),$
\begin{align}
\label{f:Upsilon'pairings} (-1)^i\langle x,y \rangle_M= \langle y,x^{**} \rangle_{\tilde{M}} =\Kl \Upsilon'_{\tilde{M}}(y), x^{**} \Kr_{\tilde{M}},
\end{align}
by Remark \ref{rem:dualePaarung}.
\end{remark}

\begin{proof}[Proof of the Theorem]
(i) is an immediate consequence Lemma \ref{lem:Zinvertible}. Now consider (ii): The first isomorphism is induced by (i).
 Up to signs, $(-)^*$ transforms $\varphi_L$ into $\psi_L$ and $\mf Z$ into $\mf Z^\iota$. Using that $\mf Z^\iota=\lambda \mf Z$   one easily verifies that also the second map is an  isomorphism. Finally, the last isomorphism stems from the identification $M^*\cong\tilde{M}$ by \cite[Cor.\ 4.5.4]{SV20}.

For the pairing on the level of cohomology groups, we want to apply (ii) of Lemma \ref{lem:strict}, for which we have to check strictness of the differentials. But this is not sufficient: in order to get perfectness of the pairings - which amounts to an algebraic duality while the functor $(-)^*$ only measures continuous duals - we also have to check that the induced topology on the cohomology groups is Hausdorff. In detail this boils down to the following reasoning: Since   by \ref{thm:CohomologyFinite} all the $H^i_{\varphi_L,\mf Z}(M)$ are finite-dimensional, we may apply Lemma \ref{lem:strict}(iii) to first conclude that $d^1$   (and trivially $d^2$) is strict. By the same reasoning for $H^i_{\varphi,\mathfrak{Z}}(\tilde{M})$ the $d^1$-differential of $K_{\Psi,{\mathfrak{Z}^\iota}}(M^*)$ is strict. Moreover, the $H^2$s are always Hausdorff by \ref{lem:strict} (iv) and we note that the $H^0$ are always Hausdorff (as they are subspaces of Hausdorff spaces). Applying \ref{lem:strict} (ii) and using that for a finite dimensional Hausdorff space the continuous and algebraic dual agree we conclude the claim for the pairings involving $H^0$ and $H^2$.
 By the strictness of $d_1$ we have $H^1_{\varphi,\mathfrak{Z}}(M)^* \cong H^1_{\varphi,\mathfrak{Z}}(\tilde{M})$ and, vice versa,  $H^1_{\varphi,\mathfrak{Z}}(\tilde{M})^* \cong H^1_{\varphi,\mathfrak{Z}}(M).$ A priori we don't know if the finite dimensional $H^1$s are Hausdorff but combining both isomorphisms we see that $(H^1_{\varphi,\mathfrak{Z}}(M)^*)^*$ has the same dimension as $H^1_{\varphi,\mathfrak{Z}}(M)$ which for a finite dimensional space can only occur, if every functional is continuous, forcing the $H^1$s to be Hausdorff, which allows us to argue analogously for the pairing of $H^1$s.
\end{proof}

\subsection{Cohomological computations in the character case}
Recall \cite[Lem.\ 4.6]{ST2} or \cite[\S 2]{Co2}
for the following.
The \textbf{Amice-Katz transform} is the  map
\[A_{-}:D(o_L,K)\to \cR_K^+,\] sending a distribution $\mu$ to
\[A_\mu(Z)=\int_{o_L}\eta(x,Z)\mu(x),\] satisfying:
\begin{enumerate}
	\item $A_-$ is a $\varphi$- and $\Gamma_L$-equivariant topological isomorphism of rings.
	\item for $z\in o_K$ with $v_p(z)>0$: $A_{\eta(x,z)\mu}(Z)=A_\mu(Z+_{LT}z),$ where $\int_{o_L}g(x)(f\cdot\mu)(x)=\int_{o_L}f(x)g(x)\mu(x)$ for any locally analytic function $f:o_L\to \mathbb{C}_p.$
	\item (multiplicativity regarding convolution) $A_{\lambda*\mu}=A_\lambda\cdot A_\mu$
	\item $A_{\mathrm{Res}_{b+\pi_L^no_L}(\mu)}=\frac{1 }{q^n}\sum_{[\pi_L^n](a)=0}\eta(-b,a)A_\mu(Z+_{LT}a)=\mathrm{Res}_{b+\pi_L^no_L}A_\mu,$ where the latter denotes the multiplication with the corresponding characteristic function.
	\item $\partial A_\mu=A_{{\Omega}x\mu}$ where $\partial=\frac{d}{dt_{LT}}=\frac{1}{\log_{LT}'}\frac{d}{dZ}=\Omega\eta(1,Z)\frac{d}{d\eta(1,Z)}.$\footnote{Note that $\partial\eta(x,Z)=x\Omega\eta(x,Z).$}
	\item $A_{d\mu}=t_{LT}A_\mu,$ where $\int_{o_L}f(x)(d\mu)(x)=\int_{o_L}f'(x)\mu(x)$ with $f'(x)=\frac{d}{dx}f.$
\end{enumerate}



\begin{lemma}\label{lem:Mellin-Amice}(Mellin transform) The natural inclusion $D(o_L^\times,K)\hookrightarrow D(o_L,K)$ combined with the Fourier isomorphism induces the map
\begin{align*}
  D(o_L^\times,K) & \xrightarrow{\; \cong \;} D(o_L,K)^{\psi_L^D=0}\cong \mathcal{O}_K(\mathfrak{X})^{{\psi_L^\mathfrak{X}} = 0} \\
  \lambda & \longmapsto \phantom{mmmmll} \lambda(\delta_1)\corresponds\lambda(\ev_1)
\end{align*}
which is a topological isomorphism of $D(o_L^\times,K)$-modules. Here $\ev_1$ denotes the map on the character variety which evaluates a character in $1$. Moreover, we have a commutative diagram
	\[\xymatrix{
		D(o_L^\times,K) \ar[rr]_{\subseteq} \ar[d]_{\mathfrak{M}}^{\cong}
		&  &    D(o_L,K) \ar[d]^{A_-} \ar@/_/[ll]_{\mathrm{Res}_{o_L^\times}}   \\
		(\cR_K^+ )^{\Psi=0}     \ar[rr]_{\subseteq}   &  &   \cR_K^+    \ar@/_/[ll]_{1-\varphi\circ\Psi}   }
	\] where $\mathfrak{M}$ denotes the Mellin transform, which by definition sends $\mu$ to \[\mu\cdot \eta(1,Z)=\int_{o_L^\times} \eta(x,Z)\mu(x),\] see \cite[\S 2.1.4, Lem.\ 2.6, Thm.\ 2.33,\S 2.2.7]{SV20}.
\end{lemma}
\begin{proof} $\mu\in D(o_L^\times)\subseteq D(o_L)$ satisfies $\mathrm{Res}_{o_L^\times}(\mu)=\mu,$ whence
	$A_\mu(Z)=\int_{o_L}\eta(x,Z)\mu(x)=\int_{o_L^\times}\eta(x,Z)\mu(x)=\mathfrak{M}(\mu).$
\end{proof}

We write $LA(o_L):=LA(o_L,K)$ for the set of locally $L$-analytic functions $\phi:o_L\to K$ endowed with the following operators:
\begin{align*}
	\varphi(\phi)(x):=&\left\{
	\begin{array}{ll}
		\phi(\frac{x}{\pi_L}), & \hbox{if $x \in \pi_Lo_L;$} \\
		0, & \hbox{otherwise.}
	\end{array}
	\right. \\
	\Psi(\phi)(x):= & \phi(\pi_Lx) \\
	\gamma(\phi)(x):= & \phi(\chi_{LT}^{-1}(\gamma)x).
\end{align*}

By \cite[Thm.\ 2.3]{Co} (for the exact sequence), \cite[Cor.\ 2.3.4]{BF} (for the surjectivity on $\cR_K^+(\delta)$), we have for all  $\delta \in \Sigma_{an} $ the following commutative diagram of $D(\Gamma_L,L)$-modules with exact rows
\begin{equation}\label{f:Colmeztransform}
	\xymatrix{
		0   \ar[r]^{ } & \cR_K^+(\delta) \ar@{->>}[d]_{\Psi-1} \ar[r]^{ } & \cR_K(\delta) \ar[d]_{\Psi-1} \ar[r]^{ } & LA(o_L)(\chi^{-1}\delta) \ar[d]_{\Psi-1} \ar[r]^{ } & 0  \\
		0   \ar[r]^{ } & \cR_K^+(\delta)   \ar[r]^{ } & \cR_K(\delta) \ar[r]^{ } & LA(o_L)(\chi^{-1}\delta)   \ar[r]^{ } & 0  }
\end{equation}
which we can also interpret as short exact sequence of complexes of $D(\Gamma_L,L)$-modules
\[\xymatrix@C=0.5cm{
	0 \ar[r] & \cT_\Psi(\cR_K^+(\delta)) \ar[rr]^{ } &&  \cT_\Psi(\cR_K(\delta)) \ar[rr]^{ } &&  \cT_\Psi( LA(o_L)(\chi^{-1}\delta)) \ar[r] & 0. }\]
with $\cT_\Psi(\cR_K^+(\delta))\cong (\cR_K^+(\delta))^{\Psi=1}[0]$ in degree zero. Here the map $\cR_K(\delta)\to LA(o_L)(\chi^{-1}\delta)$ sends $f\mathbf{e}_\delta$ to $\phi_f \mathbf{e}_{\chi^{-1}\delta}$ with\footnote{Our map is $\frac{1}{\Omega}$ times Colmez' one.}
\begin{equation}\label{f:phif}
	\phi_f(z):=Res(\eta(-z,Z)fdt_{LT})=Res(\eta(-z,Z)f(Z)g_{LT}(Z)dZ).
\end{equation}
In particular we obtain a short exact sequence
\begin{align}\label{f:psiinv}
	& \xymatrix@C=0.5cm{
		0 \ar[r] & \cR_K^+(\delta)^{\Psi=1} \ar[rr]^{ } && \cR_K(\delta)^{\Psi=1} \ar[rr]^{} && LA(o_L)(\chi^{-1}\delta)^{\Psi=1} \ar[r] & 0 }
\end{align}
and an isomorphism
\begin{align}\label{f:psicoinv}
	&  \cR_K(\delta)/{\Psi-1}\cong LA(o_L)(\chi^{-1}\delta)/{\Psi-1}.
\end{align}


Let $Pol_{\leq N}(o_L):=Pol_{\leq N}(o_L,K):=\bigoplus_{i=0}^NKz^i\subseteq LA(o_L)$ denote the polynomial functions on $o_L$. This subspace is $\Gamma_L$- and $\Psi$-stable, more precisely we have
\begin{align*}
	\Psi(z^i)&= \pi_L^i z^i  \\
	\gamma(z^i) &= \chi_{LT}^{-i}z^i.
\end{align*}
for all $i\geq 0$ and $\gamma\in\Gamma_L.$ In particular, we obtain, for $i=0,1$,
\begin{align}\label{f:PolPsi}
	H^i_\Psi(Pol_{\leq N}(o_L)(\chi^{-1}\delta))\cong \left\{
	\begin{array}{ll}
		K z^k\mathbf{e}_{\delta\chi^{-1}}, & \hbox{if $\delta(\pi_L)=\frac{\pi_L^{k+1}}{q}$ for some $0\leq k\leq N$;} \\
		0, & \hbox{otherwise.}
	\end{array}
	\right.
\end{align}
It follows that
\begin{align}\label{f:PolZPsi}
	H^j_{\mf Z}( H^i_\Psi(Pol_{\leq N}(o_L)(\chi^{-1}\delta)))\cong \left\{
	\begin{array}{ll}
		K z^k\mathbf{e}_{\delta\chi^{-1}}, & \hbox{if $\delta =x^k\chi$ for some $0\leq k\leq N$;} \\
		0, & \hbox{otherwise.}
	\end{array}
	\right.
\end{align}

\begin{lemma}\label{f:Chenevier} For $N>v_\pi(\chi^{-1}\delta(\pi))$ we have a quasi-isomorphism
	\[\cT_\Psi( LA(o_L)(\chi^{-1}\delta))\simeq \cT_\Psi( Pol_{\leq N}(o_L)(\chi^{-1}\delta))\] and an isomorphism \[Pol_{\leq N}(o_L)(\chi^{-1}\delta)^{\Psi=1}\cong Pol_{\leq N}(o_L)(\chi^{-1}\delta)/(\Psi-1)\] as $L$-vector spaces.
\end{lemma}

\begin{proof} (see \cite[Lem.\ 2.9]{chen} for the cyclotomic case, even over affinoid algebras $A$ instead of $L$).
	Use the decomposition $LA(o_L)\cong x^{N+1}LA(o_L) \oplus Pol_{\leq N}(o_L)$ and show that for $N$ as in  the assumption $\Psi-1$ is a topological  isomorphism on $x^{N+1}LA(o_L).$
\end{proof}

Similarly, regarding the $\Gamma_L$- and $\Psi$-stable submodule   $D_{N}:=D_{K,N}:=\bigoplus_{l=0}^NKt_{LT}^l\subseteq \cR_K^+$ we obtain
for $i,j\in\{0,1\}$,
\begin{align}\label{f:DPsi}
	H^i_\Psi(D_{N}(\delta))\cong \left\{
	\begin{array}{ll}
		K t_{LT}^k\mathbf{e}_{\delta}, & \hbox{if $\delta(\pi_L)= \pi_L^{-k} $ for some $0\leq k\leq N$;} \\
		0, & \hbox{otherwise,}
	\end{array}
	\right.
\end{align}
and
\begin{align}\label{f:DZPsi}
	H^j_{\mf Z}( H^i_\Psi(D_{N}(\delta)))\cong \left\{
	\begin{array}{ll}
		K t_{LT}^k\mathbf{e}_{\delta}, & \hbox{if $\delta =x^{-k}$ for some $0\leq k\leq N$;} \\
		0, & \hbox{otherwise.}
	\end{array}
	\right.
\end{align}

\begin{remark}\label{rem:compcoh}
Note that, by the same reasoning, the analogue of Lemma \ref{lem:compcoh} (ii) (but in general not (iii)) does also hold for $M$ of the form $ \cR_K^+(\delta)$ or $ LA(o_L)(\delta)$.
\end{remark}

Recall that $\Sigma_1=\{x^{-i}| i\in\mathbb{N}\},$ $\Sigma_2=\{x^{i}\chi| i\in\mathbb{N}\}$ and $\Sigma_{gen}=\Sigma_{an}\setminus (\Sigma_1\cup\Sigma_2).$

\begin{lemma}\label{lem:dim} The dimensions of the analytic cohomology groups are as follows:
	\begin{enumerate}
		\item $\dim_K H^j_{\varphi,D(\Gamma_L,K)}(\cR_K^+(\delta ))=\left\{
		\begin{array}{ll}
			0, & \hbox{$\delta\not\in\Sigma_1$;} \\
			1, & \hbox{$\delta\in\Sigma_1$, $j=0$;} \\
			2, & \hbox{$\delta\in\Sigma_1$, $j=1$;} \\
			1, & \hbox{$\delta\in\Sigma_1$, $j=2$.}
		\end{array}
		\right.$
		\item For $\delta^{-1}\not\in\Sigma_1$ we have $\dim_K H^j_{\varphi,D(\Gamma_L,K)}( LA(o_L)(\delta ))=\left\{
		\begin{array}{ll}
			0, & \hbox{$j=0$;} \\
			1, & \hbox{$j=1$;} \\
			0, & \hbox{$j=2$.}
		\end{array}
		\right.$
		\item For $\delta^{-1}\in\Sigma_1$ we have $\dim_K H^j_{\varphi,D(\Gamma_L,K)}( LA(o_L)(\delta ))=\left\{
		\begin{array}{ll}
			0, & \hbox{$j=0$;} \\
			2, & \hbox{$j=1$;} \\
			1, & \hbox{$j=2$.}
		\end{array}
		\right.$
		\item For $\delta\in \Sigma_1$ we have $\dim_K H^j_{\varphi,D(\Gamma_L,K)}(\cR_K(\delta ))=\left\{
		\begin{array}{ll}
			1, & \hbox{$j=0$;} \\
			2, & \hbox{$j=1$;} \\
			0, & \hbox{$j=2.$}
		\end{array}
		\right.$
		\item For $\delta\in \Sigma_2$ we have $\dim_K H^j_{\varphi,D(\Gamma_L,K)}(\cR_K(\delta ))=\left\{
		\begin{array}{ll}
			0, & \hbox{$j=0$;} \\
			2, & \hbox{$j=1$;} \\
			1, & \hbox{$j=2.$}
		\end{array}
		\right.$
		\item For $\delta\in \Sigma_{gen}$ we have $\dim_K H^j_{\varphi,D(\Gamma_L,K)}(\cR_K(\delta ))=\left\{
		\begin{array}{ll}
			0, & \hbox{$j=0$;} \\
			1, & \hbox{$j=1$;} \\
			0, & \hbox{$j=2.$}
		\end{array}
		\right.$
	\end{enumerate}
In particular, generic characters are precisely those with vanishing $H^0$ and $H^2.$
\end{lemma}
\begin{proof}
	By Remark \ref{rem:analCoh}  $H^\bullet_{an}(A^+,M)$ in \cite[\S 5]{Co2} coincides with $H^\bullet_{\varphi,D(\Gamma_L,K)}(M).$ Note that Colmez uses $L$ to denote a large field such as our field $K.$
\end{proof}

 It is easy to check that analogous results as in this subsection hold for modules of the form $\cR_A(\delta)$ for affinoids $A$ over K  instead of the base field $K$. The only subtlety is the appearance of non-trivial zero divisors. By imposing some additional conditions we can strengthen \ref{thm:perf} to cover the affinoid case as well.

 Note that the action of $\Gamma_L$ on $A(\delta\chi_{LT}^{\pm i})$ extends to an action of $D(\Gamma_L,K)$ by continuity. The element $\mathfrak{Z}$ acts as an $A$-linear endomorphism on $A(\delta\chi_{LT}^{\pm i})$ hence by multiplication with an element $\mathfrak{Z}(\delta\chi_{LT}^{\pm i}) \in A.$
\begin{remark}\label{rem:perf-delta}\phantom{section}
\begin{enumerate}
\item 	Let $A$ be affinoid over $K$ and let $\delta \colon L^\times \to A^\times$ be a locally $L$-analytic character. Assume that $1-\delta(\pi)\pi^i$ is not a non-trivial zero divisor in $A$ for every $i \in \ZZ$ and assume that (the image of) $\mathfrak{Z}(\delta\chi_{LT}^{\pm i}) \in A$
is not a non-trivial zero divisor in $A$ or any $A/(1-\delta(\pi)\pi^i)$\footnote{If we drop the zero divisor assumption the same proof would show that the complexes lie in $D_{perf}^-.$ If $A$ is a domain and $\delta(\pi) \in K^\times$ then $(1-\delta(\pi)\pi^i)$ is either $0$ or a unit and hence the condition on $\delta(\mathfrak{Z})$ is automatically satisfied!}. Then $\mathcal{T}_{\Psi}(M)$ is perfect as a $D(\Gamma_L,A)$-module for $M$ in \[\{\cR_A^{+}(\delta), \cR_A(\delta), LA(o_L,A)(\delta\chi^{-1}),  D_{A,N}(\delta), Pol_{\leq N}(o_L,A)(\chi^{-1}\delta) \}.\]
\item  As in Nakamura's setting we expect the statement of (i) to be true without any condition. Unfortunately, the methods of \cite[Section 5]{KPX} do not transfer to our situation directly due to the fact that \cite{KPX} makes use of the Euler characteristic formula and perfectness of the $\Psi$-complex in the \'{e}tale case. The analogues of these results are not known to us for analytic cohomology over affinoids.
\end{enumerate}
\end{remark}
\begin{proof}
 First observe that for any locally analytic character $\rho\colon L^\times \to A^\times$ the free rank one module $A(\rho)$ is perfect as a $D(U,A)$-module if $\mathfrak{Z}(\rho)$ is not a non-trivial zero divisor in $A.$ Indeed, let $\alpha:=\mathfrak{Z}(\rho) \in A.$ Then, using the assumptions on $\alpha,$ one sees that $A(\rho) \cong D(U,A)/(\mathfrak{Z}-\alpha)D(U,A)$ is perfect as a $D(U,A)$-module but then also perfect as a $D(\Gamma_L,A)$-module by \cite[Lemma 2.5]{Ste2}. Let us call a module of the form $A(\rho)$ of type $\mathcal{F}$. Now consider the sequence
	\be \label{eq:seqpsi} 0 \to \cR_A^+(\delta) \to \cR_A(\delta) \to LA(o_L,A)(\chi^{-1}\delta)\to 0.\ee We have that $\mathcal{T}_{\Psi}(LA(o_L,A)(\chi^{-1}\delta))$ is perfect by \cite[\href{https://stacks.math.columbia.edu/tag/066T}{Tag 066T}]{SP} since the inclusion of the $Pol_{\leq N}(o_L,A)(\chi^{-1}\delta)$ induces a quasi-isomorphism for $N \gg 0$ to a complex whose terms are perfect as they are finite direct sums of modules of type $\mathcal{F}$. Similarly for $D_{A,N}(\delta).$ To see that $\cR_A^+(\delta)^{\Psi=1}$ is perfect, consider the exact sequence
	$$0 \to V_1 \to (\cR_A(\delta)^+)^{\Psi=1} \xrightarrow{\varphi-1} (\cR_A(\delta)^{+})^{\Psi=0} \to V_2\to 0,$$ where $V_1,V_2$ are defined as kernel and co-kernel of the middle map. By an adaptation of \cite[Lemma 2.9 and Proposition 2.20]{chen} to our situation the kernel is of type $\mathcal{F}$ over $A$ while the cokernel is a finite direct sum of modules of type $\mathcal{F}$ over $A/(1-\delta(\pi)\pi^i)$ with varying $i.$ It suffices to see that they are perfect as $D(\Gamma_L,A)$-modules. This follows from the assumption that $(1-\delta(\pi)\pi^i)$ is not a zero divisor and hence $D(\Gamma_L,A/(1-\delta(\pi)\pi^i))$ is itself perfect as a $D(\Gamma_L,A)$-module.  It remains to see that $\cR_A^+(\delta)/(\Psi-1)$ is perfect. Again by a similar argument it is a finite direct sum of perfect $D(\Gamma_L,A/(1-\delta(\pi)^{-1}\pi^i))$ modules (the appearance of $\delta(\pi)^{-1}$ is due to using $\Psi-1$ instead of $\varphi-1$). Our assumptions ensure that $1-\delta(\pi)^{-1}\pi^i = (-\delta(\pi)^{-1}\pi^i)(1-\delta(\pi)\pi^{-i})$ is not a zero divisor and we can proceed as before. This proves the perfectness of $\mathcal{T}_{\Psi}(\cR_A(\delta)^+).$
Finally the perfectness of $\mathcal{T}_{\Psi}(\cR_A(\delta))$ follows from the exact sequence \eqref{eq:seqpsi}.
\end{proof}

\section{Bloch–Kato exponential for analytic \texorpdfstring{$(\varphi_L,\Gamma_L)$}{(phi,Gamma)}-modules}

 \label{sec:BKexponential}

Recall that for a module $M$ over a commutative
ring $R$ and $t\in R$ a non-zerodivisor, we
use the notation $M_t$ for the localisation of $M$
at the multiplicatively closed subset $\{t^n\}_{n\geq 0}$.

\subsection{$D_{dR}$ and $D_{cris}$ for analytic \texorpdfstring{$(\varphi_L,\Gamma_L)$}{(phi,Gamma)}-modules}
 In this section we will define versions of $D_{dR}$ and $D_{cris}$ for $L$-analytic $(\varphi_L, \Gamma_L)$-modules $M$. The idea is that, for an \'{e}tale  $(\varphi_L, \Gamma_L)$-module attached to a representation $V$, these versions correspond to the  identity component of the full $D_{dR}(V)$, which  arise as $(B_{dR}\otimes_L V)^{G_L}$ instead of $(B_{dR}\otimes_{\mathbb{Q}_p}  V)^{G_L}$, and similarly for $D_{cris}(V)$. The comparison between the definitions used in this article and Fontaine's classical ones is described in \cite[Section 5.2]{Por}. \\

\begin{definition}
For an $L$-analytic $(\varphi_L, \Gamma_L)$-module $M$ over $\cR_F$, we define
\begin{align*}
\bD_{\tn{dR}}(M):=\bD_{\textnormal{dif}}^{}(M)^{\Gamma_L}
\end{align*}
and
\begin{align*}
\bD_{\tn{cris}}(M):=M_{t_{LT}}^{\Gamma_L}.
\end{align*}
\end{definition}

\begin{remark}
	Let $M$ be an $L$-analytic $(\varphi_L, \Gamma_L)$-module $M$ over $\cR_L.$
	Then $\bD_{\mathrm{dR}}(M)$ and $\bD_{\mathrm{cris}}(M)$ are finite dimensional $L$-vector spaces of dimension $\leq\rk(M)$. Furthermore $\varphi_M$ induces an automorphism of $\bD_{\mathrm{cris}}(M).$
\end{remark}
\begin{proof} We first show that $\bD_{\mathrm{dR}}(M)$ is finite dimensional. By construction $D:=\bD_{\textnormal{dif}}(M)$ is a finite-dimensional $B:=\bigcup_{n\geq0}L_{n}((t_{LT}))$-semilinear representation of $\Gamma_L.$ We claim that the natural map $$B \otimes_{B^{\Gamma_L}} D^{\Gamma_L} \to D$$ is injective and $B^{\Gamma_L}=L,$ which shows $\dim_L(D^{\Gamma_L}) \leq \dim_B(D) =\rk(M)$.
We first show $B^{\Gamma_L}=L.$ Let $f = \sum a_i t_{LT}^i \in B^{\Gamma_L}.$ We conclude $a_i= \chi_{LT}(\gamma)^i\gamma(a_i)$ for every $\gamma \in \Gamma_L.$ Let $n$ be large enough such that all $a_i$ belong to $L_n.$ Then $\gamma(a_i)=a_i$ for every $\gamma \in \Gamma_n$ and we conclude that $a_i=0$ holds for every $i \neq 0.$ Finally $\gamma(a_0) \in L_n^{\Gamma_L} = L$ which proves the claim. For the injectivity we argue like in  the proof of 2.13 in \cite{FO}. Consider $L$-linearly independent vectors $v_1,\dots,v_d\in D^{\Gamma_L}$ such that $$\sum_{i=1}^d \lambda_iv_i=0$$ with some $\lambda_i \in B.$ Suppose $d \geq 2,$ $\lambda_1 \neq 0 $ and assume without loss of generality $\lambda_1=1.$ We obtain $v_1 = \gamma(v_1) = \sum_{i=2}^d- \gamma(\lambda_i)v_i.$ Arguing by induction we may assume that $v_2,\dots,v_d$ are linearly independent over $B$ and conclude $\lambda_2,\dots,\lambda_d \in B^{\Gamma_L}=L$, a contradiction. From the injectivity of $\iota_n$ according to Remark \ref{iota t} we deduce that $\dim_L(\bD_{\mathrm{cris}}(M)) \leq \dim_L(\bD_{\mathrm{dR}}(M)).$ Finally $\varphi_M$ induces an injective endomorphism of $M_{t_{LT}}$ and by a dimension argument an automorphism of $\bD_{\mathrm{cris}}(M).$
\end{proof}

\subsection{$\exp$   for analytic \texorpdfstring{$(\varphi_L,\Gamma_L)$}{(phi,Gamma)}-modules}\label{sec:exp}

Let $M$ be an $L$-analytic $(\varphi_L, \Gamma_L)$-module over $\cR_K$. By Prop. \ref{ana:to:distr}, we have an action of $D(\Gamma,K)$ on $M$. Thus we may (for some fixed n) consider the complex $K_{\varphi,\mathfrak Z_n}(M)$, which (up to sign) amounts to
\begin{align*}
K_{\varphi,\mathfrak Z_n}(M)=[M \xrightarrow{(\varphi-1, \mathfrak Z_n)}  M \oplus M \xrightarrow{(\mathfrak Z_n \oplus 1-\varphi)} M]
\end{align*}
concentrated in degree $[0,2]$. \\
On the other hand, for any $D(\Gamma,K)$-module $N$, we define
\begin{align*}
K_{\mathfrak Z_n}(N):=[N \xrightarrow{ \mathfrak Z_n}  N]
\end{align*}
concentrated in degree $[0,1]$, and denote its cohomology by $H^\bullet_{\mf Z_n}(N)$. \\
Next we want to define $K^{(\varphi)}_{\varphi,\mathfrak Z_n}(M_0)$ and $K^{(\varphi)}_{\mathfrak Z_n}(M_0)$ for $M_0 \in \{M, M_{t_{LT}}\}$. By inspecting the proof in the reference for Prop. \ref{ana:to:distr}, one sees that the action of $D(\Gamma,K)$ on $M$ preserves all the $M^{(m)}$. For $m \gg 0$, we set
\begin{align*}
\widetilde{K}_{\varphi,\mathfrak Z_n}(M_0^{(m)}):=[M_0^{(m)} \xrightarrow{(\varphi-1, \mathfrak Z_n)}  M_0^{(m+1)} \oplus M_0^{(m)} \xrightarrow{(\mathfrak Z_n \oplus 1-\varphi)} M_0^{(m+1)}]
\end{align*}
concentrated in degree $[0,2]$. Passing to the limit with respect to the transition maps induced by the canonical inclusions $M_0^{(m)} \nachinj M_0^{(m+1)}$ recovers $K_{\varphi,\mathfrak Z_n}(M_0)$, but taking the limit with respect to the transition maps induced by $\varphi \colon M_0^{(m)} \nach M_0^{(m+1)}$ produces a new complex
\begin{align*}
K^{(\varphi)}_{\varphi,\mathfrak Z_n}(M_0):= \varinjlim_{m, \varphi} \widetilde{K}_{\varphi,\mathfrak Z_n}(M_0^{(m)})
\end{align*}
whose cohomology we denote by $H^{(\varphi),\bullet}_{\varphi,\mathfrak Z_n}(M_0)$. Similarly we define
\begin{align*}
K^{(\varphi)}_{\mathfrak Z_n}(M_0):= \varinjlim_{m, \varphi} {K}_{\mathfrak Z_n}(M_0^{(m)})
\end{align*}
with cohomology groups denoted by $H^{(\varphi),\bullet}_{\mathfrak Z_n}(M_0)$.
\begin{remark}
Note that we have
\begin{align*}
\bD_{\mathrm{dR}}(M)=H^0_{\mf Z_n}(\bD_{\mathrm{dif}}^{}(M))^{\Gamma/\Gamma_n}
\end{align*}
and
\begin{align*}
\bD_{\tn{cris}}(M)=H^0_{\mf Z_n}(M_{t_{LT}})^{\Gamma/\Gamma_n}.
\end{align*}
\end{remark}
\begin{lemma}\label{2.17Nak}
For $m \gg 0$ and $M_0 \in \{M, M[1/t_{LT}]\}$, the following natural maps induced by $\varphi$ are quasi-isomorphisms:
\begin{align*}
&K_{\mf Z_n}(\bD_{\textnormal{dif},m}^{(+)}(M)) \nach K_{\mf Z_n}(\bD_{\textnormal{dif},m+1}^{(+)}(M)),\\
&K_{\mf Z_n}(M_0^{(m)}) \nach K_{\mf Z_n}(M_0^{(m+1)}) \tn{ and } \\
&\widetilde{K}_{\varphi,\mathfrak Z_n}(M_0^{(m)}) \nach \widetilde{K}_{\varphi,\mathfrak Z_n}(M_0^{(m+1)})
\end{align*}
 In particular, the maps
\begin{align*}
&K_{\mf Z_n}(\bD_{\textnormal{dif},m}^{(+)}(M)) \nach K_{\mf Z_n}(\bD_{\textnormal{dif}}^{(+)}(M)), \\
&K_{\mf Z_n}(M_0^{(m)}) \nach K^{(\varphi)}_{\mathfrak Z_n}(M_0) \tn{ and } \\
&{\widetilde K}_{\varphi,\mathfrak Z_n}(M_0^{(m)}) \nach {
 K}^{(\varphi)}_{\varphi,\mathfrak Z_n}(M_0)
\end{align*}
are quasi-isomorphisms.
\end{lemma}
\begin{proof}
We only need prove the first statement.By Lemma \ref{lem:Zinvertible} the action of $\mathfrak{Z}_n$ on $M_0^{\psi=0}$ (resp. $(M_0^{(m)})^{\psi=0}$) is invertible. Using this fact, one can conclude the proof with the same (purely formal) arguments as in the proof of \cite[Lemma 2.17]{NaANT}.
\end{proof}
\begin{lemma}\phantomsection\label{2.20Nak}
\begin{enumerate}[(i)]
\item For $m \gg 0$ and $M_0 \in \{M, M_{t_{LT}}\}$, the map
\begin{align*}
\widetilde{K}_{\varphi,\mathfrak Z_n}(M_0^{(m)}) \nach {K}_{\varphi,\mathfrak Z_n}(M_0^{})
\end{align*}
induced by the inclusion $M_0^{(m)} \nachinj M_0$ is a quasi-isomorphism.
\item In $\mathbf D^-(K)$, by composing the inverse of the isomorphism in (i) with the isomorphism $\widetilde{K}_{\varphi,\mathfrak Z_n}(M_0^{(m)}) \nach {K}^{(\varphi)}_{\varphi,\mathfrak Z_n}(M_0)$ from Lemma \ref{2.17Nak}, one obtains an isomorphism
\begin{align*}
{K}_{\varphi,\mathfrak Z_n}(M_0^{}) \iso {K}^{(\varphi)}_{\varphi,\mathfrak Z_n}(M_0)
\end{align*}
which is independent of the choice of $m \gg 0$.
\end{enumerate}
\end{lemma}
\begin{proof}
Both statements follow by purely formal arguments from Lemma \ref{2.17Nak}, just as in the proof of \cite[Lemma 2.20]{NaANT}.
\end{proof}

\begin{definition}\label{DefPfeileVorPROP21}
\begin{enumerate}[(a)]
\item By the compatibility of the maps $\iota_m$
with $\varphi\colon M^{(m)}\nach M^{(m+1)}$ and the
inclusions $\bD_{\mathrm{dif},m}^{(+)}(M)\nachinj
\bD_{\mathrm{dif},m+1}^{(+)}(M)$ as in Definition \ref{Ddif.Def}, the $\iota_m$ induce canonical morphisms\footnote{Note that obviously we have $K_{\mf Z_n}(\bD_{\textnormal{dif}}^{(+)}(M))\cong\varinjlim_m
K_{\mf Z_n}(\bD_{\textnormal{dif},m}^{(+)}(M))$.}
\begin{align*}
 K_{\mathfrak{Z}_n}^{(\varphi)}(M)\nach K_{\mathfrak{Z}_n}(\bD_{\mathrm{dif}}^+(M))\ \ \ \
 \ \textnormal{and}\ \ \ \ \
 K_{\mathfrak{Z}_n}^{(\varphi)}(M_{t_{LT}})\nach
 K_{\mathfrak{Z}_n}(\bD_{\mathrm{dif}}(M))
\end{align*}
which we will both call $\iota$. Moreover, the
inclusions $M_0^{(m)}\nachinj M_0^{(m+1)}$ induce
a map
\begin{align*}
 \Frob\colon K_{\mf Z_n}^{(\varphi)}(M_0)\nach
 K_{\mf Z_n}^{(\varphi)}(M_0).
\end{align*}
\item We construct morphisms
\begin{align*}
 f_n\colon K_{\varphi,\mathfrak{Z}_n}(M_0)\nach
 K_{\mathfrak{Z}_n}^{(\varphi)}(M_0)\ \ \ \ \
 \textnormal{and}\ \ \ \ \
 g_n\colon K_{\varphi,\mathfrak{Z}_n}(M_0)\nach
 K_{\mathfrak{Z}_n}(\mathbf D_{\mathrm{dif}}^{(+)}(M))
\end{align*}
in the following way:\\ Define $f_n$ as the composition
of the isomorphism ${K}_{\varphi,\mathfrak Z_n}(M_0^{}) \iso {K}^{(\varphi)}_{\varphi,\mathfrak Z_n}(M_0)$ from Lemma \ref{2.20Nak}(ii) with the
map $K_{\varphi,\mathfrak{Z}_n}^{(\varphi)}
(M_0)\nach K_{\mathfrak{Z}_n}^{(\varphi)}(M_0)$
obtained by taking the direct limit of the morphisms
\begin{align*}
 \begin{xy}
  \xymatrix{
    \widetilde{K}_{\varphi,\mathfrak{Z}_n}(M_0^{(
    m)}):\ar[d]& [M_0^{(m)} \ar[rr]^-{(\varphi-1, \mathfrak Z_n)}
    \ar[d]^{\id} && M_0^{(m+1)}\oplus
    M_0^{(m)}\ar[rr]^-{\mathfrak Z_n \oplus (1-\varphi)}\ar[d]^{(x,y)\auf y}&&M_0^{(m+1)}] \\
    K_{\mathfrak{Z}_n}(M_0^{(m)}): &
    [M_0^{(m)}\ar[rr]^{\mathfrak{Z}_n} && M_0^{(m)}] &&
  }
 \end{xy}
\end{align*}
Furthermore, the morphism $g_n$ is defined as
\begin{align*}
 g_n\colon K_{\varphi,\mf Z_n}(M_0)\overset {f_n}\nach
 K_{\mf Z_n}^{(\varphi)}(M_0)\overset \iota\nach
 K_{\mf Z_n}(\bD_{\mathrm{dif}}^{(+)}(M)).
\end{align*}
\end{enumerate}
\end{definition}

\begin{proposition}\label{Nak:Prop2.21}
 Consider the following diagram:
 \begin{align*}
 \begin{xy}
  \xymatrix{
   K_{\varphi,\mf Z_n}(M)\ar[dd]_\id\ar[r]^-{d_1}
   &K_{\varphi,\mf Z_n}(M_{t_{LT}})\oplus K_{\mf Z_n}(\bD_{\mathrm{dif}}^+(M))
   \ar[dd]_{f_n\oplus\id}\ar[r]^-{d_2}&
   K_{\mf Z_n}(\Ddif(M))\ar[dd]_{x\auf(0,x)}\ar[r]^-{
   +1}&
   \\ &&& \\
   K_{\varphi,\mf Z_n}(M)\ar[r]^-{d_3}&
   K_{\mf Z_n}^{(\varphi)}(M_{t_{LT}})\oplus
   K_{\mf Z_n}(\bD_{\mathrm{dif}}^+(M))\ar[r]^-{d_4}&
   K_{\mf Z_n}^{(\varphi)}(M_{t_{LT}})
   \oplus K_{\mf Z_n}(\Ddif(M))\ar[r]^-{+1}&
  }
 \end{xy}
 \end{align*}
 where the $d_i$ are given by
 \begin{align*}
  & d_1(x):=(x,g_n(x)), && d_2(x,y):=g_n(x)-y,\\
  & d_3(x):=(f_n(x),g_n(x)), && d_4(x,y):=
  (\Frob(x)-x,\iota(x)-y).
 \end{align*}
 Then the vertical map is a morphism between two distinguished triangles.
\end{proposition}
\begin{proof}
 The proof can be carried out analogously to the
 proof of
 \cite[Prop. 2.21]{NaANT}:
 We make use of the following well-known fact from
 homological algebra (see for instance
 \cite[Ex. 10.4.9]{Wei}):\\
 {\it Let $A$ be a ring and
 \begin{align*}
  0\nach X^\bullet\nach Y^\bullet\nach Z^\bullet\nach 0
 \end{align*}
 an exact sequence of complexes of $A$-modules. Then there
 exists a natural map $Z^\bullet\nach X^\bullet[1]$
 in the derived category $\mathbf{D}(A)$ such that
 \begin{align*}
  X^\bullet \nach Y^\bullet\nach Z^\bullet\overset{+1}{\nach
  } X^\bullet[1]
 \end{align*}
 is a distinguished triangle.}
 \\
 First, we show that the upper row is a distinguished triangle.
 Our goal is to replace the complexes
 $K_{\mf Z_n}(\bD_{\mathrm{dif}}^{(+)}(M))$
 by new, quasi-isomorphic complexes $\widetilde K_{\varphi,\mf Z_n}(\bD_{
 \mathrm{dif}}^{(+)}(M))$, which we define below, and
 construct an exact sequence
 \begin{align}
  0\nach K_{\varphi,\mf Z_n}(M)\nach K_{\varphi,\mf Z_n}(M_{t_{LT}})\oplus \widetilde K_{\varphi,\mf Z_n}(\bD_{\mathrm{dif}}^+(M))
   \nach
   \widetilde K_{\varphi,\mf Z_n}(\Ddif(M))\nach 0
    \label{2.21:lim7}
 \end{align}
 that will induce the upper triangle in the statement
 by the above-stated fact.
 For $k\geq 0$ and $m\gg 0$,
 we put
 \begin{align*}
  \PDif_{m,k}(M):=\prod_{\mu\geq m}t_{LT}^{-k}\cdot\bD_{\dif,\mu}^+(M)
 \end{align*}
 and denote by
 $\widetilde K_{\varphi,\mf Z_n}(t_{LT}^{-k}\cdot \bD_{\dif,m}^+
 (M))$ the complex concentrated in degree $[0,2]$:
 \begin{align*}
  \PDif_{m,k}(M)\overset{b_0}
  \nach \PDif_{m,k}(M)
  \oplus \PDif_{m+1,k}(M)
  \overset{b_1}\nach \PDif_{m+1,k}(M),
 \end{align*}
 where
 \begin{align*}
  b_0((x_\mu)_\mu):=((\mf Z_nx_\mu)_{\mu\geq m},
  (x_{\mu-1}-x_\mu)_{\mu\geq m+1}),
 \end{align*}
 and
 \begin{align*}
   b_1((x_\mu)_{\mu\geq m},(y_\mu)_{\mu\geq m+1}):=
  ((x_{\mu-1}-x_\mu)-\mf Z_n y_\mu)_{\mu\geq m+1}.
 \end{align*}
 Furthermore, let
 \begin{align*}
  \widetilde K_{\varphi,\mf Z_n}(\bD_{\dif,m}(M))
  :=\bigcup_{k\geq 0}\widetilde K_{\varphi,\mf Z_n}(t_{LT}^{
  -k}\bD_{\dif,m}^+(M)).
 \end{align*}
 We now define
 \begin{align}
  0\nach\widetilde K_{\varphi,\mf Z_n}(M^{(m)})
  \nach \widetilde K_{\varphi,\mf Z_n}(M^{(m)}_{t_{LT}})\oplus\widetilde K_{\varphi,\mf Z_n}(\bD^+_{\dif,m}(M))\nach
  \widetilde K_{\varphi,\mf Z_n}(\bD_{\dif,m}(M))
  \nach 0 \label{2.21:7}
 \end{align}
 as the sequence of complexes
 induced by applying $\widetilde K_{\varphi,\mf Z_n}(-)$ to
 \begin{align}
  0\nach M^{(m)}\overset{c_1}\nach M^{(m)}_{t_{LT}}
  \oplus\prod_{\mu\geq m}\bD_{\dif,\mu}^+(M)
  \overset{c_2}\nach \bigcup_{k\geq 0}
  \prod_{\mu\geq m}t_{LT}^{-k}\bD_{\dif,\mu}^+(M)
  \nach 0\label{2.21:5},
 \end{align}
 where
 \begin{align*}
  c_1(x):=(x,(\iota_\mu(x))_{\mu\geq m})\ \ \
  \text{and}\ \ \
  c_2(x,(y_\mu)_\mu):=(\iota_\mu(x)-y_\mu)_{\mu\geq m}.
 \end{align*}
 Down below, the sequence (\ref{2.21:lim7}) will be
 obtained as a direct limit of the sequences
 (\ref{2.21:7}).
 We claim that the sequence (\ref{2.21:5}) and hence also
 (\ref{2.21:7}) is exact. Consequently, the
 same will then hold
 for the direct limit (\ref{2.21:lim7}).

 The crucial part now is the exactness of (\ref{2.21:5}),
 which generalizes the exactness of the sequence
 (5) in the proof of \cite[Prop. 2.21]{NaANT}. The
 latter is demonstrated in \cite[Lem. 2.9]{Na}, and
 we check that the arguments carry over to our sequence
 (\ref{2.21:5}). The non-trivial statements are
 $\ker(c_2)=\im(c_1)$ and the surjectivity of $c_2$.

 The second statement can be reduced to showing that
 the map
 \begin{align}\label{c2-surj}
  M^{(m)}\nach \prod_{\mu\geq m}\bD_{\dif,\mu}^+(M)
  /t_{LT},\ x\auf (\overline{\iota_\mu(x)})_{\mu\geq
  m}
 \end{align}
 is surjective, using the fact that $M_{t_{LT}}^{(m)}
 =\bigcup_{k\geq 0}t_{LT}^{-k}M^{(m)}$ and reducing
 inductively to the case $k=1$ via d\'{e}vissage.
 Now we fix an $\cR_F^{(m)}$-basis $e_1,\ldots,e_d$
 of $M^{(m)}$,
 assuming $m$ is large enough for Proposition
 \ref{phiGammaModel} to hold. From \cite[Prop. 4.8 \& Lem. 4.9]{Be02} it follows
 that for any $\mu\geq m$, the composition
 \begin{align*}
  \cR_L^{(m)}\overset{\iota_\mu}{\nach} L_\mu[[t_{LT}]]
  \overset{t_{LT}\auf 0}{\nach} L_\mu
 \end{align*}
 induces an isomorphism $\cR_L^{(m)}/(Q_\mu)\cong
 L_\mu$ where
 $Q_\mu:=\frac{\varphi^\mu(Z)}{\varphi^{\mu-1}(Z)}$.\footnote{In case the underlying Lubin-Tate group law
 is special, then $Q_m$ is just the minimal
 polynomial of a uniformiser of $L_m/L$.}
 Therefore, using
 \cite[Lemma 4.3.6]{SV20}, we obtain on the
 level of the extension $F$ an isomorphism
 \begin{align*}
  \cR_F^{(m)}/Q_\mu\cong(\cR_L^{(m)}\widehat{\otimes}_L
  F)/Q_\mu\cong (\cR_L^{(m)}/Q_\mu)\widehat{\otimes}_L F
  \cong L_\mu\otimes_L F=F_\mu
 \end{align*}
 for $\mu\geq m$.
 As a result, we see that $(\overline{\iota_\mu(e_i)})_{i=1,\ldots,d}$ is an $F_\mu$-basis of
$\bD_{\dif,\mu}^+(M)
  /t_{LT}$ for any $\mu\geq m$. Now the surjectivity
  of (\ref{c2-surj}) is proven just as in
  \cite[Lem. 2.9]{Na}: For a family $(y_\mu)_{\mu
  \geq m}$ in the target, we write $y_\mu=
  \sum_{i=1}^da_{\mu,i}\cdot\overline{\iota_\mu(
  e_i)}$ for $\mu\geq m$. Choosing a representative
  $a_i\in\cR_F^{(m)}$ of the preimage of $(a_{\mu,i})_{\mu\geq m}$ under the natural isomorphism

  \begin{align}\label{prod iso}
   \cR_F^{(m)}/(t_{LT})\iso \prod_{\mu\geq m}
   F_\mu,\ a\auf(\overline{\iota_\mu(a)})_{\mu\geq m}
  \end{align}
  for each $i$, we obtain a preimage
  $\sum_{i=1}^d a_ie_i$ of $(y_\mu)_\mu$ under
  (\ref{c2-surj}). To see \eqref{prod iso} one uses $t_{LT}=Z\prod_{\mu\geq 1}\frac{Q_\mu}
  {\pi_L}$ as well as the fact that $t_{LT}$ and
  $\prod_{\mu\geq m}\frac{Q_\mu}{\pi_L}$ differ by a unit
  in $\cR_F^{(m)}$ since for $\nu<m$ the $Q_\nu$ are
  units as they have no zeros inside the annulus of convergence
  of $\cR_F^{(m)}$. Now \eqref{prod iso} follows via a
  projective limit argument from the isomorphisms
  $\cR_F^{(m)}/(Q_\mu)\cong F_\mu$ and the chinese remainder theorem.\\
 Concerning the first statement $\ker(c_2)=
 \im(c_1)$, one needs to show
 for any $x\in M^{(m)}_{t_{LT}}$ that if
 $\iota_\mu(x)\in \bD_{\dif,\mu}^+(M)$ for all
 $\mu\geq m$, then we have in fact $x\in M^{(m)}$.
 Writing $x=x_0\cdot t_{LT}^{-k}$ with
 $x_0\in M^{(m)}$, Remark \ref{iota t} implies
 \begin{align*}
  \iota_\mu(x)=\frac{\iota_\mu(x_0)\cdot \pi_{LT}^{\mu k}}{t_{LT}^k}.
 \end{align*}
 So the claim follows if we show that if
 $t_{LT}^k$ divides $\iota_\mu(x_0)$ in $F_\mu[[t_{LT}]]$ for all $\mu\geq m$,
 then it also divides $x_0$ in $M^{(m)}$.
 Of course, we can assume $k=1$ as well as $M=\cR_F$
 after choosing
 a basis of $M$. Then the isomorphism (\ref{prod iso})
 yields the desired result.\\ 
  \\
  Now that the exactness of (\ref{2.21:5}) is
  established, it follows by construction that the
  sequence (\ref{2.21:7}) of complexes is exact as
  well.\\
  Next we form the direct limit of the sequences
  (\ref{2.21:7}) over $m$, where the transition
  maps are the ones
  induced by the natural inclusions
  $M_0^{(m)}\nachinj M_0^{(m+1)}$ and the maps
  \begin{align*}
   a^\bullet \colon \widetilde K_{\varphi,\mf Z_n}(\bD_{\dif,m}
   ^{(+)}(M))\nach \widetilde K_{\varphi,\mf Z_n}(\bD_{\dif,m+1}
   ^{(+)}(M))
  \end{align*}
  given by "cutting off" the component at the lowest
  index. So by defining
  \begin{align*}
   \widetilde K_{\varphi,\mf Z_n}(\bD_{\dif}^{(+)}
   (M)):=\varinjlim_{m,a^\bullet}
   \widetilde K_{\varphi,\mf Z_n}(\bD_{\dif,m}^{(+)}
   (M))
  \end{align*}
  we obtain the desired
  exact sequence (\ref{2.21:lim7}).\\
  This sequence
  yields a distinguished triangle as explained in the
  beginning of the proof; in order to bring this
  triangle into the desired form, it remains to
  define suitable quasi-isomorphisms of complexes
  $K_{\mf Z_n}(\bD_{\mathrm{dif}}^{(+)}(M))\nach\widetilde K_{\varphi,\mf Z_n}(\bD_{
 \mathrm{dif}}^{(+)}(M))$, which is done in the following way:
 First, for $m\gg 0$ consider the morphisms
 $K_{\mf Z_n}(\bD_{\mathrm{dif},m}^{+}(M))\nach\widetilde K_{\varphi,\mf Z_n}(\bD_{
 \mathrm{dif},m}^{+}(M))$ defined by
 \begin{align}
  \begin{xy}
   \xymatrix{
    \bD_{\dif,m}^+(M) \ar[d]_{x\auf (x)_{\mu\geq m}}
    \ar[r]^{\mf Z_n} & \bD_{\dif,m}^+(M)\ar[d]^{
    x\auf((x)_{\mu\geq m},0)} &\\
    \prod\limits_{\mu\geq m}\bD_{\dif,\mu}^+(M) \ar[r] &
    \prod\limits_{\mu\geq m}\bD_{\dif,\mu}^+(M)
    \oplus\prod\limits_{\mu\geq m+1}\bD_{\dif,\mu}^+(M)
    \ar[r] & \prod\limits_{\mu\geq m+1}\bD_{\dif,\mu}^+(M),
   }
  \end{xy} \label{2.21:8}
 \end{align}
 noting that $b_0((x)_{\mu\geq m}) = ((\mathfrak{Z}_nx)_{\mu\geq m}, (x-x)_{\mu\geq{m+1}}) = ((\mathfrak{Z}_nx)_{\mu\geq m}, 0).$
 There are similar morphisms $K_{\mf Z_n}(\bD_{\mathrm{dif},m}(M))\nach\widetilde K_{\varphi,\mf Z_n}(\bD_{
 \mathrm{dif},m}(M))$, and one checks that they
 are all quasi-isomorphisms, using the exactness
 of the sequence
 \begin{align*}
 \begin{xy}
 \xymatrix{
  0\ar[r] & \bD_{\dif,m}^{(+)}(M)\ar[rr]^-{x\mapsto
  (x)_{\mu\geq m}} && \prod\limits_{\mu\geq m}
  \bD_{\dif,\mu}^{(+)}(M)
  \ar[rrr]^-{(x_\mu)\mapsto (x_{\mu-1}-x_\mu)_{\mu\geq
  m+1}} &&&\prod\limits_{\mu\geq m+1}
  \bD_{\dif,\mu}^{(+)}(M) \ar[r] & 0.
  }
  \end{xy}
 \end{align*}
 It is obvious that the quasi-isomorphisms
 $K_{\mf Z_n}(\bD_{\mathrm{dif},m}^{(+)}(M))\iso\widetilde K_{\varphi,\mf Z_n}(\bD_{
 \mathrm{dif},m}^{(+)}(M))$ are
 compatible with the transition maps, induced
 by the inclusions $\bD_{\dif,m}^{(+)}(M)
 \nachinj\bD_{\dif,m+1}^{(+)}(M)$ on the left and
 given by the $a^\bullet$
 on the right, so they induce a quasi-isomorphism
 \begin{align*}
  K_{\mf Z_n}(\bD_{\mathrm{dif}}^{(+)}(M))\iso\widetilde K_{\varphi,\mf Z_n}(\bD_{
 \mathrm{dif}}^{(+)}(M)).
 \end{align*}
 Putting everything together, and inspecting the
 explicit definitions of the morphisms involved,
 we get that the upper row of the diagram in the
 statement is in fact a distinguished triangle.\\
 \\
 To demonstrate that the second row is also a
 distinguished triangle, we start again with forming
 a certain direct limit of the
 exact sequences (\ref{2.21:7}) of complexes.
 But this time, instead of the $a^\bullet$ from
 above, we define morphisms
 \begin{align*}
  (a')^\bullet\colon \widetilde K_{\varphi,\mf Z_n}
  (\bD_{\dif,m}^{(+)}(M))\nach
  \widetilde K_{\varphi,\mf Z_n}
  (\bD_{\dif,m+1}^{(+)}(M))
 \end{align*}
 given by shifting $(x_\mu)_{\mu\geq m}\auf
 (x_{\mu-1})_{\mu\geq m+1}$ instead of cutting off.
 Then let
 \begin{align*}
  \widetilde K_{\varphi,\mf Z_n}^{(\varphi)}
  (\bD_{\dif}^{(+)}(M)):=
  \varinjlim_{m,(a')^\bullet}
  \widetilde K_{\varphi,\mf Z_n}
  (\bD_{\dif,m}^{(+)}(M)).
 \end{align*}
 Furthermore, note that the quasi-isomorphisms
 $K_{\mf Z_n}(\bD_{\mathrm{dif},m}^{(+)}(M))\iso\widetilde K_{\varphi,\mf Z_n}(\bD_{
 \mathrm{dif},m}^{(+)}(M))$ from (\ref{2.21:8})
 also form a morphism of directed systems if we
 use the $(a')^\bullet$ instead of the $a^\bullet$
 as transition maps on the right, so they yield a quasi-isomorphism
 \begin{align}
  K_{\mf Z_n}(\bD_{\dif}^{(+)}(M))\iso
  \widetilde K_{\varphi,\mf Z_n}^{(\varphi)}(
  \bD_{\dif}^{(+)}(M)). \label{qis'}
 \end{align}
 After these preparations, we consider the
 chain of quasi-isomorphisms
 \begin{align}
  \label{qis}&\widetilde K_{\varphi,\mf Z_n}(M^{(m)})\cong
   \cone\left(
  \widetilde K_{\varphi,\mf Z_n}(M^{(m)}_{t_{LT}})\oplus\widetilde K_{\varphi,\mf Z_n}(\bD^+_{\dif,m}(M))\nach
  \widetilde K_{\varphi,\mf Z_n}(\bD_{\dif,m}(M))\right)[-1]\\
   &\cong
   \cone\left(K_{\mf Z_n}(M^{(m)}_{t_{LT}})\oplus \widetilde K_{\varphi,\mf Z_n}(\bD^+_{\dif,m}(M))\nach
   K_{\mf Z_n}(M^{(m+1)}_{t_{LT}})\oplus \widetilde K_{\varphi,\mf Z_n}(\bD_{\dif,m}(M))\right)[-1],\notag
 \end{align}
 where the first one follows from
 applying the fact from
 homological algebra stated at the beginning of the
 proof to
 the sequence (\ref{2.21:7}) and the second one is
 formally
 obtained by the identity
 \begin{align*}
  \widetilde K_{\varphi,\mf Z_n}(M^{(m)}_{t_{LT}})
  =\cone\left(K_{\mf Z_n}(M^{(m)}_{t_{LT}})\overset{1-\varphi}{\nach
  }K_{\mf Z_n}(M^{(m+1)}_{t_{LT}})\right)[-1].
 \end{align*}
 Taking the direct limit of the quasi-isomorphisms
 (\ref{qis}) with respect to the transition maps
 $(a')^\bullet$ and the morphisms induced by
 $\varphi\colon M_0^{(m)}\nachinj
 M_0^{(m+1)}$, and applying the quasi-isomorphism
 $K_{\varphi,\mf Z_n}(M)
 \cong K_{\varphi,\mf Z_n}^{(\varphi)}(M)$ from
 Lemma \ref{2.20Nak}(ii) to the left-hand side
 and (\ref{qis'}) to right-hand side,
 we obtain the distinguished triangle
  \begin{align*}
  K_{\varphi,\mf Z_n}(M)\nach
   K_{\mf Z_n}^{(\varphi)}(M_{t_{LT}})\oplus
   K_{\mf Z_n}(\bD_{\mathrm{dif}}^+(M))
   \nach
   K_{\mf Z_n}^{(\varphi)}(M_{t_{LT}})
   \oplus K_{\mf Z_n}(\Ddif(M))\overset {+1} \nach
 \end{align*}
 which is the bottom row in the statement of the
 proposition.
\end{proof}

We define
\begin{align*}
\mathbf{D}^{(n)}_{\tn{dR}}(M):= H^{0}_{\mathfrak Z_n}(\mathbf{D}_{\textnormal{dif}}(M)) \quad \tn{ and } \quad \mathbf{D}^{(n)}_{\tn{cris}}(M):=H^{0}_{\mathfrak Z_n}(M_{t_{LT}}).
\end{align*}
For $m \gg 0$, the map $$\varphi \colon H^{0}_{\mathfrak Z_n}(M^{(m)}_{t_{LT}}) \nach H^{0}_{\mathfrak Z_n}(M^{(m+1)}_{t_{LT}})$$ is an isomorphism by Lemma \ref{2.17Nak}. Moreover, the inclusions $H^{0}_{\mathfrak Z_n}(M^{(m)}_{t_{LT}}) \nachinj \mathbf{D}^{(n)}_{\tn{cris}}(M)$ are isomorphisms by a result analogous to \cite[Lemma 2.18]{NaANT} which is formally deduced from Lemma \ref{2.17Nak} {(and the fact that the cohomologies are finite-dimensional). Thus the above $\varphi$ can be viewed as an automorphism
\begin{align*}
\varphi \colon \mathbf{D}^{(n)}_{\tn{cris}}(M) \iso \mathbf{D}^{(n)}_{\tn{cris}}(M).
\end{align*}
Next we construct two isomorphisms $j_1, j_2 \colon \mathbf{D}^{(n)}_{\tn{cris}}(M) \iso  H^{(\varphi),0}_{\mathfrak Z_n}(M_{t_{LT}})$ making the diagram
\[
\begin{tikzcd}
\mathbf{D}^{(n)}_{\tn{cris}}(M) \arrow[d, "j_1"] \arrow[r, "1-\varphi"] & \mathbf{D}^{(n)}_{\tn{cris}}(M) \arrow[d, "j_2"] \\
H^{(\varphi),0}_{\mathfrak Z_n}(M_{t_{LT}}) \arrow[r, "\Frob - \id"]                & H^{(\varphi),0}_{\mathfrak Z_n}(M_{t_{LT}})
\end{tikzcd}
\]
commute, where $\Frob$ is induced by the $\Frob$ in Definition \ref{DefPfeileVorPROP21}(a). Let
\begin{align*}
j_1 \colon \mathbf{D}^{(n)}_{\tn{cris}}(M)=H^{0}_{\mathfrak Z_n}(M^{(m)}_{t_{LT}}) \xrightarrow{\varphi} H^{0}_{\mathfrak Z_n}(M^{(m+1)}_{t_{LT}}) \iso H^{(\varphi),0}_{\mathfrak Z_n}(M_{t_{LT}})
\end{align*}
where the last map is an isomorphism by Lemma \ref{2.17Nak}. Note that $j_1$ is independent of the choice of $m \gg 0$.
Finally, we set $$j_2 \colon \mathbf{D}^{(n)}_{\tn{cris}}(M) \xrightarrow{j_1} H^{(\varphi),0}_{\mathfrak Z_n}(M_{t_{LT}}) \xrightarrow{\Frob} H^{(\varphi),0}_{\mathfrak Z_n}(M_{t_{LT}}).$$
Additionally, we define
\begin{align*}
\mf i \colon \mathbf{D}^{(n)}_{\tn{cris}}(M) \xrightarrow{j_1} H^{(\varphi),0}_{\mathfrak Z_n}(M_{t_{LT}}) \xrightarrow{\iota} \mathbf{D}^{(n)}_{\tn{dR}}(M)
\end{align*}
where $\iota$ is induced by the $\iota$ in Definition \ref{DefPfeileVorPROP21}(a).
\begin{definition}
Denote by
\begin{align*}
\expnM \colon \mathbf{D}^{(n)}_{\tn{dR}}(M) \nach H^{1}_{\varphi,\mathfrak Z_n}(M)
\end{align*}
and
\begin{align*}
\expnfM \colon \mathbf{D}^{(n)}_{\tn{cris}}(M) \xrightarrow{j_2} H^{(\varphi),0}_{\mathfrak Z_n}(M_{t_{LT}}) \nach H^{1}_{\varphi,\mathfrak Z_n}(M)
\end{align*}
the boundary maps obtained by taking cohomology of the exact triangles in Proposition \ref{Nak:Prop2.21}.
\end{definition}
Set
\begin{align*}
H^{1}_{\varphi,\mathfrak Z_n}(M)_e:=\tn{Im}(\mathbf{D}^{(n)}_{\tn{dR}}(M) \xrightarrow{\expnM} H^{1}_{\varphi,\mathfrak Z_n}(M))
\end{align*}
and
\begin{align*}
H^{1}_{\varphi,\mathfrak Z_n}(M)_f:=\tn{Im}(\mathbf{D}^{(n)}_{\tn{cris}}(M) \oplus \mathbf{D}^{(n)}_{\tn{dR}}(M) \xrightarrow{\expnfM + \expnM} H^{1}_{\varphi,\mathfrak Z_n}(M))
\end{align*}
and
\begin{align*}
t^{(n)}_M:=\mathbf{D}^{(n)}_{\tn{dR}}(M)/\mathbf{D}^{(n)}_{\tn{dR}}(M)^0 \quad \tn{ where } \quad \mathbf{D}^{(n)}_{\tn{dR}}(M)^0:=H^{0}_{\mathfrak Z_n}(\mathbf{D}^+_{\textnormal{dif}}(M)).
\end{align*}
Then Proposition \ref{Nak:Prop2.21} yields the following diagram with exact rows
\begin{equation}\label{f:exp-sequence}
  \begin{tikzcd}
0 \arrow[r] & H^{0}_{\varphi,\mathfrak Z_n}(M) \arrow[r, "x \mapsto x"] \arrow[d, "\id"] & \mathbf{D}^{(n)}_{\tn{cris}}(M)^{\varphi=1} \arrow[r, "x \mapsto \overline{\mf i(x)}"] \arrow[d, "x \mapsto x"] & t^{(n)}_M \arrow[r, "\expnM"] \arrow[d, "{x \mapsto (0,x)}"] & H^{1}_{\varphi,\mathfrak Z_n}(M)_e \arrow[r] \arrow[d, "x \mapsto x"] & 0 \\
0 \arrow[r] & H^{0}_{\varphi,\mathfrak Z_n}(M) \arrow[r, "x \mapsto x"]                  & \mathbf{D}^{(n)}_{\tn{cris}}(M) \arrow[r, "d_5"]                                                    & \mathbf{D}^{(n)}_{\tn{cris}}(M) \oplus t^{(n)}_M \arrow[r, "d_6"]                          & H^{1}_{\varphi,\mathfrak Z_n}(M)_f \arrow[r]                         & 0
\end{tikzcd}
\end{equation}
 where
 \begin{align*}
 d_5(x,y)=((1-\varphi)x, \overline{\mf i (x)}) \quad \tn{ and } \quad d_6=\expnfM + \expnM.
 \end{align*}
 For later calculations we state  the analogue (in the LT-setting) of the $\mathfrak{f}$-version (with $\mathfrak{f}\in\{\Psi,\varphi\}$) of the explicit formula for $\exp_{M} $  and $\exp_{f,M} $ in \cite[Prop.\ 2.23]{NaANT} and leave the straight forward adaption of the proof to the reader.
  \begin{proposition}\label{prop:nak17-2.23}
  \begin{enumerate}
    \item For $x\in \mathbf{D}^{(n)}_{\tn{dR}}(M) $ choose $\tilde{x}\in M_{t_\mathrm{LT}}^{(n')}$ for some sufficiently large $n'\geq n$ such that
    \[\iota_m(\tilde{x})-x\in \bD_{
 \mathrm{dif},m}^{+}(M)\]for any $m\geq n'.$ Then we have
 \begin{align*}
   \exp_M^{(n)}(x)= & [(\mathfrak{f}-1)  \tilde{x}, {\mf Z_n} \tilde{x}]\in H^{1}_{\mathfrak{f},\mathfrak Z_n}(M).
 \end{align*}
    \item  For $x\in \mathbf{D}^{(n)}_{\tn{dR}}(M) $ choose $\tilde{x}\in M_{t_\mathrm{LT}}^{(n')}$ for some sufficiently large $n'\geq n$ such that
    \[\iota_m(\tilde{x}) \in \bD_{
 \mathrm{dif},n'}^{+}(M)\]
 and
 \[\iota_{n'+k}(\tilde{x})-\sum_{l=1}^{k}\iota_{n'+l}(\varphi^{n'}(x))\in \bD_{
 \mathrm{dif},n'+k}^{+}(M)\]
 for any $k\geq 1.$
  Then we have
 \begin{align*}
   \exp_{f,M}^{(n)}(x)= & [(\varphi-1)  \tilde{x}+\varphi^{n'}(x), {\mf Z_n} \tilde{x}]\in H^{1}_{\varphi,\mathfrak Z_n}(M).
 \end{align*}
  \end{enumerate}
 \end{proposition}

\subsection{Derivatives of measures}

In cyclotomic Iwasawa theory the constant $\log(\chi_{cyc}(\gamma))$ shows up at various places (see \cite{NaANT}) in order to make constructions independent of the choice of a topological generator $\gamma$ of $\Gamma_{\Qp}.$ Since we have replaced the element $\gamma-1$ by ${\mf Z_n}$ we again have to check the dependence on this choice.
As our computations below show, the constant $\Omega$ plays a role in normalisation and seems conceptually new at a first glance since in the case $L=\QQ_p$ one can take $\Omega=1.$ But recall that $\Omega$ is only unique up to units in $o_L,$ hence in the cyclotomic case one could just as well take any element of $\ZZ_p^\times.$
Comparing \eqref{f:CTr} with \cite[Proposition 5.2]{NaANT} we see that we should take $\Omega_{\QQ_p} = \operatorname{log}_0(\chi(\gamma))^{-1},$ where $\log_0(a):= \log(a)/p^{v_p(\log(a))}.$
We first generalize the derivative of a measure from the cyclotomic case (e.g.\ \cite[\S 2.1]{LVZ15}) in a {\it naive} way:

By $\langle-\rangle:o_L^\times \to 1+\pi_Lo_L$ we denote the projection which is induced by the Teichm\"{u}ller character $\omega:k_L^\times \to o_L^\times.$
Fix $m_0>\frac{e}{p-1}$ and $m_1\geq 0$ such that $\log_p(o_L^\times)\subseteq \pi_L^{-m_1}o_L$. Then, for $s\in  \pi_L^mo_L$ with $m:=m_0+m_1$ the map
\[\langle-\rangle^s:o_L^\times \to 1+\pi_L^{m_0}o_L, x\mapsto \langle x\rangle^s:=\exp(s \log(x))\]
is well-defined.
For $\lambda\in D(\Gamma_L,K)$ and $f\in C^{an}(\Gamma_L,K)$ we define
\[\cL_\lambda(f,s):=\lambda(f \langle\chi_{LT}\rangle^s)\;\;\;\mbox{ ($\cL_\lambda(-,s)$ in $D( \Gamma_L,K)$ for fixed $s$)}\]
and
\[\cL_\lambda'(f):= \lim_{0\neq s\rightarrow 0}\frac{\cL_\lambda(f,s)-\cL_\lambda(f,0)}{s}\in D( \Gamma_L,K).\]
This limit exists  and we have
\begin{equation}\label{f:derivativedistribution}
  \cL_\lambda'(f)=\lambda(\log(\chi_{LT})f)
\end{equation}
using that $\lim_{0\neq s\rightarrow 0}\frac{\langle\chi_{LT}(\gamma)\rangle^s-1}{s}=\log(\chi_{LT}(\gamma)).$
As an example one easily  sees using Remark \ref{rem:Zproperties} that the expression
\begin{equation}\label{f:derivative}
  \frac{\cL_{{\mf Z_n}}'(\mathbf{1})}{\pi_L^n}
\end{equation}
is independent of $n.$

For $ D(\Gamma_n)$ as at the beginning of subsection \ref{sec:exp} there is another way of attaching such a derivative better adapted to the Lubin-Tate situation as follows:

By \cite[\S 3, Thm.\ 3.6]{ST} the characters of $\Gamma_n$ are all of the form $\psi_z(\gamma):=\kappa_z(\ell_n(\gamma))$ using their notation. For small $\gamma$ we have
\[\psi_z(\gamma)=\exp(\Omega\ell_n(\gamma)\log_{LT}(z))\]
and for $z=\exp_{LT}(\frac{\pi_L^n}{\Omega})$ the characters $\psi_z$ and $ \chi_{LT}$ coincide on an open subgroup of $\Gamma_L.$
For $\lambda\in D(\Gamma_n)$ and $f\in C^{an}(\Gamma_n,K)$ we may define
\[\cL\mathcal{T}_\lambda(f,z):=\lambda(f \psi_z)\;\;\;\mbox{ ($\cL\mathcal{T}_\lambda(-,z)$ in $D( \Gamma_n,K)$ for fixed $z$)}\]
and
\[\cL\mathcal{T}_\lambda'(f):= \lim_{0\neq z\rightarrow 0}\frac{\cL\mathcal{T}_\lambda(f,z)-\cL\mathcal{T}_\lambda(f,0)}{z}\in D( \Gamma_n,K).\]
This limit exists  and we have
\begin{equation}\label{f:derivativedistributionLT}
  \cL\mathcal{T}_\lambda'(f)=\frac{\Omega}{\pi_L^n}\lambda(\log(\chi_{LT})f)
\end{equation}
using that $\lim_{0\neq z\rightarrow 0}\frac{\psi_z(\gamma)-1}{z}=\frac{\Omega}{\pi_L^n}\log(\chi_{LT}(\gamma))$ as $g_{LT}(0)=1$ by \eqref{f:tLT}.

We conclude this discussion by considering again $\lambda={\mf Z_n}$ and  the trivial character $f= \mathbf{1}.$ Then $\cL\mathcal{T}_{\mf Z_n}(\mathbf{1},z)={\mf Z_n}(\psi_z)= z$ by \cite[Lem.\ 4.6]{ST}, whence $\cL\mathcal{T}_{\mf Z_n}(\mathbf{1},0)=0$ and \eqref{f:derivativedistributionLT} becomes
\begin{equation}\label{f:LTderiv}
   1= \cL\mathcal{T}_{\mf Z_n}'(\mathbf{1})=\frac{\Omega}{\pi_L^n}{\mf Z_n}(\log(\chi_{LT})) \mbox{ and } \frac{\cL_{{\mf Z_n}}'(\mathbf{1})}{\pi_L^n}=\frac{1}{\Omega} .
\end{equation}

\subsection{The dual exponential map $\exp^*$}\label{sec:dualexp}
 Let $M$ be a free $L$-analytic $(\varphi_L, \Gamma_L)$-module over $\cR_K$. We say that $M$ is \textbf{de Rham} if
 the $B = \bigcup_m K_m((t_{LT}))$-module $\mathbf{D}_{\dif}(M)$ is trivial as a $(B,\Gamma_L)$-module. By Galois descent (technically in the form of \cite[\href{https://stacks.math.columbia.edu/tag/0CDR}{Tag 0CDR}]{SP} for the Scheme $X = \operatorname{Spec}(K)$) this is equivalent to $\mathbf{D}_{\dif}(M)$ being trivial as a $(B,\Gamma_n)$-module for some $n$. Indeed, in this case $\mathbf{D}_{dR}^{(n)}(M)$ is a $K\otimes_{L}L_n$-module with a semi-linear $\Gamma_L/\Gamma_n= \operatorname{Gal}(L_n/L)$-action (which is trivial on $K$). Note that we have $\dim_K \mathbf{D}^{(n)}_{\tn{dR}}(M)=[L_n : L] \cdot r_M$ where $r_M$ is the rank of $M$ over $\cR_K$ and that $\mathbf{D}_{dR}^{(n)}(M)$ is in fact free as a $K \otimes_L L_n$-module. We denote by $\delta_{2,M}$ the connecting homomorphism $$H^1_{\mf Z_n}(\bD_{\dif}(M)) \to H^2_{\varphi,\mf Z_n}(M)$$ obtained from the sequence in \ref{Nak:Prop2.21}. We define
 \begin{equation}\label{f:DefCg}
 C_g({\mf Z_n}):=\cL_{{\mf Z_n}}'(\mathbf{1})={\mf Z_n}(\log(\chi_{LT})) =\frac{\pi_L^n}{\Omega}
  \end{equation}    for the trivial character $\mathbf{1}.$ We stress that this is compatible with Nakamura's definition when specializing to the cyclotomic situation.
 	\begin{lemma}\label{lem:g_M}
 		Let $M$ be de Rham. Then the natural map $$\left(\bigcup_m K_m((t_{LT}))\right)\otimes_{K_n} \bD_{dR}^{(n)}(M)\to \bD_{\dif}(M)$$ is an isomorphism and the induced map
 		$$g_{M}^{(n)}:\bD_{dR}^{(n)}(M) \to H^1_{\mf Z_n}(\bD_{\dif}(M)),\;\;x \mapsto C_g({\mf Z_n})(1 \otimes x) $$ is an isomorphism.  The inverse is induced by sending $f\otimes d\in K_m((t_{LT}))\otimes_{K_n} \bD_{dR}^{(n)}(M)$ to $C_g({\mf Z_n})^{-1}\frac{1}{[K_m:K_n]}Tr_{K_m/K_n}(f_{|t_{LT}=0})d,$ where by abuse of notation (although $t_{LT}$ gets inverted!) we denote by $ f_{|t_{LT}=0}$ the constant term of $f$ with respect to $t_{LT}$.
 	\end{lemma}
 \begin{proof} The first part follows immediately from the definition and implies that $\bD_{\dif}(M)$ is isomorphic to the trivial $B$-semi-linear $\Gamma_n$-representation. For the second statement it thus suffices to consider the rank $1$ case and prove the statement for $B$ itself, namely that the natural map $B^{\Gamma_n} \to B \to B/\mf Z_n$ is an isomorphism. Because the $\Gamma$-action respects the direct product decomposition $K_m((t_{LT})) \cong \prod_{k \in \ZZ} K_mt^{k}_{LT}$ and $B^{\Gamma_n} = K_n$ it suffices to show that any Laurent series, whose constant term vanishes, lies in the image of $\mf Z_n$   and that there is an exact sequence of the form
 \[\xymatrix@C=0.5cm{
   0 \ar[r] & K_n \ar[rr]^{ } && K_m \ar[rr]^{\mf Z_n} && K_m \ar[rr]^{\mathrm{Tr}} && K_n \ar[r] & 0 }\] with $\mathrm{Tr}=\frac{1}{[K_m:K_n]}\mathrm{Tr}_{K_m/K_n} $by Wedderburn theory. Using the product decomposition it suffices to treat the monomials $at_{LT}^k$ with some a $\in K_m.$ Taking $1\neq\gamma \in \Gamma_m$ we obtain  $\gamma(a)=a$ and $\gamma(t_{LT})= \chi_{LT}(\gamma)t_{LT}.$ By construction $(\gamma-1)(t_{LT}) = \pi_L^lu t_{LT}$ for some $l \in \ZZ, u\in o_L^{\times}$ and hence $(\gamma-1)(at_{LT}\pi_L^{-l}u^{-1}) = a\pi_L^{-l}u^{-1}((\gamma-1)(t_{LT})) = at_{LT}.$ Since $\delta_{\gamma}-1$ is divisible by $\mf Z_n$ in $D(\Gamma_n,K),$ we conclude that $at^{k}_{LT}$ lies in the image of $\mf Z_n.$
 \end{proof}
Note that $\mathbf{D}^{(n)}_{\dR}(M),$ for a $(\varphi_L,\Gamma_L)$-module $M,$ carries a natural filtration given by \[\Fil^i \mathbf{D}^{(n)}_{\dR}(M) = \mathbf{D}^{(n)}_{\dR}(M)\cap t_{LT}^i\mathbf{D}^{+}_{\dif}(M).\]
\begin{definition}\label{def:dualexp}
	Let $M$ be a de Rham $(\varphi_L,\Gamma_L)$-module over $\cR_K$. We define the dual exponential map as the composite $$H^1_{\varphi,\mf Z_n}(M) \to H^1_{\mf Z_n} (\bD_{\dif}^+(M)) \to H^1_{\mf Z_n}(\bD_{\dif}(M))\xrightarrow{(g_M^{(n)})^{-1}}D_{\dR}^{(n)}(M).$$ Where the first map is given by mapping $[x,y]$ to $[\iota_\mu(y)]$ with $\mu \gg 0.$
	Its image is contained in $\Fil^0(D_{\dR}^{(n)}(M))$ and we thus obtain a map $$\exp^{*,(n)}_{\tilde{M}}: H^1_{\varphi,\mf Z_n}(M)\to \Fil^0(D_{\dR}^{(n)}(M)).$$ We define
\[\exp^{*}_{\tilde{M}}:H^1_{\varphi,D(\Gamma_L)}(M)\to \Fil^0(D_{\dR} (M))\]
 by taking $\Gamma_L$-invariants of $\exp^{*,(n)}_{\tilde{M}}$, which is independent of the choice of $n.$ Indeed, as shown in \cite[Lem.\ 3.2.7]{Ste1} the restriction map
 \[H^1_{\varphi,\mf Z_n}(M) \to H^1_{\varphi,\mf Z_m}(M), [x,y]\mapsto [x,\mathfrak{Q}_{m-n}y], \]
 for $m\geq n$ induces an isomorphism after taking $\Gamma_L$-invariants, where   $\mathfrak{Q}_{m-n}({\mf Z_n}):=\frac{{\mf Z_m}}{{\mf Z_n}}=\frac{\varphi_L^{m-n}({\mf Z_n})}{{\mf Z_n}}$ with $\mathfrak{Q}_{m-n}(0)=\pi_L^{m-n}$ and we have $C_g({\mf Z_m})= \pi_L^{m-n}C_g({\mf Z_n})$ by \eqref{f:derivative}.

\end{definition}
\begin{definition}
	We define a pairing $$\cup_{\dif}\colon H^0_{\mf Z_n}(\bD_{\dif}(M_1)) \times H^1_{\mf Z_n}(\bD_{\dif}(M_2)) \to H^1_{\mf Z_n}(\bD_{\dif}(M_1 \otimes_{\cR_K}M_2))$$ given by $(x,y)\mapsto [x \otimes y].$
	Furthermore we define
	$$\langle-,-\rangle_{\dif}:H^0_{\mf Z_n}(\bD_{\dif}(M)) \times H^1_{\mf Z_n}(\bD_{\dif}(\tilde{M})) \xrightarrow{\cup_{\dif}} H^1_{\mf Z_n}(\bD_{\dif}(M \otimes_{\cR_K}\tilde{M})) \to K
$$
as composite of $\cup_{\dif}$ with
\[H^1_{\mf Z_n}(\bD_{\dif}(M \otimes_{\cR_K}\tilde{M})) \xrightarrow{\ev} H^1_{\mf Z_n}(\bD_{\dif}(\Omega^1))\xrightarrow{(g^{(n)}_{\Omega^1})^{-1}} \bD_{\dR}^{(n)}(\Omega^1)\cong K_n \xrightarrow{\frac{1}{[K_n:K]}Tr_{K_n/K}} K\]
using 	that $\bD_{\dR}^{(n)}(\Omega^1) \cong K_n$ via $t_{LT}^{-1}e \mapsto 1,$ where $e$ corresponds to $1$ in $\Omega^1\cong \cR_K(\delta)$.
	We further define $\langle-,- \rangle_{\dR}$ via the composite
	$$\bD_{\dR}^{(n)}(M) \times \bD_{\dR}^{(n)}(\tilde{M}) \to \bD_{\dR}^{(n)}(M \otimes \tilde{M}) \xrightarrow{ev} \bD_{\dR}^{(n)}(\Omega^1) \cong K_n \xrightarrow{\frac{1}{[K_n:K]}Tr_{K_n/K}}K.$$
\end{definition}
\begin{remark}\label{rem:dualitydeRham}
	The pairing $\langle -,- \rangle_{\dR}:=\langle -,- \rangle_{\dR,M}$  $$\bD_{\dR}^{(n)}(M) \times \bD_{\dR}^{(n)}(\tilde{M}) \to K$$is perfect if $M$ is de Rham  and induces a perfect pairing $$\bD_{\dR}^{(n)}(M)/\Fil^0\bD_{\dR}^{(n)}(M) \times \Fil^0\bD_{\dR}^{(n)}(\tilde{M}) \to K.$$
\end{remark}
\begin{proof}
Let us abbreviate $V:= \bD_{\dR}^{(n)}(M), G = \Gamma_n$ and $B = K_n((t_{LT})).$ For a suitable $r_n$ we have that $V = (B\otimes_{\iota_n}M^{(n)})^G = (B\otimes_{K_n}V)^G$ by definition.  We will first show that the pairing on the level of $K_n$ is perfect. Observe that base change to $B$ provides us with an injection $\operatorname{Hom}_{K_n}(V,K_n) \hookrightarrow \operatorname{Hom}_B(B\otimes_{K_n} V,B).$ The target can be endowed with a $G$ action  by $(g\lambda)(x) = g\lambda(g^{-1}x)$ and, because the action on $V$ is trivial, we see that the image of the above map is precisely the set of $G$-invariant elements. Indeed, since $B^G = K_n$, a linear form $\lambda$ which is fixed by $g$ has to map elements of the form $1 \otimes v$ into $B^G = K_n$ and hence restricts to an element of $\Hom_{K_n}(V,K_n).$
The perfectness now follows from $\Hom_B(B \otimes_{K_n} V,B) = \Hom_B(B \otimes_{\iota_n} M,B) \cong \Hom_B(B \otimes_{\iota_n} M,B(\chi_{LT}))\cong B \otimes_{\iota_n}\tilde{M}$ by taking $G$-invariants, using that the evaluation pairing commutes with base change.
Note that we used that $B = B(\chi_{LT})$ as $\Gamma_L$-modules (since $B^\times$ contains $t_{LT}$) and that $\mathbf{D}_{\dif}(M(\chi_{LT})) \cong \mathbf{D}_{\dif}(M(\chi))$ as they have ``the same'' $\Gamma_L$-action.
To conclude perfectness of the $K$-valued pairing, it suffices to show that the pairing is non-degenerate on one side.
Set $W = \operatorname{Hom}_{K_n}(V,K_n).$ Let $V':=\Hom_K(V,K)$ which we view as a $K_n$-module in the obvious way;
 we endow
$W$ and $V'$ with a $\Gamma_L/\Gamma_n$-action
via $\gamma f(-)=\gamma f(\gamma^{-1}-)$.
By the above perfectness at the level of $K_n$ it thus suffices to show that  the map $W \to V'$ given by $w \mapsto \operatorname{Tr}_{K_n/K}(w(-))$ is injective.
One easily checks that it is compatible with the $K_n$- and $\Gamma_L/\Gamma_n$-structure on $V'.$ We thus have constructed a $\Gamma_L/\Gamma_n$-semilinear map between free $K_n$-modules of the same rank. By Galois descent it suffices to show that it is injective on $\Gamma_L/\Gamma_n$-invariant elements. Suppose $w \in W^{\Gamma_L}$ satisfies $\operatorname{Tr}_{K_n/K}(w(v))=0$ for all $v \in V.$
This means that the image of the map $w\colon V \to K_n$ is contained in the kernel of the trace map. For any $x\in V^{\Gamma_L}$, we obtain $w(x) \in K \cap \operatorname{ker}(Tr)=0$ by the $\Gamma_L
$-equivariance of $w$. Thus $w$ is trivial on $\Gamma_L$-invariants and by Galois descent trivial, because $V$ is generated by $\Gamma_L$-invariant elements, which implies $w=0.$

For the second statement observe first that $\Fil^0(\mathbf{D}^{(n)}_{dR}(\Omega^1))=0$ and $\Fil^{-1}(\mathbf{D}_{dR}^{(n)}(\Omega^1)))=K_n.$ Hence $\Fil^{0}(\mathbf{D}^{(n)}_{dR}(M))$ is contained in the subspace orthogonal to $\operatorname{Fil}^0\mathbf{D}_{dR}^{(n)}(\tilde{M})$. In order to see that this inclusion is an equality, it suffices to show that the canonical bijective morphism of filtered vectorspaces $\mathbf{D}^{(n)}_{\dR}(M) \otimes \mathbf{D}^{(n)}_{\dR}(\tilde{M}) \cong \mathbf{D}_{\dR}^{(n)}(M\otimes\tilde{M})$ is in fact an isomorphism. This is not entirely trivial and can be achieved by an analogue of \cite[Proposition 6.3.3]{BC}. As in their case one reduces to the corresponding statement about graded objects and finally to the corresponding statement of rank one objects (which is clear in our case as well).
\end{proof}

\begin{lemma}
\label{lem:ResTr}
The diagram
\begin{equation}\label{f:g-trace}
   \xymatrix{
     D_{\dR}^{(n)}(\Omega^1) \ar[d]_{g^{(n)}_{\Omega^1}} \ar@{=}[r]^{ } & K_n \frac{\mathbf{e}_\chi}{t_{LT}}  \ar[rr]^{ \frac{a}{t_{LT}}\mathbf{e}_\chi\mapsto a } && K_n \ar[d]^{\frac{1}{[K_n:K]}Tr_{K_n/K}} \\
     H^1_{\mf Z_n}(\bD_{\dif}(\Omega^1)) \ar[r]^{-\delta_2} & H^2_{\varphi,\mf Z_n}(\Omega^1) \ar[rr]^(0.65){Tr=C_{Tr}({\mf Z_n})Res} && K }
\end{equation}
is commutative.
\end{lemma}

\begin{proof}
Given any $a\in K_n$ (in the right upper corner of the diagram) we first have to calculate $\delta_2(1\otimes \frac{a}{t_{LT}}\mathbf{e}_\chi)$ as $x:=1\otimes \frac{a}{t_{LT}}\mathbf{e}_\chi\in K_n[[t_{LT}]][\frac{1}{t_{LT}}]\mathbf{e}_\chi=\bD_{\dif,n}(\Omega^1)$ represents - up to a constant -  the image of $\frac{a}{t_{LT}}\mathbf{e}_\chi\in D_{\dR}^{(n)}(\Omega^1)  $ under $g^{(n)}_{\Omega^1}.$  In order to calculate the transition map $\delta_2$ we use an analogue of  \cite[Lem.\ 2.12(2)]{Na}, which is an  easy snake-lemma application to \ref{Nak:Prop2.21}: Assume that $x$ belongs to $\bD_{\dif,k}(\Omega^1)$ for some $ k\geq 0.$ For  any element $\tilde{x}\in \frac{1}{t_{LT}}\Omega^{1,(k)}=\frac{1}{t_{LT}}\cR(\chi)^{(k)}$ such that
\[\iota_m(\tilde{x})-\mathrm{can}_{k,m}(x) \in \bD_{\dif,m}^+(\Omega^1)\]   (using the notation of Definition \ref{Ddif.Def})
for all $m\geq k,$ we then have $\delta_2([x])=[(\varphi-1)\tilde{x}]\in H^2_{\varphi,\mf Z_n}(\Omega^1).$

We construct $\tilde{x}$ as follows. Consider the isomorphism
    \begin{align}\label{prod iso0}
   \cR_K^{+}/(t_{LT})\iso \prod_{\mu\geq 0}
  K_\mu,\ a\auf(\overline{\iota_\mu(a)})_{\mu\geq 0}
  \end{align}
analogous to \eqref{prod iso} and let $f$ be an element in $ \cR_F^{+}$, whose class in the left hand side corresponds  to the tuple $(a_\mu)_{\mu\geq 0}$ with
\[a_\mu:=\left\{
         \begin{array}{ll}
           \frac{a}{q^{\mu-k}\pi_L^k}, & \hbox{if $\mu\geq k$;} \\
           \frac{\mathrm{Tr}_{K_k/K_\mu}(a)}{\pi_L^k}, & \hbox{if $k\geq\mu\geq 0$;}
         \end{array}
       \right.
\]
on the right hand side. Note that the operator $\Psi$ on $\cR^+_K$ induces the map
\begin{align*}
 \Psi:\prod_{\mu\geq 0}  K_\mu &\to \prod_{\mu\geq 0}
  K_\mu, \\
(x_\mu)&\mapsto (q^{-1}x_0+q^{-1}\mathrm{Tr}_{K_1/K_0}(x_1),q^{-1}\mathrm{Tr}_{K_2/K_1}(x_2),\ldots, q^{-1}\mathrm{Tr}_{K_{\mu+1}/K_\mu}(x_{\mu+1}),\ldots ).
\end{align*}
Moreover, $(x_\mu)_\mu$ satisfies $q\Psi((x_\mu)_\mu)= (x_\mu)_\mu$ if and only $\mathrm{Tr}_{K_1/K_0}(x_1)=0$ and $\mathrm{Tr}_{K_{m+1}/K_{m}}(x_{m+1})=x_{m}$ for all $m\geq 1.$ In particular, $q\Psi((a_\mu)_\mu)= (a_\mu)_\mu$ if and only if $ \mathrm{Tr}_{K_{k}/K_{0}}(a)=0.$
We now set $\tilde{x}=\frac{f}{t_{LT}}\mathbf{e}_\chi\in\frac{1}{t_{LT}}\cR(\chi)^{(k)}$ and check that, for $m\geq k$,
\begin{align*}
  \iota_m(\tilde{x}) & \equiv \frac{\pi_L^mf(u_m)}{t_{LT}}\mathbf{e}_\chi \; \mod \bD_{\dif,m}^+(\cR(\chi))  \\
    & = \frac{\pi_L^m a_m}{t_{LT}}\mathbf{e}_\chi\\
 & = \frac{\pi_L^{m-k} a }{q^{m-k}t_{LT}}\mathbf{e}_\chi\\
    & = \mathrm{can}_{k,m}(x)
\end{align*}
as required, i.e., $\delta_2(1\otimes \frac{a}{t_{LT}}\mathbf{e}_\chi)=[(\varphi-1)(\frac{f}{t_{LT}}\mathbf{e}_\chi)]=[\left(\frac{\varphi(f)}{q}-f\right)\frac{1}{t_{LT}}\mathbf{e}_\chi]$. Since
\begin{align*}
  \left(\frac{\varphi(f)}{q}-f\right)(u_m)= & \frac{ f(u_{m-1})}{q}-f(u_m)= \frac{ a_{m-1} }{q}-a_m=0
\end{align*}
for all $m>k$, we conclude from \eqref{prod iso0} that $\frac{\varphi(f)}{q}-f\in \left(\prod_{m>k}^{\infty}\frac{Q_m}{\pi_L} \right)\cR^+_K, $ whence
\begin{equation}\label{f:}
  \left(\frac{\varphi(f)}{q}-f\right)\frac{1}{t_{LT}}\in \frac{\left(\prod_{m>k}^{\infty}\frac{Q_m}{\pi_L} \right)}{t_{LT}}\cR^+_K=\frac{1}{Z\prod_{\mu\geq 1}^k\frac{Q_\mu}
 {\pi_L}}\cR^+_K
\end{equation}
using $t_{LT}=Z\prod_{\mu\geq 1}\frac{Q_\mu}
 {\pi_L}.$
Since all involved maps are $K$-linear and $K_n=K\oplus \ker(\mathrm{Tr}_{K_n/K})$ it suffices to check the commutativity in the two cases $a\in K,$ i.e., $k=0$, or $\mathrm{Tr}_{K_n/K}(a)=0,$ i.e., $q\Psi((a_\mu)_{\mu\geq 0})  = (a_\mu)_{\mu\geq 0}.$

If $k=0$,  the element $\frac{a}{t_{LT}}\mathbf{e}_\chi $  is sent via the lower composite to
\begin{align*}
 C_{Tr}({\mf Z_n})Res\circ (-\delta_2)\circ g_{\Omega^1}(\frac{a}{t_{LT}}\mathbf{e}_\chi ) & =-C_g({\mf Z_n}) C_{Tr}({\mf Z_n})Res\circ \delta_2(1\otimes \frac{a}{t_{LT}}\mathbf{e}_\chi) \\
    & =-\frac{q}{q-1} Res (\left(\frac{\varphi(f)}{q}-f\right)\frac{1}{t_{LT}}dt_{LT}) \\
   & =-\frac{q}{q-1} \left(\frac{\varphi(f)}{q}-f\right)_{|Z=0}\\
& =-\frac{q}{q-1}  \left(\frac{f(0)}{q}-f(0) \right) =a
\end{align*}
where we use for the second equality  the definition \eqref{f:CTr} and for the third equality \eqref{f:}.
Thus the claim follows on the subspace $D_{\dR}^{(n)}(\Omega^1)^{\Gamma_L}, $  because $ \frac{1}{[K_n:K]}Tr_{K_n/K}(a)    =a$.

If $q\Psi((a_\mu)_{\mu\geq 0})  = (a_\mu)_{\mu\geq 0}$, i.e., $q\Psi(f\mod t_{LT})  \equiv f\mod t_{LT}$, it follows from the surjectivity of $\frac{1}{\pi_L}\Psi-\frac{1}{q}$ on $\cR^+$ by \cite[Cor.\ 2.3.4]{BF} and the commutative diagram
with exact rows
\begin{equation*}
	\xymatrix{
		0   \ar[r]^{ } & \cR_K^+  \ar@{->>}[d]_{\frac{1}{\pi_K}\Psi-\frac{1}{q}} \ar[r]^{t_{LT} } & \cR_K^+ \ar[d]_{\Psi-\frac{1}{q}} \ar[r]^{ } & \cR_K^+/t_{LT}\cR_K^+ \ar[d]_{\Psi-\frac{1}{q}} \ar[r]^{ } & 0  \\
		0   \ar[r]^{ } & \cR_K^+    \ar[r]^{ t_{LT}} & \cR_K^+ \ar[r]^{ } & \cR_K^+/t_{LT}\cR_K^+  \ar[r]^{ } & 0  }
\end{equation*}
that we may assume without loss of generality that $f$ also satisfies $q\Psi(f)=f,$ whence we obtain $\Psi(\frac{\varphi(f)}{q}-f)=0.$ Using the identity $Res(\Psi(f)dt_{LT})=\chi(\pi_L)Res(fdt_{LT})$ from \cite[Prop.\ 1.5]{Co2} we conclude that $Res (\left(\frac{\varphi(f)}{q}-f\right)\frac{1}{t_{LT}}dt_{LT})$ vanishes, from which the commutativity follows also in this case by a similar calculation as above.
\end{proof}

\begin{lemma} \label{lem:expadjoint} Let $z \in H^0_{\mf Z_n}(\bD_{\dif}(M)),[x,y] \in H^1_{\varphi, \mf Z_n}(\tilde{M}), a \in H^0_{\varphi,\mf Z_n}(M)$ and $[b] \in H^2_{\varphi,\mf Z_n}(\tilde{M}).$
Using $\langle-,-\rangle:=\langle-,-\rangle_M$ as before to denote the pairing $$H^i_{\varphi,\mf Z}(M)\times H^{2-i}_{\varphi,\mf Z}(\tilde{M}) \to K$$ obtained from \ref{lem:explicitPairing} we have $$\langle \exp_M^{(n)}(z),[x,y]\rangle=\langle z,[\iota_n(y)]\rangle_{\dif} $$ and
$$\langle a, \delta_{2,M}([b]) \rangle = \langle  \iota_{n}(a),[b] \rangle_{\dif}.$$

\end{lemma}
\begin{proof}

Let $z$ be in $\bD_{\dR}^{(n)}(M)$ and $[x,y]\in H^1_{\varphi, \mf Z_n}(\tilde{M})$ with $x\in \tilde{M}^{(n)}, y\in \tilde{M}^{(n+1)}.$ Then
we have
\[\langle z,[\iota_n(y)]\rangle_{\dif}= \frac{1}{[K_n:K]}Tr_{K_n/K} \circ{ g_{\Omega}^{-1}\circ}\mathrm{ev}\Big( [z\otimes \iota_n(y)]\Big)=-Tr\circ\delta_2\Big( [z\otimes \iota_n(y)]\Big)\] by \eqref{f:g-trace} and, by the same snake-lemma application in order to calculate the transition map $\delta_2$ induced by Proposition \ref{Nak:Prop2.21} (compare with \cite[Lem.\ 2.12(2)]{Na}),
\begin{align*}
 -Tr\Big( \delta_2([z\otimes \iota_n(y)])\Big)  & =-C_{Tr}({\mf Z_n})Res\Big([(\varphi-1)(\tilde{z}\otimes y)]\Big) \\
    & =-C_{Tr}({\mf Z_n})Res\Big([(\varphi-1)(\tilde{z})\otimes\varphi( y)+\tilde{z}\otimes {\mf Z}_n x]\Big)\\
    & =-C_{Tr}({\mf Z_n})Res\Big( \varphi(y)\big((\varphi-1)(\tilde{z})\big)+ \big({\mf Z}_n x \big)(\tilde{z})  \Big),
\end{align*}
where in the second equality we have used the co-boundary condition $(\varphi-1)(y)={\mf Z}_nx.$ Moreover, $\tilde{z}\in M^{(n)}_{t_{LT}}$ is an element with the property that $\iota_m(\tilde{z})-z$ belongs to $\bD_{\dif,m}^+(M)$ for all $m\geq n$, the existence of which is granted by the exactness of \eqref{2.21:5}, whence $\iota_m(\tilde{z}\otimes y)-z\otimes\iota_m(y)\in \bD_{\dif,m}^+(M\otimes \tilde{M})$ for all $m\geq n$.

On the other hand we have by a straightforward analogue of \cite[Lem.\ 2.12 (1)]{Na} for the first, the formula  in Theorem \ref{lem:explicitPairing} for the second equality and \eqref{f:lambda} for the third equality
\begin{align*}
  \langle\exp_M^{(n)}(z), [x,y]\rangle & =\langle [(\varphi-1)(\tilde{z}), {\mf Z}_n\tilde{z} ],[x,y]\rangle\\
    & =-C_{Tr}({\mf Z_n})Res\Big( \varphi(y)\big((\varphi-1)(\tilde{z})\big)+(\lambda^\iota x)({\mf Z}_n\tilde{z})   \Big) \\
  &=-C_{Tr}({\mf Z_n})Res\Big( \varphi(y)\big((\varphi-1)(\tilde{z})\big)+ \big({\mf Z}_n x \big)(\tilde{z})  \Big),
\end{align*}
which agrees with the above formula. We leave the easy proof of the second identity to the reader.
\end{proof}
%
\begin{proposition} \label{prop:NakIw2.16}
	Let $M$ be de Rham. Let $x \in \bD_{\dR}^{(n)}(M)/\Fil^0\bD_{\dR}^{(n)}(M)$ and $y \in H^1_{\varphi,\mf Z_n}(\tilde{M}).$ We have $$\langle \exp_{M}^{(n)}(x),y \rangle_M=   \langle x,\exp^{*,(n)}_{M}(y) \rangle_{\dR},$$
	i.e., $\exp^{(n)}_{M}$ is adjoint to $ \exp_{M}^{*,(n)}.$
\end{proposition}
\begin{proof}
	This is a formal consequence of Lemma \ref{lem:expadjoint} after plugging in the definition \ref{def:dualexp} of $\exp_{M}^{*,(n)}.$
\end{proof}

Only for the purpose of the next lemma (needed in the proof of the subsequent proposition) we introduce the notation
$H^i_{\mathrm{mix}}(N)$ as the $i$-th cohomology of the complex $ K_{\mf Z_n}^{(\varphi)}(N_{t_{LT}})
   \oplus K_{\mf Z_n}(\Ddif(N))  $ of the bottom right in Proposition \ref{Nak:Prop2.21}.
We define a pairing $$\cup_{\mathrm{mix}}\colon H^0_{\mathrm{mix}}(M_1) \times H^1_{\mathrm{mix}}({M}_2) \to H^1_{\mathrm{mix}}(M_1 \otimes_{\cR_K}M_2)$$ given by $(x,y)\mapsto [x \otimes y].$
	Furthermore, we set
	$$\langle-,-\rangle_{\mathrm{mix}}:H^0_{\mathrm{mix}}(M) \times H^1_{\mathrm{mix}}(\tilde{M}) \to H^1_{\mathrm{mix}}(M \otimes_{\cR_K}\tilde{M})) \xrightarrow{\ev} H^1_{\mathrm{mix}}(\Omega) .$$ Finally, by
\[G:H^{1}_{\varphi,\mf Z_n}(\tilde{M})\to H^1(K_{\mf Z_n}^{(\varphi)}(\tilde {M}_{t_{LT}})\oplus
   K_{\mf Z_n}(\bD_{\mathrm{dif}}^+(\tilde{M)}) ) \to H^1_{\mathrm{mix}}(\tilde{M})\] we denote the composite $H^1(d_7)\circ H^1(d_3),$ where $d_7:K_{\mf Z_n}^{(\varphi)}(\tilde{M}_{t_{LT}})\oplus
   K_{\mf Z_n}(\bD_{\mathrm{dif}}^+(\tilde{M})) \to
   K_{\mf Z_n}^{(\varphi)}(\tilde{M}_{t_{LT}})
   \oplus K_{\mf Z_n}(\Ddif(\tilde{M}))$ sends $(x,y)$ to itself using the natural inclusion $  \bD_{\mathrm{dif}}^+(\tilde{M}) \hookrightarrow\Ddif(\tilde{M})$.
Then the next Lemma is formally analogous to Lemma \ref{lem:expadjoint}, thus we leave the details to the interested reader.
\begin{lemma}\label{lem:expadjointmix}
The following diagram is commutative
\[ \xymatrix{
H^0_{\mathrm{mix}}(M)\ar[d]^{\expnfM + \expnM}\ar@{}[r]|(0.65){\times} & { \phantom{m} }H^1_{\mathrm{mix}}(\tilde{M}) \ar[rr]^(0.5){\langle-,-\rangle_{\mathrm{mix}}} &&  H^1_{\mathrm{mix}}(\Omega^1)\ar[d]^{\partial_2}\\
   H^1_{\varphi,\mf Z_n}(M) \ar@{}[r]|(0.65){\times} &  H^{1}_{\varphi,\mf Z_n}(\tilde{M}) \ar[u]_{-G} \ar[rr]^(0.5){\langle-,-\rangle} &  & H^2_{\varphi,\mf Z_n}(\Omega^1).}\]
\end{lemma}

\begin{proposition}\label{prop:orthogonal}
Let $M $ be a trianguline $L$-analytic $(\varphi_L,\Gamma_L)$-module over $\cR_K$ which is de Rham. Then $H^1_{\varphi,\mf Z_n}(M)^{\Gamma_L}_f$ is the orthogonal complement of $H^1_{\varphi,\mf Z_n}(\tilde{M})^{\Gamma_L}_f$ with respect to the duality pairing $\langle -,- \rangle_M.$
\end{proposition}

\begin{proof}
 Analogous to \cite[Prop.\ 2.24]{NaANT}: Replacing the sequence (13) in (loc.\ cit.) by \eqref{f:exp-sequence}, using the  Euler-Poincar\'{e} formula \ref{rem:EP-formula} as well as duality \ref{lem:explicitPairing} and the de Rham property of $M$ one shows that $\dim_KH^1_{\varphi,\mf Z_n}(M)^{\Gamma_L}_f+\dim_K H^1_{\varphi,\mf Z_n}(\tilde{M})^{\Gamma_L}_f=\dim_K H^1_{\varphi,\mf Z_n}({M})^{\Gamma_L}.$ Therefore it suffices to show that $\langle x,y \rangle=0$ for all $x\in  H^1_{\varphi,\mf Z_n}(M)^{\Gamma_L}_f $ and $y\in H^1_{\varphi,\mf Z_n}(\tilde{M})^{\Gamma_L}_f  .$ This is accomplished by Lemma \ref{lem:expadjointmix}, because $G(y)=0$ since $y\in\ker H^1(d_3)$ by assumption.
\end{proof}

\section{\texorpdfstring{$\epsilon$}{Epsilon}-constants}\label{sec:varepsilon-constants}
Let $E$ be a field of characteristic zero containing $\mu_{p^\infty}$,
$\psi_0:=\psi_\xi: \qp\to E^\times$ the character (with kernel $\zp$) attached to a fixed compatible system $\xi=(\xi_n)_{n\geq 1}$ of $p$-power roots of unity via $\psi_0(\frac{1}{p^n})=\xi_n$.

Similarly, we may define for the compatible system $u=(u_n)_{n\in\mathbb{N}}\in T_\pi$ (and a choice of generator $t_0'$ of $\TLT'$) the character $\psi_u:=\psi_{u,t_0'}:L\to E^\times,\; \frac{x}{\pi_L^n}\mapsto \eta_{t_0'}(x,u_n).$

But there is another (canonical) choice:
$\psi_L:=\psi_0\circ \mathrm{Tr}_{L/\qp}: L\to E^\times$ is a locally constant character (with kernel the inverse of the different ideal $\mathfrak{D}_{L/\qp} $).

\begin{remark}\label{rem:psiu}
The character $\psi_L$ factorizes over $o_L.$ Hence, by \eqref{f:Tate} there exists $a=a(t_0',u)\in o_L$ such that the following diagram commutes
\[\xymatrix{
  L/o_L \ar[d]_{x\mapsto u\otimes x}^{\cong} \ar[r]^{Tr_{L/\mathbb{Q}_p}} & \mathbb{Q}_p/\mathbb{Z}_p  \ar[r]^{\psi_{\eta_{t_0'}(1,u)}} & \mu(p)  \\
  \TLT\otimes_{o_L}L/o_L  \ar[rr]^{a\cdot} &&   \TLT\otimes_{o_L}L/o_L   \ar[u]^{  \eta_{t_0'} }. }\]
  Here $\eta_{t_0'}(1,u):=(\eta_{t_0'}(1,u_n))_n$ is a generator of $\zp(1)$, again by \eqref{f:Tate}. In particular, for the choice $\xi=\eta_{t_0'}(1,u)$ we obtain
\begin{equation}\label{f:psiLpsiu}
  \psi_L(x)=\psi_u(ax)
\end{equation}
for all $x\in L.$ It is clear that $a$ is a generator of the  different ideal $\mathfrak{D}_{L/\qp} $.
\end{remark}

Let $dx$ be the Haar measure on $L$ such that $\int_{o_L}dx=1.$ Let $\psi \colon L \to E^\times$ be a non-trivial character which kills an open subgroup of $L$.

For a finite-dimensional $E$-linear representation $D$ of the Weil-group $W_L:=W(\bar{L}/L)$ which is locally constant (i.e., the image of the inertia group is finite) we have local constants
\[\epsilon(D,\xi):=\epsilon_E(L,D,\psi,dx)\in E^\times,\]
see \cite{deligne} or \cite{tate} and \cite[\S 2.2]{Daigle&Flach16}.

If $\dim_E D=1$ corresponding to a locally constant homomorphism $\delta:L^\times \to E^\times$ via local class field theory (see section \ref{sec:deRham} for the normalisation we choose), i.e., $D=E(\delta),$ then
\begin{align} \label{eq:constantsrank1}
  \epsilon_E(L,D,\psi_L,dx) &= \delta(\pi_L)^{a(\delta)} q^{n({\psi_L})} \sum_{i\in(o_L/\pi_L^{a(\delta)})^\times}\delta(i)^{-1} \psi_L(\frac{i}{\pi_L^{a(\delta)}}),\\\label{eq:constantsrank1u}
 \epsilon_E(L,D,\psi_u,dx) &= \delta(\pi_L)^{a(\delta)} q^{n(\psi_u)} \sum_{i\in(o_L/\pi_L^{a(\delta)})^\times}\delta(i)^{-1} \eta({i},u_{a(\delta)}).
\end{align}

Here $n(\psi)$ denotes the largest integer $n$ such that $\pi_L^{-n}o_L\subseteq \ker{\psi}$, $a(\delta)$  denotes the \textbf{conductor} of $\delta,$ ($0$ if $\delta$ is unramified, the smallest positive integer $m$ such that $1+\pi_L^mo_L\subseteq \ker\delta$,
if $\delta$ is ramified).
If $W= (D,N)$ is a Weil Deligne representation of $W_L$ with monodromy operator $N$ and underlying Weil group
 representation $D$ we modify, following Nakamura, its $\epsilon$-constant by the factor \be \label{eq:adjustedconstant} \epsilon(W):= \epsilon(D)\det(-\operatorname{Frob}\mid (D/D^{N=0})^{I_L}),\ee where $I_L$ denotes the inertia subgroup. Both definitions agree if $N$ acts as $0$ on $D.$

\begin{remark}\label{rem:npsiunull}
	$\psi_u(y)$ is independent of the choice of $n$ such that $y=x/\pi_L^n$ and $n(\psi_u)=0.$
\end{remark}
\begin{proof}
	The independence follows inductively from $\eta(\pi_Lx,u_{n}) = \eta(x,\varphi_L(u_{n})) = \eta(x,u_{n-1}).$ On the one hand, by definition $o_L \subset \ker \psi_u.$ On the other hand by §1, Fact 2 in \cite{boxall}, using that $u_1$ is a non-zero $\pi_L$-torsion point we may find $a \in o_L$ such that $\psi_u(a/\pi_L)=\eta(a,u_1)$ is a primitive $p$-th root of unity. This proves that $\pi_L^{-1}o_L$ is not contained in $\ker(\psi_u).$ We conclude $n(\psi_u)=0.$
\end{proof}
How do the epsilon-constants for the two choices $\psi_L$ and $\psi_u$ compare? The first choice behaves well under induction: there is  a constant $\lambda\in E$ depending on $L/\qp$, the choices of Haar measures $dx_{\qp}, dx_L  $ and the choice of $\psi_0,$ such that\footnote{As $\epsilon_E$ is inductive with regard to virtual representations of dimension $0,$ one concludes that $\lambda=\frac{\epsilon_E(\qp,\mathrm{Ind}^L_{\qp}\delta_{triv},\psi_0,dx_{\qp})}{\epsilon_E(L,\delta_{triv},\psi_L,dx_L)}$ for the trivial representation $\delta_{triv}.$}
\[ \epsilon_E(\qp,\mathrm{Ind}^L_{\qp}\delta,\psi_0,dx_{\qp}) =\lambda  \epsilon_E(L,\delta,\psi_L,dx_L) \] for all locally constant characters $\delta:L^\times \to E^\times$ (see \cite{BB} or \cite[(5.6)]{deligne}).\\
The second choice is obviously better adapted to the Lubin-Tate situation. By \eqref{f:psiLpsiu} there exists $a\in L^\times$ such that $\psi_L( x)=\psi_u(ax)$. Moreover, one knows that $\epsilon_E(L,\delta,\psi_L(bx),dx_L)=\frac{ \delta(b) }{ |b|}\epsilon_E(L,\delta,\psi_L,dx_L)$ by \cite[(3.2.3) or (3.4.4)]{tate} for all $b\in L^\times$.
Combining the above we get the following:
\begin{remark}\label{rem:epsind}
There exists $\lambda \in E$  depending on $L/\qp$, the choices of Haar measures $dx_{\qp}, dx_L  $ and the choice of $\psi_0$, as well as $a \in L^\times$ depending on $\psi_0$ and $u$ such that
\begin{equation}\label{f:comppsiu}
   \epsilon_E(\qp,\mathrm{Ind}^L_{\qp}\delta,\psi_0,dx_{\qp})=\lambda \frac{ \delta(a) }{ |a|}\epsilon_E(L,\delta,\psi_u,dx_L)
\end{equation}
for all $\delta.$

\end{remark}

If we start with a Haar measure $dx$ of $L$, then the dual Haar measure $\hat{dx}$ with respect to the duality induced by $\psi$, i.e.,
\[L\times L \to \mu_{p^\infty}\subseteq E^\times , (x,y)\mapsto \psi (xy),\]
is the unique Haar measure such that $f(x)=\hat{\hat{f}}(-x)=\int_L \hat{f}(y)\psi(-xy)\hat{dx}(y)$ holds for all test functions in $L^1(L),$ where
\[\hat{f}(y):=\int_L f(x)\psi(xy)dx(x)\]
denotes the Fourier transform of $f.$ Especially for $f=\mathbf{1}_{\pi_L^{-n(\psi)}o_L}$ we obtain:
\begin{align*}
\widehat{\mathbf{1}_{\pi_L^{-n(\psi)}o_L}}(y)= \left(\int_{\pi_L^{-n(\psi)}o_L} dx\right)\mathbf{1}_{o_L}(y),
\end{align*}
whence
\begin{align*}
\widehat{\widehat{\mathbf{1}_{\pi_L^{-n(\psi)}o_L}}}(-x)&=\left(\int_{\pi_L^{-n(\psi)}o_L} dx\right)\int_{o_L} \psi(-xy)\hat{dx}(y)\\
&=\left(\int_{\pi_L^{-n(\psi)}o_L} dx\right)\left(\int_{o_L}\hat{dx}(y)\right)  \mathbf{1}_{_{\pi_L^{-n(\psi)}o_L}}(x),
\end{align*}
i.e., $\int_{o_L}\hat{dx}(y)=\frac{1}{q^{n(\psi)}}   $ and $\hat{dx}=\frac{1}{q^{n(\psi)}}dx.$

 From \cite[(3.4.7)]{tate} we obtain
\begin{equation}\label{f:epsduality}
  \epsilon(L,\delta,\psi,dx)\epsilon(L,\delta^{-1}|-|,\psi(-x),\hat{dx})=1
\end{equation}
and similarly for higher rank representations $D$ instead of $\delta.$
Since by (3.2.2/3) in (loc.\ cit.) we have $\epsilon(L,\delta,\psi,rdx)=r\epsilon(L,\delta,dx)$ for $r>0$ and $\epsilon(L,\delta,\psi(ax),dx)=\delta(a)|a|^{-1}\epsilon(L,\delta,\psi(x),dx)$, we conclude that
\begin{equation}\label{f:epsduality2}
  \epsilon(L,\delta,\psi,dx)\epsilon(L,\delta^{-1}|-|,\psi(x),dx)=\delta(-1)q^{n(\psi)}.
\end{equation}
Moreover, by (3.4.5) in (loc.\ cit.) it holds that
\begin{equation}\label{f:epsabsolutevalue}
   \epsilon(L,\delta^{-1}|-|,\psi(x),dx)=q^{-a(\delta)-n(\psi)} \epsilon(L,\delta^{-1},\psi(x),dx)=|\pi_L^{a(\delta)+n(\psi)}| \epsilon(L,\delta^{-1},\psi(x),dx).
\end{equation}

\section{Epsilon-isomorphisms - the statement of the conjecture}

\label{sec:EpsilonStatement}

\subsection{Determinant functor}\label{sec:Det}
	
Let $R$ be a commutative ring. A graded invertible $R$-module is a pair $(\cL,r),$ where $\cL$ is an invertible $R$-module and $r: \Spec(R)\to \ZZ$ is a locally constant function. We define the category $\cP_R$ of graded invertible $R$-modules by setting $\operatorname{Mor}((\cL_1,r),(\cL_2,s)):=\operatorname{Isom}_R(\cL_1,\cL_2)$ if $r=s$ and empty otherwise.
We further define $$(\cL_1,r)\cdot (\cL_2,s):=(\cL_1,r)\otimes (\cL_2,s):= (\cL_1 \otimes \cL_2,r+s)$$ for each pair of objects and we identify $(\cL_1,r)\otimes (\cL_2,s)$ with $(\cL_2,s)\otimes (\cL_1,r)$ via the morphism induced by $l_1\otimes l_2 \mapsto (-1)^{r+s}l_2 \otimes l_1.$ We denote by $\mathbf{1}_R$ the object $(R,0),$ which acts as a unit with respect to the tensor product and we remark that every object $(\cL,r)$ has an inverse given by $(\cL^{-1},-r),$ where $\cL^{-1}$ denotes the $R$-dual of $\cL.$ For a ring morphism $R \to S$ and $(\cL,r) \in \cP_R$ we set $(\cL,r)_{S}:=(\cL \otimes S,r^*),$ where $r^*$ denotes the pullback of $r$ along $R \to S.$  An isomorphism $\mathbf{1}_R \to \cL$ is called a \textbf{trivialisation} of $\cL.$ Let $\cP_{fg}(R)$ be the category of finitely generated projective $R$-modules and let $(\cP_{fg}(R),is)$ be its core, i.e. the subcategory consisting of the same objects with isomorphisms as morphisms. We have a functor
	\begin{align*} \d_R\colon(\cP_{fg}(R),is) &\to \cP_R \\
	P &\mapsto (\det P, \operatorname{rank}_R(P)),
	\end{align*} where $\det P$ denotes the highest exterior power of $P.$ Note that $\d_R$ is compatible with short exact sequences and base change in the sense that given an exact sequence $0 \to P_1 \to P_2 \to P_3\to 0$ the natural isomorphism $$\det P_1 \otimes \det P_3\cong \det P_2$$ induces an isomorphism $$\d_RP_1 \cdot \d_RP_3 \cong \d_RP_2.$$ Moreover, for a morphism of rings $R \to S$ we have $\d_R(P)_{S} = \d_S(P\otimes S).$ This functor can be extended to the category $(\cC^p(R),qis)$ of bounded complexes in $\cP(R)$ with quasi-isomorphisms as morphisms. On the level of objects this extension can be described as follows: Let $C^\bullet \in \cC^p(R)$ then $$\d_R(C^\bullet):= \bigotimes_{i \in \ZZ} \d_R(C^i)^{{-1}^i}.$$
This functor is again compatible with exact sequences and if $C^\bullet$ is acyclic, then the quasi isomorphism $0 \to C^\bullet$ induces a trivialisation of $\d_R(C^\bullet)$ that we take as an identification. One can show that $\d_R$ factorises over $(\cD^b_{\mathrm{perf}}(R),qis)$, the image of the category of bounded complexes of finitely generated projective modules in the derived category with quasi isomorphisms as morphisms. If a complex $C^\bullet$ is cohomologically perfect meaning that $H^i(C^\bullet)$ considered as a complex concentrated in degree $0$ is in $\cD^b_{\mathrm{perf}}(R)$ for all $i$, then we have a canonical isomorphism $$\d_R(C^\bullet)=\bigotimes \d_R(H^i(C^\bullet))^{(-1)^i},$$ that we take as an identification. This extension is further compatible with duality and base change in the following sense:
There exist canonical isomorphisms $$\d_R(R\Hom_R(C^\bullet,R)) \cong \d_R(C^\bullet)^{-1}$$ and $$\d_S(S\otimes^\mathbb{L}_R(C^\bullet)) \cong \d_R(C^\bullet)_S.$$

\subsection{Fundamental lines}

%

Let $M$ be a  $(\varphi_L,\Gamma_L)$-module over $\cR_A,$ where $A$ is an affinoid algebra over $K.$ We assume that $M$ satisfies the following technical condition:
\begin{equation}
	\label{cond:charactertype}
	\text{There exist } \mathcal{L} \in \operatorname{Pic}(A) \text{ and } \delta=  \delta_{\det M} \in \Sigma_{an}(A)  \text{ such that} \operatorname{det}_{\cR_A}M \cong \mathcal{L} \otimes_A \cR_A(\delta),
\end{equation}
where $\det M$ denotes the highest exterior power of $M.$ Clearly $\det M$ is always a module of rank $1$ and the technical condition is asking $\det M$ to be of character type up to a twist on the base. The full subcategory of $(\varphi_L,\Gamma_L)$-modules satisfying the above contains all modules that arise as a base change from $\cR_L$  by \cite[Proposition 1.9]{FX} and furthermore contains all trianguline modules (even with $\mathcal{L}=A$).
If $M$ satisfies the above condition the isomorphism class of $\mathcal{L}$ and the character $\delta$ are uniquely determined. Furthermore $\mathcal{L}$ can be identified with the subset $$\mathcal{L}_A(M):= \{x \in\det M \mid \varphi_L(x)= \delta_{\det M}(\pi_L)x, \gamma x = \delta_{\det M}(\gamma)x\}$$
by sending $l\in\mathcal{L}$ to $l\otimes \mathbf{e}_\delta\in \mathcal{L} \otimes_A \cR_A(\delta).$
\begin{definition}
	Let $M$ be an $L$-analytic $(\varphi_L,\Gamma_L)$-module of rank $r_M$ over $\cR_A$ satisfying \eqref{cond:charactertype}. Write $\det(M) = \mathcal{L} \otimes \cR_A( \delta_{\det M}).$  We define
	$$\Delta_{1,A}(M):= \d_{A[\Gamma_L/U]}( K_{\varphi_L,D(U,K)} (M)) \otimes_{A[\Gamma_L/U]}A , $$
	using Remark \eqref{rem:CohomologyRefinement}, and
	$$\Delta_{2,A}(M):= \bigg(\{x \in \det M \mid \varphi_L(x)= \delta_{\det M}(\pi_L)x, \gamma x = \delta_{\det M}(\gamma)x\}, -\chi_{A[\Gamma_L/U]}(K_{\varphi,D(U,K)}(M))\bigg),$$
	i.e., the underlying line bundle of $\Delta_{2,A}$ is $\mathcal{L}$ which has a canonical  $(\varphi_L,\Gamma_L)$-action given by $\delta_{\det M},$. We also write $\mathcal{L}(\delta_{\det M})$ if we wish to emphasize the action.
\end{definition}

\begin{remark}\label{rem:DeltaVonRdelta}
We have
\begin{align*}
\{x \in \cR_K(\delta) \mid \varphi_L(x)= \delta(\pi_L)x, \gamma x = \delta(\gamma)x\}=\cR_K^{\varphi=1,\Gamma_L}\mathbf{e}_\delta= K\mathbf{e}_\delta\cong K
\end{align*}
whence $\Delta_{2,K}(\cR_K(\delta))=(K\mathbf{e}_\delta,1)\cong(K,1)$ using Remark \ref{rem:EP-formula}.
\end{remark}

\begin{proposition}
	$\Delta_{1,A}(M)$ and $\Delta_{2,A}(M)$ are well-defined graded invertible modules and
	$$\Delta_A(M):=\Delta_{1,A}(M) \cdot \Delta_{2,A}(M)$$ satisfies the following properties
	\begin{enumerate}
		\item For any continuous map of affinoid algebras $A \to B$ induces a canonical isomorphism $$\Delta_A(M)\otimes_AB \cong \Delta_B(M\hat{\otimes}_AB).$$
		\item $\Delta_A(M)$ is multiplicative in short exact sequences.
		\item  $ \Delta_A(M) \cong \Delta_A(\tilde{M})^*\otimes (A(\chi^{r_M}),0).$
	\end{enumerate}
\end{proposition}
\begin{proof}
Compatibility with base change can be checked for $\Delta_{i}:=\Delta_{i,A}$ individually. For $\Delta_{1}$ it follows from Theorem \ref{thm:CohomologyFinite} and for $\Delta_2$ it is clear.
	The compatibility with short exact sequences can also be checked individually for $\Delta_{i}.$ For $i=2$ it follows from the corresponding statement for determinants and for $i=1$ it follows from the fact that a short exact sequence of $(\varphi_L,U)$-modules induces a short exact sequence of the complexes $K_{f,D(U,K)}.$
  The quasi-isomorphism $\digamma(M):{K_{\varphi,{\mf Z}}(M)}\cong K_{\varphi,\mf Z}(\tilde{M} )^*[-2]$  induced from  \eqref{f:dualitycomplex} by identifying $\tilde{\tilde{M}}\cong M$   gives an isomorphism $\Delta_{?,1}(M)\cong\Delta_{?,1}(\tilde{M})^*$ while the isomorphism  $\Delta_{A,2}(M)\cong\Delta_{A,2}(\tilde{M})^*\otimes (A(\chi^{r_M}),0)$ arises as follows: First observe that $\tilde{M}$ satisfies \eqref{cond:charactertype}, if $M$ does, and since $\tilde{M} = A(\chi) \otimes_A M^*$ one sees that $\det(\tilde{M}) \cong A(\chi^{\operatorname{rk}(M)}) \otimes \det(M^*).$ Hence we see $\Delta_{A,2}(\tilde{M}) = \Delta_{A,2}(M^*) \otimes A(\chi^r).$ A small calculation shows $\Delta_{A,2}(M^*)= \Delta_{A,2}(M)^*,$ hence the claim.

\end{proof}

\begin{definition}\label{def:quasi-Stein}
	Let $X$ be a  rigid analytic space over $K.$ Given a family of $(\varphi_L,\Gamma_L)$-modules $M$ over $\mathcal{O}_X,$ i.e., a compatible collection of $(\varphi_L,\Gamma_L)$-modules $M_A$ over $\cR_A$ for every affinoid $\operatorname{Sp}(A)\subseteq X$, we define $\Delta_X(M)$ as the global sections of the line-bundle   $\mathfrak{D}_X(M)$ defined by $\operatorname{Sp}(A) \mapsto \Delta_{A}(M_A).$
	If $X$ is quasi-Stein covered by an increasing union $X_n$ of affinoids we also have
	$\Delta_X(M)[0]=\operatorname{R}\Gamma(X,\mathfrak{D}_X(M))=\operatorname{Rlim}(\Delta_{X_n}(M_{X_n})) \cong \lim \Delta_{X_n}(M_{X_n})[0]$ by Theorem B for quasi-Stein spaces.  We have analogous definitions and statements   for $\Delta_{i,X}(M)$and $\mathfrak{D}_{i,X},$ $i=1,2$  respectively. $\mathfrak{D}_X,\mathfrak{D}_{i,X}$ are graded invertible $\mathcal{O}_X$-modules by definition.
\end{definition}
 A word of caution is in order. A priori the $\Delta_{i,X}(M)$ are not necessarily graded invertible $\mathcal{O}_X(X)$-modules because the global sections do not have to be finitely generated over $\mathcal{O}_X.$ In our applications (in section \ref{sec:RankoneCase}) we will have $\Delta_{2,X} = \mathcal{O}_X(X)$ and will be in a position to apply the subsequent remark in order to conclude that $\Delta_{1,X}$ is an invertible $\mathcal{O}_X(X)$-module.

\begin{remark} \label{rem:formalismRlimDet}
	Let $ X = \bigcup X_n$ be a quasi-Stein space. Let $C_n^\bullet$ be a family of perfect complexes of $\mathcal{O}_{X}(X_n)$-modules together with quasi-isomorphisms $\mathcal{O}_{X}(X_{n-1}) \otimes^{\mathbb{L}}_{\mathcal{O}_{X}(X_n)} C^\bullet_n \simeq C^\bullet_{n-1}.$ Assume that there exists a perfect complex $C^\bullet$ of $\mathcal{O}_X(X)$-modules (in the ring-theoretic sense \footnote{Here one has to make a distinction between a perfect complex of $\mathcal{O}(X)$-modules and a perfect complex of sheaves of $\mathcal{O}_X$-modules, i.e., a complex whose restriction to each $\mathcal{O}_{X_n}$ is perfect. One can show that $C$ is isomorphic to $\mathbf{Rlim}C_n.$ Hence this remark could be restated to require $\mathbf{Rlim}C_n$ to be perfect.}) such that $\mathcal{O}_{X}(X_n) \otimes^{\mathbb{L}}_{\mathcal{O}_{X}(X)} C^\bullet\simeq C^\bullet_{n}$.\\
 Then we have $\d_{\mathcal{O}_{X}(X_n)}(C^\bullet_n)\cong \mathcal{O}_{X}(X_n)\otimes \d_{\mathcal{O}_X(X)}(C^\bullet)$.
Furthermore $\d_{\mathcal{O}_X(X)}(C^\bullet)$ is coadmissible, i.e., $\d_{\mathcal{O}_X(X)}(C^\bullet) = \varprojlim_n \d_{\mathcal{O}_{X}(X_n)}(C_n^\bullet).$
\end{remark}
\begin{proof}
	The proof is formal using that determinant functors commute with derived tensor products and $\mathcal{O}_{X}(X_n) \to \mathcal{O}_{X}(X_{n-1})$ is flat together with the fact that $\d(C^\bullet)$ is a rank one projective module over $\mathcal{O}_X(X)$ and hence coadmissible by \cite[Corollary 3.4]{ST}.
\end{proof}

\subsection{Statement}
 We expect that the results in section \ref{sec_analcoho} extend to affinoids (where only stated or proven over fields) and to all analytic $(\varphi_L,\Gamma_L)$-modules (where only stated for rank one or trianguline ones), explicitly this refers to Remark \ref{rem:EP-formula} and Theorems \ref{thm:perf}, \ref{thm:def}, \ref{lem:explicitPairing}. Hence we state the conjecture below in this level of generality.

\begin{conjecture}\label{conj}
	Choose a compatible system $u=(u_n)$ of $[\pi_L^n]$-torsion points of the Lubin-Tate group and a generator $t'_0$ of $T_\pi'$. Let $A$ be an affinoid algebra over $K,$ a complete field extension of $L$ containing $L^{ab}.$ For each $L$-analytic $(\varphi_L,\Gamma_L)$-module $M$ over $\cR_A$  satisfying condition \eqref{cond:charactertype}
	there exists a unique trivialisation $$\varepsilon_{A,u}(M): \u_A \xrightarrow{\cong}\Delta_A(M)$$ satisfying the following axioms:
	\begin{enumerate}
		\item For any affinoid algebra $B$ over $A$ we have $$\varepsilon_{A,u}(M)\otimes_A \id_B = \varepsilon_{B,u}(M\hat{\otimes}_AB)$$ under the canonical isomorphism $\Delta_A(M)\otimes_AB \cong \Delta_B(M\hat{\otimes}_AB).$
		\item $\varepsilon_{A,u}$ is multiplicative in short exact sequences.
		\item For any $a \in  o_L^\times$ we have $$\varepsilon_{A,a \cdot u}(M)=\delta_{\det M}(a)\varepsilon_{A,u}.$$
		\item $\varepsilon_{A,u}(M)$ is compatible with duality in the sense that
		$$\varepsilon_{A,u}(\tilde{M})^* \otimes h(\chi^{r_M}) = (-1)^{\dim_KH^0(M)}\Omega_{t'_0}^{-r_{M}} \varepsilon_{A,-u}(M)$$ under the natural isomorphisms $\mathbf{1}_A \cong \mathbf{1}_A \otimes \mathbf{1}_A$ and $\Delta(M)\cong \Delta(\tilde{M})^* \otimes (A(r_M),0),$ where $h(\chi^{r_M})  \colon A(r_M) \to A$ maps $e_{\chi^{r_M}}$ to $1.$
		\item For $L=\Qp$, $\pi_L=p$ and $u=(\zeta_{p^n}-1)_n$ the trivialisation coincides with that of Nakamura, in the sense of Proposition \ref{prop:precisecomp}.
		\item Let $F/L$ be a finite subextension of $K,$ $M_0$ be a de Rham $(\varphi_L,\Gamma_L)$-module over $\cR_F$ and $M = K\hat{\otimes}_FM_0.$ Then $$\varepsilon_{K,u}(M)=\varepsilon^{dR}_{F,u}(M_0).$$
	\end{enumerate}
\end{conjecture}

\begin{remark}
	\label{rem:compNakamura}
	\begin{enumerate}
	
	\item The occurrence of the power of $\Omega$ in the compatibility with duality (iv) is a conceptually new phenomenon in our conjecture, see also Proposition \ref{prop:duality}.
\item Due to the equivalence of categories stated in \cite[Thm.\ 3.16]{BSX} there is an analogous conjecture for   $L$-analytic $(\varphi_L,\Gamma_L)$-modules over the character variety, i.e., by replacing the usual Robba ring $\cR_K=\cR_K(\mathbf{B})$ (attached to the open unit ball $\mathbf{B}$) by  the Robba ring $\cR_K(\mathfrak{X}_{o_L} )$ of the character variety $\mathfrak{X}_{o_L}$ attached to the group $o_L$, see \cite[\S 2.4]{BSX}   or \cite[\S 4.3.6]{SV20}. In this situation, we expect that the conditions concerning $K$ can be weakened and perhaps the descent to   $L$ (or any finite extension of it) instead of the huge field $K$ should be feasible, compare with  Thm.\ 4.3.23 in (loc.\ cit.). Moreover, due to \cite[Lem.\ 4.3.25]{SV20} there should be no occurrence of $\Omega$! We will pursue this in future work.
\item   The assumption that $K$ contains $\widehat{L^{ab}}$ can be dropped in the case that $L=\QQ_p$ as the period $\Omega_{\QQ_p}$ can be taken to be any element in  $\ZZ_p^\times.$ In order to specialise our construction to Nakamura's one has to make more specific choices. Fixing an element $\gamma\in \Gamma,$ whose image in $\Gamma/\Gamma_{p-\text{power-torsion}}$ is a topological generator implicitly determines the period as $\Omega_{\QQ_p}= \log_0(\chi_{cyc}(\gamma))^{-1}.$ But this would not necessarily be compatible with Nakamura's variant of the de Rham isomorphism, since his variant does not involve any period. Instead one should choose a $\gamma$ such that $\log_0(\chi_{cyc}(\gamma))=1.$
This defect is due to the fact our variant of the exponential map involves the period $\Omega$ as part of its definition and hence so does our de Rham isomorphism.
This is not a contradiction to the uniqueness of the $\varepsilon$-isomorphisms in question. Indeed in the rank one case, we can see the $\varepsilon$-isomorphism is determined by its behaviour at de Rham points. If $\Omega_{\QQ_p} \neq 1$ then our variant asks for a different behaviour at these de Rham points thus leading to a different result.
\end{enumerate}
\end{remark}

\subsection{The de Rham case} \label{sec:deRham}
	In this section we explain how to attach a Weil-Deligne Representation to an $L$-analytic de Rham $(\varphi_L,\Gamma_L)$-module over $\cR_L$ in order to define the de Rham epsilon-constants.  We denote by $\mathbf{B}_{?}$ for $?  \in \{\cris,\dR,\mathrm{st}\}$ Fontaine's usual period rings. Without difficulty this construction can be generalised to $(\varphi_L,\Gamma_L)$-modules over $F \otimes_L \cR_L$ for a finite extension $F$ with trivial action. In order to keep notation light we will assume without loss of generality $F=L.$ We write  $\mathbf{B}_{e,LT}=\tilde{\cR}[1/t_{LT}]^{\varphi_L=1}.$
	We will make use of the equivalence of categories between $L$-analytic $(\varphi_L,\Gamma_L)$-modules and $L$-analytic $B$-pairs originally suggested in \cite[Remark 10.3]{Be16} and detailed in \cite[Theorem 5.5]{Poy}.
	A priori these results  are only applicable to $E$-linear representations of $G_L,$ where $E$ denotes a Galois closure of $L/\QQ_p.$ If we start with an analytic $(\varphi_L,\Gamma_L)$-module $M$ over $\cR_L$ then by \cite{Poy} we can attach to $E \otimes_L M$ a $B$-pair (called $B_{\id}$-pair in (loc.\ cit.)), i.e., a pair consisting of a finite free $E\otimes_L\mathbf{B}_{dR}^+$-module $W^+_{dR,\id,E}$ with a $\mathbf{B}_{dR}^+$-semi-linear (and $E$-linear) $G_L$-action and a finite free $\mathbf{B}_{e,LT,E}:= E\otimes_L\mathbf{B}_{e,LT}$-module $W^{LT}_{\id,E}$ with semi-linear $G_L$-action together with an isomorphism after base change to $\mathbf{B}_{dR}.$
	By Galois descent, taking invariants with respect to the $G(E/L)$-action (acting via the first tensor factor) provides us with a $B$-pair $W(M):=(W_{dR}^+(M),W_e(M))$ over $(\mathbf{B}_{dR}^+,\mathbf{B}_{e,LT}).$
The ring $\mathbf{B}_{e,LT}$ can be viewed as a subring of $\mathbf{B}_{\cris,L}.$ Indeed, since $\varphi(t_{LT}) = \pi_L t_{LT}$ it suffices to consider elements of $\tilde{\cR}$ satisfying $\varphi(x) = \pi_L^jx$ for some $j \in \ZZ,$ which by Frobenius regularisation are already contained in $\tilde{\cR}^+$ (cf. \cite[Proposition 3.2]{Be02} in the cyclotomic case, and a similar result holds for ramified Witt-vectors as well (cf. \cite[Satz 3.19]{Ste3})). The ring $\tilde{\cR}^+$ is a subring of $\mathbf{B}_{\cris,L}.$
We call a $B$-pair $(W_{dR}^+,W_e)$ \textbf{de Rham} if $W_{dR}^+[1/t_{LT}]$ admits a $G_L$-invariant basis. One can show, that this is equivalent to the corresponding $(\varphi_L,\Gamma_L)$-module being de Rham (cf. \cite[Section 3.2, Proposition 3.7]{Por} for a proof in the \'{e}tale case). Note that our notion of de Rham coincides with $L$-de Rham in loc.\ cit.). Consider for $F/L$ finite the  vector space $$D_{st}(M_{\mid F}):= (\mathbf{B}_{\mathrm{st}} \otimes_{\mathbf{B}_{e}}W_e(M))^{G_F}$$
over the maximal unramified subextension $F'$ of $F/L$. We  define $D_{\pst}(M)$ as their colimit over all $F/L$ finite.
By a standard argument (cf. proof of Theorem 2.13 Part (1) in \cite{FO}), each $F'$-vector space $D_{st}(M_{\mid F})$ is of dimension $\leq \rk M$ and $D_{pst}(M)$ is hence an $L^{nr}$-vector  space of dimension $\leq \rk M.$ We say that $M$ is \textbf{potentially semi-stable} if this dimension is precisely $\rk M$ or, equivalently, if there exists a finite extension $F/L$ such that $D_{st}(M_{\mid F})$ is an $F'$-vector space of dimension $\rk M.$
The $p$-adic monodromy theorem also holds for $B$-pairs in the cyclotomic case and there is an obvious $L$-analytic analogue providing us with the following (see \cite[Corollary 3.10]{Por} for a treatment in the \'{e}tale case).
\begin{remark}
	$M$ is de Rham if and only if $M$ is potentially semi-stable.
\end{remark}
 Note that $D_{\pst}(M)$ naturally has a semi-linear $G_L$-action and inherits from $\mathbf{B}_{\mathrm{st},L}=\mathbf{B}_{\mathrm{st}}\otimes_{L_0}L$ an action of $\varphi_q$ and the monodromy operator $N$ satisfying $N\varphi_q = q\varphi_qN.$

 We now explain how to modify this action in order to obtain an $L^{nr}$-linear representation of the Weil group $W_L.$ By local class field theory the maximal abelian extension $L^{ab}$ of $L$ is given by the composite $L^{nr}L_\infty$ and $L^{nr}\cap L_\infty=L.$ Consider the reciprocity map $$\operatorname{rec}_L: L^\times \to \Gal(L^{ab}/L),$$ which by our convention sends  $\pi_L$ to the geometric Frobenius on $L^{nr}.$ This induces an isomorphism $L^{\times} \cong W_L^{ab} \cong \varphi_L^{\ZZ} \times \Gamma_L.$ We denote by $\bar{\phantom{m}}:W_L\twoheadrightarrow W_L^{ab}$ the canonical surjection and define a linearised action of $W_L$ on $D_{pst}(M)$ by setting
$$\rho_{lin}(g)(x) := \varphi_q^{v_\pi(\operatorname{rec}^{-1}(\bar{g}))}(\rho_{semi-lin}(g)(x)),$$   where $\rho_{semi-lin}$ denotes the action we considered previously. For $a \in L^{nr}$, we then have
 \begin{align*}\rho_{lin}(g)(ax) &= \varphi_q^{v_\pi(\operatorname{rec}^{-1}(\bar{g}))}\big(\rho_{semi-lin}(g)(ax)\big)\\
 	&=\varphi_q^{v_\pi((\operatorname{rec}^{-1}(\bar{g}))}\big(\rho_{semi-lin}(g)(a)\rho_{semi-lin}(g)(x)\big)\\
 	&= \varphi_q^{v_\pi(\operatorname{rec}^{-1}(\bar{g}))}\big((\varphi_L^{-1})^{v_\pi(\operatorname{rec}^{-1}(\bar{g}))}(a)\big)\cdot \varphi_q^{v_\pi(\operatorname{rec}^{-1}(\bar{g}))}\big(\rho_{semi-lin}(g)(x)\big)\\
 	&= a \rho_{lin}(g)(x).
 \end{align*}

By passing to the base change $D_{pst}(M)\otimes_LL_{\infty}$ (with trivial action on $L_\infty$) we are finally able to define $W(M):= (D_{pst}(M)\otimes_LL_{\infty}, \rho_{lin}, N)$ which is an $L^{ab}$-linear Weil-Deligne representation (Note that since $D_{pst}(M)$ can be written as a base extension of some $D_{st}(M_{\mid F})$, the action of the inertia group $I_F$ is discrete and because $I_F$ is open in $W_L$ the action of $W_L$ is discrete.)

\begin{example}\label{rem:ExChar}
The linearized Weil-Deligne representation $W:=W(\cR_K(\delta))$ with $\delta=\delta_\mathrm{lc}x^k $is given by the character $ \delta_W= \delta_\mathrm{lc}  \delta^{un}_{\pi_L^{-k}}:L^\times\to (L^{ab})^\times$ via class field theory sending $\pi_L$ to the geometric Frobenius.
In particular,
\begin{equation}\label{f:changeOfCharacter}
  (\delta_W)_{|o_L^\times}=\delta_{|o_L^\times}(x^{-k})_{|o_L^\times}.
\end{equation}
\end{example}
%
\begin{proof}
	For the convenience of the reader we give a proof using $B$-pairs.
	Let $\mathbf{e}_{\delta}$ be the obvious basis of $\cR_L(\delta)$ and write $\delta(\pi) = \delta_{lc}(\pi) \pi^k =\pi^l\alpha$ with $\alpha \in o_L^\times.$ We can find $a \in L^{nr} \subset \tilde{\cR}$ such that $\varphi_q(a) = \alpha a$ and hence $y:=\frac{1}{a t_{LT}^l}\otimes  \mathbf{e}_{\delta} \in W_e(\cR_L(\delta)) = (\tilde{\cR}[1/t_{LT}] \otimes_{\cR} \cR_L(\delta))^{\varphi_L=1}.$ Note that $G_L$ acts diagonally on $W_e(\cR_L(\delta)),$ where the action on $\cR_L(\delta)$ is given via the quotient $\Gamma_L.$ Let $F$ be a field extension of $L$ such that $\delta_{lc}$ is trivial when restricted to the image of $G_F$ in $\Gamma_L.$ Then the action of $g \in G_F$ is given by $g(y) = a/g(a) \chi_{LT}(g)^{k-l}y$ and hence $z:=  t_{LT}^{l-k}a \otimes y$ is a basis of $(B_{st}\otimes W_e(\cR_L(\delta)))^{G_F}.$
	Write $\delta = x^k \delta^{un}_{\pi_L^{l-k}} \delta^{un}_{\alpha} \rho$ where $\rho$ is a locally constant character   with $\rho(\pi)=1$ and $\rho(\gamma) = \delta_{lc}(\gamma)$ for $\gamma \in \Gamma_L.$ In this representation it is clear that the residual (non-linearised) action of $G_L$ is given by $gz = \rho(g)z$ and hence the linearised action is given by
\begin{align*}
	\rho(g)\varphi_q^{v_\pi(\operatorname{rec}^{-1}g)}(z)& = (\pi^{l-k}\alpha)^{v_\pi(\operatorname{rec}^{-1}g)}\rho(g)z \\
	& = \delta_{lc}(\pi)^{v_\pi(\operatorname{rec}^{-1}g)} (\pi^{-k})^{v_\pi(\operatorname{rec}^{-1}g)}\rho(g) = \delta_{lc} \delta^{un}_{\pi^{-k}}(g)z.
\end{align*}
\end{proof}

\subsubsection{Equivariant de Rham epsilon constants}

	For   a de Rham $(\varphi_L,\Gamma_L)$-module  $M$ over $\cR_L$ we would like to define the epsilon constant of $M$ to be the $\varepsilon$-constant associated to $W(M)$ $$\varepsilon(M,\psi,dx):=\varepsilon_{L^{ab}}(L,W(M),\psi,dx)$$ defined in section \ref{sec:varepsilon-constants} using the adjustment \eqref{eq:adjustedconstant}.
In the cyclotomic case (take for simplicity $L = K = \QQ_p$), these $\varepsilon$-constants can be viewed as elements of $L_n = \QQ_p(\zeta_{p^n}).$ In our case the constants are defined using $p$-power roots of unity which are "built" from the LT-torsion points using the power series $\eta(-,T).$ The problem we run into is that, contrary to the classical case, we can not assume that $L_n$ contains the $p$-power roots of unity.

Suppose $K$ contains $L^{ab}.$ Then it makes sense to view $\varepsilon(M,\psi,dx)$ as an element of $K,$ but by our convention that $K$ carries the trivial $\Gamma_L$-action, we do not have $\gamma(\varepsilon(M,\psi_u,dx)) = \varepsilon(M,\psi_{\gamma(u)},dx),$ which we will need for technical reasons below in \eqref{f:fMu},in the form of Remark \ref{rem:equivariance}. Roughly speaking we would like to define the $\varepsilon$-constants as elements of $L_n\otimes_L K$ with $n$ large enough, such that the definition of the epsilon constants ``involves only'' the $\pi_L^n$-division points of the Lubin-Tate group. We make this concept precise via the following equivariant construction.

\begin{definition} \label{def:equivariantconstant} Suppose the complete subfield $K$ of $\CC_p$ contains $L^{ab}$ and
let $W$  be a Weil-Deligne representation of $W_L$ with coefficients in $K.$ Building on the $\varepsilon$-constants  defined in section \ref{sec:varepsilon-constants} with $E=K$ we   define the \textbf{$\Gamma_L$-equivariant $\varepsilon$-constant}
 $$\tilde{\varepsilon}(W,u,dx):= (\varepsilon_{K}(L,W,\psi_{\hat{\tau}(u)},dx) )_{\tau},$$  for   lifts $\hat{\tau}$ of $\tau$ to $\Gamma_L$,  viewed  as an element of\footnote{If $\delta$ takes values in a finite extension $F$ of $L$ and $W=W(\cR_F(\delta))$, then  as an element of $$\prod_{\tau\colon L_n\to K; \sigma\colon F\to K} K\cong \prod_{\tau\colon L_n\to K}F\otimes_L K \cong F\otimes_L L_n\otimes_LK\subseteq F\otimes_L L_\infty\otimes_L K$$ assuming $F\subseteq K$ for the first isomorphism. Also the $\sigma$ should be involved as $W(M)_\sigma$ in the defining tuple then.} $$\prod_{\tau\colon L_n\to K} K \cong L_n\otimes_LK\subseteq L_\infty\otimes_L K$$
via the canonical isomorphism, where $n \gg 0$ is large enough such that the $\varepsilon$-constant can be defined in terms of characters of conductor $\leq n$ according to Deligne's (inductive) construction principle: In the rank one case, i.e., in the case of a locally constant character $\delta \colon L^\times \to K^\times$, one can take $n \geq a(\delta).$ In general, the definition of the $\varepsilon$-constant involves multiple such characters defined over finite extensions of $L$ (cf. \cite[p. 536, Equation 4.2.1]{deligne}) and one has to choose $n$ greater than the supremum of all appearing conductors.
\end{definition}
\begin{remark}\label{rem:equivariance}
The $\varepsilon$-constant $\tilde{\varepsilon}(W,u,dx)$
is well-defined, i.e. its definition above is independent
of the choices of the lifts $\hat{\tau}$.
Furthermore,
with respect to the $\Gamma_L$-action on $L_n \otimes_LK$ via the left tensor factor we have
	$$\gamma(\tilde{\varepsilon}(W,u,dx)) = \tilde{\varepsilon}(W,\gamma(u),dx) = (1 \otimes \delta_{\det W}(\chi_{LT}(\gamma)))\tilde{\varepsilon}(W,u,dx).$$
	
\end{remark}
\begin{proof}
	Without loss of generality we can assume $W$ is of rank one corresponding to a locally constant character $\delta \colon L^\times \to K^\times$ due to Deligne's   construction principle.
 First of all we note that $\tilde{\varepsilon}$ is well-defined since $ u_{a(\delta)} \in L_n$ by assumption. Because the natural isomorphism $L_n \otimes_L K \cong \prod_{\tau\colon L_n \to K}K$ maps $u \otimes 1$ to $(\tau(u))_\tau,$ we can see that $\tilde{\varepsilon}$ is obtained by replacing in \eqref{eq:constantsrank1u} the elements $\eta(a,u_{a(\delta)})$ by the series $\eta(a,T)$ evaluated at the element $(u_{a(\delta)} \otimes 1)$, i.e., by $\sum_{i \in \NN_0} (1 \otimes a_i) (u_{a(\delta)} \otimes 1)^i,$ where $\eta(a,T) = \sum a_iT^i$ (this expression converges with respect to the tensor product topology). The formula for the $\gamma$-action can be read off from \eqref{eq:constantsrank1u}.
\end{proof}

\begin{definition}
	For   a de Rham $(\varphi_L,\Gamma_L)$-module  $M$ over $\cR_L$ we define the epsilon constant of $M$ to be the $\Gamma_L$-equivariant $\varepsilon$-constant associated to $W(M)$ $$\tilde{\varepsilon}(M,u,dx):=\tilde{\varepsilon}(W(M),u,dx).$$ We usually omit $dx$ from the notation and write \[\tilde{\varepsilon}(M,u):=\tilde{\varepsilon}(M,\psi_u,dx).\]
\end{definition}

\begin{remark}\label{f:epsilondual} Let $dx$ be the self dual Haar measure with respect to $\psi_u,$ then
 \[\tilde{\varepsilon}(M,-u,dx)\tilde{\varepsilon}(\tilde{M},u,dx)=1\]
\end{remark}	
\begin{proof}
 In order to apply \eqref{f:epsduality} we  check that  we have an isomorphism $\mathbf{D}_{\mathrm{pst}}(\tilde{M})\cong \mathbf{D}_{\mathrm{pst}}({M})^*(|x|)$. Using the usual  functorialities it suffices to check that $\mathbf{D}_{\mathrm{pst}}(\Omega^1)\cong L^{nr}(|x|)$,  which is a special case of Example \ref{rem:ExChar}.
%
The proof of the other required equation
 \[\det(-\varphi|\mathbf{D}_{\mathrm{st}}({M})/\mathbf{D}_{\mathrm{cris}}({M}))\det(-\varphi|\mathbf{D}_{\mathrm{st}}(\tilde{M})/\mathbf{D}_{\mathrm{cris}}(\tilde{M}))=1\]
is then also standard, see e.g.\ \cite[claim 5 in proof of Prop.\ 2.2.20]{da}.
\end{proof}

We now describe how our construction relates to the \'{e}tale and the cyclotomic case.
The comparison of $\varepsilon$-constants involves a number of choices and we will only give an informal comparison of the constructions presented here and the ones from \cite{NaANT} - by which we mean that we give a comparison up to constants that only depend on $L/\QQ_p.$
There are two avenues to be considered. On the one hand, we can specialise our constructions to the cyclotomic case $L=\QQ_p,$ taking $u_n = \zeta_{p^n}-1$ and $\Omega=1.$ Because $\zeta_{p^n} = 1+ u_n = \eta(1,(u_n))$ in this case our construction specialises to Nakamura's, more precisely, our $\tilde{\varepsilon}$ is equal to $\varepsilon_{\text{Na}} \otimes 1$ viewed as an element of $L_\infty \otimes_{\QQ_p} K,$ where $\varepsilon_{Na}$ denotes the constant from \cite[Section 3C]{NaANT}. Indeed the elements $\eta(a,u_n \otimes 1) = \zeta^a_{p^n} \otimes 1$ appearing in \ref{rem:equivariance} lie inside $L_n \otimes_{\QQ_p} \QQ_p.$

On the other hand, we take the induction of an $L$-linear $G_L$-representation $V$ and treat it as an $L$-linear representation of $G_{\QQ_p}$. For the moment let us assume $V \in \Rep_L{G_L}$ is semi-stable and $L$-analytic and set $X:= \Ind_{L/\QQ_p}V.$ Let $\QQ_p\subseteq L_0 \subseteq L$ be the maximal unramified subextension.
We can decompose $$(\mathbf{B}_{st}\otimes_{\QQ_p}V)^{G_L} \cong \prod_{\tau\colon L_0 \to \overline{\QQ_p}}(\mathbf{B}_{st}\otimes_{L_0,\tau}V)^{G_L}$$
and have a similar decomposition for $(\mathbf{B}_{st} \otimes_{\QQ_p} X)^{G_{\QQ_p}}.$
The epsilon constants of the induction (given suitable choices of additive characters) are related by explicit constants independent of $V$ (see \eqref{f:comppsiu}). Ignoring these,
the $\varepsilon$-constants defined by Nakamura are the product of the $\varepsilon$-constants of each component in the sense that he attaches to $X$ a tuple  $(W_\tau)_\tau$ of $W_{\QQ_p}$-representations to which he attaches a tuple $(\varepsilon(W_\tau))_\tau$ (cf. \cite[p.359]{NaANT} for details) of constants living over $\QQ_{p}(\zeta_{p^{\infty}})\otimes_{\QQ_p}L$. In contrast we attach (informally speaking) to the $\tau=\id$ component a constant $\varepsilon(W_{id}).$ As we can not assume that $L_\infty$ contains the $p$-th roots of unity, an analogous construction involving $L_\infty$ does not work in the obvious sense and taking the base change to $L^{ab}$ with $G_{L}$ acting naturally on $L^{ab}$ does not provide us with the Galois action needed to make the constructions in \ref{rem:equivariance} work. By assuming $L^{ab} \subset K$ we can make sense of the elements $\eta(a,(u_n \otimes 1)) \in L_n \otimes_LK,$ which allow us to define $\tilde{\varepsilon}(W_{id})$ with the desired technical properties now living over $L_n \otimes_LK= \prod_{\sigma \in \Hom_L(L_n,K)} K$ for $n \gg 0$ (note that the index set of the product is different in comparison to Nakamura's situation). By projecting to the $\sigma = \id$ component we can recover  Nakamura's $\varepsilon(W_{\id})$  and our constant $\tilde{\varepsilon(W_{\id})}=(\varepsilon(W_{\id},\psi_{\hat{\sigma}(u)},dx))_\sigma$ should be informally thought of as $(\sigma(\varepsilon(W_{\id})))_{\sigma}$, which is not well-defined as $\sigma$ does not act on $K$ but only on $L_n.$

  The fact that $V$ is semi-stable and $L$-analytic forces each non-identity component to be potentially unramified (since they are semi-stable with Hodge-Tate weights $0$). If they are even unramified, all $\varepsilon$-constants at non-identity components would be $1$ and both methods give comparable $\varepsilon$-constants (more precisely, at $\sigma=\id$ they would be the same up to explicit constants independent of $V$).
If the action on the non-identity components is only potentially unramified, we cannot assume that the $\varepsilon$-constants at the non-identity embeddings are $1.$ In particular these embeddings contribute to the $\varepsilon$-constant attached to $V$ by Nakamura in a way that can not be captured by only considering the identity embedding.

%

%

\subsubsection{The de Rham epsilon-isomorphism}

For each de Rham $(\varphi,\Gamma_L)$-module $M$ over $\cR_K$ which arises as a base change of a $(\varphi_L,\Gamma_L)$-module $M_0$ over $\cR_F$ for some finite extension $F/L$, and for each generator $u$ of $T_{\pi }$  we are going to define a trivialization
\begin{equation}\label{f:epsdR}
  \varepsilon^{dR}_{F,u}(M_0):\mathbf{1}_K\xrightarrow{\cong} \Delta_K(M)
\end{equation}
as product of three terms
\[\varepsilon^{dR}_{F,u}(M_0):=\Gamma(M)\cdot \Theta(M)\cdot \Theta_{dR,u}(M_0)\]
where
\begin{align*}
  \Theta(M)\colon & \mathbf{1}\xrightarrow{\cong}  \Delta_{K,1}(M)\d_K(  \bD_{\tn{dR}}(M)  ),\\
  \Theta_{dR,u}(M_0)\colon & \d_K(  \bD_{\tn{dR}}(M)  ) \xrightarrow{\cong}\Delta_{K,2}(M),\\
  \Gamma(M) & \in K^\times.
\end{align*}
To keep notation light and consistent with the previous subsection we will, without loss of generality, restrict ourselves to the case $L=F.$
Firstly, we define $\Gamma(M)$, which depends only on the Hodge--Tate weights of $M$. For $r \in \ZZ$ let
  \[ n(r) = \dim_K \gr^{-r} \bD_{\tn{dR}}(M) ,\]
  so $n(r)$ is the multiplicity of $r$ as a Hodge--Tate weight of $M$. We adopt the convention in this paper that the Hodge--Tate weight of the cyclotomic character is $1$. We define
  \[ \Gamma^*(r) := \begin{cases} (r-1)! & \text{if $r > 0$,}\\ \frac{(-1)^r}{(-r)!} & \text{if $r \le 0$,}\end{cases}\]
  the leading coefficient of the Taylor series of $\Gamma(s)$ at $s = r$.
  Then we set
  \[ \Gamma(M) := \prod_{r \in \ZZ} (\Omega^r\Gamma^*(r))^{-n(r)}.\footnote{  $\Gamma^*(k)$ in \cite{NaANT} has been replaced by   $\Omega^k\Gamma^*(k)$ in our setting. }
\]

Secondly, $\Theta(M)$ is obtained by applying
the determinant functor to the following exact sequence \begin{equation}\label{f:longExactexp}
\begin{split} &\xymatrix@C=0.5cm{
0\ar[r]^{ } & {H^{0}_{\varphi,\mathfrak Z_n}(M)^{\Gamma_L} } \ar[r]^{ } & \mathbf{D}^{}_{\tn{cris}}(M) \ar[r]^{ } &
\mathbf{D}^{}_{\tn{cris}}(M)\oplus t_M \ar[rr]^<(0.3){ } && {H^{1}_{\varphi,\mathfrak Z_n}(M)^{\Gamma_L}} \ar[r]^{ } &     }\\
&\xymatrix@C=0.5cm{ \ar[rr]^<(0.1){  } && \mathbf{D}^{}_{\tn{cris}}(\tilde{M})^*\oplus
\mathbf{D}^{}_{\tn{dR}}({M})^0 \ar[r]^{ } & \mathbf{D}^{}_{\tn{cris}}(\tilde{M})^* \ar[r]^{ } & {\H^2_{\varphi,\mathfrak Z_n}(M)^{\Gamma_L} } \ar[r]^{ } & 0,  }
\end{split}
\end{equation}
 which arises from joining the bottom exact sequence of \eqref{f:exp-sequence} with the dual of the same sequence applied to $\tilde{M}$
  by local duality
  $\langle -,-\rangle_{{M}}$
 in \ref{prop:orthogonal} and using Remark \ref{rem:dualitydeRham}, upon
  \begin{enumerate}
  \item using the tautological exact sequence $ \xymatrix@C=0.5cm{
    0 \ar[r] & \mathbf{D}^{}_{\tn{dR}}({M})^0\ar[rr]^{ } && \mathbf{D}^{}_{\tn{dR}}({M}) \ar[rr]^{} && t_M \ar[r] & 0 }$ as well as de Rham duality in the form \[\mathbf{D}^{}_{\tn{dR}}({M})^0\xrightarrow{\cong} t_{\tilde{M}}^*, \; x\mapsto\{\bar{y}\mapsto [y,x]_{dR,\tilde{M}}\},\]
    and
  \item identifying each time the two instances of $\mathbf{D}^{}_{\tn{cris}}(\tilde{M})^*$ and $\mathbf{D}^{}_{\tn{cris}}({M})$, respectively, by the identity.
  \end{enumerate}

 Thirdly - here comes the reason why we use a model $M_0$ -  $\Theta_{dR,u}(M_0):=f_{M_0,u}^{-1}$ is defined by the analogue of \cite[Lem.\ 3.4]{NaANT}  which   - using Remark \ref{rem:equivariance} and \eqref{f:changeOfCharacter} -     induces an isomorphism $f_{M_0,u}:\Delta_{K,2}(M)\xrightarrow{\cong} d_K( \mathbf{D}^{}_{\tn{dR}}({M}))$ from the map   (taking into account Remark \ref{rem:DdifBasechange})
 \begin{align}\label{f:fMu}
   \mathcal{L}_K(M) \to \mathbf{D}^{}_{\tn{dif},n}(\mathrm{det}_{\cR_K}M)=&K_n((t_{LT}))\otimes_{\iota_n,\cR_K^{(n)}}(\mathrm{det}_{\cR_K}M)^{(n)} \\ \notag
   & x\mapsto  \big(\tilde{\varepsilon}(M_0,u)^{-1}  \cdot  \frac{1}{t_{LT}^{h_M}}\big) \otimes \varphi^n(x)
 \end{align}
 for sufficiently large $n$ such that  the equivariant constant $\tilde{\varepsilon}(W(M_0),u,dx)$ from Definition \ref{def:equivariantconstant} lies in $L_n \otimes K\subseteq K_n((t_{LT}))$, where $h_M$ denotes the Hodge-Tate weight of $\det M.$  One easily checks independence of the choice of a model $M_0$ - the reason why we use $M_0$ in the notation is to indicate that we need a model to define these objects.
  Note that \eqref{f:fMu} depends on $u$ in two ways. On the one hand via $\tilde{\varepsilon}$ and on the other hand due to the explicit appearance of $t_{LT}$ which, as pointed out in \ref{rem:hiddendependence}, depends on the choice of $u.$ An analogous computation to \cite[Remark 3.5]{NaANT} shows that $f_{M_0,au} = \delta_{\det_{\cR_K}M}(a)^{-1} f_{M_0,u}$ for $a \in o_L^\times.$

\begin{proposition}[Properties (ii) and (iv) for $ \varepsilon^{dR}_{L,u}(M_0)$]\label{prop:duality}
\phantom{mmmmmmmmmmmmmmmmmmmmmmmmmmmmmmmmmmmmmmmmm}
\begin{enumerate}
  \item For any exact sequence $\xymatrix@C=0.5cm{
    0 \ar[r] & M_1 \ar[rr]^{ } && M_2 \ar[rr]^{ } && M_3 \ar[r] & 0 }$, we have
    \[\varepsilon^{dR}_{L,u}(M_{2,0})=\varepsilon^{dR}_{L,u}(M_{1,0})\otimes\varepsilon^{dR}_{L,u}(M_{3,0})\]
    under the canonical isomorphism $\Delta_K(M_2)\cong \Delta_K(M_1)\otimes \Delta_K(M_3).$
  \item The following diagram of isomorphisms commutes
  \[\xymatrix{
    \Delta_K(M) \ar[rr]^(0.4){\mathrm{can}} && \Delta_K(\tilde{M})^*\otimes (K(r_M),0) \ar[d]^{\varepsilon^{dR}_{L,u}(\tilde{M}_0)^*\otimes h(\chi^{r_M})} \\
    \mathbf{1}_K  \ar[u]^{\varepsilon^{dR}_{L,-u}({M_0}) }\ar[rr]^{(-1)^{\dim_K \H^0(M)}\Omega^{-r_M}\mathrm{can}} && \mathbf{1}_K\otimes\mathbf{1}_K   ,}\]
    where $h(\chi^r):(K(\chi^r),0)\to \mathbf{1}_K$ sends $\mathbf{e}_r$ to $1.$
\end{enumerate}
\end{proposition}
\begin{proof}
Analogous to \cite[Lem.\ 3.7]{NaANT}, but with some differences. Due to the period $\Omega$ in the definition of $\Gamma(M)$ we now obtain
\begin{equation}\label{f:Gammadual}
   \Gamma(M)\Gamma(\tilde{M})=\Omega^{-r_M}(-1)^{h_M+\dim_K t_M}
\end{equation}
instead of (27) in (loc.\ cit.). By definition  the second part of the long exact sequence \eqref{f:longExactexp} for $\tilde{M}$ is given by the commutativity of the following diagram with exact rows
\begin{equation}\label{f:part2tildeM}
  \xymatrix{
    0  \ar[r]^{ } & H^1(\tilde{M})/ H^1(\tilde{M})_f\ar[d]_{\digamma_{/f}^1(\tilde{M})} \ar[r]^{ } & \mathbf{D}^{}_{\tn{cris}}({M})^*\bigoplus\mathbf{D}^{}_{\tn{dR}}(\tilde{M})^0 \ar[d]_{ } \ar[r]^{ } &\mathbf{D}^{}_{\tn{cris}}({M})^* \ar[d]_{ } \ar[r]^{ } & \H^2_{\varphi,\mathfrak Z_n}(M) \ar[d]_{\digamma^2(\tilde{M})} \ar[r]^{ } & 0   \\
    0 \ar[r]^{ } & \H^1(M)_f^* \ar[r]^{(\exp_{f,M}\bigoplus\exp_M)^* } & \mathbf{D}^{}_{\tn{cris}}({M})^*\bigoplus t_M^* \ar[r]^{ } &  \mathbf{D}^{}_{\tn{cris}}({M})^* \ar[r]^{ } & \H^0(M)^* \ar[r]^{ } & 0   }
\end{equation}
where we have identified $M\cong\tilde{\tilde{M}}$ and abbreviated $\H^i_{\varphi,\mathfrak Z_n}(N)^{\Gamma_L}$ by $\H^i(N)$. Moreover the maps $\digamma_{/f}^1(\tilde{M})$, $\digamma^2(\tilde{M}) $  and similarly $\digamma_{f}^1(\tilde{M}): \H^1(\tilde{M})_f\to  (H^1({M})/ H^1({M})_f)^*$ are induced from the complex isomorphism ${K_{\varphi,{\mf Z}}(\tilde{M})}\cong K_{\varphi,\mf Z}(M )^*[-2]$ from \eqref{f:dualitycomplex}. Taking duals gives the following commutative diagram with exact rows
\begin{equation}\label{f:dualpart2tildeM}
   \xymatrix{
     0  \ar[r]^{ } & \H^0(M) \ar[d]_{\digamma^2(\tilde{M})^*} \ar[r]^{ } & \mathbf{D}^{}_{\tn{cris}}({M})  \ar[d]_{ } \ar[r]^{ } &  \mathbf{D}^{}_{\tn{cris}}({M}) \bigoplus t_M\ar[d]_{ } \ar[r]^{ \exp_{f,M}\bigoplus\exp_M} & \H^1(M)_f \ar[d]_{\digamma_{/f}^1(\tilde{M})^*} \ar[r]^{ } & 0 \\
     0 \ar[r]^{ } & \H^2(\tilde{M})^* \ar[r]^{ } & \mathbf{D}^{}_{\tn{cris}}({M}) \ar[r]^{ } & \mathbf{D}^{}_{\tn{cris}}({M})\bigoplus(\mathbf{D}^{}_{\tn{dR}}(\tilde{M})^0)^* \ar[r]^{ } & (H^1(\tilde{M})/ H^1(\tilde{M})_f)^* \ar[r]^{ } & 0   }
\end{equation}
Upon noting  that $\digamma^2(\tilde{M})^*=\digamma^0({M}) $ while $\digamma^1(\tilde{M})^*=-\digamma^1({M}) $, whence also $\digamma^1_{/f}(\tilde{M})^*=-\digamma^1_f({M}) $, we obtain the modified commutative diagram with exact rows
\begin{equation}\label{f:dualpart2tildeMmod}
   \xymatrix{
     0  \ar[r]^{ } & \H^0(M) \ar[d]_{-\digamma^0({M})} \ar[r]^{ } & \mathbf{D}^{}_{\tn{cris}}({M})  \ar[d]_{-\id } \ar[r]^{ } &  \mathbf{D}^{}_{\tn{cris}}({M}) \bigoplus t_M\ar[d]_{-\id\bigoplus -can } \ar[r]^{ \exp_{f,M}\bigoplus\exp_M} & \H^1(M)_f \ar[d]_{\digamma_{f}^1({M})} \ar[r]^{ } & 0 \\
     0 \ar[r]^{ } & \H^2(\tilde{M})^* \ar[r]^{ } & \mathbf{D}^{}_{\tn{cris}}({M}) \ar[r]^{ } & \mathbf{D}^{}_{\tn{cris}}({M})\bigoplus(\mathbf{D}^{}_{\tn{dR}}(\tilde{M})^0)^* \ar[r]^{ } & (H^1(\tilde{M})/ H^1(\tilde{M})_f)^* \ar[r]^{ } & 0   }
\end{equation}
Combining this diagram with the analogue of diagram \eqref{f:part2tildeM} for $M$ instead of $\tilde{M}$ we obtain the commutative diagram
\begin{equation}\label{f:Thetadual}
   \xymatrix{
     \mathbf{1}_K \ar[d]_{(-1)^{\dim_Kt_M+\dim_K\H^0(M)}} \ar[r]^(0.3){\Theta_M} & \Delta_{K,1}(M)\otimes d_K(\mathbf{D}^{}_{\tn{dR}}({M})) \ar[d]^{\digamma(M)\otimes\mathrm{can}} \\
     \mathbf{1}_K  & \Delta_{K,1}(\tilde{M})^*\otimes d_K(\mathbf{D}_{\tn{dR}}(\tilde{M}))^* \ar[l]_(0.7){\Theta_{\tilde{M}}^*}.}
\end{equation}
Finally one has the commutative diagram
\begin{equation}\label{f:delta2dual}
   \xymatrix{
     d_K(\mathbf{D}^{}_{\tn{dR}}({M}))  \ar[d]_{(-1)^{h_M}\mathrm{can}} \ar[rr]^{ \Theta_{dR,-u}(M_0)} &&  \Delta_{K,2}(M) \ar[d]^{\mathrm{can}} \\
       d_K(\mathbf{D}_{\tn{dR}}(\tilde{M}))^*\otimes \mathbf{1}_K  & & \Delta_{K,2}(\tilde{M})^*\otimes(K(\chi^{r_M}),0)  \ar[ll]_{ \Theta_{dR,u}(\tilde{M}_0)^*\otimes h(\chi^{r_M})} ,}
\end{equation}
because of Remark \ref{f:epsilondual} and since changing $u$ to $-u$ requires the change $t_{LT}$ to $-t_{LT}$ (compare with \cite[Rem.\ 3.5]{NaANT} which applies analogously here) in the definition of $f_{M_0,u}$ above.
Then (ii) follows from  \eqref{f:Gammadual},\eqref{f:Thetadual} and \eqref{f:delta2dual}
while the proof of (i) is literally the same as in (loc.\ cit.).
\end{proof}

\begin{remark}\label{rem:dRchangeofu}
As in \cite[Rem.\ 3.5]{NaANT} one shows property (iii) for $ \varepsilon^{dR}_{L,u}(M_0)$ using Remark \ref{rem:equivariance}:
\[\varepsilon^{dR}_{L,au}(M_0)=\delta_{\det_{\cR_K}(M)}(a)\varepsilon^{dR}_{L,u}(M_0)\]
for all $a\in o_L^\times.$
\end{remark}

\section{Epsilon-isomorphisms for (Lubin-Tate deformations of) rank one modules}
\label{sec:RankoneCase}
 In order to construct the Epsilon-isomorphism for rank one modules $M$ in $\mathfrak{M}^{an}(K)$ we shall construct it on the level of the deformation $\mathbf{Dfm}(M)$ of $M$ (introduced in \S \ref{sec:Iwasawa}) and descend the results to $M.$ As this deformation lives over the character variety $\mathfrak{X}_{\Gamma_L}$ (base changed to $K$) of the locally $L$-analytic group $\Gamma_L$,
we can use density arguments to deduce many of its properties  just from its de Rham points.

\begin{definition}
Using that the complexes $C_n^\bullet:=K_{\Psi,D(\Gamma_L,K)}(\mathbf{Dfm}_n(M))$ are perfect by Theorem \ref{thm:CohomologyFinite} (1) we  can apply our definition $$\Delta_{1,X_n}(\mathbf{Dfm}_n(M)):=\d_{D_{r_n}(\Gamma_L,K)}(K_{\Psi,D(\Gamma_L,K)}(\mathbf{Dfm}_n(M))),$$ which defines a (graded) line-bundle on $\mathfrak{X}_{\Gamma_L}$ by (2) of the same theorem, with global sections
\[\Delta_{1,\mathfrak{X}_{\Gamma_L}}(\mathbf{Dfm}(M))=\varprojlim_n \Delta_{1,X_n}(\mathbf{Dfm}_n(M)).\]
From the proof of Theorem \ref{thm:def} we know that for the derived limit $C^\bullet$ and for every $n$,
\[\mathcal{O}_{\mathfrak{X}_{\Gamma_L}}(X_n)\otimes^\mathbb{L}_{\mathcal{O}_{\mathfrak{X}_{\Gamma_L}}(\mathfrak{X}_{\Gamma_L})}C^\bullet\cong C_n^\bullet\] in $D(\mathcal{O}_{\mathfrak{X}_{\Gamma_L}}(X_n)).$
Hence, by Definition \ref{def:quasi-Stein}, Remark \ref{rem:formalismRlimDet} and again Theorem \ref{thm:def} together with Remark \ref{rem:normalization} we obtain
\begin{align}\label{f:detdfm}
\Delta_{1,\mathfrak{X}_{\Gamma_L}}(\mathbf{Dfm}(M))\cong\d_{D(\Gamma_L,K)}(R\Gamma_{\Psi,D(\Gamma_L,K)}(\mathfrak{X}_{\Gamma_L}, \mathbf{Dfm}(M)))\cong \d_{D(\Gamma_L,K)}(\cT_{\Psi}(M)).
\end{align}
Furthermore,
\begin{align}\label{f:detdfm2}
\Delta_{2,\mathfrak{X}_{\Gamma_L}}(\mathbf{Dfm}(M))=\varprojlim_n \Delta_{2,X_n}(\mathbf{Dfm}_n(M)) \cong\varprojlim_n (\mathcal{O}_{\mathfrak{X}_{\Gamma_L}}(X_n),1)=(D(\Gamma_L,K),1).
\end{align}
\end{definition}

We survey some preliminary results that allow us to construct an isomorphism \[\Delta_{2,\mathfrak{X}_{\Gamma_L}}(\mathbf{Dfm}(M))\cong \d_{D(\Gamma_L,K)}(\cT_{\psi}(M))^{-1}.\]
Let $\delta \in \Sigma_{an}.$ Using  $(\cR_K^+(\delta))^{\Psi=0}=(\cR_K^+)^{\Psi=0}(\delta)$ combined with Lemma \ref{lem:Mellin-Amice} and since $\varphi(\mathbf{e}_\delta)$ differs from $\mathbf{e}_\delta$ only by a scalar in $K^\times,$ we can take $\eta(1,Z)\mathbf{e}_\delta$ as a $D(\Gamma_L,K)$-basis of  $(\cR_K^+(\delta))^{\Psi=0},$ which gives rise to the Mellin isomorphism
 \begin{align}\label{f:Mellingdelta}
   \mathfrak{M}_\delta: & D(\Gamma_L,K)\cong (\cR_K^+(\delta))^{\Psi=0},\; \lambda\mapsto \lambda(\eta(1,Z)\mathbf{e}_\delta).
 \end{align}

It turns out that for technical reasons (more precisely, in order to obtain the commutative diagram \eqref{f:Nakamura} below),
we have to renormalize the Mellin isomorphism by inserting the operator $\sigma_{-1}\in \Gamma_L$ with $\chi_{\mathrm{LT}}(\sigma_{-1})=-1:$
 \begin{equation}\label{f:Mellinsigmaminus1}
 \mathfrak{M}_\delta\circ\sigma_{-1}: D(\Gamma_L,K)\cong (\cR_K^+(\delta))^{\Psi=0}, \lambda \mapsto  \lambda (\sigma_{-1}(\eta(1,Z)\mathbf{e}_\delta)).
\end{equation}

\begin{remark}\label{rem:perf}\phantom{section}
\begin{enumerate}
  \item  The complexes $\cT_\Psi( LA(o_L)(\chi^{-1}\delta)) \cong\cT_\Psi( \cR_K(\delta)/\cR_K(\delta)^+)$,  $\cT_\Psi(\cR_K(\delta)) $ and\linebreak   $\cT_\Psi(\cR_K^+(\delta)),$ are all perfect complexes of $D(\Gamma_L,K)$-modules. Indeed, by  Lemma \ref{f:Chenevier}, the cohomology groups of $\cT_\Psi( LA(o_L)(\chi^{-1}\delta))$ are finite-dimensional $K$-vector spaces, whence perfect as $D(\Gamma_L,K)$-modules by \cite[Lem.\ 3.7]{Ste2} (with $r=0$ and using the Fourier-isomorphism). Then \cite[\href{https://stacks.math.columbia.edu/tag/066U}{Tag 066U}]{SP} implies that $\cT_\Psi( LA(o_L)(\chi^{-1}\delta))$ belongs to $\mathbf{D}^b_\mathrm{perf}(D(\Gamma_L,K))$. Since  $\cT_\Psi(\cR_K(\delta)) $ is in $\mathbf{D}^b_\mathrm{perf}(D(\Gamma_L,K))$ by Theorem \ref{thm:perf}, so is $\cT_\Psi(\cR_K^+(\delta))$ as the third complex in an obvious exact triangle with the previous ones. The same holds  for $\cT_\Psi(D_N(\delta))$  and  $\cT_\Psi(\cR_L^+(\delta)/D_N(\delta)) $ for similar reasons.
  \item Since over affinoids $A$ the analogous conclusion of \cite[Lem.\ 3.7]{Ste2} - i.e., that a $D(\Gamma_L,A)$-module, which is finitely generated as an $A$-module, is perfect - is not available, we are not sure whether the construction below also carries over  to families directly. It certainly does, if $\cR_A(\delta)\in\mathfrak{M}^{an}(A)$    satisfies the conditions of Remark \ref{rem:perf-delta}.
\end{enumerate}
\end{remark}

 \begin{lemma} \label{lem:isokette} Let $\delta \in \Sigma_{an}$ and let $M=\cR_K(\delta)$ be the associated $(\varphi_L,\Gamma_L)$-module of rank one. We denote by $M^+$ the submodule $\cR_K^+(\delta).$
		 We have the following isomorphisms in $\cP_{D(\Gamma_L,K)}:$
		\begin{enumerate}

				\item $\d_{D(\Gamma_L,K)}(\cT_\Psi(M)) \cong \d_{D(\Gamma_L,K)}(\cT_\Psi(M^+))$ induced by the canonical inclusion $M^+ \subset M$ and the trivialisation of $\d_{D(\Gamma_L,K)}(\cT_\Psi(M/M^+))$ from Lemma \ref{f:Chenevier}.
			\item $\d_{D(\Gamma_L,K)}(\cT_\Psi(M^+)) \xrightarrow{\cong} \d_{D(\Gamma_L,K)}([M^+ \xrightarrow{\Psi}M^+])$ induced by $(1-\varphi_L,\id)$ and the trivialization of   $d_{D(\Gamma_L,K)}(\cT_\Psi(D_N(\delta)))$.
			\item $(D(\Gamma_L,K),1) \cong (\d_{D(\Gamma_L,K)}[M^+ \xrightarrow{\Psi}M^+])^{-1}$ induced by identifying $\ker(\Psi)$ with $D(\Gamma_L,K) $ via $\mathfrak{M}_\delta\circ\sigma_{-1} .$
					\end{enumerate}
	Chaining these together gives an isomorphism $\d_{D(\Gamma_L,K)}(\cT_\Psi(M))^{-1} \cong (D(\Gamma_L,K)(\delta),1).$
	\end{lemma}
\begin{proof}
	The first statement follows since the short exact sequence $0 \to M^+ \to M \to M/M^+ \to 0$ induces a short exact sequence of complexes. For the second statement we use that by \cite[Lem.\ 5.1]{FX} we have a commutative diagram with exact rows
	\[\xymatrix{
		0  \ar[r]^{ } &  ( \cR_K^+(\delta))^{\Psi=1} \ar[d]_{1- \varphi}^{\cong} \ar[r]^{ } &  \cR_K^+(\delta) \ar[d]_{1- \varphi }^{\cong} \ar[r]^{\Psi-1 } &  \cR_K^+(\delta) \ar@{=}[d]_{ } \ar[r]^{ } & 0 \\
		0 \ar[r]^{ } &  \cR_K^+(\delta)^{\Psi=0} \ar[r]^{ } &  \cR_K^+(\delta) \ar[r]^{\Psi} &  \cR_K^+(\delta) \ar[r]^{ } & 0,   }\]
	which induces a quasi-isomorphism between the complexes, if $\delta(\pi_L)\neq \pi_L^{-i}$ for all $i\in\mathbb{N}.$ Otherwise, kernel and cokernel of $\cR_K^+(\delta)\xrightarrow{1-\varphi} \cR_K^+(\delta)   $ are isomorphic to $Kt_{LT}^i$ and can be trivialized by each other when taking determinants (formally this is achieved by replacing $\cR_K^+(\delta)$ by $\cR_K^+(\delta)/D_N(\delta)$ and then trivializing the determinant of $\cT_\Psi(D_N(\delta))$ as in \cite[(40),(44) in \S 4.1]{NaANT}). For the third statement we first remark that the complex $M^+ \xrightarrow{\Psi} M^+$ (concentrated in degrees $1,2$) is cohomologically perfect by Lemma \ref{lem:Mellin-Amice} -  using  $(\cR_K^+(\delta))^{\Psi=0}=(\cR_K^+)^{\Psi=0})(\delta)$  -  because on the one hand $\Psi$ is surjective and on the other hand its kernel is free over $D(\Gamma_L,K)$ by \eqref{f:Mellinsigmaminus1}. Therefore the determinant of  $M^+ \xrightarrow{\Psi} M^+$ is equal to $(D(\Gamma_L,K),1)^{-1}.$
\end{proof}


From Lemma \ref{lem:isokette} we obtain  finally  an isomorphism (cf.\ \cite[Def.\ 4.1]{NaANT})
\[\Theta(\delta):\mathbf{d}_{D(\Gamma_L,K)}(\cT_\Psi( \cR_K(\delta)))^{-1}\cong
\mathbf{d}_{D(\Gamma_L,K)}(D(\Gamma_L,K))\cong\Delta_{2,\mathfrak{X}_{\Gamma_L}}(\mathbf{Dfm}(\cR_K(\delta))) \] which in turn induces an isomorphism over $K$
\begin{align}\label{f:epsilonDfm}
\notag  \varepsilon_{D(\Gamma_L,K),u}(\mathbf{Dfm}(\cR_K(\delta)):\mathbf{1}_{D(\Gamma_L,K)} & \xrightarrow{\mathrm{can}} \mathbf{d}_{D(\Gamma_L,K)}(\cT_\Psi( \cR_K(\delta)))\mathbf{d}_{D(\Gamma_L,K)}(\cT_\Psi( \cR_K(\delta)))^{-1}\\
    & \xrightarrow{\id \otimes \Theta(\delta)} \mathbf{d}_{D(\Gamma_L,K)}(\cT_\Psi( \cR_K(\delta)))\Delta_{2,\mathfrak{X}_{\Gamma_L}}(\mathbf{Dfm}(\cR_K(\delta))) \\
\notag    & \xrightarrow{=} \Delta_{\mathfrak{X}_{\Gamma_L}}(\mathbf{Dfm}(\cR_K(\delta))).
\end{align}
 Note that the map \eqref{f:Mellinsigmaminus1} depends implicitly on $u.$  If we consider instead of $\cR_L$ the isomorphic subring $\cR_L(Z_u)$ of $\widetilde{\cR}_L,$ as pointed out in Remark \ref{rem:hiddendependence} , then for $a=\chi_{LT}(\gamma_a) \in o_L^\times$ we have $Z_{au} = [a](Z_u)$ and thus we get a commutative diagram
\begin{equation}
	\label{eq:dependenceonu}
	\xymatrix{
	D(\Gamma_L,K) \ar[d]_{\delta_{\gamma_a^{-1}}\cdot} \ar[r]^{\mathfrak{M}_{\delta,u}} & (\cR_K^{+}(\delta))^{\Psi=0} \ar[d]^{\delta(a)^{-1}\cdot} \\
	D(\Gamma_L,K)  \ar[r]^{\mathfrak{M}_{\delta,au}} &(\cR_K^{+}(\delta))^{\Psi=0}.}
\end{equation}
Indeed, we have \begin{align*}
                  \mathfrak{M}_{\delta,au}(\delta_{\gamma_a^{-1}}\lambda)  & =\lambda\bigg(\delta_{\gamma_a^{-1}}\big(\eta(1,Z_{au})\big)\delta_{\gamma_a^{-1}}\mathbf{e}_\delta\bigg)  \\
                    &= \lambda\bigg( \eta(a^{-1},[a](Z_{u}))\delta(a)^{-1}\mathbf{e}_\delta\bigg)\\
                    & = \delta(a)^{-1} \cdot \bigg(\lambda\big( \eta(1,Z_u)\mathbf{e}_\delta\big)\bigg).
                \end{align*}


 Concerning the descent, we have to distinguish the following two ways.

\begin{remark}\label{rem:descent} Let $\delta: \Gamma_L \to K^\times$ be an $L$-analytic character.
Mapping a Dirac distribution $\gamma$ to $\delta(\gamma)\mathbf{e}_{\delta}$ induces a surjection of $D(\Gamma_L,K)$-modules $$p_{\delta}\colon D(\Gamma_L,K) \to K\mathbf{e}_{\delta}.$$ Alternatively we may equip $D(\Gamma_L,K)$ with the $\Gamma_L$-action $\gamma \eta = [\gamma^{-1}]\eta$,
denoting the resulting $\Gamma_L$-module by
$D(\Gamma_L,K)^\iota$, and map $\gamma$ to $\delta(\gamma^{-1})\mathbf{e}_\delta$ to obtain a surjection of $D(\Gamma_L,K)$-modules $$f_{\delta}:D(\Gamma_L,K)^{\iota} \to K\mathbf{e}_{\delta}.$$
\end{remark}
\begin{proof}
	Since $\delta$ is analytic $K(\delta)=K\mathbf{e}_{\delta}$ comes equipped with a $D(\Gamma_L,K)$-module structure extending the $K[\Gamma_L]$-module structure. The map $p_\delta$ is surjective because $1$ is mapped to a $K$-basis $\mathbf{e}_\delta$ and $D(\Gamma_L,K)$-linear by construction. The second statement follows analogously since the inverted action is also $L$-analytic.
\end{proof}

Now, for the descent we   observe that, if $f_{\delta_0}:D(\Gamma_L,K)\to K$ arises from a character $\delta_0:o_L^\times\subseteq L^\times\to K^\times$ interpreted as character of $\Gamma_L$,
we have the following:
 \begin{lemma}\label{lem:decentDelta}
The isomorphism \eqref{f:detdfm} induces the canonical isomorphism
\[\mathrm{sp}_{\delta_0}:\Delta_{{\mathfrak{X}_{\Gamma_L}}}(\mathbf{Dfm}(\cR_K(\delta)))\otimes_{D(\Gamma_L,K),f_{\delta_0}}K\cong \Delta_{K}(\cR_K(\delta\delta_0)) \]
taking the normalisation from Remark \ref{rem:normalization} into account, compare with \cite[(34), p. 370]{NaANT}.
\end{lemma}
\begin{proof} We show this isomorphism for each part of $\Delta$ separately:
\begin{align}\label{f:derivedtensor}
\notag\Delta_{1,\mathfrak{X}_{\Gamma_L}}(\mathbf{Dfm}(M))\otimes_{D(\Gamma_L,K),f_{\delta_0}}K&\cong d_{K[\Gamma_L/U]}(\mathcal{T}_\psi(M(\delta_0))\otimes^\mathbb{L}_{D(\Gamma_L,K)}D(\Gamma_L/U))\otimes_{K[\Gamma_L/U]}K\\
&\cong d_{K[\Gamma_L/U]}(\mathcal{T}_\psi(M(\delta_0))\otimes^\mathbb{L}_{D(U)}K)\otimes_{K[\Gamma_L/U]}K\\
\notag&\cong d_{K[\Gamma_L/U]}( K_{\Psi_L,{\mf Z}} (M(\delta_0)))\otimes_{K[\Gamma_L/U]}K\\
\notag&\cong d_{K[\Gamma_L/U]}( K_{\varphi_L,{\mf Z}} (M(\delta_0)))\otimes_{K[\Gamma_L/U]}K=\Delta_{1,K}(M(\delta_0))
\end{align}
and
\begin{align}\label{f:Delta2}
 \Delta_{2,\mathfrak{X}_{\Gamma_L}}(\mathbf{Dfm}(M))\otimes_{D(\Gamma_L,K),f_{\delta_0}}K&\cong \Delta_{2,K}(M)\otimes_K D(\Gamma_L,K)\otimes_{D(\Gamma_L,K),f_{\delta_0}}K\\
\notag &\cong \Delta_{2,K}(M(\delta_0))=(K\mathbf{e}_{\delta\delta_0},1)
\end{align}
 using Remark \ref{rem:DeltaVonRdelta}.
\end{proof}

With these preparations we are now able to state the main result of this article.

\begin{theorem}[Local $\varepsilon$-conjecture for Lubin-Tate deformations of rank one modules]\label{thm:LTDeform}
Let $F'/L$ be a finite subextension of $K$ and $M$ be a rank one analytic $(\varphi_L,\Gamma_L)$-module over $\cR_{F'}$ and denote by $M_K$ the completed base change $M\hat{\otimes}_{F'} K$. Then the isomorphism
  \[\varepsilon_{D(\Gamma_L,K),u}(\mathbf{Dfm}(M_K)):\mathbf{1}_{D(\Gamma_L,K)}\xrightarrow{\cong} \Delta_{\mathfrak{X}_{\Gamma_L}}(\mathbf{Dfm}(M_K))\]
  induces for every  $L$-analytic character $\vartheta:\Gamma_L\to F^\times$ with finite intermediate extension $F'\subseteq F\subseteq K$ such that $M_K(\vartheta)$ is de Rham the following commutative diagram
\begin{equation}\label{f:LTDeformatoiondecent}
  \xymatrix{
    \mathbf{1}_{D(\Gamma_L,K)}\otimes_{D(\Gamma_L,K),f_\vartheta} K \ar[d]_{\varepsilon_{D(\Gamma_L,K),u}(\mathbf{Dfm}(M_K))\otimes\id_{K}} \ar[r]^-{\mathrm{can} } & \mathbf{1}_K \ar[d]^{\varepsilon^{dR}_{F,u}(\cR_F(\delta\vartheta))} \\
    \Delta_{\mathfrak{X}_{\Gamma_L}}(\mathbf{Dfm}(M_K))\otimes_{D(\Gamma_L,K),f_\vartheta} K  \ar[r]^-{\mathrm{sp}_\vartheta} & \Delta_{K}(M_K(\vartheta)), }
\end{equation}
where the notation $f_\vartheta$ has been defined in Remark \ref{rem:descent} and the specialisation isomorphism $\mathrm{sp}_\vartheta$ is explained in Lemma \ref{lem:decentDelta} above.
  Moreover, $\varepsilon_{D(\Gamma_L,K),u}(\mathbf{Dfm}(M_K)) $ is uniquely determined by this property.
\end{theorem}
The uniqueness follows from the considerations in Appendix \ref{App:density} while the specialisation property will be proved in subsection \ref{sec:descent} below.

Note that the isomorphism $\varepsilon_{D(\Gamma_L,K)}(\mathbf{Dfm}(\cR_K(\delta)) $ does  not literally fit into Conjecture \ref{conj}, because $D(\Gamma_L,K)$ is not an affinoid algebra over $K$. But for any    morphism of rigid analytic spaces  $f:Sp(A)\to \mathfrak{X}_\Gamma$ with
an affinoid algebra $A$ (e.g.\ $D_n$) over $K$ it induces the isomorphism
\begin{align*}
\varepsilon_{A}(f^*\mathbf{Dfm}(\cR_K(\delta))(Sp(A)))&:=\varepsilon_{D(\Gamma_L,K)}(\mathbf{Dfm}(\cR_K(\delta))\otimes_{D(\Gamma_L,K)}A)
\colon
\\
\mathbf{1}_{A}&\xrightarrow{\cong} \Delta_{A}(f^*\mathbf{Dfm}(\cR_K(\delta))(Sp(A)))
\end{align*}
which provides instances of the conjectured type. Note that for   the inclusion $f:Sp(D_n)\hookrightarrow\mathfrak{X}_\Gamma$ we obtain \[f^*\mathbf{Dfm}(\cR_K(\delta))(Sp(D_n))\cong \mathbf{Dfm}_n(\cR_K(\delta)). \]


\subsection{Property (v)}
Specialization  to the case considered by Nakamura requires some special care, because we used a different definition of $\varepsilon$-constants. As discussed in \ref{rem:compNakamura} the assumption that $K$ contains $L^{ab}$ can be dropped since $L_\infty$ contains the $p$-power roots of unity. We can thus even assume $K=\QQ_p$ in the construction of the de Rham $\varepsilon$-constants. Similarly we can take $\Omega=1$ and hence do not need any special assumptions on $K$ in order to make use of $p$-adic Fourier theory.
\begin{proposition} \label{prop:precisecomp}
	Assume $L= \QQ_p,$ assume $\pi_L=p.$ take $u_n = \zeta_{p^n}-1$ for a compatible system of $p$-power roots of unity and choose a $\gamma \in \Gamma_{\QQ_p},$ which is a topological generator of the torsion-free part, such that $\log_0(\chi_{cyc}(\gamma)) = 1,$ and take $\Omega_{\QQ_p} = 1.$ Then, if one assumes $K=\QQ_p$, our construction agrees with the one in \cite{NaANT}.
\end{proposition}
\begin{proof}
	Note that by a density argument and by property (vi) it suffices to see that the constructions in the de Rham case coincide.
	The condition of $L$-analyticity is automatic, if $L=\QQ_p.$
We remark that the complex $K^\bullet_{\varphi,\mathfrak{Z}}(-)$ considered by us specialises to a variant of the usual Herr-complex as we can take $\mathfrak{Z}=\gamma-1$, but there is a small difference to \cite[Definition 2.10]{NaANT}. The order of $\varphi-1$ and $\gamma-1$ is exchanged (which poses no problem), Nakamura uses a topological generator $\gamma_{Na}$ of $\Gamma/\Delta,$ with $\Delta = \Gamma_{p-\text{power-torsion}}$ while we use a generator of the free part. In the case $p=2$ the terms of Nakamura's complex are $M^{\Delta}.$ In this case our choice of $\gamma$ is a valid choice for the variant in (loc.cit.) while in the case $p\neq 2$ we can arrange that $\gamma_{Na}^{p-1} = \gamma.$ In both cases the torsion subgroup $\Delta' \subseteq \Gamma$ is a split subgroup and taking $\Delta'$-invariants is exact in characteristic $0.$ Let $U = \overline{\langle\gamma\rangle} \subseteq \Gamma.$ For $p=2$ we have $\Delta = \Delta'$ and
plugging in the isomorphism $M^\Delta \cong M/\Delta$ and $\Delta \cong \Gamma/U$ we see that our $\QQ_p\otimes_{\QQ_p[\Delta]}K_{\varphi,\mathfrak{Z}}(M)$ is canonically isomorphic to the complex considered in (loc.cit.).
For $p\neq 2$ we can consider instead the natural map of complexes
\[ [M \xrightarrow{\gamma_{Na}-1} M] \to [M \xrightarrow{\gamma-1} M] \] given by $m \mapsto \frac{1}{p-1} \sum_{g \in \Delta'}gm$ in both degrees, which induces a quasi-isomorphism onto the $\Delta'$-(co)invariants of the right-hand side and induces a corresponding quasi-isomorphism of the Herr-complexes by taking $\varphi-1$-cones.
We can thus conclude that the fundamental lines are canonically isomorphic to the ones considered by Nakamura. Similarly the exponential maps are the same. Because $\pi_L=p=q$ we see that the character $\chi$ is just $\chi_{cyc}$ and the duality pairing $\langle-,- \rangle_{\tilde{M}}$ from section \ref{sec:duality} is the 	pairing used by Nakamura. In \eqref{f:longExactexp} we use $\langle-,-\rangle_M$ which by the same reasoning corresponds to the pairing used by Nakamura, namely the duality pairing for $M_0 = \tilde{M}$ . The assumptions on $\gamma$ and $\Omega$ avoid the problem discussed in \ref{rem:compNakamura} (ii) concerning normalisation factors and the appearance of $\Omega$ in the $\Gamma$-factor. Finally, the series $\eta(1,Z)$ is just $1 + Z$ and we can view $\eta(1,(u_n \otimes 1))$ appearing in the construction of the equivariant $\varepsilon$-constants as an element of $L_\infty,$ in fact we have $\eta(1,u_n\otimes 1) = \zeta_{p^n}$ under the isomorphism $L_\infty \otimes_{\QQ_p} \QQ_p \cong L_\infty.$ Combining all of the above shows that our  $\varepsilon$-constants constructed in the de Rham case agree with those in Nakamura's work.
\end{proof}

\subsection{Property (i)}
For all $f:A\to A'$, such that we are able to construct the $\varepsilon$-isomorphism as above for $A$ and $A'$, the base change property (i) with respect to $f:A\to A'$ obviously holds by construction.

\subsection{Property (iii)}
We can rephrase the diagram \eqref{eq:dependenceonu} to the following commutative diagram for any $a\in o_L^\times$
\begin{equation*}
	\xymatrix{
		D(\Gamma_L,K)(\delta) \ar[d]_{[a^{-1}]} \ar[r]^{\mathfrak{M}_{\delta,au}} & (\cR_K^{+}(\delta))^{\Psi=0} \ar[d]^{\id} \\
		D(\Gamma_L,K)(\delta)  \ar[r]^{\mathfrak{M}_{\delta,u}} &(\cR_K^{+}(\delta))^{\Psi=0}},
\end{equation*}
where $[a]$ acts on $D(\Gamma_L,K)(\delta)$ as $\delta_{\gamma_a}^{-1}\cdot\delta(a)$ (here $\delta_{\gamma_a}$ denotes the dirac distribution attached to $\gamma_a\in\Gamma_L$ with $\chi_{LT}(\gamma_a)=a$). Note that the action on $N:=\mathbf{Dfm}(\cR_K(\delta))$ with respect to the basis $1 \otimes e_\delta$ is given precisely by the character $\overline{\delta}\colon \Gamma_L \to D(\Gamma_L,K)^\times; \gamma \mapsto (\delta_{\gamma})^{-1}\delta(\chi_{LT}(\gamma))$ and hence property (iii) follows from the above diagram by specialising along $D(\Gamma_L,K)(\delta) \to K(\delta).$

\subsection{Descent}\label{sec:descent}

For $\delta\in \Sigma_{an}(F)$ with  $F$ a finite extension of $L$, we consider the decomposition $\delta=\delta^\mathrm{un}\delta_0$ as in section \ref{sec:rankonecharacters} and  define on the basis of \eqref{f:epsilonDfm}
\[\varepsilon_{K,u}(\cR(\delta)):\mathbf{1}_{K}\xrightarrow{\cong} \Delta_K(\cR(\delta)) \]
as $\varepsilon_{D(\Gamma_L,K)}(\mathbf{Dfm}(\cR_K(\delta^\mathrm{un}))\otimes_{D(\Gamma_L,K),f_{\delta_0}}K)
$  followed by the isomorphism from Lemma \ref{lem:decentDelta}.
In order to make this definition more explicit we have to understand the isomorphism $\bar{\Theta}(\delta):=\Theta(\delta^\mathrm{un})\otimes_{D(\Gamma_L,K),f_{\delta_0}}L$, which we will consider as an isomorphism
\[\bar{\Theta}(\delta):\bigotimes_{i=0}^2 \d_{K[\Gamma_L/U]}(H^i_{\Psi_L,{\mf Z}}(\cR_K(\delta)))^{(-1)^{i+1}}\otimes_{K[\Gamma_L/U]}K\cong
(K\mathbf{e}_\delta,1) \]
by using \eqref{f:Delta2} and the inverse of the natural isomorphism
\begin{align*}
\d_{D(\Gamma_L,K)}(\cT_\psi( \cR_L(\delta)))\otimes_{D(\Gamma_L,K),f_{\delta_0}}K&\cong \d_{K[\Gamma_L/U]}( K_{\Psi_L,{\mf Z}} (M(\delta_0)))\otimes_{K[\Gamma_L/U]}K\\
&\cong \bigotimes_{i=0}^2 \d_{K[\Gamma_L/U]}(H^i_{\Psi_L,{\mf Z}}(\cR_K(\delta)))^{(-1)^{i}}\otimes_{K[\Gamma_L/U]}K
\end{align*}
induced from \eqref{f:derivedtensor} using properties of the determinant functor from section \ref{sec:Det}.

From the exact sequences \eqref{f:Colmeztransform}, \eqref{f:psiinv}, \eqref{f:psicoinv} we derive the following exact sequences and isomorphisms:
\begin{align}\label{f:longExact}
  0\to &H^0_{\Psi,{\mf Z}}(\cR_K^+(\delta))\to H^0_{\Psi,{\mf Z}}(\cR_K (\delta))\to  H^0_{\Psi,{\mf Z}}(LA(o_L)(\chi^{-1}\delta))\to \\
\notag    &  H^1_{{\mf Z}}(H^0_\Psi(\cR_K^+(\delta)))\to H^1_{{\mf Z}}(H^0_\Psi(\cR_K(\delta)))\to H^1_{{\mf Z}}(H^0_\Psi(LA(o_L)(\chi^{-1}\delta)))\to 0,
\end{align}
\begin{align}\label{f:cohR+}
 & H^2_{\Psi,{\mf Z}}(\cR_K^+(\delta))=H^0_{{\mf Z}}(H^1_\Psi(\cR_K^+(\delta)))=0   \\ & \label{f:H1R+} H^1_{\Psi,{\mf Z}}(\cR_K^+(\delta))\cong  H^1_{{\mf Z}}(H^0_\Psi(\cR_K^+(\delta)))\cong \cR_K^+(\delta)^{\Psi=1}/Z,
\end{align}
\begin{align}\label{f:H0H1}
  H^0_{{\mf Z}}(H^1_\Psi(\cR_K(\delta)))&\cong H^0_{{\mf Z}}(H^1_\Psi(LA(o_L)(\chi^{-1}\delta))) \\ H^2_{\Psi,{\mf Z}}(\cR_K (\delta))&\cong H^2_{\Psi,{\mf Z}}(LA(o_L)(\chi^{-1}\delta)),
\end{align}
\begin{equation}\label{f:spectralsequ}
\xymatrix@C=0.5cm{
  0 \ar[r] & H^1_{{\mf Z}}(H^0_\Psi(\cR_K(\delta))) \ar[rr]^{ } && H^1_{\Psi,{\mf Z}}(\cR_K (\delta)) \ar[rr]^{ } && H^0_{{\mf Z}}(H^1_\Psi(\cR_K(\delta))) \ar[r] & 0. }
\end{equation}

For the descent it is useful to recall that the determinant functor $\mathbf{d}_?$ commutes with taking the derived tensor product  $-\otimes^\mathbb{L}_{D(\Gamma_L,K),f_{\delta_0}}K$. E.g.\  the additivity on short exact sequences above turns into the additivity on the associated long exact sequences of cohomology groups below. Finally, the determinant functor commutes with attached spectral sequences by \cite{ven-det}.
%


\subsection{Verification of the conditions (iv),   (vi).}

 In this subsection, we prove
the condition (iv)   using density arguments in the process of verifying the
condition (vi).
Indeed, it suffices to prove (vi) as the duality statement for de Rham characters was shown in \ref{prop:duality} and by Zariski density of the de Rham characters (see Corollary \ref{cor:density}) the validity of property (iv) holds in general once we establish (vi), i.e., the interpolation property in the de Rham case. We follow the strategy of Nakamura and consider first a generic $L$-analytic de Rham character of weight $k.$ The case $k\leq0$ boils down to Proposition \ref{prop:explicit}. The case $k\geq 1$ is treated in Proposition \ref{prop:explicit2}. The remaining so-called exceptional case is treated in Section \ref{sec:exceptional}.

\subsubsection{Twisting}\label{sec:twisting}
 We define the operator $\partial\colon\cR_K\to\cR_K, f\mapsto   \frac{1}{\log'_{LT}}\frac{df}{dZ}=\frac{df}{ dt_{LT}},$ and the residuum map\linebreak
  $Res\colon\cR_K\to K, f\mapsto res(fdt_{LT})$ with $res(\sum_{i\in\mathbb{Z}} a_iZ^idZ)=a_{-1}.$ Extending theses maps coefficientwise, i.e., applying it to $f$ in $f\mathbf{e}_\delta$ and using  \cite[Lem.\ 2.11, 2.12]{FX} we obtain  an exact sequence\footnote{This sequence already exists over $L$ instead of $K!$}
\begin{equation}\label{f:partial}
   \xymatrix@C=0.5cm{
    0 \ar[r] & K(\delta) \ar[rr]^{ } && \cR_K(\delta) \ar[rr]^{{\partial}} && \cR_K(x\delta) \ar[rr]^{Res} && K(\delta|x|^{-1}) \ar[r] & 0. }
\end{equation}

It is well-known that the partial operator $\partial\colon \cR_K\to\cR_K$ is related to twisting, see e.g. \cite[\S 4.3.9]{SV20}:\footnote{  Here  $\Omega$ is required!}

\begin{equation}\label{f:twist}
  \xymatrix{
   D(\Gamma_L,K) \ar[d]_{Tw_{\chi_{LT}}} \ar[r]^-{\mathfrak{M}_{}} & (\cR_K^+)^{\psi_L=0} \ar[d]^{\frac{1}{\Omega}\partial}_{\cong} \\
   D(\Gamma_L,K) \ar[r]^-{\mathfrak{M}_{}} &(\cR_K^+)^{\psi_L=0} .  }
\end{equation}
Here, for a locally $L$-analytic character $\rho : \Gamma_L \rightarrow K^\times$ we  denote   by
\begin{equation*}
  Tw_\rho : D(G,K) \xrightarrow{\;\cong\;} D(G,K) \ ,
\end{equation*}
the isomorphism which on  Dirac distributions satisfies  $Tw_\rho(\delta_g) = \rho(g) \delta_g$.

Using for $\mathbf{d}_{D(\Gamma_L,K)}\cT_\Psi(K(\delta'))$, $\delta'=\delta,\delta|x|^{-1},$ the trivialization by identity, the operator $\partial$ induces via the above exact sequence the isomorphism
\[\partial:\Delta_{1,\mathfrak{X}_{\Gamma_L}}(\mathbf{Dfm}(\cR_K(\delta)))\xrightarrow{\cong}\Delta_{1,\mathfrak{X}_{\Gamma_L}}(\mathbf{Dfm}(\cR_K(x\delta))),  \]
which also descends   to an isomorphism
\[\partial:\Delta_{1,K}( \cR_K(\delta))\xrightarrow{\cong}\Delta_{1,K}(\cR_K(x\delta)). \]
Moreover, we have isomorphisms
\[\partial:\Delta_{2,\mathfrak{X}_{\Gamma_L}}(\mathbf{Dfm}(\cR_K(\delta)))\xrightarrow{\cong}\Delta_{2,\mathfrak{X}_{\Gamma_L}}(\mathbf{Dfm}(\cR_K(x\delta))),  \]
and
\[\partial:\Delta_{2,K}( \cR_K(\delta))\xrightarrow{\cong}\Delta_{2,K}(\cR_K(x\delta)). \]
by sending $f\mathbf{e}_{\delta}$ to $\frac{-1}{\Omega}f\mathbf{e}_{x\delta}$.
Altogether we obtain an isomorphism
\[\partial:\Delta_{\mathfrak{X}_{\Gamma_L}}(\mathbf{Dfm}(\cR_K(\delta)))\xrightarrow{\cong}\Delta_{\mathfrak{X}_{\Gamma_L}}(\mathbf{Dfm}(\cR_K(x\delta))),  \]
which also descends   to an isomorphism
\[\partial:\Delta_K( \cR_K(\delta))\xrightarrow{\cong}\Delta_K(\cR_K(x\delta)). \]
Using diagram \eqref{f:twist} and the definition of $\epsilon_{D(\Gamma_L,K)}(\mathbf{Dfm}(\cR_K(\delta))$ and $ \epsilon_{L}(\cR_K(\delta))$ respectively, we conclude the following
\begin{proposition}\label{prop:twisteps} If $\delta\neq  \mathbf{1}, |x|$, then there are canonical equalities
 \[\partial\circ\epsilon_{D(\Gamma_L,K)}(\mathbf{Dfm}(\cR_K(\delta))= \epsilon_{D(\Gamma_L,K)}(\mathbf{Dfm}(\cR_K(x\delta)) \mbox{ and }\partial \circ \epsilon_{L}(\cR_K(\delta))= \epsilon_{L}(\cR_K(x\delta)).\]
\end{proposition}
\begin{proof} Since the second statement follows by descent from the first one, we only have to consider the case of the deformation following the construction in Lemma \ref{lem:isokette} step by step. Regarding  \ref{lem:isokette}(i) we observe that  the operator $\partial$ restricts to an operator $\cR_K^+(\delta)\to \cR_K^+(x\delta) $ while it induces the operator $LA(o_L)(\chi^{-1}\delta) \to LA(o_L)(\chi^{-1}x\delta), \phi\mathbf{e}_{\chi^{-1}\delta}\mapsto \Omega x \phi\mathbf{e}_{\chi^{-1}x\delta},$ which can easily be derived from Remark \ref{rem:FColmezFormeln} (v) combined with the exactness of \eqref{f:partial}. The compatibility with $1-\varphi_L$ in  \ref{lem:isokette}(ii) is a consequence of  Remark \ref{rem:FColmezFormeln} (i).  Finally,  the compatibility of $\partial$ with $\mathfrak{M}_\delta\circ\sigma_{-1} $  in  \ref{lem:isokette}(iii) follows from diagram \eqref{f:twist} together with  the $\sigma_{-1}$ in the definition of \eqref{f:Mellinsigmaminus1} using \ref{rem:FColmezFormeln} (ii). Combining both yields the factor $-\Omega$ which cancels against the factor in the definition of $\partial_{|\Delta_2}.$ One can check that the twisting construction is compatible with the various trivializations involved.
%
%
\end{proof}


\begin{proposition}\label{prop:twistepsdR} Let $\delta\in\Sigma^{an}(F)$ with $F/L$ finite
such that $\cR_F(\delta)$ is a de Rham $(\varphi,\Gamma_L)$-module with Hodge-Tate weight different form zero. Then we have the equality
\[\partial\circ \varepsilon^{dR}_{F,u}(\cR_F(\delta))= \varepsilon^{dR}_{F,u}(\cR_F(x\delta)).\]
\end{proposition}
\begin{proof}
The proof is analogous to that of \cite[4.14]{NaANT} upon noting that $\Gamma^*(k)$ has to be replaced by $\Omega^k\Gamma^*(k)$.
\end{proof}

Since $\varepsilon^{dR}$ and $\varepsilon$ are compatible with respect to $\partial$ by the above propositions,  it can be used to transport the validity of the Conjecture between  characters $\delta x $ and $\delta.$

\subsubsection{Generic case}

This subsection has been inspired by \cite[4B1]{NaANT} and \cite{V-Kato}. In this subsection $U=\Gamma_n$ and ${\mf Z}={\mf Z_n}$ for an appropriate sufficiently large $n\gg0,$ which might be adapted to the specific situation. This is possible because due to our normalisations the constructions and the  factorization of the descent over $K[\Gamma_L/\Gamma_n]$ are independent of $n$, see  Lemma \ref{lem:decentDelta}, \eqref{f:derivedtensor},  Definition \ref{def:dualexp} and Remark \ref{rem:normalization}.

\begin{lemma}\label{lem:genericvan}
For $\delta\in \Sigma_{gen}{(F)}$ we have
\begin{align}\label{f:genricvanishing}
  H^i_{\Psi,{\mf Z}}(LA(o_L)(\chi^{-1}\delta))&=H^i_{\Psi,{\mf Z}}(Pol_{\leq N}(o_L)){ = H^i_{\Psi,{\mf Z}}(D_N(\delta))=0,} \\\notag
  H^1_{{\mf Z}}(H^0_{\Psi}(LA(o_L)(\chi^{-1}\delta)))
  &=H^0_{{\mf Z}}(H^1_{\Psi}(LA(o_L)(\chi^{-1}\delta)))=0
\end{align}
for all $i$ and $N\geq 0,$ and
\begin{equation}\label{f:genericvanishingRR+}
  H^i_{\Psi,{\mf Z}}(\cR_K^+(\delta))=H^i_{\Psi,{\mf Z}}(\cR_K(\delta))=0
\end{equation}
for $i\neq 1,$ and
\begin{equation}\label{f:isoR+R}
  H^1_{\Psi,{\mf Z}}(\cR_K^+(\delta))\cong H^1_{{\mf Z}}(H^0_\Psi(\cR_K^+(\delta)))\cong  H^1_{\Psi,{\mf Z}}(\cR_K(\delta)).
\end{equation}

\end{lemma}
\begin{proof}
The first claim  follows from  \eqref{f:PolZPsi},\eqref{f:DZPsi} and Lemma \ref{f:Chenevier}.  The second follows from Lemma \ref{lem:compcoh} (combined with Remark \ref{rem:compcoh}) and \ref{lem:dim}  (there for $\Gamma_L$ instead of $U$, but applying the result to all twists by characters of the finite group $\Gamma_L/U$ also implies the statement concerning $U$) combined with \eqref{f:longExact} and \eqref{f:cohR+}. The last assertion follows from  the previous ones combined with \eqref{f:spectralsequ},\eqref{f:H0H1},\eqref{f:H1R+}.
\end{proof}

By construction according to Lemma \ref{lem:isokette} and using  Lemma \ref{lem:genericvan} we see that $\bar{\Theta}(\delta)$ arises - upon taking determinants and descending further by $-\otimes^\mathbb{L}_{K[\Gamma_L/U]}K$ - by the composite of

\begin{enumerate}
\item (the inverse of) the isomorphism  $H^1_{{\mf Z}}(H^0_\Psi(\cR_K^+(\delta)))\cong  H^1_{\Psi,{\mf Z}}(\cR_K(\delta))$    together with the trivializations of $\mathbf{d}_{D(\Gamma_L,K)}(\cT_\Psi( LA(o_L)(\chi^{-1}\delta)))$ and $ \mathbf{d}_{D(\Gamma_L,K)}(\cT_\Psi( Pol_{\leq N}(o_L)(\chi^{-1}\delta)))$,
\item $H^1_{{\mf Z}}(H^0_\Psi(\cR_K^+(\delta)))\cong H^1_{{\mf Z}}(\cR_K^+(\delta)^{\Psi=0})$ induced by $1-\varphi,$   together with the trivialisation of kernel and cokernel of $\cR_K^+(\delta)\xrightarrow{1-\varphi} \cR_K^+(\delta)   $ - each isomorphic to  $Kt_{LT}^i$  - respectively with the trivialization of   $\mathbf{d}_{D(\Gamma_L,K)}(\cT_\Psi(D_N(\delta)))$ and
\item $H^1_{{\mf Z}}(\cR_K^+(\delta)^{\Psi=0})\xrightarrow{C_{Tr}({\mf Z_n})}(\cR_K^+(\delta)^{\Psi=0})_U \cong D(\Gamma_L,K)_U\cong K[\Gamma_L/U]$ up to choosing basis elements and using the Mellin transform $\mathfrak{M}_\delta\circ \sigma_{-1}.$
\end{enumerate}


Altogether - up to the isomorphism    $H^1_{{\mf Z}}(H^0_\Psi(\cR_K^+(\delta))) \cong H^1_{\Psi,{\mf Z}}(\cR_K(\delta)),$ $[x]\mapsto [(0, x)]$  - this amounts to
\begin{align}\label{f:bartheta}
   H^1_{{\mf Z}}(H^0_\Psi(\cR_K^+(\delta)))\xrightarrow{1-\varphi} H^1_{{\mf Z}}(\cR_K^+(\delta)^{\Psi=0})\cong  D(\Gamma_L,K)_U\mathbf{e}_\delta\cong K[\Gamma_L/U]\mathbf{e}_\delta.
\end{align}
For the remainder of the section we assume in addition that $\delta$ is de Rham.
We have to compare \eqref{f:bartheta} with
\begin{equation}\label{f:expdual+}
H^1_{{\mf Z}}(H^0_\Psi(\cR_K^+(\delta))) \cong H^1_{\Psi,{\mf Z}}(\cR_K(\delta))\xrightarrow{\exp^{*,(n)}_{\cR_K(\delta)^*}}\mathbf{D}_{\mathrm{dR}}^{(n)}(\cR_K(\delta)). \end{equation}
By the commutativity of the upper square in the (second) diagram of Lemma \ref{lem:epsilon}
one immediately sees that a class $[A_\mu\mathbf{e}_\delta]$ is mapped under \eqref{f:bartheta} to $pr_{\Gamma_n}(Tw_{\delta^{-1}}(\mathrm{Res}_{o_L^\times}(\mu)))\mathbf{e}_\delta$ while under \eqref{f:expdual+} to $\iota_n(A_\mu\mathbf{e}_\delta  )_{|t_{LT}=0}=\theta\circ\iota_n(A_\mu\mathbf{e}_\delta )$ by Definition \ref{def:dualexp} combined with Lemma \ref{lem:g_M}. Recall that $\theta $ was defined above Definition \ref{def:iota}.
Consider the $K[\Gamma_L/\Gamma_n] $-linear map
\begin{equation}
\Sigma:K[\Gamma_L/\Gamma_n] \cong \mathbf{D}_{\mathrm{dR}}^{(n)}(\cR_K(\delta)),
\end{equation}
whose $\rho$-component, for $\rho$ running through the characters of $G_n:=\Gamma_L/\Gamma_n$, is given as the $K$-linear map
\begin{equation}
\mathfrak{e}_\rho\Sigma:K\cong K \mathfrak{e}_\rho\to  \mathfrak{e}_\rho \mathbf{D}_{\mathrm{dR}}^{(n)}(\cR_K(\delta))\cong \mathbf{D}_{\mathrm{dR}} (\cR_K(\delta\rho^{-1})), 1\mapsto \mathfrak{C}(\delta\rho^{-1})\frac{1}{t_{LT}^k}\mathbf{e}_{\delta\rho^{-1}},
\end{equation}
upon noting that   $\mathbf{D}_{\mathrm{dR}} (\cR_K(\delta'))=(L_\infty\otimes_LK\frac{1}{t_{LT}^k}\mathbf{e}_{\delta'})^{\Gamma_L}.$
Here $\mathfrak{e}_\rho:=\frac{1}{|G_n|}\sum_{g\in G_n} \rho(g^{-1})g\in K[G_n]$ denotes the idempotent attached to  $\rho$ satisfying $g\mathfrak{e}_\rho=\rho(g)\mathfrak{e}_\rho$ for all $g\in G_n$, while for an analytic character { $\delta':L^\times \to (F')^\times$   (of weight $k\leq 0$) we set
\begin{equation}\label{f:DeffrakC}
\mathfrak{C}(\delta'):=\frac{(-\Omega)^k}{(-k)!}\left\{
                                             \begin{array}{ll}
                                              { \tilde{\varepsilon}(\cR_{K}(\delta'),u)^{-1},} & \hbox{if $a(\delta')\neq 0$;} \\
                                               \frac{\det(1-q^{-1}\varphi^{-1}|\bD_{\tn{cris}}(\cR_K(\delta')))}{\det(1-\varphi|\bD_{\tn{cris}}(\cR_K(\delta')))}, & \hbox{otherwise}
                                             \end{array}
                                           \right.
\end{equation}
in $L_n\otimes_LK.$}   Unravelling the definition of $\varepsilon^{dR}_{L,u}$   and using Proposition \ref{prop:NakIw2.16} one easily sees   that part (vi) of  Conjecture \ref{conj} is equivalent, for $k\leq 0,$ to the next
\begin{proposition}[Explicit reciprocity formula]\label{prop:explicit} { Let $\delta=\delta_{lc}x^k$ be de Rham.}
For $k\leq 0,$ the following diagram is commutative:\footnote{The factor $C_{Tr}({\mf Z_n})^{-1}$ in the left vertical map takes  \eqref{f:normalization}, i.e., (iii) above into account.}
\begin{equation}
\label{f:Nakamura}
\xymatrix{
  (\cR^+_K(\delta))^{\Psi=1} \ar[dd]_{x\mapsto [(0,C_{Tr}({\mf Z_n})^{-1} x)] }\ar[ddr]^(0.7){C_{Tr}({\mf Z_n})^{-1}\iota_n}\ar[ddrr]^{\frac{q-1}{q}\theta\circ\iota_n} \ar[r]^{1-\varphi } & (\cR^+_K(\delta))^{\Psi=0}  & & D(\Gamma_L,K) \ar[ll]_{\mathfrak{M}_{\delta}\circ \sigma_{-1}}^{\cong}\ar[dd]^{pr_{\Gamma_n}} \\
      &   &  & \\
  H^1_{\Psi,{\mf Z_n}}(\cR_K(\delta))\ar@<0ex>`d[r]`[rr]_{\exp^{*,(n)}}[rr] \ar[r]^{can} & H^1_{\mf Z_n}(D_{dif})  & H^0_{\mf Z_n}(D_{dif})=D_{dR}^{(n)}\cong L_n\otimes_L D_{dR} \ar[l]_(0.6){\cong }^(0.65){g^{(n)}_{\cR_K(\delta)}} &K[\Gamma_L/U],\ar[l]_(0.3){\Sigma} }
\end{equation}
i.e., a  class $[A_\mu \mathbf{e}_\delta]\in H^1_{\mf Z_n}(H^0_\Psi(\cR_K^+(\delta)))^{\Gamma_L} \cong H^1_{\Psi,\mf Z_n}(\cR_K(\delta))^{\Gamma_L},$  is mapped under $\exp^{*}$ to
\[\mathfrak{C}(\delta )p_{\delta^{-1}}( \mu))\frac{1}{t_{LT}^k}\mathbf{e}_{\delta }=\mathfrak{C}(\delta )\int_{ o_L^\times}\delta(x)^{-1} \mu(x)\frac{1}{t_{LT}^k}\mathbf{e}_{\delta }.\]
\end{proposition}


The left hand triangle in \eqref{f:Nakamura} is induced by the commutative diagrams
\begin{equation}\label{}
  \xymatrix{
    \cR^+_K(\delta) \ar@{^(->}[d]_{ } \ar[r]^-{\iota_n} & D_{dif}(\cR_K(\delta))=K_\infty((t_{LT}))\mathbf{e}_\delta \phantom{mmmm} \\
    \cR^{(n)}_K(\delta) \ar[r]^-{\iota_n} & \phantom{}K_n((t_{LT}))\mathbf{e}_\delta   \ar[u]^{ } }
\end{equation}
and
\begin{equation}\label{}
   \xymatrix{
     (\cR^+_K(\delta))^{\Psi=1}/{\mf Z_n} \ar[d]_{\cong } \ar[r]^{\iota_n} &  D_{dif}(\cR_K(\delta))/{\mf Z_n} \ar[d]^{ } \\
     H^1_{\Psi,{\mf Z_n}}(\cR_K(\delta)) \ar[r]^{can} & H^1_{\mf Z_n}(D_{dif}(\cR_K(\delta)) .  }
\end{equation}

The middle triangle is commutative by Lemma \ref{lem:g_M} upon recalling that
\[C_g({\mf Z_n})C_{Tr}({\mf Z_n})=\frac{q}{q-1}\]
by \eqref{f:CTr}, while
 the commutativity of the right upper triangle of diagram \eqref{f:Nakamura}
follows for $k=0$ from the (lower rectangle of the) following lemma (applied to each $\rho$-component) which explains how $\epsilon$-constants show up naturally  in the descent procedure (cf.\ with \cite[Lem.\ 4.9/Cor.\ 4.10]{BB} and \cite[Prop.\ 4.11]{NaANT} in the cyclotomic situation):

\begin{lemma} \label{lem:epsilon} Let $\delta$ be a locally constant character.
Then the following diagram is commutative:
\footnote{Here, the notation of a map $f(-)\mathbf{e}_\delta$ means that $d\mathbf{e}_\delta$ or $d$ is sent to $f(d)\mathbf{e}_\delta.$}
\[\xymatrix{
  (D(o_L,K)\mathbf{e}_\delta)^{\Psi=1} \ar[d]_{A_{(-)}\mathbf{e}_\delta} \ar[rrrr]^{\mathrm{Res}_{o_L^\times}(-)\mathbf{e}_\delta} & && & D(\Gamma_L,K)\mathbf{e}_\delta \ar[d]_{\mathfrak{M}(-)\mathbf{e}_\delta }^{\cong} \\
  (\cR^+_K(\delta))^{\Psi=1} \ar[d]_{\frac{q-1}{q}\frac{1}{[L_n:L]}\mathrm{Tr}_{K_n/K}\circ \theta\circ \iota_n} \ar[rrrr]^{ (1-\varphi)\circ \Psi= 1-\varphi } && & &(\cR^+_K(\delta))^{\Psi=0} && D(\Gamma_L,K)\ar[llu]_{\phantom{mmm}\delta(-1)Tw_\delta(-)\mathbf{e}_\delta}^{\cong}\ar[d]^{p_{\mathbf{1}}(-)  }\ar[ll]_{\mathfrak{M}_\delta\circ\sigma_{-1}  }^{\cong}  \\
 H^0(\Gamma_L,  L_n\otimes_L K \mathbf{e}_\delta ) 
  & & &&&& K . \ar[llllll]_{  \mathfrak{C}(\delta)\mathbf{e}_\delta }  }\]
\end{lemma}

\begin{proof}
The commutativity of the upper rectangle  in this diagram is an immediate consequence of Lemma \ref{lem:Mellin-Amice}, that of the triangle is immediate from the definitions, while that for the lower part  is obviously equivalent to the   commutativity of the outer diagram
\[\xymatrix{
     (D(o_L,K)\mathbf{e}_\delta)^{\Psi=1}\ar[d]_{\frac{q-1}{q}\frac{1}{[L_n:L]}\mathrm{Tr}_{K_n/K}\circ \theta\circ \iota_n(A_{(-)}\mathbf{e}_\delta)} \ar[rrrr]^{\mathrm{Res}_{o_L^\times}(-)\mathbf{e}_\delta } && & &D(\Gamma_L,K)\mathbf{e}_\delta  \ar[d]^{p_{\delta^{-1}}(-)\mathbf{e}_\delta}\\
 L_n\otimes_L K \mathbf{e}_\delta 
  & &&&  K \mathbf{e}_\delta, \ar[llll]_{ \delta(-1) \mathfrak{C}(\delta) }}\]
where $p_{\delta}(\mu):=\int_{ o_L^\times}\delta(x)  \mu(x)$ denotes the evaluation at a character $\delta. $
In order to check this, assume $p\neq 2$ (the case $p=2$ can be dealt with similarly as in the proof of \cite[Pro.\ 4.11]{NaANT}) and first assume that $n:=a(\delta)\geq 1.$ Then we have
\begin{align*}
\mathrm{Tr}_{K_n/K}\circ \theta\circ \iota_n(A_\mu \mathbf{e}_\delta)&=
\sum_{i\in (o_L/\pi_L^n)^\times}\sigma_i\left( \theta\circ \iota_n(A_\mu \mathbf{e}_\delta)\right)\\
&=\sum_{i\in (o_L/\pi_L^n)^\times}\sigma_i\left(   \iota_n(A_\mu \mathbf{e}_\delta)_{|t_{LT}=0}\right)\\
&=\sum_{i\in (o_L/\pi_L^n)^\times}\sigma_i\left(    A_\mu( u_n\otimes1) \varphi^{-n}(\mathbf{e}_\delta) \right)\;\;\;\; \mbox{  in $L_n\otimes K\mathbf{e}_\delta$}\\
&=\bigg(\sum_{i\in (o_L/\pi_L^n)^\times}\delta(i)   A_\mu( \sigma_i\tau(u_n)\otimes1) \bigg)_\tau \varphi^{-n}(\mathbf{e}_\delta) \;\;\;\; \mbox{  in $\prod_\tau K\mathbf{e}_\delta$}\\
&=\bigg(\frac{1}{\delta(\pi_L)^n}\sum_{i\in (o_L/\pi_L^n)^\times}\delta(i)     \int_{o_L}\eta(x,\sigma_i\tau(u_n))\mu(x) \bigg)_\tau \mathbf{e}_\delta \\
&=\bigg(\frac{1}{\delta(\pi_L)^n}\sum_{i\in (o_L/\pi_L^n)^\times} \delta(i)\left(    \int_{o_L}\eta(xi,\tau(u_n))\mu(x) \bigg)_\tau \mathbf{e}_\delta \right)\\
&=\bigg(\frac{1}{\delta(\pi_L)^n}\sum_{i\in (o_L/\pi_L^n)^\times} \delta(i)\sum_{j\in o_L/\pi_L^n}     \eta(ji,\tau(u_n)) \int_{j+\pi_L^n o_L}\mu(x) \bigg)_\tau\mathbf{e}_\delta \\
&=\bigg(\frac{1}{\delta(\pi_L)^n}\sum_{j\in o_L/\pi_L^n}\sum_{i\in (o_L/\pi_L^n)^\times} \delta(i)     \eta(ji,\tau(u_n)) \int_{j+\pi_L^n o_L}\mu(x) \bigg)_\tau\mathbf{e}_\delta \\
&\stackrel{ (\ast)}{=}\bigg( \frac{1}{\delta(\pi_L)^n}\sum_{j\in (o_L/\pi_L^n)^\times}\sum_{i\in (o_L/\pi_L^n)^\times} \delta(i)    \eta(ji,\tau(u_n)) \int_{j+\pi_L^n o_L}\mu(x) \bigg)_\tau\mathbf{e}_\delta \\
&=\bigg(\frac{1}{\delta(\pi_L)^n}\sum_{j\in (o_L/\pi_L^n)^\times} \sum_{i'\in (o_L/\pi_L^n)^\times} \delta(i'j^{-1})    \eta(i',\tau(u_n)) \int_{j+\pi_L^n o_L}\mu(x) \bigg)_\tau\mathbf{e}_\delta \\
&=\bigg(\frac{1}{\delta(\pi_L)^n}\left(\sum_{i\in (o_L/\pi_L^n)^\times}\delta(i)\eta(i,\tau(u_n))\right)\sum_{j\in (o_L/\pi_L^n)^\times}  \delta(j^{-1})    \int_{j+\pi_L^n o_L}\mu(x)\bigg)_\tau \mathbf{e}_\delta \\
&=\bigg(q^{-n(\psi_u)} \epsilon_K(L,K(\delta^{-1}),\psi_{\hat{\tau}u},dx)  \int_{ o_L^\times}\delta(x)^{-1} \mu(x) \bigg)_\tau \mathbf{e}_\delta \\
&=\bigg(q^{-n(\psi_u)} \epsilon_K(L,K(\delta^{-1}),\psi_{\hat{\tau}u},dx) p_{\delta^{-1}}({\mathrm{Res}_{o_L^\times}\mu} ) \bigg)_\tau\mathbf{e}_\delta \\
&=\bigg(q^{a(\delta) }\epsilon_K(L,\delta^{-1}|-|,\psi(x)_{\hat{\tau}u},dx)p_{\delta^{-1}}({\mathrm{Res}_{o_L^\times}\mu} ) \bigg)_\tau \mathbf{e}_\delta \\
&=\bigg(\frac{\delta(-1)q^{a(\delta)+n(\psi_{\hat{\tau}u})}}{ \epsilon_K(L,\delta,\psi_{\hat{\tau}u},dx)}p_{\delta^{-1}}({\mathrm{Res}_{o_L^\times}\mu} ) \bigg)_\tau\mathbf{e}_\delta.\\
\end{align*}
In the two last equalities we used \eqref{f:epsduality2} and \eqref{f:epsabsolutevalue}. Moreover, the equation $(\ast)$ requires part (i) of the next lemma. Finally, by Remark \ref{rem:npsiunull} we have $n(\psi_u)=0,$ whence the result in this case as $[K_n:K]=q^{n-1}(q-1)$ {upon comparing with \eqref{f:DeffrakC}, Example \ref{rem:ExChar} and Definition \ref{def:equivariantconstant}.}

Now we consider the case $a(\delta)=0$   and obtain
\begin{align*}
\mathrm{Tr}_{K_1/K}\circ \theta\circ \iota_1(A_\mu \mathbf{e}_\delta)&=
\sum_{i\in (o_L/\pi_L)^\times}\sigma_i\left( \theta\circ \iota_1(A_\mu \mathbf{e}_\delta)\right)\\
&=\sum_{i\in (o_L/\pi_L)^\times}\sigma_i\left(   \iota_1(A_\mu \mathbf{e}_\delta)_{|t_{LT}=0}\right)\\
&=\sum_{i\in (o_L/\pi_L)^\times}\sigma_i\left(    A_\mu(u_1\otimes 1) \varphi^{-1}(\mathbf{e}_\delta) \right)\\
&=\left(\frac{1}{\delta(\pi_L)}\sum_{i\in (o_L/\pi_L)^\times}    \int_{o_L}\eta(x,\sigma_i\tau(u_1))\mu(x) \right)_\tau \mathbf{e}_\delta \\
&=\left(\frac{1}{\delta(\pi_L)}\sum_{i\in (o_L/\pi_L)^\times}     \int_{o_L}\eta(xi,\tau(u_1))\mu(x)\right)_\tau \mathbf{e}_\delta \\
&=\left(\frac{1}{\delta(\pi_L)}\sum_{i\in (o_L/\pi_L)^\times} \sum_{j\in o_L/\pi_L}     \eta(ji,\tau(u_1)) \int_{j+\pi_L o_L}\mu(x)\right)_\tau \mathbf{e}_\delta \\
&=\left(\frac{1}{\delta(\pi_L)}\sum_{j\in o_L/\pi_L}\sum_{i\in (o_L/\pi_L)^\times}   \eta(ji,\tau(u_1)) \int_{j+\pi_L o_L}\mu(x) \right)_\tau\mathbf{e}_\delta \\
&\stackrel{ (\ast)}{=}\left(\frac{1}{\delta(\pi_L)}\left( (q-1)\int_{\pi_L o_L}\mu(x)-\int_{ o_L^\times}\mu(x)\right) \right)_\tau\mathbf{e}_\delta     \\
&=\frac{1}{\delta(\pi_L)}\left( (q-1)\frac{\delta(\pi_L)}{1-\delta(\pi_L)}-1\right) \int_{ o_L^\times}\mu(x) \mathbf{e}_\delta \\
&=q \frac{1-\frac{1}{q\delta(\pi_L)}}{1-\delta(\pi_L)} p_{\delta^{-1}}({\mathrm{Res}_{o_L^\times}\mu} ) \mathbf{e}_\delta,
\end{align*}
where the fact that $\delta(i)=1$ for all $i\in o_L^\times$ by assumption is used in the fourth and last equality, while part (ii) from the next Lemma is the justification for  the equality $(\ast)$ . The second last equality can be derived from the observation that the condition $A_\mu\mathbf{e}_\delta\in \cR_K^+(\delta)^{\Psi=1}$  implies that $\Psi(A_\mu)=\delta(\pi_L)A_\mu$ by the product formula, whence
\begin{align*}
 \int_{\pi_L o_L}\mu(x) & =\int_{ o_L}\Psi(\mu)(x)= \delta(\pi_L)\int_{ o_L}\mu(x) =\delta(\pi_L) \left( \int_{ o_L^\times}\mu(x)+\int_{\pi_L o_L}\mu(x)  \right).
\end{align*}
It follows that $\int_{\pi_L o_L}\mu(x) =\frac{\delta( \pi_L)}{1-\delta(\pi_L)}\int_{ o_L^\times}\mu(x).$
\end{proof}

\begin{lemma} Assuming   $n=a(\delta)\geq 1$ we have for all $j\in  o_L/\pi_L^n$
\begin{enumerate}
\item $\sum_{i\in (o_L/\pi_L^n)^\times} \delta(i)\eta(ij,u_n)=0$ if $\pi_L|j,$
\item $ \sum_{i\in (o_L/\pi_L)^\times}   \eta(ji,u_1)=\left\{
                                                        \begin{array}{ll}
                                                          q-1, & \hbox{if $\pi_L|j$;} \\
                                                          -1, & \hbox{otherwise.}
                                                        \end{array}
                                                      \right.
 $
\end{enumerate}
\end{lemma}
\begin{proof} If $\pi_L$ divides $j,$ then $\eta(ji,u_1)=1$ for all $i$ and both  statements (for $n=1$) follow by a character sum argument (Note that the assumption $n \geq 1$ asserts that $\delta$ is not trivial). Otherwise the claim (ii) follows from the character formula $\sum_{o_L/\pi_L}\eta(i,u_1)=0$ while for
	 (i)  we may assume $n \geq 2.$ We first show
\be \label{eq:charsum}\sum_{i \in (o_L/\pi_L^n)^{\times}, i \equiv r \mod j'} \delta(i) =0\ee
 for every $r\in  (o_L/\pi_L^n)^\times$ and every proper divisor $j' \mid \pi_L^n.$ By shifting it suffices to consider $r =1.$ In this case we are looking at $$\sum_{i \in H} \delta(i),$$ where $H = \operatorname{ker}(o_L/\pi_L^n)^\times \to (o_L/(j'))^\times.$ This character sum can only be different from zero if $\delta$ is trivial on the subgroup $H,$ contradicting the minimality of $n.$
Without loss of generality assume that $v_{\pi_L}(j)\leq n$, whence $\pi_L^{n}/j$ belongs to $o_L.$
Now let $\mathfrak{R}$ be a system of representatives of $(o_L/(\pi_L^{n}/j))^{\times}$ inside $(o_L/\pi_L^n)^\times$ and rewrite $$\sum_{i\in (o_L/\pi_L^n)^\times} \delta(i)\eta(ij,u_n)= \sum_{r \in \mathfrak{R}}\left( \eta(jr,u_n) \left(\sum_{i, i \equiv r \mod \pi_L^n/j}\delta(i)\right) \right) =0$$ by \eqref{eq:charsum}
applied to $j' = \pi_L^n /j,$ using that $\eta(ji,u_n) = \eta(jr,u_n)$ if $i \equiv r \mod \pi_L^n/j.$

\end{proof}
Proposition \ref{prop:explicit} for $k<0$ we will be reduced to the case $k=0$ by a twisting argument based on the previous subsection \ref{sec:twisting}. Similarly the cases $k>1$ of the following proposition will also be reduced to the case $k=1$. But first we have to slightly modify our notation.
Consider the $K[\Gamma_L/\Gamma_n] $-linear map
\begin{equation}
\Sigma':K[\Gamma_L/\Gamma_n] \cong \mathbf{D}_{\mathrm{dR}}^{(n)}(\cR_K(\delta)),
\end{equation}
whose $\rho$-component, for $\rho$ running through the characters of $G_n:=\Gamma_L/\Gamma_n$, is given as the $K$-linear map
\begin{equation}
\mathfrak{e}_\rho\Sigma':K\cong K \mathfrak{e}_\rho\to  \mathfrak{e}_\rho \mathbf{D}_{\mathrm{dR}}^{(n)}(\cR_K(\delta))\cong \mathbf{D}_{\mathrm{dR}} (\cR_K(\delta\rho^{-1})), 1\mapsto \mathfrak{C}'(\delta\rho^{-1})\frac{1}{t_{LT}^k}\mathbf{e}_{\delta\rho^{-1}}
\end{equation}
with
\begin{equation}
\mathfrak{C}'(\delta'):=\Omega^k(k-1)!\left\{
                                             \begin{array}{ll}
                                             {  \tilde{\varepsilon}(\cR_{K}(\delta'),u)^{-1},} & \hbox{if $a(\delta')\neq 0$;} \\
                                               \frac{\det(1-q^{-1}\varphi^{-1}|\bD_{\tn{cris}}(\cR_K(\delta')))}{\det(1-\varphi|\bD_{\tn{cris}}(\cR_K(\delta')))}, & \hbox{otherwise,}
                                             \end{array}
                                           \right.
\end{equation}
{in $L_n\otimes_LK.$}

\begin{proposition}[Explicit reciprocity formula]\label{prop:explicit2} { Let $\delta=\delta_{lc}x^k$ be de Rham. }
For $k\geq 1,$ the following diagram is commutative:
\begin{equation}
\label{f:Nakamura2}
\xymatrix{
  (\cR^+_K(\delta))^{\Psi=1} \ar[dd]^{x\mapsto [(0,C_{Tr}({\mf Z_n})^{-1} x)] }  \ar[r]^{1-\varphi } & (\cR^+_K(\delta))^{\Psi=0}  & & D(\Gamma_L,K) \ar[ll]_{\mathfrak{M}_{\delta}\circ\sigma_{-1}}^{\cong}\ar[dd]^{pr_{\Gamma_n}} \\
      &   &  & \\
  H^1_{\Psi,{\mf Z_n}}(\cR_K(\delta)) &   & D_{dR}^{(n)}(\cR_K(\delta))\cong L_n\otimes D_{dR}(\cR_K(\delta)) \ar[ll]_{\exp^{(n)}} &K[\Gamma_L/U],\ar[l]_(0.3){\Sigma'} }
\end{equation}
i.e., a  class $[A_\mu \mathbf{e}_\delta]\in H^1_{\mf Z_n}(H^0_\Psi(\cR_K^+(\delta)))^{\Gamma_L} \cong H^1_{\Psi,\mf Z_n}(\cR_K(\delta))^{\Gamma_L},$  is mapped under $\exp^{-1}_{\cR_K(\delta)}$ to
\[\mathfrak{C}'(\delta )p_{\delta^{-1}}( \mu))\frac{1}{t_{LT}^k}\mathbf{e}_{\delta }=\mathfrak{C}'(\delta )\int_{ o_L^\times}\delta(x)^{-1} \mu(x)\frac{1}{t_{LT}^k}\mathbf{e}_{\delta }.\]
\end{proposition}
\begin{proof}
 As mentioned earlier - by the twisting technique - we only have to show the case $k=1$ (i.e. $\delta=\tilde{\delta}x$).   We will show the commutativity of the following diagram
 \[\xymatrix{
   H^1_{\Psi,{\mf Z_n}}(\cR_K(\tilde{\delta})) \ar[d]_{\partial} \ar[rr]^-{\exp^{*,(n)}_{\cR_K(\tilde{\delta}^{-1}x|x|)}} & & \mathbf{D}_{\mathrm{dR}}^{(n)}(\cR_K(\tilde{\delta}))=({ L_\infty\otimes_LK}\mathbf{e}_{\tilde{\delta}})^{\Gamma_n} \ar[d]^{a\mathbf{e}_{\tilde{\delta}}\mapsto \frac{a}{t_{LT}}\mathbf{e}_{\delta}} \\
   H^1_{\Psi,{\mf Z_n}}(\cR_K({\delta}))  & & \mathbf{D}_{\mathrm{dR}}^{(n)}(\cR_K({\delta}))=({ L_\infty\otimes_LK}\mathbf{e}_{{\delta}})^{\Gamma_n} \ar[ll]_-{\exp_{\cR_K({\delta})}}   }
\]
on the image of $(\cR_K^+(\tilde{\delta}))^{\Psi=1}$ in $H^1_{\Psi,{\mf Z_n}}(\cR_K(\tilde{\delta})) $, which together with the diagram \eqref{f:twist} implies the desired formula by comparing the cases $k=1$ and $k=0.$
To this end assume
\[\exp^{*,(n)}_{\cR_K(\tilde{\delta}^{-1}x|x|)}([0,f\mathbf{e}_{\tilde{\delta}}])=\alpha \mathbf{e}_{\tilde{\delta}}\]
with $f\mathbf{e}_{\tilde{\delta}}\in (\cR^+_K(\tilde{\delta}))^{\Psi=1}$. Then it follows from definition \ref{def:dualexp} in combination with Lemma \ref{lem:g_M} that
\[C_g({\mf Z_n})^{-1}[\iota_n(f\mathbf{e}_{\tilde{\delta}})]=[\alpha \mathbf{e}_{\tilde{\delta}}]\in  H^1_{\mf Z_n} (\bD_{\dif}^+(\cR_K(\tilde{\delta}))\]
for sufficiently large $n\geq 1,$ i.e., there exists $y_n\in \bD_{\dif}^+(\cR_K(\tilde{\delta}) )$ such that
\begin{equation}\label{f:condition1}
 C_g({\mf Z_n})^{-1} \iota_n(f\mathbf{e}_{\tilde{\delta}})-\alpha \mathbf{e}_{\tilde{\delta}}={\mf Z_n}y_n.
\end{equation}
By Remark \ref{rem:Zproperties} the element $\nabla \in \Lie(\Gamma_n)$ is divisible by $\mathfrak{Z}_n$ in $D(\Gamma_n,K)$ and the quotient  $\frac{\nabla}{\mathfrak{Z}_n}$ corresponds to $\frac{\Omega}{\pi_L^n}\frac{\log_{LT}(Z)}{Z}$ under the Fourier-LT-isomorphism, which takes the value $\frac{\Omega}{\pi_L^n}=C_g(\mathfrak{Z}_n)$ at $Z=0$ (cf. \eqref{f:DefCg}).

We wish to apply   the $\Psi$-version of Proposition \ref{prop:nak17-2.23} (1) for $\exp_{\cR_K({\delta})}^{(n)} $    with $\tilde{x}=\frac{\nabla}{\mf Z_n}\big(\frac{f}{t_{LT}}\mathbf{e}_\delta\big)$ and $x=\frac{\alpha}{t_{LT}}\mathbf{e}_{\delta}$, which would tell us that
\begin{align*}
  \exp_{\cR_K({\delta})}^{(n)}(\frac{\alpha}{t_{LT}}\mathbf{e}_{\delta}) & =[(\Psi-1)\frac{\nabla}{\mf Z_n}\big(\frac{f}{t_{LT}}\mathbf{e}_\delta\big), {\mf Z_n} \frac{\nabla}{\mf Z_n}\big(\frac{f}{t_{LT}}\mathbf{e}_\delta\big)]\\
    & =[ \frac{\nabla}{\mf Z_n}(\Psi-1)\big(\frac{f}{t_{LT}}\mathbf{e}_\delta\big), {\nabla}\big(\frac{f}{t_{LT}}\mathbf{e}_\delta\big) ] \\
   &  =[0,\partial(f) \mathbf{e}_\delta  ],
\end{align*}
whence the claim. Here, for the last equality we used the formula (iv) of Remark \ref{rem:FColmezFormeln}.
\[\nabla\big(f\frac{1}{t_{LT}}\mathbf{e}_\delta\big)=\big((\nabla+\omega_{\chi_{LT}^{-1}\delta})f\big)\frac{1}{t_{LT}}\mathbf{e}_\delta=t_{LT}\partial(f)\frac{1}{t_{LT}}\mathbf{e}_\delta=\partial(f)\mathbf{e}_\delta\]
noting that $\cR_K(\tilde{\delta})\xrightarrow{\cong}\frac{1}{t_{LT}}\cR_K(\delta), f\mathbf{e}_{\tilde{\delta}}\mapsto \frac{f}{t_{LT}}\mathbf{e}_{{\delta}},$ is an isomorphism of $(\varphi,\Gamma_L)$-modules and that the Hodge-Tate weight of $\omega_{\chi_{LT}^{-1}\delta} $ vanishes.

Thus it remains to verify the assumption of Proposition \ref{prop:nak17-2.23} (1), i.e., $\iota_m(\tilde{x})-x\in \mathbf{D}^+_{\mathrm{dif},m}(\cR_K(\delta))\cong t_{LT}\mathbf{D}^+_{\mathrm{dif},m}(\cR_K(\tilde{\delta}))$ for all $m\geq n.$

From \eqref{f:condition1} and Remark \ref{rem:FColmezFormeln} (iii)   we obtain the equality
\begin{align}\label{f:inductionbeginning}
 C_g({\mf Z_n})^{-1} \iota_n\big(\frac{\nabla}{\mf Z_n}(  f\mathbf{e}_{\tilde{\delta}}) \big) & = C_g({\mf Z_n})^{-1}\alpha \mathbf{e}_{\tilde{\delta}}+\nabla(y_n)   \in C_g({\mf Z_n})^{-1}\alpha \mathbf{e}_{\tilde{\delta}} + t_{LT}\mathbf{D}^+_{\mathrm{dif},m}(\cR_K(\tilde{\delta})).
\end{align}
Using that $(1-\varphi)(f\mathbf{e}_{\tilde{\delta}})={\mf Z_n}\beta$ for some $\beta\in \cR_K(\tilde{\delta})^{\Psi=0}$ by Lemma \ref{lem:Zinvertible}, we conclude for any $m\geq n+1,$
\begin{align*}
   \iota_m\big(\frac{\nabla}{\mf Z_n}(  f\mathbf{e}_{\tilde{\delta}}) \big)&- \iota_{m-1}\big(\frac{\nabla}{\mf Z_n}(  f\mathbf{e}_{\tilde{\delta}}) \big)  \\
   & =  \iota_m\big((1-\varphi)\frac{\nabla}{\mf Z_n}(  f\mathbf{e}_{\tilde{\delta}}) \big)\\
    & =\iota_m\big(\frac{\nabla}{\mf Z_n}( (1-\varphi) f\mathbf{e}_{\tilde{\delta}}) \big)\\
    & =\iota_m\big(\frac{\nabla}{\mf Z_n}( {\mf Z_n}\beta) \big)=\iota_m(\nabla(\beta))\in t_{LT}\mathbf{D}^+_{\mathrm{dif},m}(\cR_K(\tilde{\delta}))
\end{align*}
In particular, we obtain
\[\iota_m\big(\frac{\nabla}{\mf Z_n}(  f\mathbf{e}_{\tilde{\delta}}) \big)- \iota_{n}\big(\frac{\nabla}{\mf Z_n}(  f\mathbf{e}_{\tilde{\delta}}) \big)\in  t_{LT}\mathbf{D}^+_{\mathrm{dif},m}(\cR_K(\tilde{\delta}))\] for any $m\geq n$
by induction. This finishes the proof.
\end{proof}

By an analogous density argument (using the results from Appendix \ref{App:density}) as in \cite[Cor.\ 4.17]{NaANT} the Propositions \ref{prop:explicit} and \ref{prop:explicit2} imply that   $\epsilon_{K}(\cR_K(\delta)):\mathbf{1}_{K}\xrightarrow{\cong} \Delta_K(\cR_K(\delta))$ satisfies conditions (iii), (iv) of Conjecture \ref{conj} for any analytic character $\delta,$ i.e., for any rank one analytic $(\varphi,\Gamma_L)$-module.

\subsubsection{Exceptional case}
\label{sec:exceptional}
This subsection has been inspired by \cite[4B2]{NaANT} and \cite[\S 2.5]{V-Kato}.

By observing that the character $x^0$ is dual to $\chi=x|x|$ with respect to the pairing in Theorem \ref{lem:explicitPairing} and upon applying compatibility with this duality  \ref{prop:duality} as well as with twisting according to Propositions \ref{prop:twisteps} and \ref{prop:twistepsdR} one easily  reduces the verification of condition (vi) in the {\it exceptional case}, i.e., $\delta$ being of the form $x^{-i}$ or $x^i\chi=x^{i+1}|x|$ for $i\in\mathbb{N}$ (recall $0\in\mathbb{N}$), to the case of $\delta=\chi=x|x|.$

First we are going to describe $\bar{\Theta}(\delta).$ To this aim note that the natural inclusion $Kz^0=Pol_{\leq 0}(o_L)\hookrightarrow LA(o_L)$,  which is a splitting of the projection sending $\phi$ to $\phi(0),$ induces a quasi-isomorphism
\begin{equation}\label{f:constantfunctionsquasi}
  K_{\Psi,{\mf Z}}(Kz^0)\hookrightarrow K_{\Psi,{\mf Z}}(LA(o_L))
\end{equation}
by Lemma \ref{f:Chenevier}.
%

The long exact $H^i_{\Psi, D(\Gamma_L,K)}$-sequence attached to
\eqref{f:Colmeztransform} together with \eqref{f:cohR+}, \eqref{f:constantfunctionsquasi} induces for dimension reasons (compare with Lemma \ref{lem:dim} (v)) an isomorphism
\begin{equation}
\label{f:H1Rchi}\alpha_1: H^1_{\Psi, D(\Gamma_L,K)}(\cR_K(\chi))\cong H^1_{\Psi, D(\Gamma_L,K)}(LA(o_L))\cong H^1_{\Psi, D(\Gamma_L,K)}(Kz^0)\cong Kz^0\oplus Kz^0,
\end{equation}
which - induced by the evaluation at $0$ of the Colmez transform given by \eqref{f:phif}  - sends $[f_1\mathbf{e}_\chi,f_2 \mathbf{e}_\chi]$ to
\begin{equation}\label{f:isoH1explizit}
  \left( Res( f_1(Z)g_{LT}(Z)dZ)z^0,  Res( f_2(Z)g_{LT}(Z)dZ)z^0 \right)
\end{equation}
as well as
\begin{equation}
\label{f:H2Rchi}\alpha_2: H^2_{\Psi, D(\Gamma_L,K)}(\cR_K(\chi))\cong H^2_{\Psi, D(\Gamma_L,K)}(LA(o_L))\cong H^2_{\Psi, D(\Gamma_L,K)}(Kz^0)\cong Kz^0 ,
\end{equation}
which  sends $[f\mathbf{e}_\chi]$ to
\begin{equation}\label{f:isoH2explizit}
   Res( f(Z)g_{LT}(Z)dZ)z^0.
\end{equation}
Finally, again as part of the long exact $H^i_{\Psi, D(\Gamma_L,K)}$-sequence attached to
\eqref{f:Colmeztransform}, we have an isomorphism
\begin{equation}\label{f:cohH1R+exceptional}
 \alpha_0:  (Kz^0 \cong) H^0_{\Psi, D(\Gamma_L,K)}(Kz^0)\cong H^1_{\Psi, D(\Gamma_L,K)}(LA(o_L))\cong H^1_{\Psi, D(\Gamma_L,K)}(\cR_K^+(\chi))\cong H^1_{ D(\Gamma_L,K)} ( H^0_{\Psi}(\cR_K^+(\chi))).
\end{equation}

But note that {\it in contrast to the generic case} the canonical map
\begin{equation}\label{f:exceptionalR+R}
  H^1_{\Psi,D(\Gamma_L,K)}(\cR_K^+(\delta))\cong H^1_{D(\Gamma_L,K)}(H^0_\Psi(\cR_K^+(\delta)))\to  H^1_{D(\Gamma_L,K)}(H^0_\Psi(\cR_K^+(\delta)))\hookrightarrow  H^1_{\Psi,D(\Gamma_L,K)}(\cR_K(\delta))
\end{equation}
is the zero map, which can be seen by using \eqref{f:longExact}, \eqref{f:spectralsequ} and counting dimensions. Moreover, we have $ H^0_{\Psi,D(\Gamma_L,K)}(\cR_K^+(\chi))=H^0_{\Psi,D(\Gamma_L,K)}(\cR_K(\chi))=0= H^1_{D(U)}(H^1_\Psi(\cR_K^+(\delta)))=H^2_{\Psi,D(\Gamma_L,K)}(\cR_K^+(\chi))$ by  Lemma \ref{lem:dim} (v) and \eqref{f:cohR+} as well as
 $H^0_{\Psi, D(\Gamma_L,K)}(LA(o_L))\cong H^0_{\Psi, D(\Gamma_L,K)}(Kx^0)\cong Kx^0$ (cf.\   \eqref{f:PolZPsi}.

Altogether it follows that the isomorphism
\[\bar{\Theta}(\chi):\bigotimes_{i=0}^2 \d_{K[\Gamma_L/U]}(H^i_{\Psi_L,{\mf Z}}(\cR_K(\chi)))^{(-1)^{i+1}}\otimes_{K[\Gamma_L/U]}K\cong
(K\mathbf{e}_\chi,1) \]
coincides with the composite
\begin{align}\label{f:exceptionalTheta}
 \notag \bigotimes_{i=1}^2 \d_{K }(H^i_{\Psi_L,D(\Gamma_L,K)}&(\cR_K(\chi)))^{(-1)^{i+1}}    \\
 & \xrightarrow{\alpha}  \bigotimes_{i=0}^2 \d_{K }(H^i_{\Psi_L,D(\Gamma_L,K)}(Kx^0))^{(-1)^{i+1}} \otimes \d_{K}(H^1_{D(\Gamma_L,K)}(H^0_{\Psi}(\cR_K^+(\chi))) )  \\
   \notag & \xrightarrow{\beta\otimes \id} \d_{K }(H^1_{D(\Gamma_L,K)}(H^0_{\Psi}(\cR_K^+(\chi))) ) \xrightarrow{\varrho} (K\mathbf{e}_\chi,1),
\end{align}
 where $\alpha$ is induced by $\alpha_i,$ for $i=0,1,2,$ and $\beta$ is   the canonical isomorphism
 \begin{align}\label{f:exceptionalKz0}
   \bigotimes_{i=0}^2 \d_{K }(H^i_{\Psi_L,D(\Gamma_L,K)}&(Kx^0))^{(-1)^{i+1}}\cong \mathbf{1}_{K}
 \end{align}
which stems from the base change of the trivialisation of $\d_{D(\Gamma_L,K)}(\cT_\Psi(LA(o_L)))$ from (i) of Lemma \ref{lem:isokette}. Finally, $\varrho$ is induced from \eqref{f:bartheta}, i.e., by
\begin{enumerate}
 \item $H^1_{{\mf Z}}(H^0_\Psi(\cR_K^+(\chi)))\cong H^1_{{\mf Z}}(\cR_K^+(\chi)^{\Psi=0})$ induced by $1-\varphi,$  and
\item $H^1_{{\mf Z}}(\cR_K^+(\chi)^{\Psi=0})\xrightarrow{C_{Tr}({\mf Z_n})}(\cR_K^+(\delta)^{\Psi=0})_U \cong D(\Gamma_L,K)_U\cong K[\Gamma_L/U]\mathbf{e}_\chi $ using $(\mathfrak{M}_\chi\circ\sigma_{-1})^{-1}.$
\end{enumerate}

%

Consider the basis $\tilde{f}_0:=z^0,$  $(\tilde{f}_{1,1}:=(z^0,0), \tilde{f}_{1,2}:=(0,z^0))$ and $\tilde{f}_2=z^0$ of $H^0_{\Psi,D(\Gamma_L,K)}(Kx^0),$   $H^1_{\Psi,D(\Gamma_L,K)}(Kx^0)$ and $H^2_{\Psi,D(\Gamma_L,K)}(Kx^0),$ respectively. Then, analogously to  \cite[Lem.\ 4.19]{NaANT} one easily checks that
\begin{equation}\label{f:betatrivial}
  \beta(\tilde{f}_0^*\otimes(\tilde{f}_{1,1}\wedge \tilde{f}_{1,2}  )\otimes \tilde{f}_2^*)=1.
\end{equation}
where $\tilde{f}_i^*$ denotes the dual basis of $\tilde{f}_i$ for $i=0,2.$ So it remains to study the effect of $\varrho.$

 In order to calculate the effect of $\alpha_0$ consider the Coleman power series $g:=g_{\iota(u),u}(T)$\footnote{For $L=\qp$,$\pi_L=p$ odd and $LT=\widehat{\mathbb{G}}_m$ one has $g(Z)=Z$ as $\mathcal{N}(Z)=Z$ in that case. We do not know whether $\prod_{a\in LT_1}(a+_{LT}Z)=(-1)^{v_2(p)}\varphi(Z)$ holds in general? If so, this would have simplified the proof of \cite[Lem.\ 2.5]{SV15}. Moreover, it would simplify the argument here considerably as the use of the reciprocity law is quite a heavy argument. The statement is true in the case that $\varphi(Z)$ is a monic polynomial by the following argument, which was explained to us by Laurent Berger: Observe that the monic degree $q$ polynomial $h(T):=\varphi(T)-\varphi(Z)  \in \operatorname{Quot}(o_{\CC_p}\llbracket Z\rrbracket)[T]$ vanishes precisely at the $a+_{LT}Z, a \in LT_1$ and hence $h(T) = \prod(T-(a+_{LT}Z)).$ Comparing the constant coefficients yields the claim.   } in the notation of \cite[Theorem 2.2]{SV15}, where we consider \[\iota(u)=(u_n \mod u_1)_n\in\varprojlim_n o_{L_n}/u_1 o_{L_n}\cong\varprojlim_{n,Norm} o_{L_n}\] as an element of $\varprojlim_n L_n^\times$,   the group of units of the corresponding field of  norms  $ \mathbf{E}_L$ (cf.\ \cite[Lem.\ 1.4]{KR}).
 \begin{remark}\label{rem:g} The element $\frac{\partial g}{g}$
 \begin{enumerate}
   \item belongs to $\cR_K^{\Psi=\frac{\pi_L}{q}}$ and
   \item satisfies $Res(\frac{\partial g}{g}dt_{LT})=1$.
 \end{enumerate}
 \end{remark}

 \begin{proof}
 By the last sentence of section 2 of \cite{SV15}  the term  $\frac{\partial g}{g}$ belongs to $\cR_K^{\Psi=\frac{\pi_L}{q}}$. By the explicit reciprocity law Prop.\ 6.3 in (loc.\ cit.) we obtain\footnote{ Note the opposite normalisation of the reciprocity map in (loc.\ cit.).}
 \begin{align}\label{f:Res-formula}
   Res(\frac{\partial g}{g}dt_{LT})= Res(\frac{dg}{g} )=\partial_{\varphi}(1)(\mathrm{rec}_{\mathbf{E}_L}(\iota(u))^{-1})=1.
 \end{align}
 Indeed, under the reciprocity map $\mathrm{rec} _{\mathbf{E}_L} $  the inverse of the uniformiser $\iota(u) $ is sent to the Frobenius (lift) $\varphi_q$ of $\mathbf{A},$ whence the cocycle $\partial_{\varphi}(1)$, which is given by sending $h\in H$ to $ha-a$ for some $a\in \mathbf{A}$ with $\varphi_q(a)-a=1,$ sends $\mathrm{rec}_{\mathbf{E}_L}(\iota(u)) $ to $1$ tautologically.
 \end{proof}
The following Lemma should be compared to \cite[Lem.\ 4.20]{NaANT} and \cite[Lem.\ 2.9]{V-Kato}.
\begin{lemma}\label{lem:gamma}
The isomorphism
\[ H^0_{\Psi, D(\Gamma_L,K)}(Kx^0) \xrightarrow{\alpha_0}  H^1_{ D(\Gamma_L,K)} ( H^0_{\Psi}(\cR_K^+(\chi)))  \xrightarrow{\varrho} K\mathbf{e}_\chi \]
sends $\tilde{f}_0$ to  $-\Omega\frac{q-1}{q}\left({\mf Z}(\log (g(T)))\right)_{|T=0}\mathbf{e}_\chi.$
\end{lemma}

\begin{proof}
 By Remark \ref{rem:g}    we obtain an element $\frac{\partial g}{g}\mathbf{e}_\chi\in\cR_K(\chi)^{\Psi=1}$ which lifts $\tilde{f}_0$ under the Coleman transform \eqref{f:phif}. Thus $\alpha_0(\tilde{f}_0)$ is represented by \[{\mf Z}\left(\frac{\partial g}{g}\mathbf{e}_\chi  \right)={\mf Z}\left(\partial \log (g)\mathbf{e}_\chi  \right)=\partial \left({\mf Z}\log g   \right)\mathbf{e}_\chi\] by \eqref{f:twist}.
 It is mapped into $H^1_{D(\Gamma_L,K)}(\cR_K^+(\chi)^{\Psi=0})$ to the class of
 \begin{align*}
  (1-\varphi){\mf Z}\left(\frac{\partial g}{g}\mathbf{e}_\chi  \right)&={\mf Z}(1-\varphi)\left(\partial \log (g)\mathbf{e}_\chi  \right)\in \cR_K^+(\chi)^{\Psi=0}\\
  &={\mf Z}\left(\partial (1-\frac{\varphi}{q})(\log (g))\mathbf{e}_\chi  \right)\\
  &=\left(Tw_{\chi_{LT}}({\mf Z})\partial (1-\frac{\varphi}{q})(\log (g))\right)\mathbf{e}_\chi  \\
  &=\left(\partial {\mf Z}(1-\frac{\varphi}{q})(\log (g))\right)\mathbf{e}_\chi  \\
  &=\partial \left((1-\frac{\varphi}{q}){\mf Z}(\log (g))\mathbf{e}_{|x|} \right). \\
 \end{align*}
 Now we use the commutative diagram
\begin{equation*}
   \xymatrix{
     K \ar[d]_{-\Omega}  & D(\Gamma_L,K) \ar[d]_{-\Omega}\ar[l]_{p_\mathbf{1}} \ar[r]^{\mathfrak{M}_{|x|}} & \cR^r_K(|x|)^{\Psi=0} \ar[d]^{\partial} \\
     K   & D(\Gamma_L,K) \ar[l]_{p_\mathbf{1}}\ar[r]^{\mathfrak{M}_{\chi}} & \cR^r_K(\chi)^{\Psi=0}   }
\end{equation*}
 to conclude by observing that the evaluation at $\mathbf{1}$ corresponds to setting $T=0$ and that $\left((1-\frac{\varphi}{q}){\mf Z}(\log (g(T)))\right)_{|T=0}=\frac{q-1}{q}\left({\mf Z}(\log (g(T)))\right)_{|T=0}$.
\end{proof}

\begin{remark}
The map $v\mapsto (1-\frac{\varphi}{q})(\log (g_{v,u}(T)))$ generalizes Coleman's map as used in Kato's proof of the classical rank one case, cf.\ \cite[(2.5)]{V-Kato}.
\end{remark}

\begin{lemma} \label{lem:Cg} With the notation in the proof of Lemma \ref{lem:gamma} we have
\[\left({\mf Z}_n(\log (g(T)))\right)_{|T=0}=\cL_{{\mf Z_n}}'(\mathbf{1})= C_g({\mf Z}_n).\]
\end{lemma}
\begin{proof}
Note that $g\in To_L[[T]].$ Writing $g=\sum_{i\geq 1} a_iT^i$ we see that for any $\gamma\in\Gamma_n$ we have
\begin{align*}
  ((\gamma-1)\log(g(T)))_{|T=0} &=\log(\frac{g([\chi_{LT}(\gamma)](T))}{g(T)})_{|T=0} \\
  &=\log(\left(\frac{\sum_{i\geq 1}a_i\frac{([\chi_{LT}(\gamma)](T))^i}{T}}{\sum_{i\geq 1} a_iT^{i-1}}\right)_{|T=0} )\\
  &=\log(\frac{ a_1 \chi_{LT}(\gamma) }{  a_1 } ) =\log(  \chi_{LT}(\gamma)   ).
\end{align*}
It follows that for elements $\lambda=\sum_i a_i(\gamma_i-1)$ in the $K$-span $S$ of $\gamma-1,$ $\gamma\in \Gamma_n\setminus{\{1\}}$, in $D(\Gamma_n,K)$ we have
\begin{align*}
  (\lambda\log(g(T)))_{|T=0} & =\sum_i a_i \log(  \chi_{LT}(\gamma_i)   )=\lambda(\log(  \chi_{LT}))= \cL_{{\lambda}}'(\mathbf{1}),
\end{align*}
because $\sum_i a_i(\gamma_i-1)=\sum_i a_i\gamma_i - (\sum_i a_i)1$ and $\log(\chi_{LT}(1))=0.$ Since ${\mf Z}_n$ belongs to the closure of $S$ the claim follows by continuity.
\end{proof}

Now we define a basis $(f_{1,1},f_{1,2})$ of $ H^1_{\Psi_L,D(\Gamma_L,K)}(\cR_K(\chi)) $ and $f_2$ of $H^2_{\Psi_L,D(\Gamma_L,K)}(\cR_K(\chi)) $ via\footnote{Nakamura adds here the factor $\pm \frac{p}{(p-1)\log(\chi_{cyc}(\gamma))}$ in front of the $\alpha_i!$}
\begin{equation}\label{f:basisHchi}
  \alpha_1(f_{1,i})=\tilde{f}_{1,i} \mbox{ for } i=1,2 \mbox{ and } \alpha_2(f_2)=\tilde{f}_2.
\end{equation}

Combining \eqref{f:exceptionalTheta}, \eqref{f:betatrivial} and \eqref{f:basisHchi} with Lemmata \ref{lem:gamma}, \ref{lem:Cg} we obtain
\begin{corollary}\label{cor:thetabar}
 $\bar{\Theta}(\chi)((f_{1,1}\wedge f_{1,2})\otimes f_2^*)=  -\Omega\frac{q-1}{q}C_g({\mf Z}_n)C_{Tr}({\mf Z}_n)\mathbf{e}_\chi=-\Omega\ \mathbf{e}_\chi.$
\end{corollary}

Now we shall compare this to the {\it de Rham } $\epsilon$-isomorphism, i.e., mainly to the map $\Theta(\cR_K(\chi))$, because
\begin{equation}\label{f:Gammafactor}
  \Gamma(\cR_K(\chi))=\Omega^{-1}
\end{equation} and  $\Theta_{F,dR,u}(\cR_F(\chi)):  \d_K(  \bD_{\tn{dR}}(\cR(\chi))  ) \xrightarrow{\cong}\Delta_{K,2}(\cR_K(\chi))$ corresponds to the isomorphism
\begin{equation}\label{f:calL}
  \mathcal{L}_K(\cR_K(\chi))=K\mathbf{e}_\chi\xrightarrow{\cong}\mathbf{D}^{}_{\tn{dR}}({\cR_K(\chi)})=K\frac{1}{t_{LT}}\mathbf{e}_\chi, a\mathbf{e}_\chi\mapsto \frac{a}{t_{LT}}\mathbf{e}_\chi
\end{equation}
as $\tilde{\varepsilon}(\cR_K(\chi),u)=1$ due to $\chi$ being crystalline.

By the long exact sequence \eqref{f:longExactexp} the map $\Theta(\cR_K(\chi))$ is induced from  the following isomorphisms and exact sequences
\begin{equation}\label{f:D_crisphi}
  \mathbf{D}^{}_{\tn{cris}}(\cR_K(\chi))\xrightarrow{1-\varphi_L}\mathbf{D}^{}_{\tn{cris}}(\cR_K(\chi)),i.e.,\;\;K\frac{1}{t_{LT}}\mathbf{e}_\chi\xrightarrow{1-\frac{1}{q}}K\frac{1}{t_{LT}}\mathbf{e}_\chi,
\end{equation}
\begin{equation}\label{f:shortexactexp}
\mathbf{D}^{}_{\tn{dR}}({\cR_K(\chi)})\xrightarrow{\exp_{\cR_K(\chi)}}H^{1}_{\varphi,\mathfrak Z_n}(\cR_K(\chi))^{\Gamma_L}_f  \xrightarrow[\cong]{\Upsilon'_f} H^{1}_{\Psi,\mathfrak Z_n}(\cR_K(\chi))^{\Gamma_L}_f
\end{equation}
 (with $\Upsilon'_f$ induced by $\Upsilon'$ in Remark \ref{rem:asymmetricII})
\begin{equation}\label{f:dualexpexplizit}
    H^{1}_{\Psi,\mathfrak Z_n}(\cR_K(\chi))^{\Gamma_L}/H^{1}_{\Psi,\mathfrak Z_n}(\cR_K(\chi))^{\Gamma_L}_f\xrightarrow{\bar{x}\mapsto \{y\mapsto \langle(\Upsilon')^{-1}(x),y\rangle_{\cR(\chi)}\} }(H^{1}_{\varphi,\mathfrak Z_n}(\cR_K)^{\Gamma_L}_f)^*\xrightarrow{(\exp_{f,\cR_K})^*}\mathbf{D}_{\tn{cris}}(\cR_K)^{*}
\end{equation}
and
\begin{equation}\label{f:H2Dcris}
 \mathbf{D}_{\tn{cris}}(\cR_K)^{*}\cong \H^2_{\Psi,\mathfrak Z_n}(\cR_K(\chi))^{\Gamma_L} \xrightarrow{\bar{x}\mapsto \{y\mapsto \langle (\Upsilon')^{-1}(x),y \rangle_{\cR(\chi)}\}} (H^{0}_{\varphi,\mathfrak Z_n}(\cR_K)^{\Gamma_L} )^* ,
\end{equation}
which is dual to the natural isomorphism $H^{0}_{\varphi,\mathfrak Z_n}(\cR_K)^{\Gamma_L} \cong \mathbf{D}_{\tn{cris}}(\cR_K), 1\mapsto d_0:=1\in K=\mathbf{D}_{\tn{cris}}(\cR_K).$  We define basis $e_0$ and $(e_{1,1},e_{1,2})$ of  $H^{0}_{\varphi,\mathfrak Z_n}(\cR_K)^{\Gamma_L} $ and $H^{1}_{\varphi,\mathfrak Z_n}(\cR_K)^{\Gamma_L},$ respectively, as follows:\footnote{In order to normalize $e_{1,2}$, i.e., to make it independent of the choice of $\mf Z_n$, one would need the factor $C_{Tr}({\mf Z_n})$ from \eqref{f:CTr}, but for our calculations this is not needed. Since in our choice for the generalized  Herr complex the order of the operators $Z$ and $\varphi-1$ (or $\Psi-1$) is the opposite compared to  Nakamura's version, our indexing of the basis elements differs from Nakamura's!}
\[e_0:=1\in\cR_K,\;\;\;  e_{1,1}:=[(1,0)],\;\;\; e_{1,2}:=[(0,1)]. \]
\begin{lemma}\label{lem:explizitexp}
 \begin{enumerate}
   \item $\exp_{f,\cR_K}(d_0)=e_{1,1}$
   \item $ \Upsilon'_f\circ\exp_{\cR_K(\chi)}(t_{LT}^{-1}\mathbf{e}_\chi)=\frac{q-1}{q}f_{1,1}$
   \item 
     Using   the pairing $$\Kl-,-\Kr_{\cR_K(\chi)}\colon H^i_{\Psi,\mf Z_n}(\cR_K(\chi))\times H^{2-i}_{\varphi,\mf Z_n}(\cR_K) \to K$$   from  Remark \ref{rem:asymmetricII} we have
                               \begin{align*}
                                  \Kl f_{1,2},e_{1,1}\Kr_{\cR_K(\chi)}=1, & \phantom{m}\Kl f_{1,1},e_{1,1}\Kr_{\cR_K(\chi)}=0,\\
                                 \Kl f_{1,2},e_{1,2}\Kr_{\cR_K(\chi)}=0, &  \phantom{m} \Kl f_{1,1},e_{1,2}\Kr_{\cR_K(\chi)}=1,\\
                                 \Kl f_2,e_0\Kr_{\cR_K(\chi)}= -1  .&
                               \end{align*}
   \item      $(\exp_{f,\cR_K})^*(\langle (\Upsilon')^{-1}(f_{1,2}),-\rangle_{\cR(\chi)})=d_0^*\in  \mathbf{D}_{\tn{cris}}(\cR_K)^{*}$,
   where the pairing  $\langle-,-\rangle_{\cR(\chi)}$ had been introduced in Theorem \ref{lem:explicitPairing}.
 \end{enumerate}
\end{lemma}

\begin{proof}
(i) follows from Proposition \ref{prop:nak17-2.23}(ii) by taking $\tilde{x}=1.$ For (ii) we apply Proposition \ref{prop:nak17-2.23} (i)   with $\tilde{x}=\frac{f}{t_{LT}}\mathbf{e}_\chi\in\cR_K(\chi)[\frac{1}{t_{LT}}]$, where $f$ lies in $\cR_K^+$ such that $f(u_n)=\frac{1}{\pi_L^n}$ for any $n\geq 0.$ The existence of such $f$  follows from the analogue of \eqref{prod iso} over the ring $\cR_K^+$
\[\cR^+_K/t_{LT}\cong \prod_{n\geq 0} L(u_n), \bar{f}\mapsto (f(u_n))_{n\geq 0}.\]
Moreover, $\tilde{x}$ satisfies
\[\iota_n(\tilde{x})- \tilde{x}\in \mathbf{D}^+_{\mathrm{dif},n}(\cR_K(\chi))\]
for all $n\geq 1,$ because
\[\iota_n(\frac{f}{t_{LT}}\mathbf{e}_\chi)\equiv \pi_L^n\frac{f(u_n)}{t_ {LT}}  \mathbf{e}_\chi=\frac{1}{t_ {LT}}  \mathbf{e}_\chi \mod\mathbf{D}^+_{\mathrm{dif},n}(\cR_K(\chi)) \]
by Remark \ref{iota t}. Therefore the conditions of Proposition \ref{prop:nak17-2.23}(i) are satisfied and hence we conclude that
\[\exp_{\cR_K(\chi)}(t_{LT}^{-1}\mathbf{e}_\chi)=[(\varphi-1)(\tilde{x}),{\mf Z}_n(\tilde{x})  ]\in H^1_{\varphi,{\mf Z}_n}(\cR_K(\chi))\]
and
\[\Upsilon'_f\circ\exp_{\cR_K(\chi)}(t_{LT}^{-1}\mathbf{e}_\chi)=[-\Psi(\varphi-1)(\tilde{x}),{\mf Z}_n(\tilde{x})  ]\in H^1_{\Psi,{\mf Z}_n}(\cR_K(\chi)).\]

Hence it suffices to show that
\begin{equation}\label{f:ResZn}
Res\left(\frac{{\mf Z}_nf}{t_{LT}}d\log_{LT}\right)=0
\end{equation}
and
\begin{align}\label{f:ResII}
   Res_Z\left(\frac{(\frac{\varphi}{q}-1)f}{t_{LT}}d\log_{LT}\right) & = -\frac{q-1}{q},
\end{align}
because $-\Psi(\varphi-1)(\frac{f}{t_{LT}}\mathbf{e}_\chi)=-\frac{q}{\pi_L}\Psi\bigg(\frac{(\frac{\varphi}{q}-1)f}{t_{LT}}\bigg)\mathbf{e}_\chi$ and
\begin{align*}
 Res_Z\left(\Psi\bigg(\frac{(\frac{\varphi}{q}-1)f}{t_{LT}}\bigg)d\log_{LT}\right)=\frac{\pi_L}{q}Res_Z\left(\frac{(\frac{\varphi}{q}-1)f}{t_{LT}}d\log_{LT}\right)
\end{align*}
by \cite[Lem.\ 4.5.1 (iv)]{SV20}.
For \eqref{f:ResZn} one shows first  the analogous statement for $\gamma-1$, $ \gamma\in\Gamma_L$, instead of ${\mf Z}_n$ by similar arguments and then concludes by continuity. For \eqref{f:ResII} we calculate
\begin{align*}
   Res_Z\left(\frac{(\frac{\varphi}{q}-1)f}{t_{LT}}d\log_{LT}\right) & = Res_{t_{LT}}\left(\frac{(\frac{\varphi}{q}-1)f}{t_{LT}}dt_{LT}\right)  \\
   &=  Res_{t_{LT}}\left(\frac{\frac{\varphi}{q}(f)}{t_{LT}}dt_{LT}\right)-Res_{t_{LT}}\left(\frac{ f}{t_{LT}}dt_{LT}\right) \\
   &=  \left(\frac{\varphi(f)}{q} \right)_{|t_{LT}=0}- \left(  f \right)_{|t_{LT}=0} \\
   &=  \left(\frac{\varphi(f)}{q} \right)_{|Z=0}- \left(  f \right)_{|Z=0} \\
    & = \frac{1-q}{q}f(0)=-\frac{q-1}{q},
\end{align*}
where for the first equation the reasoning is as follows: since $f(0)\neq0$ the expression on the left-hand side has a simple pole at $Z=0.$ The residue formula for simple poles gives us $Res_Z\left(\frac{(\frac{\varphi}{q}-1)f}{t_{LT}}g_{LT}dZ \right)= \left(\frac{(\frac{\varphi}{q}-1)f(Z)}{g_{LT}(Z)}g_{LT}(Z)\right)_{\mid Z=0} = Res_{t_{LT}}\left(\frac{(\frac{\varphi}{q}-1)f}{t_{LT}}dt_{LT}\right)$.

(iii) follows by direct computation using the formulae of Remark \ref{rem:asymmetricII}:
\[\Kl f_{1,2},e_{1,1}\Kr_{\cR_K(\chi)}=-Res((\lambda^\iota1)(\frac{\partial g}{g}\mathbf{e}_\chi) )=-p_\mathbf{1}(\lambda^\iota)Res(\frac{\partial g}{g} (\lambda d\log_{Lt}))=1\] upon noting that $f_{1,2}= [(0, \frac{\partial g}{g}\mathbf{e}_\chi)]$ by the   proof and with notation of Lemma \ref{lem:gamma}.

 Take $f_{1,1}=[(\lambda_1\mathbf{e}_\chi,\lambda_2 \mathbf{e}_\chi)]$. Then
\begin{align*}
  \Kl f_{1,1}, e_{1,1}\Kr_{\cR_K(\chi)}&=-Res((\lambda^\iota1)(\lambda_2\mathbf{e}_\chi))=-p_\mathbf{1}(\lambda^\iota)Res(\lambda_2d\log_{LT})=0
\end{align*}
and, for $\lambda$ satisfying $[-1]({\mf Z_n})=\mf Z_n^\iota=\lambda \mf Z_n$,
\begin{align*}
  \Kl f_{1,1}, e_{1,2}\Kr_{\cR_K(\chi)}&= Res(1(\lambda_1\mathbf{e}_\chi))= Res(\lambda_1d\log_{LT})=1
\end{align*}
by definition of $f_{1,1}.$
Finally, writing $f_2=\lambda_3 \mathbf{e}_\chi $ we have
\begin{align*}
   \Kl f_2,e_0\Kr_{\cR_K(\chi)}= & Res(1(\lambda f_2))= Res((\lambda^\iota1)( f_2))=p_\mathbf{1}(\lambda^\iota)Res(\lambda_3 d\log_{LT})=-1
\end{align*}
(iv) follows from (i) and (iii) using   \eqref{f:Upsilon'pairings}, i.e.,   $\langle   (\Upsilon')^{-1}(f_{1,2}),-\rangle_{\cR(\chi)} =\Kl f_{1,2},- \Kr_{\cR_K(\chi)}$.
\end{proof}

Combining the previous lemma with \eqref{f:D_crisphi}, \eqref{f:shortexactexp}, \eqref{f:dualexpexplizit}, \eqref{f:H2Dcris} we obtain
\begin{corollary}
  $\Theta(\cR_K(\chi))((f_{1,1}\wedge f_{1,2})\otimes f_2^*)= -\frac{1}{ t_{LT}}\mathbf{e}_\chi\in K \frac{1}{ t_{LT}}\mathbf{e}_\chi=\mathbf{D}^{}_{\tn{dR}}({\cR_K(\chi)}) $.
\end{corollary}

Together with Corollary \ref{cor:thetabar} and \eqref{f:Gammafactor},\eqref{f:calL} this proves property (vi) for the exceptional case.

%

\newpage

\appendix

\section{Density Argument} \label{App:density}

When verifying that $\varepsilon_{K,u}(M)$ satisfies a given property we frequently require a density argument. This is formally justified as follows:
We reinterpret a given property as a commutative diagram in the category of graded line bundles (hence involving only isomorphisms). E.g. for the property (iii) of Conjecture \ref{conj} one takes the diagram \[\xymatrix{
	\mathbf{1}_K \ar[d]_{\id} \ar[r]^{\varepsilon_{K,u}(M)} &  \Delta(M) \ar[d]^{\delta_{\det M}(a)}  \\
\mathbf{1}_K  \ar[r]^{\varepsilon_{K,a\cdot u}(M)} &  \Delta(M). }\] The commutativity of this diagram for a given property $P$ can be reinterpreted as the automorphism $P(M)$ of $\mathbf{1}_K,$ which is obtained by going around the diagram, taking the constant value $1.$
In this section we will construct a (reduced) rigid analytic space $\mathcal{T}_{an}$ over the normal hull $E$ of $L,$ whose $K'$-points parameterise $L$-analytic $(\varphi_L,\Gamma_L)$-modules attached to characters $\delta\colon L \to K'^{\times}$ such that the map $M \mapsto P(M)$ is a map of rigid analytic spaces $$\mathcal{T}_{an} \to \mathbb{G}_m^{an}.$$ This is the same thing as a
 a global section of $\mathcal{O}_{\mathcal{T}_{an}}^{\times}.$
Since $\mathcal{T}_{an}$ is reduced the vanishing of a global section  (in our case $M\mapsto P(M)-1$)   can be checked on a Zariski dense subset (essentially by definition as we will see below).
Let $\mathcal{W}$ (resp. $\mathcal{T}$)  be the rigid spaces representing the functors
$$\mathcal{W}(X) = \operatorname{Hom}_{cts}(o_L^\times,\Gamma(X,\mathcal{O}_X)^\times)$$
resp. $$\mathcal{T}(X) = \operatorname{Hom}_{cts}(L^\times, \Gamma(X,\mathcal{O}_X)^\times).$$
For the representability of the first functor see  \cite[Lemma 2]{BuzzardAutomorphic}.
The representability of the second functor can be seen by fixing a uniformiser, which provides us with an isomorphism $\mathcal{T} \cong \mathbb{G}_m \times \mathcal{W}$ and we denote by $\mathcal{W}_{an}$ the subspace of locally $L$-analytic characters inside $\mathcal{W},$ where by convention we call a character \textbf{locally analytic} if the composite with the restriction map to $\mathcal{O}(Y)$ is $L$-analytic for any affinoid $Y \subset X$ (this makes sense because $\mathcal{O}(Y)$ is a Banach space). Similarly we define
$\mathcal{T}_{an}.$ Since $o_L^\times$ is open in $L^\times$, we conclude that  a character $\delta  \in \mathcal{T}$ is $L$-analytic if and only if its projection to $\mathcal{W}$ is $L$-analytic. Analogously we get an isomorphism (depending on the choice of uniformiser) $ \mathcal{T}_{an} \cong \mathbb{G}_m \times \mathcal{W}_{an}.$
The representability of $\mathcal{W}_{an}$ is shown in \cite[Proposition 6.4.5]{Eme}. Recall that a character is locally $L$-analytic if and only if its differential at $1$ is $L$-linear. A character $\delta\colon o_L^\times \to \Gamma(X,\mathcal{O}_X^\times)$ (with $X$ affinoid over a Galois closure $E$ of $L$) can be written (locally around $1$) as $$\delta(x) = \sum_{\mathbf{n} \in \NN^\Sigma} a_\mathbf{n} (x-1)^{\mathbf{n}},$$ with some $a_{\mathbf{n}} \in \Gamma(X,\mathcal{O}_X),$ where $\Sigma$ is the set of $\QQ_p$-homomorphisms $\sigma\colon L \to E$ and $(x-1)^{\mathbf{n}}$ is defined as $\prod_{\sigma \in \Sigma} \sigma(x-1)^{n_\sigma}$ where $\mathbf{n} = (n_{\sigma})_{\sigma}.$
The partial derivative at $x=1$ in the direction of $\sigma \in \Sigma$, i.e., the coefficient $a_{e_\sigma}$ of the power series at the $\sigma$-unit-vector is called the $\sigma$\textbf{-part} of the generalised Hodge-Tate weight of $\delta.$

\begin{remark}
	A character $\delta \in \mathcal{W}$ is $L$-analytic if and only if $a_{e_\sigma}=0$ for every $\sigma \neq \id.$
\end{remark}
\begin{proof}This is essentially \cite[Remark 2.7]{Be16}. Note that a character is $L$-analytic if and only if $1$ is an $L$-analytic vector for the corresponding representation. (Loc.\ cit.)  uses the logarithm as a chart around $1 \in o_L^\times$ rather than the map  $x \mapsto x-1.$ Since $\log_{\mathbb{G}_m}(T) = T + \dots$ the coefficients in total degree $1$ are unaffected by the change of charts. This means that our $a_{e_{\sigma}}$ agrees with $\nabla_{\sigma}(1)$ in (loc.\ cit.).
\end{proof}
Recall (cf. \cite[Chapter 3]{bellaichechenevier}) that a subset $Z$ of a rigid analytic space $X$ is called \textbf{Zariski-dense} if the only reduced analytic subset containing $Z$ is $X_{red}.$ For a reduced Stein space this is equivalent to requiring that an analytic function vanishing along $Z$ is identically zero.
An illustrating example is the set $p\mathbb{N} \subset \mathbb{B}^{[0,1)}.$ It is Zariski dense because a function vanishing on $p\mathbb{N}$ has infinitely many zeroes inside the affinoids $\mathbb{B}^{[0,r]}$ and thus vanishes identically along an admissible cover. For $n\gg 0 $ the group $U=\Gamma_n$ of $n$-units is an open subgroup of $o_L^\times$ isomorphic to $o_L.$ Recall that by \cite{ST2} the corresponding character variety  $\mathfrak{X}:=\mathfrak{X}_{\Gamma_n}$ is a smooth one-dimensional quasi-Stein space. For such spaces it is known (cf.\ \cite[Section 1.1]{BSX}) that the divisor of an analytic function has finite support in every affinoid subdomain and a similar argument as before shows, that a set having infinite intersection with infinitely many members of a given increasing family of affinoids $(\mathfrak{X}_m)_{m \in \NN}$ covering $\mathfrak{X}$ is automatically Zariski dense. Note that we have a canonical restriction map $$\mathcal{W}_{an} \to \mathfrak{X},$$ which is finite and flat.

\begin{theorem} If $e <p-1$ then
	the set $W_{int}=\{ x^d \mid d \in \NN\}$ is Zariski dense in $\mathcal{W}_{an}.$
	If $e \geq p-1$ we have that the set
	$\{ \chi \in \mathcal{W}_{an} \mid \chi_{\mid U}=x^d \}$ is Zariski dense in  $\mathcal{W}_{an}.$
\end{theorem}
\begin{proof}
 We first consider the restriction of $x^d$ to the subgroup $\Gamma_n$ as above.
 Recall that $\mathfrak{X}$ is covered by the neighbourhoods $\mathfrak{X}(r)$ consisting of characters taking values inside the disc $|z-1|\leq r$ Using the fact that for any element of $x \in o_L^\times$ we know that $x^{q-1}$ is a $1$-unit and $x^{(q-1)p^N}$ for $N \gg0$ is close to $1$ we conclude that $x^{(q-1)p^m}$ for $m \geq N$ are an infinite family of distinct points inside $\mathfrak{X}(r)$ for $N \gg0.$
  If $e(L/\QQ_p)<p-1$ we can decompose $o_L^\times \cong \kappa^\times \times (1+\pi_Lo_L).$ This allows us to cover $\mathcal{W}_{an}$ by sets of the form $\omega^j\mathfrak{X}(r),$ with $\omega$ the composition of the projection mod $\pi$ and the Teichm\"{u}ller character. Since the powers of $x$ intersect every $\omega^j$-component infinitely many times we can conclude from the preceding reasoning, that $W_{int}$ is Zariski-dense.
 In the general case we consider the finite flat restriction map. Passing to affinoids we first observe that Zariski density inside an affinoid $\operatorname{Sp}(A)$ in the sense above is equivalent to Zariski density in the scheme $\operatorname{Spec}(A)$ since affinoids are Jacobson. Furthermore, because affinoids are noetherian, we can conclude that the restriction of the map $\rho\colon \mathcal{W}_{an} \to \mathfrak{X}$ to a suitable family of affinoids is finitely presented and flat (in the ring-theoretic sense) and hence (universally) open with respect to the Zariski topology. The claim follows from the preceding density statement because openness implies that the preimage of a dense subset of $\mathfrak{X}$ is dense inside $\mathcal{W}_{an}.$ Arguing as in the first part, we can show that the image of $W_{int}$ inside $\mathfrak{X}$ is dense and hence also $\rho^{-1}(\rho(W_{int})) = \{ \chi \in \mathcal{W}_{an} \mid \chi_{\mid U}=x^d \}.$
 \end{proof}
 \begin{remark} Let $F'/L$ be a finite subextension of $K$ and fix some $\delta \in \Sigma_{an}(F').$
 	The map $M \mapsto P(M)$ (for a given property $P$) corresponds to a unique section of $\Gamma(\mathcal{W}_{an}, \mathcal{O}_{\mathcal{W}_{an}})^\times,$ where we identify $\mathcal{W}_{an}$ with the space of analytic twists of $\delta.$
 \end{remark}
\begin{proof}
	We will explain the argument for property (iii). The other properties are treated similarly.
	We consider the isomorphism   \[\varepsilon_{D(\Gamma_L,K),u}(\mathbf{Dfm}(M_K)):\mathbf{1}_{D(\Gamma_L,K)}\xrightarrow{\cong} \Delta_{\mathfrak{X}_{\Gamma_L}}(\mathbf{Dfm}(M_K))\] from Theorem \ref{thm:LTDeform}.
 The validity of property (iii) amounts to the commutativity of the diagram
\[\xymatrix{
	\mathbf{1}_{D(\Gamma_L,K)} \ar[d]_{\id} \ar[r]^{\varepsilon_{D(\Gamma_L,K),u}} & \Delta_{\mathfrak{X}_{\Gamma_L}}(\mathbf{Dfm}(M_K)) \ar[d]^{\delta_{\det \mathbf{Dfm}(M)}(a)}  \\
	\mathbf{1}_{D(\Gamma_L,K)}  \ar[r]^{\varepsilon_{D(\Gamma_L,K),au}} &  \Delta_{\mathfrak{X}_{\Gamma_L}}(\mathbf{Dfm}(M_K)). }\]
Since all arrows are isomorphisms, going around the diagram clockwise (starting at $\mathbf{1}_{D(\Gamma_L,K)}$) amounts to an automorphism of $D(\Gamma_L,K),$ or in other words, an invertible global section of $\mathfrak{X}_{\Gamma_L}.$ The isomorphism $\Gamma_L \cong o_L^\times$ induces an isomorphism $\mathfrak{X}_{\Gamma_L} \cong \mathcal{W}_{an}.$ Hence we get an invertible global section $\mathcal{P}$ of $\mathcal{W}_{an}.$ This allows us to interpret the validity of property (iii) for every twist $M_K(\delta)$ with $\delta \in \mathcal{W}_{an}$ as the section  $\mathcal{P}$ of $\mathcal{W}_{an}$ constructed above specialising to $1$ at every such $\delta.$
\end{proof}
\begin{corollary}\label{cor:density}
	The set $$S:=\{(\lambda,\delta_0) \in \mathcal{T}_{an}(K')= \mathbb{G}_m(K') \times \mathcal{W}_{an}(K')\mid K'/K \text{ finite,} \delta_\lambda\delta_0 \text{ generic and } (\delta_0)_{\mid U} \text{ de Rham }\}$$ is Zariski dense.
\end{corollary}
\begin{proof} Note that the set of de Rham characters contains the set of characters which restrict to a power of $x$ on $U$ and is hence dense in $\mathcal{W}_{an}.$ As a conclusion the analogously defined set without the genericity condition is dense.
For every $d$ there is precisely one $\lambda$ such that $\delta_\lambda x^d$ is non-generic. It is not difficult to see that the set $S$ remains dense.
\end{proof}
To restrict some considerations to $(\varphi_L,\Gamma_L)$-modules arising as a base change from a finite extension of $L$ we introduce the following notion.
\begin{definition}
A character $\rho \colon o_L^\times  \to \CC_p$ is called \textbf{classical}, if its image is contained in $\overline{\QQ_p}.$ Analogously a character $L^\times \to \CC_p$ is called \textbf{classical}, if it takes values in $\overline{\QQ_p}.$
\end{definition}
\begin{remark}~
		\begin{enumerate}
		\item The image of a classical character $\rho\colon o_L^\times \to \CC_p$ is contained in some finite extension $F$ of $\QQ_p.$
		\item A character is classical if and only if its restriction to some open subgroup $U$ takes values inside $\overline{\QQ_p}.$
	\end{enumerate}
\end{remark}
\begin{proof}
 Since $o_L^\times $ is topologically finitely generated we can see that the image of some set of topological generators is contained inside $F^\times$ for a suitable finite extension $F$ of $\QQ_p.$
Moreover, by compactness of $o_L^\times$, its image is contained inside the maximal compact subgroup $o_F^\times\subseteq F^\times$.
Now suppose $\rho(U)\subset F$ for some open subgroup $U\subset o_L^\times .$ Let $\gamma \in o_L^\times ,$ then $\gamma^{[o_L^\times :U]} \in U$ and hence $\rho(\gamma)$ is algebraic over $F.$ Setting $F' = F(\rho(\gamma), \gamma \in \mathfrak{R})$ for a system of representatives $\mathfrak{R} \subset o_L^\times $ of $o_L^\times /U$ we can see that the image of $\rho$ is contained in $F'.$
\end{proof}

\begin{lemma}
	Let $\delta \colon L^\times \to K^\times$ be a  de Rham ($L$-analytic)   character, i.e., such that $\cR_K(\delta)$ is de Rham  in the sense of section \ref{sec:dualexp}. Then $\delta(o_L^\times)\subset\overline{\QQ_p}.$
\end{lemma}
\begin{proof} Let $n \gg0.$
Note that $\mathbf{D}_{\mathrm{dif},n}(\cR_K(\delta))$ embeds $\Gamma_L$-equivariantly into $\prod_{j \in \ZZ}(L_n\otimes_L K t_{LT}^j(\delta))$ and the latter is $\Gamma_n$-equivariantly isomorphic to $\prod_{j\in \ZZ} (\prod_{l=1}^{[L_n:L]} K t_{LT}^j(\delta)).$ The de Rham condition hence forces that $\delta$ agrees with $\chi_{LT}^j$ for some (unique) $j$ when restricted to $\Gamma_n.$ As a consequence the restriction $\delta_{\mid_{o_L^\times}}$ is classical.
\end{proof}

\begin{remark} \label{rem:deRhamchar}
The proof of the previous Lemma shows that any de Rham $L$-analytic character $\delta \colon L^\times \to K^\times$ is of the form
\[\delta=\delta_{\mathrm{lc}}x^k\]
for some $k\in\mathbb{Z}$ and some locally constant character $\delta_\mathrm{lc}:L^\times\to K^\times.$ Vice versa any character of this form is obviously de Rham $L$-analytic.
\end{remark}

\begin{corollary}Using the notation from \ref{cor:density}
	the subset $S'$ of $S$ consisting of classical points is Zariski dense.
\end{corollary}
\begin{proof}
This follows from the following easy observations: The set of characters whose restriction to $U$ is of the form $x^d$ is classical and the subset of $\mathbb{G}_m$ defined by $\overline{\QQ_p}^\times$ both are Zariski dense.
\end{proof}
%

\newpage

Competing interests: The authors declare none

\let\stdthebibliography\thebibliography
\let\stdendthebibliography\endthebibliography
\renewenvironment*{thebibliography}[1]{%
	\stdthebibliography{RJRC23}}
{\stdendthebibliography}

\bibliographystyle{abbrv}

\bibliography{rkoneJIMJ}

\end{document}